\documentclass{article}

\usepackage{quartics-layout}
\usepackage{quartics}

\usepackage[
	style=alphabetic,
	giveninits=true,
]{biblatex}
\title{Rational cohomology of the space of equations of GIT-stable plane quartic curves}
\author{Maksym Dolgikh \\ Institut f\"ur Mathematik, Freie Universit\"at Berlin}

\begin{document}

\maketitle

\begin{abstract}
	We study the space of equations of GIT-stable quartic curves in the complex projective plane, that is, of plane quartics with neither a triple point nor a tacnode.
	We compute the rational cohomology of this space as a graded vector space, together with its mixed Hodge structure.
	The computation uses the Gorinov--Vassiliev method of conical resolutions: we construct a conical resolution of the discriminant and determine the resulting spectral sequence, including all of its differentials, from which the result follows.
	This result is the basis of two follow-up papers, which determine the cohomology rings of this space, of the space of GIT-stable quartic curves and of their moduli space.
\end{abstract}

\section{Introduction}
\label{sec:intro}

Spaces of non-singular sections of algebraic vector bundles have been extensively studied in the literature, but they remain largely enigmatic.
Their rational cohomology is known in a range of degrees \cite{Tom14, Aum25}, and completely only in a few cases of low degree \cite{Vass99, Gor05}.
If one is instead interested in spaces of sections with mild singularities, further complications arise.
For example, if one is interested in nodal or GIT-stable projective hypersurfaces, the singular loci become subvarieties of some flag variety and not of the ambient projective space.
The case of nodal plane cubics and quartics was treated by Berdnikov, Gorinov and Konovalov \cite{BerGorKon23}.
In this paper we consider GIT-stable quartic curves in $ \C P^2 $, and calculate the cohomology of the space of equations of such curves.

Let $ \C[x_0,x_1,x_2]_4 $ denote the $15$-dimensional complex vector space of homogeneous polynomials of degree $4$.
By \cite[Proposition 4.2]{mfk94}, a quartic in $\C P^2$ is \emph{not} GIT-stable if and only if it has a triple point or a tacnode (see also \cite[Lemma 1.4]{Art09}).

\begin{definition}
	Let $ Q = \V(F) \subset \C P^2 $ be a curve corresponding to a polynomial $ F $, and let $ p \in \C P^2 $ be a point.
	\begin{itemize}[leftmargin=*]

		\item The point $ p \in \C P^2 $ is an \emph{$ m $-point} if in an affine chart centred at $ p $, the homogeneous components $ F_k $ of $ F $ of degree $ k $ are zero for all $ k < m $.

		\item A line $ l = \V(L) \subset \C P^2 $ is an \emph{$ m $-fold tangent at $ p $} if in an affine chart centred at $ p $, the homogeneous components $ F_k $ of $ F $ of degree $ k $ are multiples of $ L $ for all $ k < m $.

		\item A \emph{cusp} is a double point of a curve with only one $3$-fold tangent line, i.e. a point whose projective tangent cone consists of a single point.
		A \emph{(classical) tacnode} is a cusp whose $3$-fold tangent is also a $4$-fold tangent.
	\end{itemize}
\end{definition}

Let $\Sigma \subset \C[x_0,x_1,x_2]_4$ be the discriminant consisting of polynomials corresponding to quartics with a triple point or a tacnode, and set
\begin{equation*}
	\UU \coloneq \C[x_0,x_1,x_2]_4 \setminus \Sigma .
\end{equation*}
Thus $ \UU $ is the space of equations of GIT-stable quartic curves, that is, of the quartics whose only singular points are ordinary nodes and ordinary cusps.
A tacnode has a distinguished tangent line, and in \autoref{sec:tacnodes} the singular points that define $ \Sigma $ are recorded as pairs of a point and a line through it, which is how the variety of flags in $ \C P^2 $ enters.

For a space $ X $, we denote by $ H^*(X) $ its rational singular cohomology, and by $ \overline H_*(X) $ its Borel--Moore homology.
Let $\Q(w)$ denote Tate mixed Hodge structure on $\Q$ of weight $-2w$.

The main result of this paper is the computation of the graded vector space $ H^*(\UU) $ and its mixed Hodge structure.
\begin{alphatheorem}[\autoref{thm:cohom_of_equations}]
	\label{thm:intro_a}
	As a graded vector space, together with its mixed Hodge structure, the singular cohomology $ H^*(\UU) $ of the space $ \UU $ of equations of GIT-stable quartics is given by the table below.
	\begin{equation*}
		\resizebox{0.95\textwidth}{!}{$
			\begin{array}{c|ccccccccc}
				n      & 0     & 5      & 7                 & 9      & 10     & 12                & 14     & 15     & 17      \\ \hline
				H^n(\UU) & \Q(0) & \Q(-3) & \Q(-4)^{\oplus 2} & \Q(-5) & \Q(-6) & \Q(-7)^{\oplus 3} & \Q(-8) & \Q(-9) & \Q(-10)
			\end{array}
		$}
	\end{equation*}
\end{alphatheorem}

It is instructive to compare this table with the two cases in which the answer was known, those of smooth and of nodal quartics, whose spaces of equations are open subsets of $ \UU $.
Their Poincar\'e polynomials are $ (1+t)(1+t^3)(1+t^5)(1+t^6) $ \cite{Vass99, PetSte03} and $ (1+t^3)(1+t^5)(1+t^7) $ \cite[Theorem 2]{BerGorKon23}.
In both, the factor $ (1+t^3)(1+t^5) $ is the Poincar\'e polynomial of $ \SL_3(\C) $, which acts on these spaces by linear changes of coordinates: the orbit maps induce a surjection onto $ H^*(\SL_3(\C)) $ \cite{PetSte03, Kon22}.
The table of \autoref{thm:intro_a} has no such factor, since $ H^3(\UU) = 0 $.
The lowest class now lies in $ H^5(\UU) $.
It is the linking class of $ \Sigma $, which has complex codimension three in $ \C[x_0,x_1,x_2]_4 $.
The orbit map sends it to the generator of degree $ 5 $ of $ H^*(\GL_3(\C)) $, and every other class of positive degree to zero, as will be shown in a follow-up paper.
In particular, $ H^*(\UU) $ is not the tensor product of the cohomology of the group and the cohomology of a quotient, and how the cohomology of $ \UU $ is related to that of its quotients has to be worked out separately.

This is done in two follow-up papers, for which \autoref{thm:intro_a} is the starting point.
The first of them determines the cohomology ring of the moduli space $ \UU/\GL_3(\C) $ of GIT-stable quartics, which is a quotient of $ H^*(\BGL_3(\C)) = \Q[s_1, s_2, s_3] $ by an explicit ideal, and then the ring structure of $ H^*(\UU) $ itself.
The second determines the cohomology ring of the space $ \UU/\C^{\times} $ of GIT-stable quartic curves, the morphisms induced by the natural quotient maps between $ \UU $, $ \UU/\C^{\times} $ and $ \UU/\GL_3(\C) $ and by the orbit maps, and the morphisms induced by the inclusion of the space of nodal quartics into the space of GIT-stable quartics.
The three papers are based on the author's PhD thesis, written at Freie Universit\"at Berlin.
Each of them can be read on its own: the later papers restate what they use from the earlier ones.

\subsection*{The method and its applications}

At the heart of the proof of \autoref{thm:intro_a} is a technique, introduced by V. Vassiliev \cite{Vass92}, known as the method of conical resolutions.
By Alexander duality, the cohomology of $ \UU = \C[x_0,x_1,x_2]_4 \setminus \Sigma $ can be determined from the Borel--Moore homology of $ \Sigma $.
The latter is described as the limit of the spectral sequence of a certain filtered space, a conical resolution of $ \Sigma $.
This filtered space is typically more tractable than the one we started with.

This technique is adapted in \cite{Vass99, Gor05} to the problem of finding the cohomology of ``non-singular'' projective algebraic hypersurfaces, and in this form it is also known as the Gorinov--Vassiliev method.
Here ``non-singular'' is defined axiomatically and includes in particular non-singular in the usual sense as well as GIT-stable, nodal or GIT-semistable.
O. Tommasi \cite{Tom05} further augmented the method with a way of finding mixed Hodge structures on the cohomology.
Further tools and techniques are presented in \cite{BerGorKon23}, and we use the method as it is set up in \cite[Section 1]{BerGorKon23}.
An overview of the method, illustrated by the example of smooth plane cubics, is given in \autoref{sec:v1_method}.

The method has been used to compute cohomology of some examples which are difficult to handle otherwise, such as the space of non-singular curves of degree $ 5 $ in $ \C P^2 $ \cite{Gor05}, the moduli space of genus $ 4 $ curves \cite{Tom05}, the moduli space of stable complex genus $ 4 $ curves \cite{BerTom07}, and the moduli space of non-singular complex projective curves of genus $ 3 $ with two marked points \cite{Tom07}.
In recent years, the method was used to study the moduli space of trigonal curves of genus $ 5 $ \cite{Zhe22} and of general genus \cite{Zhe24}, the moduli space of non-singular genus $ 4 $ curves with one marked point \cite{WonZhe25}, and the moduli spaces of hyperelliptic curves of a fixed genus on a fixed Hirzebruch surface \cite{BerZhe25}.
Other applications are \cite{Tom10, PagTom13, Tom14}.

The method can be applied to a broad class of other spaces.
For example, one could use it to compute the cohomology of the spaces associated with GIT-semistable quartics in $ \C P^2 $ and cubics in $ \C P^3 $.
To the best of our knowledge, the cohomology of the spaces of equations and of hypersurfaces remains unknown in both cases.
The method can also be applied to find the cohomology of spaces associated with non-singular real quartics in $ \PP^2 $ and real cubics in $ \PP^3 $.
In the real case, the phrase ``non-singular'' can have two meanings: we may require that the corresponding complex hypersurface should have no singularities, or that it should have no real singularities.
Both versions of the problem are interesting, although the first one seems more tractable.

\subsection*{Strategy of the proof}

The proof of \autoref{thm:intro_a} proceeds in two main steps.

First, using the method of conical resolutions, we construct a Borel--Moore spectral sequence that converges to the Borel--Moore homology of $ \Sigma $, and compute its first page, given in \autoref{table:main_sp_seq_p1}.
In particular, we construct a conical resolution of $ \Sigma $, which has the same Borel--Moore homology as $ \Sigma $, and a filtration on it that gives us the spectral sequence.
The construction starts from thirteen configuration spaces (\autoref{table:configs_def}), which contain the singular loci of all quartics in $ \Sigma $.
We verify that they satisfy the conditions required for the method (\autoref{sec:v1_conditions}), compute for each configuration the dimension of the space of quartics that are singular along it (\autoref{sec:v1_dimensions}), and from this obtain the first page (\autoref{sec:v1_first_page}).

Next, we determine all the differentials in the spectral sequence, and therefore its limiting term, given in \autoref{table:main_sp_seq_pinf}.
It turns out that every differential that can be non-trivial is indeed non-trivial.
To show this, we express the differentials through cap products with Euler classes (\autoref{sec:diff_vs_leq}, following \cite[Lemma 18]{BerGorKon23}), and compute these Euler classes using the framework of homogeneous vector bundles introduced in \autoref{sec:hvb}.
The differentials themselves are computed in \autoref{sec:diff_2_to_1}, \autoref{sec:diff_9_to_8} and \autoref{sec:diff_13_to_10}.

Both steps are written out in full detail, in the hope that this will also be of use to readers who wish to apply the method to other spaces.

\section{The method of conical resolutions}
\label{sec:v1_method}

We begin with an overview of the general method, and then illustrate it using the non-singular cubic curves in $ \C P^2 $ as an example.

\subsection{General overview}

The type of problems to which this method applies is as follows.
Suppose $ \Pi $ is a finite-dimensional vector space over a field $ k $, and elements of $ \Pi $ parameterise geometric objects of some kind.
These objects are typically subdivided into singular $ \Sigma \subset \Pi $ and non-singular (generic) $ \Pi \setminus \Sigma $.
The subset $ \Sigma $ of all singular elements is a \emph{generalised discriminant}.
The problem that can be solved with the method is the computation of cohomology groups of the space $ \Pi \setminus \Sigma $ of non-singular objects.

It was observed by V. I. Arnold that sometimes it is more convenient to reduce computation of cohomology groups of the space of generic objects to computing the Borel--Moore homology groups of its complement, the space of singular objects.
Indeed, if our field $ k $ is either $ \mathbb{R} $ or $ \mathbb{C} $, and if $ \Sigma \subset \Pi $ is good enough (e.g. a sub-polyhedron of some triangulation of $ \Pi $), then by Alexander duality we have an additive isomorphism
\begin{equation*}
	\widetilde{H}^*(\Pi\setminus \Sigma) = \overline{H}_{\dim_{\mathbb{R}}(\Pi) - * - 1}(\Sigma)
\end{equation*}
where on the right we have the Borel--Moore homology groups.

To compute these Borel--Moore homology groups we can then use the {\em method of conical resolutions}, introduced by Victor Vassiliev \cite{Vass92}.
To this end, an auxiliary space $ \sigma $, called a \emph{resolution} of $ \Sigma $.
By construction, there is a proper map $ \pi \colon \sigma \to \Sigma $, that induces an isomorphism of the Borel--Moore homology groups, and a natural filtration
\begin{equation*}
	\emptyset \subset \sigma_1 \subset \ldots \subset \sigma_N = \sigma
\end{equation*}
such that the difference $ \sigma_i \setminus \sigma_{i - 1} $ of any two consecutive terms is a fibre bundle over a configuration space of some kind.
We use these fibre bundles to calculate all $ \overline{H}_*(\sigma_i \setminus \sigma_{i - 1}) $.
We then combine the results to get the first page of the spectral sequence associated with this filtration.

The resolution and the filtration on it are both constructed using an ordered family $ X_1, \ldots X_N $ of configuration spaces that satisfy certain conditions.
The configuration spaces must be compact non-empty subspaces of some compact metric space $ M $.
Choosing $ M $ and $ X_i $ is a bit of an art, but making a good choice can simplify the spectral sequence as a result.

To build this family, one first needs to define a subset $ \Sing(f) \subset M $ for each $ f \in \Sigma $.
It has to satisfy certain conditions that require, roughly speaking, that $ L(K) \coloneq \{ f \in \Pi \ |\ \Sing(f) \supset K \} $ is a vector subspace.
These singular loci are then classified to produce the family $ X_i $, which must itself satisfy certain requirements necessary for the method to work.
Note that $ \Sigma = \{ f \in \Pi \ | \ \Sing(f) \neq \emptyset \} $.

These configuration spaces are then used to construct a self-join (symmetric join).
Let $ X = \bigcup X_i $.
Then an \emph{N-th self-join} $ X^{*N} $ is, roughly speaking, a simplicial complex spanned by the points of $ X $ as vertices (which are also compact subspaces of $ M $).
More formally, let $ X $ be a closed subset of a finite-dimensional Euclidean space and let $ N $ be a natural number.
Let us call a proper embedding $ \iota \colon X \to \mathbb{R}^{\Omega} $, $ \Omega < \infty $, \emph{$N$-generic} if every $ (N - 1) $-simplex with vertices in $ \iota (X) $ is non-degenerate and the intersection of any two such simplices is their common face.
In particular, the intersection is empty if the sets of the vertices are disjoint.
The space $ X^{*N} $ is then defined as the union of all $ (N - 1) $-simplices in $ \mathbb{R}^{\Omega} $ with all vertices in $ \iota (X) $ for an arbitrary $N$-generic $ \iota $.
For more details, consult \cite{BerGorKon23}.

The conical resolution is then constructed using some subspaces of the self-join.
First, let us call a simplex $ \Delta \subset X^{*N} $ \emph{coherent} if its vertices form a chain, i.e. for each pair of vertices $ K_1, K_2 $ of $ \Delta $, either $ K_1 \subset K_2 $ or $ K_2 \subset K_1 $ as subsets of $ M $.
The \emph{principal vertex} of a coherent simplex is defined as the vertex that interpreted a set contains all the others.
Let $ \widetilde{ \Lambda } $ be the union of all coherent simplices.
For $ K \in X $, let us denote by $ \widetilde{ \Lambda }(K) $ the union of all coherent simplices with principal vertex $ K $.
Then the \emph{conical resolution} $ \sigma $ of $ \Sigma $ is the subspace of $ \Sigma \times \widetilde{ \Lambda } $ that is formed by all pairs $ (f, x) $ such that $ f \in \Sigma $, $ x \in \widetilde{ \Lambda }(\Sing(f)) $.

Let $ \Phi_0 \coloneq \emptyset $ and $ \Phi_i \coloneq \bigcup_{j \leq i} \bigcup_{K \in X_j} \widetilde{ \Lambda }(K) $ for $ i \in \{1, \ldots , N \} $.
Then we have a filtration on $ \widetilde{ \Lambda } $:
\[
\emptyset = \Phi_0 \subset \Phi_1 \subset \ldots \subset \Phi_N = \widetilde{ \Lambda } .
\]
There exist natural projections $ \pi \colon \sigma \to \Sigma $ and $ p \colon \sigma \to \widetilde{ \Lambda } $.
We introduce a filtration on $ \sigma $ by setting $ \sigma_i = p^{-1}(\Phi_i) $.

Although we have outlined the method, there are many technical details in it that provide tools for practical calculations.
Among many other things, the details concern computations of $ \overline{H}_*(\sigma_i \setminus \sigma_{i - 1}) $.
Further, if our field $ k $ is $ \mathbb{C} $ and the space $ \Sigma $ of singular objects is affine algebraic, then the machinery we have allows us to also find the mixed Hodge structure on the cohomology.
Many other shortcuts and techniques are present in \cite{BerGorKon23}.

\subsection{Example: smooth projective plane cubics}

We illustrate the above definitions with the classical example of smooth cubic curves in $\C P^2$.
These correspond to elliptic curves, whose moduli space is well understood.

Consider the space $ \Pi \coloneq \C [x_0, x_1, x_2]_3 $ of homogeneous degree $ 3 $ polynomials in $3$ variables.
Let $ M \coloneq \mathbb{C}P^2 $ be the complex projective plane.
For $ f \in \Pi $, let $ \Sing(f) $ be the set of $ x \in \mathbb{C}P^2 $ such that $ \dif f(x) = 0 $, i.e. the set of singular points of $ f $ in the usual sense.
Let $ \Sigma $ and $ \UU \coloneq \Pi \setminus \Sigma $ be as above.
One possible choice of configuration spaces is:
\begin{enumerate}
	\item $ X_1 = B(\mathbb{C}P^{2} , 1) = \mathbb{C}P^{2}$
	\item $ X_2 = B(\mathbb{C}P^{2} , 2) $
	\item $ X_3 = \{ l \ |\ l \subset \mathbb{C}P^2 \text{ is a line} \} $
	\item $ X_4 = \{ \{x_1, x_2, x_3\} \subset \mathbb{C}P^2 \ | \ x_1, x_2, x_3 \text{ are not collinear} \} $
	\item $ X_5 = \{ \mathbb{C}P^2 \} $
\end{enumerate}
These configuration spaces satisfy conditions 1-8 from \cite[Condition list 1]{BerGorKon23}.
For example, we could insert the space ${ \{ \{x_1, x_2, x_3\} \subset \mathbb{C}P^2 \ | \ x_1, x_2, x_3 \text{ are collinear} \} }$ of collinear triples between $X_2$ and $X_3$.
The resulting family of configuration spaces would in addition satisfy condition 5+ from the same list.
We will shortly see that rationally, the resulting spectral sequence is essentially the same.
In general, there is some degree of freedom in choosing a collection of configuration spaces for any given problem, and the choice one makes may simplify or complicate the computations, and improve or impair convergence.

Starting from the spaces $X_1,\ldots, X_5$ we construct a space $\sigma$ and a filtration $\varnothing \subset \sigma_1\subset\cdots\subset \sigma_5$ on it as sketched above.

\begin{lemma*}
	\begin{enumerate}[(1), leftmargin=*]
		\item The filtration $ (\sigma_i) $ corresponding to the configuration spaces $ X_i $ gives rise to a spectral sequence whose first page $ E^1 $ is given in \autoref{table:intro_sp_seq_p1}.

		\item Moreover, all differentials vanish, so $ E^1 = E^{\infty} $.
	\end{enumerate}
\end{lemma*}
\begin{table}[!htbp]
	\centering
	\caption{$ E^{1} = E^{\infty} $}
	\label{table:intro_sp_seq_p1}
	\renewcommand{\arraystretch}{1.2}
	\begin{tabular}{|r|*{5}{>{\centering\arraybackslash}m{0.7cm}|}}
		\hline
		17 & $ \Q(9)$ &         &   &         &            \\
		\hline
		16 &          &         &   &         &            \\
		\hline
		15 & $\Q(8)$  &         &   &         &            \\
		\hline
		14 &          &         &   &         &            \\
		\hline
		13 & $\Q(7)$  & $\Q(7)$ &   &         &            \\
		\hline
		12 &          &         &   &         &            \\
		\hline
		11 &          & $\Q(6)$ &   &         &            \\
		\hline
		10 &          &         &   &         &            \\
		\hline
		9  &          & $\Q(5)$ &   &         &            \\
		\hline
		8  &          &         &   &         &            \\
		\hline
		7  &          &         &   &         &            \\
		\hline
		6  &          &         &   & $\Q(4)$ &            \\
		\hline
		& 1        & 2       & 3 & 4       & 5          \\
		\hline
	\end{tabular}
\end{table}

This allows us to compute the Borel--Moore homology of $ \Sigma $, to which the spectral sequence converges.

\begin{proposition*}
	As a graded vector space, together with its mixed Hodge structure, the Borel--Moore homology of $ \Sigma $ is given by the table below.
	\begin{equation*}
		{
			\begin{array}{c|ccccccc}
				n                     & 10     & 11      & 13      & 14    & 15    & 16    & 18     \\
				\hline
				\overline H_n(\Sigma) & \Q(4)  & \Q(5)   & \Q(6)   & \Q(7) & \Q(7) & \Q(8) & \Q(9)
			\end{array}
		}
	\end{equation*}
\end{proposition*}

By Alexander duality, we obtain the cohomology of the space $ \UU $ of non-singular cubic curves in $ \C P^2 $.

\begin{proposition*}
	As a graded vector space, together with its mixed Hodge structure, the cohomology of $ \UU $ is given by the table below.
	\begin{equation*}
		{
			\begin{array}{c|ccccccccc}
				n      & 0     & 1      & 3       & 4      & 5      & 6       & 8      & 9         \\
				\hline
				H^n(\UU) & \Q(0) & \Q(-1) & \Q(-2)  & \Q(-3) & \Q(-3) & \Q(-4)  & \Q(-5) & \Q(-6)
			\end{array}
		}
	\end{equation*}
\end{proposition*}

This proposition shows that $ H^*(\UU) = H^*(\GL_3(\C)) $.

\begin{corollary*}
	\begin{enumerate}[(1), leftmargin=*]
		\item As a graded vector space, together with its mixed Hodge structure, the cohomology of the space of non-singular cubic curves $ \UU / \C^{\times} $ is given by the table below.
		\begin{equation*}
			{
				\begin{array}{c|ccccc}
					n                  & 0      & 3       &  5     & 8          \\
					\hline
					H^n(\UU/\C^{\times}) & \Q(0)  & \Q(-2)  & \Q(-3) & \Q(-5)
				\end{array}
			}
		\end{equation*}

		\item The cohomology of the moduli space $ \UU/\GL_3(\C) $ is trivial, i.e.
		\[
		H^*(\UU/\GL_3(\C)) = H^*(pt) .
		\]
	\end{enumerate}
\end{corollary*}

By the corollary, $ H^*(\UU/\C^{\times}) = H^*(\PGL_3(\C)) = H^*(\SL_3(\C)) $, and $ H^*(\UU/\GL_3(\C)) = H^*(\mathbb{A}^1) $.
This agrees with the classical description of isomorphism classes of elliptic curves via the $ j $-invariant: the coarse moduli space $ M_{1,1} $ of elliptic curves is isomorphic to the affine line $ \mathbb{A}^1_{\C} $.

To demonstrate how the method can be applied in practice, we provide a proof of part $(1)$ of the above key lemma.
This illustrates ideas and techniques used in the method.

\begin{proof}{\it (Lemma)}

	Throughout the proof, we use the explicit form of $ \overline H_*(B(\C P^k, n); \pm \Q) $, where $B(X,n)$ denotes the configuration space of unordered $n$-tuples of distinct points in $X$.
	See \cite[Lemma 2]{Vass99}.

	We also use the following two fibre bundle structures.
	First, $ \sigma_i \setminus \sigma_{i-1} $ is a vector bundle over $ \Phi_i \setminus \Phi_{i-1} $ with fibre $ L(K) \simeq \C^{d_i} $, where $ K \in X_i $ and $ d_i $ is the same for all such $ K $.
	Second, $ \Phi_i \setminus \Phi_{i-1} $ is a fibre bundle over $ X_i $ with fibre a certain subspace of $ \widetilde{\Lambda}(K) $.
	See \autoref{sec:defs_and_core_tools}.

	Throughout the proof we use notation $ V = [\overset{n_1}{\Q(w_{n_1})}, \ldots , \overset{n_k}{\Q(w_{n_k})}] $ to mean that the graded vector space $ V $, together with its MHS structure, is equal to $ \Q(w_{n_i}) $ in degree $ n_i $, and is equal to zero in degrees that are not listed.

	{\boldmath $ E^1_{1,*} = \overline H_{*+1}(\sigma_{1}) $.}

	The space $ \sigma_{1} $ is the total space of a fibre bundle with base $ X_1 = \C P^2 $ and generic fibre $ \C^7 $.
	As a result, \[
	E^1_{1,*} = \overset{14}{\Q(7)} \otimes \overline H_*(\C P^2) = [
	\overset{13}{\Q(7)},
	\overset{15}{\Q(8)},
	\overset{17}{\Q(9)}
	] .
	\]

	{\boldmath $ E^1_{2,*} = \overline H_{*+2}(\sigma_{2}\setminus \sigma_{1}) $.}

	Similarly, $ \sigma_{2}\setminus \sigma_{1} $ is a fibre bundle over the base $ X_2 = B(\C P^2, 2) $ with fibre $ \C^4 \times I $, where $ I \sim \R $ is an open interval.
	Interchanging the two points in a configuration from $ X_2 $ reverses the orientation of the open interval.
	This gives rise to a non-trivial local system on $ X_2 = B(\C P^2, 2) $:
	\[
	E^1_{2,*} = \left(\overset{8}{\Q(4)} \otimes \overset{1}{\Q(0)}\right) \otimes \overline H_*(B(\C P^2, 2); \pm \Q) = [
	\overset{9}{\Q(5)},
	\overset{11}{\Q(6)},
	\overset{13}{\Q(7)}
	] .
	\]

	{\boldmath $ E^1_{3,*} = \overline H_{*+3}(\sigma_{3}\setminus \sigma_{2}) $.}

	The space $ \sigma_{3}\setminus \sigma_{2} $ is a fibre bundle over the base $ X_3 = \C P^{2\vee} $ with fibre $ \C^3 \times \mathring{\Cone}((\C P^1)^{*2}) $, where the second self-join $ (\C P^1)^{*2} $ has the rational homology of a point (\cite{Vass99}).
	As a result, by \autoref{lem:tool_cone}, the Borel--Moore homology of the fibre is zero, and therefore $ E^1_{3,*} = 0 $ .

	{\boldmath $ E^1_{4,*} = \overline H_{*+4}(\sigma_{4}\setminus \sigma_{3}) $.}

	The space $ \sigma_{4}\setminus \sigma_{3} $ is a fibre bundle over the base $ X_4 $ with fibre $ \C^1 \times \mathring\Delta^2 $, where the open triangle $ \mathring\Delta^2 $ is spanned by the points $ \{ p, q, r \} \in X_4 $.
	Interchanging any two of these points results in a change of orientation on $ \mathring\Delta^2 $, which gives rise to a non-trivial local system on $ X_4 $.
	Similarly to \cite[Lemma 4]{Vass99}, the Borel--Moore homology \mbox{$ \overline H_*(X_4; \pm \Q) = \overline  H_*(B(\C P^2, 3); \pm \Q) = \overset{6}{\Q(3)} $} is concentrated in degree $6$.
	Taking into account the degree shift,
	\[
	E^1_{4,*} = \left(\overset{2}{\Q(1)} \otimes \overset{2}{\Q(0)}\right) \otimes \overset{2}{\Q(3)} = \overset{6}{\Q(4)} .
	\]

	{\boldmath $ E^1_{5,*} = \overline H_{*+5}(\sigma_{5}\setminus \sigma_{4}) $.}

	The space $ \sigma_{5}\setminus \sigma_{4} $ is a fibre bundle over $ X_5 = pt $ with fibre the open cone over the union
	\[
	\partial \Lambda (\C P^2) \coloneq \left(\bigcup\limits_{L \in X_3} \widetilde \Lambda(L)\right) \bigcup \left( \bigcup \limits_{K \in X_4} \widetilde \Lambda(K) \right) .
	\]
	If the space $ \partial \Lambda (\C P^2)$ has the rational homology of a point, then the entire fifth column is zero.

	To compute the homology of the union, consider the filtration
	\[
	\emptyset \subset f_0 \subset f_1 \subset f_2 = \partial \Lambda (\C P^2),
	\]
	where
	\[
	f_0 = \bigcup_{p\in\C P^2} \widetilde \Lambda(p) = \C P^2,\quad
	f_1 = \bigcup_{L \in X_3} \widetilde \Lambda(L),\quad
	f_2 = f_1 \cup \bigcup_{K \in X_4} \widetilde \Lambda(K).
	\]
	Let $ e^* $ be the associated Borel--Moore homological spectral sequence.
	Its first page $ e^1 $ is given in \autoref{table:intro_sp_seq2_p1}.
	We claim that all potentially non-zero differentials in $ e^1 $ are non-zero, hence $ E^1_{5,*} = 0 $.

	\begin{table}[!htbp]
		\centering
		\caption{$ e^{1} $}
		\label{table:intro_sp_seq2_p1}
	\renewcommand{\arraystretch}{1.2}
		\begin{tabular}{|r|*{3}{>{\centering\arraybackslash}m{0.7cm}|}}
			\hline
			6 &          & $\Q(3)$ &   $\Q(3)$  \\
			\hline
			5 &          &         &             \\
			\hline
			4 & $ \Q(2)$ & $\Q(2)$ &               \\
			\hline
			3 &          &         &            \\
			\hline
			2 & $\Q(1)$  & $\Q(1)$ &           \\
			\hline
			1 &          &         &             \\
			\hline
			0 & $\Q(0)$  &         &          \\
			\hline
			& 0        & 1       & 2         \\
			\hline
		\end{tabular}
	\end{table}

	Indeed, if some class in $ E^1_{5,*} = \widetilde H_{*+4}(\partial \Lambda (\C P^2)) $ was non-zero, then $ \overline H_n(\Sigma) \neq 0 $ for some $ n \leq 8 $.
	By Alexander duality, $ H^{19-n}(\UU) \neq 0 $ with $ 19-n \geq 11 $.
	This contradicts the fact that $ \UU $ is an affine variety of complex dimension $10$, so $ H^i(\UU)=0 $ for $ i \geq 11 $.
\end{proof}

If, instead, we choose to insert the space of collinear triples between $X_2$ and $X_3$, then essentially identical calculations yield the spectral sequence given in \autoref{table:intro_sp_seq1_p1_alt}, which converges to the same Borel--Moore homology.

\begin{table}[!htbp]
	\centering
	\caption{$ E^{1} = E^{\infty} $}
	\label{table:intro_sp_seq1_p1_alt}
	\renewcommand{\arraystretch}{1.2}
	\begin{tabular}{|r|*{6}{>{\centering\arraybackslash}m{0.7cm}|}}
		\hline
		17 & $ \Q(9)$ &         &   & &         &            \\
		\hline
		16 &          &         &   & &         &            \\
		\hline
		15 & $\Q(8)$  &         &   & &         &            \\
		\hline
		14 &          &         &   & &         &            \\
		\hline
		13 & $\Q(7)$  & $\Q(7)$ &   & &         &            \\
		\hline
		12 &          &         &   & &         &            \\
		\hline
		11 &          & $\Q(6)$ &   & &         &            \\
		\hline
		10 &          &         &   & &         &            \\
		\hline
		9  &          & $\Q(5)$ &   & &         &            \\
		\hline
		8  &          &         &   & &         &            \\
		\hline
		7  &          &         &   & &         &            \\
		\hline
		6  &          &         &   & &         &            \\
		\hline
		5  &          &         &   & & $\Q(4)$ &            \\
		\hline
		& 1        & 2       & 3 &4 & 5       & 6          \\
		\hline
	\end{tabular}
\end{table}

\section{Singularities relevant to GIT stability}
\label{sec:main_spse}

\subsection{Tacnodes of quartic curves}
\label{sec:tacnodes}

The goal of this work is to compute the rational cohomology groups of the space of homogeneous equations of GIT-stable quartic curves in $ \C P^2 $.
It follows from the results of Mumford \cite{mfk94} that a quartic curve in $ \C P^2 $ is not GIT-stable if and only if it has a triple point or a tacnode.

\begin{definition}
	For a curve $ Q = \V(F) \subset \C P^n $, a point $ p \in Q \subset \C P^n $ is an \emph{$ m $-point} or a \emph{point of multiplicity $ m $} if in an affine chart centred at $ p $, the homogeneous components $ F_k $ of degree $ k $ of $ F $ are zero for all $ k < m $.
\end{definition}

This definition is equivalent to saying that the multiplicity $ m_p(Q) $ of $ p $ on $ Q $ is at least $ m $.
It does not depend on an affine chart.
Moreover, it is invariant under changes of affine coordinates.
The linear factors of $ F_{m_p(Q)} $ correspond to tangent lines to $ Q $ at $ p $.

\begin{remark}
	If we additionally require that $ F_m \neq 0 $, we get the usual definition of multiplicity of a point on a curve.
	Why do we omit it in our definition?
	The unstable points form a Zariski-closed subset.
	Thus, if it contains all curves with a point of multiplicity $ m $ in the classical sense, then it also contains all curves with a point of multiplicity more than $ m $.
	Indeed, by removing the condition $ F_m \neq 0 $ we pass to the closure of the set of curves with a point of multiplicity $ m $.

	Let $ V $ be the set of polynomials with an $ m + 1 $-point $ p $, i.e. $ V = \{ f \ |\ f_0 = \ldots = f_m = 0 \} $.
	Let $ W $ be the set of polynomials with the point $ p $ being an $ m $-point, i.e. $ W = \{ f \ |\ f_0 = \ldots = f_{m-1} = 0 \} $.
	Then $ V \subset W $, as $ V  = \{ f \in W\ |\ f_m = 0 \} $.
	For any point $ f \in V $ there is a sequence $ g_n \in W \setminus V $ of polynomials such that $ \lim g_n = f $.
	Indeed, take $ (g_n)_i = f_i $ for $ i \neq m $ and choose homogeneous $ \lim_{n}(g_n)_m = 0 $.
	Thus, $ V \subset \Cl(W \setminus V) $.
	Hence if $ W $ are included in the set of GIT-stable polynomials, then its closure must also be contained in it, and therefore $ V $ is also there.
\end{remark}

\begin{definition}
	For a curve $ Q = \V(F) \subset \C P^n $ and a point $ p \in \C P^n $, a line $ l = \V(L) \subset \C P^n $ is a \emph{$ m $-fold tangent at $ p $} if in an affine chart centred at $ p $, the homogeneous components $ F_k $ of degree $ k $ of $ F $ are multiples of $ L $ for all $ k < m $.
\end{definition}

This definition can also be restated in terms of the local ring at $ p $.
In these terms, it means that the intersection number $ \mu_p(l, Q) = \dim \mathscr{O}_p / \langle L, F \rangle  $ of the line $ l $ and the quartic $ Q $ at $ p $ is at least $ m $.

\begin{definition}
	A \emph{cusp} is a double point of a curve with only one $3$-fold tangent line, i.e. a point whose projective tangent cone consists of a single point.
	A \emph{(classical) tacnode} is a cusp whose $3$-fold tangent is also a $4$-fold tangent.
\end{definition}

Let $ \Pi $ denote the vector space of homogeneous polynomials of degree four in three variables over the field of complex numbers, and let $ \Sigma \subset \Pi $ be the subset consisting of those equations whose corresponding quartic curves are not GIT-stable.
We are interested in the cohomology of $ \Pi \setminus \Sigma $.

To check whether a given equation corresponds to a quartic curve that is not GIT-stable, it is convenient to introduce a new type of singularity, generalising the triple point and the tacnode.
Furthermore, we will want this singularity to carry information about both the singular point of a quartic and a particular tangent line through it, i.e. it should be an element of the incidence variety
\begin{equation*}
	P \coloneq \PP T\C P^2 = \left \{ (p,l)\in \C P^2\times \C P^{2\vee}\ |\ p\in l \right \}.
\end{equation*}
In case of a triple point $ p $ of a quartic curve $ Q $, we want $ l $ to be a $4$-fold tangent to $ Q $ at $ p $, of which there may be one, two, three, or infinitely many.
In case of a (classical) tacnode $ p \in Q $, we want $ l $ to be the unique $3$-fold tangent that is also a $4$-fold tangent.

In fact, this is the full description of the singularity that we wish to define.
From now on, we call such a pair $ (p, l) \in P $ a \emph{tacnode} of a projective quartic curve $ Q $, and the line $ l $ is called the \emph{direction} of the tacnode.
This generalises the notion defined in the Introduction, and the two do not agree: a classical tacnode $ p $ gives exactly one tacnode $ (p,l) $ in the present sense, namely the one whose direction is its unique $3$-fold tangent, whereas a triple point is a tacnode in the present sense and is not a classical one.
Whenever the distinction matters we write \emph{classical tacnode} for the notion of the Introduction.

For a quartic $ Q = \V(F) \subset \C P^2 $ given by a homogeneous equation $ F \in \Pi $, denote the set of tacnodes of $ Q $ (in the above sense) by $ \Sing(Q) $ or $ \Sing(F) $.
It is a subset of $ P $.
\label{pg:16}

Now $ \Sigma \subset \Pi $ can be thought of as the space of polynomials $ F \in \Pi $ with $ \Sing(F) \neq \emptyset $.

We also need to prove Lemma 1 from \cite{BerGorKon23} for our case.
We omit the proof, but it is important to note that the lemma, among other things, asserts that our sets $ \Sing(F) $ behave similarly to how the usual sets of singular points of a curve behave, that $ \Sigma \subset \Pi $ is closed, and that for each $ K \subset P $, the set $ L(K) \coloneq \left \{ F \in \Pi \ |\ K \subset \Sing(F) \right \} $ is a vector space.
\label{pg:17}

\begin{theorem}\label{thm:tacnode_calc}
	Let $ Q $ be a quartic curve, $ (p, l) \in P $ be an element of the incidence variety, $ U $ be an affine chart $ U \subset \C P^2 $ containing $ p $ with fixed affine coordinates on it.
	Let the intersection $ Q \cap U $ be given by a polynomial equation $ f(x,y) = 0 $.
	If
	\begin{equation*}
		f(x,y) = f_0(x,y) + f_1(x,y) + \ldots + f_4(x,y)
	\end{equation*}
	is the Taylor expansion of $ f $ around $ p $ (or more formally the image of $ p $ in $ \mathbb{A_C}^2 $ under the ``choice of affine coordinates on $ U $'' map), then $ (p,l) \in \Sing(Q) $ (that is, $ (p,l) $ is a tacnode of $ Q $ in the generalised sense fixed above) if and only if the following conditions hold:
	\begin{enumerate}[(1)]
		\item $ f_0(x,y) = 0 $;
		\item $ f_1(x,y) = 0 $;
		\item $ f_2(x,y) = \lambda (ax + by + c)^2 $, where $ l \cap U $ is given by $ ax + by + c = 0 $, $ \lambda \in \C $;
		\item $ f_3(l \cap U) = 0 $.
	\end{enumerate}
	Here $ \lambda = 0 $ occurs exactly when $ p $ is a triple point of $ Q $, and $ \lambda \neq 0 $ exactly when $ p $ is a classical tacnode with direction $ l $.
	The conditions therefore do \emph{not} characterise classical tacnodes.
	Both cases must be admitted, since a quartic of either kind lies in $ \Sigma $ and so fails to be GIT-stable.
\end{theorem}

\begin{proof}

	The first condition simply asserts that $ p $ is a point of $ Q $.

	The second condition means that all partial derivatives of $ f $ are zero at $ p $.
	If the first condition holds, then by the Euler's theorem for homogeneous functions, the homogenisation of $ f $ also has all of its partial derivatives equal to zero at $ p $.

	If $ \lambda $ in the third condition is zero, then $ p $ is a triple point of $ Q $, and vice versa.
	The fourth condition then just asserts that $ l $ is a $4$-fold tangent.

	If $ \lambda $ is not zero, then $ l $ is the only $3$-fold tangent to $ Q $ at $ p $, i.e. $ p $ is a cusp.
	The fourth condition then just asserts that it is also a $4$-fold tangent, i.e. the point $ p\in Q $ is a tacnode in the classical sense.
\end{proof}

\begin{example}
	Let $ Q $ be given locally by $ 0 = (x^2 + y^2)^2 - y(y^2 - 3x^2) = x^4 + 2 x^2 y^2 + y^4 + 3 x^2 y - y^3 $.
	Then $ (0,0) $ is a triple point.
	The lines $ y = 0 $, $ y = \pm \sqrt{3}x $ are $4$-fold tangents.
\end{example}

\begin{example}
	$ 0 = (y^2 + x^2 - 3 y)^2 - 4 y^2 (2 - y) = x^4 + 2 x^2 y^2  + y^4 - 6 x^2 y - 2 y^3 + y^2 $ has only one tacnode, namely point $ (0,0) $ with direction $ y = 0 $.
\end{example}
Tacnodes are preserved by the action of $ \PGL_3(\C) $.
Therefore, the conditions can be checked in any affine chart.
It is often convenient to assume $ p = (0:0:1) \in \C P^2 $, $ l = \{ (x:0:z) \} \subset \C P^2 $, and $ U = \{ (x:y:1) \} $.

Let a quartic $ Q = \V(F) $ be given by the polynomial $ F $.
Let $ p $ and $ l $ be as above.
Conditions $(1)$ and $(2)$ together are equivalent to $ \dif F(p) = 0 $, and if it holds, then $(3)$ is equivalent to $ L \subset \ker (\text{Hess}_p F) $, where $ L \subset \C^3 $ is a $2$-plane such that  $ l = \PP(L) $.

Furthermore, it is sufficient to check the conditions in each of the charts of the standard affine cover of $ \C P^2 $ to find $ \Sing(Q) $ of a given quartic $ Q $.

\subsection{Configurations}

Let $ 2^P $ denote the space of \emph{non-empty compact subsets} of $ P $.
We equip $ 2^P $ with the Hausdorff metric constructed using a metric on $ P $ that induces the standard topology.
Elements of $ 2^P $ are called \emph{configurations}.

We apply a general method of calculating the rational cohomology groups of a certain class of hypersurfaces in $ \C P^n $.
To do so, we first define a sequence $ X_1, \ldots, X_N \subset 2^P $ of \emph{configuration spaces}, where $ N = 13 $.
This is most conveniently done by describing the general form of a configuration in each $ X_i $.

The configuration spaces that we chose are listed in \autoref{table:configs_def}.

\begin{table}
	\centering
	\caption{Configurations}
	\label{table:configs_def}
	\resizebox{\textwidth}{!}{
		\renewcommand{\arraystretch}{1.5}
		\begin{tabular}{|c|c|c|}
			\hline
			$ i $ & A generic configuration $K\in X_i$                                                                                   & The picture of $K \in X_i$ \\
			\hline

			1     & $ \{ (p, l) \} $ a singleton in $ P $.                                                                               & $ \configlonetac $         \\
			\hline

			2     & $ \{ (p,l_1), (p,l_2) \} $, where $ l_1 \neq l_2 $, an unordered pair of elements of $ P $.                          & $ \configsmallcross $      \\
			\hline

			3     & $ \{ (p, l_1), (p, l_2), (p, l_3) \} $, where $ l_i \neq l_j $ for $ i \neq j $, an unordered triple.                & $ \configtriple $          \\
			\hline

			4     & $ \{ (p,l) \in P \mid l\text{ is a fixed line}\} = l\times\{l \}$                                                    & $ \configline $            \\
			\hline

			5     & $ \{ (p,l) \in P \mid p\text{ is a fixed point}\} = \{ p \} \times \PP(T_p\C P^2) $                                  & $ \configcirc $            \\
			\hline

			6     & $ (l_1\times \{ l_1 \}) \cup \{(p, l_2) \} $, where $ p = l_1 \cap l_2 $.                                            & $ \configlinetac $         \\
			\hline

			7     & $ (\{ p \}\times \PP(T_p\C P^2)) \cup ( l\times \{ l \}) $, where $ p\in l $.                                        & $ \configlinecirc $        \\
			\hline

			8     & $ \{(p, l_p), (q,l_q) \} $, where $ p, q $ and $ l_p \cap l_q $ are 3 distinct points.                               & $ \configparallelpair $    \\
			\hline

			9     & $ ( l\times \{ l \}) \cup \{(p,l_p),  (q,l_q)\} $, where $p, q \in l $ and $ l_p \cap l_q $ are 3 distinct points.   & $ \configlinepair $        \\
			\hline

			10    & $ (l_1\times \{ l_1 \}) \cup (l_2\times \{ l_2 \}) \cup (\{p\} \times \PP(T_p\C P^2)) $, where $ p = l_1 \cap l_2 $. & $ \configcross{} $         \\
			\hline

			11    & $ l \times \C P^{2\vee} $                                                                                            & $ \configfatline $         \\
			\hline

			12    & $ \PP(TC) $, where $ C $ is a fixed smooth conic.                                                                    & $ \configconic $           \\
			\hline

			13    & $ P = \PP T\C P^2 $                                                                                                  &                            \\
			\hline
		\end{tabular}
	}
\end{table}

\section{Verification of the conditions required for the method}
\label{sec:v1_conditions}

The following propositions verify that the conditions in \cite[Condition list 1]{BerGorKon23} are satisfied, allowing us to apply the method.

\begin{proposition}[Condition 1]
	For every polynomial $ F \in \Sigma $, the set $ \Sing(F) $ is contained in $ X_i $ for some $ i $.
\end{proposition}

\begin{proof}

	We analyse the possible decompositions of a quartic curve $ Q $, starting with the most degenerate cases and proceeding towards a smooth quartic.
	In each case we determine the corresponding configuration of tacnodes by analysing the curve in affine charts.

	\textbf{Case 1. Suppose $ Q = \V(F) $ is a union of four distinct lines.}

	\begin{enumerate}[leftmargin=*, label=(\alph*)]

		\item If the lines intersect in six distinct points, then each intersection point is a double points with two distinct $3$-fold tangent lines (the components passing through it).
		Hence there are no tacnodes.

		\begin{minipage}{0.5\textwidth}
			\centering
			\begin{tikzpicture}[scale=0.66]
				\draw (-1,0) -- (2,0);
				\draw (-1,1.5) -- (2.5,-1);
				\draw (-1,2.5) -- (1.5,-1);
				\draw (0,-1) -- (0,2);
			\end{tikzpicture}
		\end{minipage}
		\begin{minipage}{0.5\textwidth}
			\centering
			$ \emptyset $
		\end{minipage}

		\item If the lines intersect in four distinct points, then exactly one of the intersection points is a triple point.
		It has three $4$-fold tangents -- the components of $ Q $ passing through it.
		Thus $ \Sing(Q) \in X_3$.

		\begin{minipage}{0.5\textwidth}
			\centering
			\begin{tikzpicture}[scale=0.66]
				\draw[fill] (0,0) circle (1.5pt) node[above left]{$p$};
				\draw (-1,0) -- (2,0);
				\draw (-0.5,-0.5) -- (1.5,1.5);
				\draw (-1,2) -- (2,-1);
				\draw (0,-1) -- (0,2);
			\end{tikzpicture}
		\end{minipage}
		\begin{minipage}{0.5\textwidth}
			\centering
			\begin{tikzpicture}[scale=0.66]
				\draw[fill] (0,0) circle (1.5pt) node[above left]{$p$};
				\draw (-0.5,0) -- (0.5,0);
				\draw (-0.707/2,-0.707/2) -- (0.707/2,0.707/2);
				\draw (0,-0.5) -- (0,0.5);
			\end{tikzpicture}
		\end{minipage}

		\item If all four lines intersect in the same point, then this point is a triple point, and each line through it is a direction.
		Thus $ \Sing(Q) \in X_5 $.

		\begin{minipage}{0.5\textwidth}
			\centering
			\begin{tikzpicture}[scale=0.66]
				\draw[fill] (0,0) circle (1.5pt);
				\draw (-1.5,0) -- (1.5,0);
				\draw (-1,-1) -- (1,1);
				\draw (-1,1) -- (1,-1);
				\draw (0,-1) -- (0,1);
			\end{tikzpicture}
		\end{minipage}
		\begin{minipage}{0.5\textwidth}
			\centering
			\begin{tikzpicture}
				\draw (0,0) circle (2.5pt);
			\end{tikzpicture}
		\end{minipage}

	\end{enumerate}

	\textbf{Case 2. Suppose $ Q $ is a union of three distinct lines.}

	\begin{enumerate}[leftmargin=*, label=(\alph*)]
		\item If the three lines intersect in three distinct points, then the two intersection points lying on the double line are triple points of $ Q $, and the two components containing them are $4$-fold tangents.
		Furthermore, any other point on the double line has the unique $3$-fold tangent -- namely the line itself, which is also a $4$-fold tangent.
		Thus $ \Sing(Q) \in X_{9}$.

		\begin{minipage}{0.5\textwidth}
			\centering
			\begin{tikzpicture}[scale=0.66]
				\draw[fill] (-1,0) circle (1.5pt) node[above left]{$p$};
				\draw[fill] (1,0) circle (1.5pt) node[above right]{$q$};
				\draw (-2,0) -- (2,0);
				\draw (-1,-1) -- (-1,1);
				\draw (1,-1) -- (1,1);
			\end{tikzpicture}
		\end{minipage}
		\begin{minipage}{0.5\textwidth}
			\centering
			\begin{tikzpicture}[scale=0.66]
				\draw[fill] (-1,0) circle (1.5pt) node[above left]{$p$};
				\draw[fill] (1,0) circle (1.5pt) node[above right]{$q$};
				\draw (-2,0) -- (2,0);
				\draw (-1,-0.5) -- (-1,0.5);
				\draw (1,-0.5) -- (1,0.5);
			\end{tikzpicture}
		\end{minipage}

		\item If the lines intersect in a single point, then this point is a triple point and also a $4$-point, so every line through it is a direction.
		By the same argument as in the previous case, the double line itself is a direction for any point on it.
		Thus $ \Sing(Q) \in X_7$.

		\begin{minipage}{0.5\textwidth}
			\centering
			\begin{tikzpicture}[scale=0.66]
				\draw[fill] (0,0) circle (1.5pt);
				\draw (-1.5,0) -- (1.5,0);
				\draw (-0.707,-0.707) -- (0.707,0.707);
				\draw (-0.707,0.707) -- (0.707,-0.707);
			\end{tikzpicture}
		\end{minipage}
		\begin{minipage}{0.5\textwidth}
			\centering
			\begin{tikzpicture}[scale=0.66]
				\draw (-1.5,0) -- (1.5,0);
				\draw (0,0) circle (2.5pt);
			\end{tikzpicture}
		\end{minipage}
	\end{enumerate}

	\textbf{Case 3. Suppose $ Q $ is a union of two lines, one with multiplicity $3$.}

	Each point of the triple line is a triple point, and has a unique $4$-fold tangent -- the line itself.
	The intersection point is a $4$-point, and therefore every line through it is a direction.
	Thus $ \Sing(Q) \in X_7$.

	\begin{minipage}{0.5\textwidth}
		\centering
		\begin{tikzpicture}[scale=0.66]
			\draw[fill] (0,0) circle (1.5pt) node[above left]{$q$};
			\draw (-1.5,0) -- (1.5,0);
			\draw (0,-1) -- (0,1);
		\end{tikzpicture}
	\end{minipage}
	\begin{minipage}{0.5\textwidth}
		\centering
		\begin{tikzpicture}[scale=0.66]
			\draw (-1.5,0) -- (1.5,0);
			\draw (0,0) circle (2.5pt) node[above left]{$q$};
		\end{tikzpicture}
	\end{minipage}

	\textbf{Case 4. Suppose $ Q $ is a union of two double lines.}

	By the same argument as above, each point on one of the lines has that line as a direction, and for the intersection point, every line through it is a direction.
	Thus $ \Sing(Q) \in X_{10}$.

	\begin{minipage}{0.5\textwidth}
		\centering
		\begin{tikzpicture}[scale=0.66]
			\draw[fill] (0,0) circle (1.5pt) node[above left]{$q$};
			\draw (-1.5,0) -- (1.5,0);
			\draw (0,-1) -- (0,1);
		\end{tikzpicture}
	\end{minipage}
	\begin{minipage}{0.5\textwidth}
		\centering
		\begin{tikzpicture}[scale=0.66]
			\draw (-1.5,0) -- (1.5,0);
			\draw (0,-1) -- (0,1);
			\draw (0,0) circle (2.5pt) node[above left]{$q$};
		\end{tikzpicture}
	\end{minipage}

	\textbf{Case 5. Suppose $ Q $ is a quadruple line.}

	In this case each point is a $4$-point, so every line through the point is a direction.
	Thus $ \Sing(Q) \in X_{11}$.

	\begin{minipage}{0.5\textwidth}
		\centering
		\begin{tikzpicture}[scale=0.66]
			\draw (-1.5,0) -- (1.5,0);
		\end{tikzpicture}
	\end{minipage}\hfill
	\begin{minipage}{0.5\textwidth}
		\centering
		$\configfatline[1.0]$
	\end{minipage}

	This covers all quartic curves that split into four lines.

	\textbf{Case 6. Suppose $ Q = C \cup l \cup l' $ is a union of a smooth conic and two lines.}

	\begin{enumerate}[leftmargin=*, label=(\alph*)]
		\item If the two lines intersect the conic in four distinct points, then the same argument as in Case 1(a) applies: the two $3$-fold tangents at a singular point $ p $ are the linear component of $ Q $ containing $ p $ and the tangent line to the conic at $ p $.

		\begin{minipage}{0.5\textwidth}
			\centering
			\begin{tikzpicture}[scale=0.66]
				\draw (-0.5,1) -- (2.5,1);
				\draw (1,-0.5) -- (1,2.5);
				\draw (1,1) circle (1cm);
			\end{tikzpicture}
		\end{minipage}\hfill
		\begin{minipage}{0.5\textwidth}
			\centering
			$ \emptyset $
		\end{minipage}

		\item If the lines intersect the conic in exactly three points, none of which is the intersection point of the two lines, then one of the lines must be tangent to the conic.
		By the same argument as above, only the point $ p $ of tangency is relevant.
		Obviously, the point $ p $ is not a triple point.
		However, the tangent line $ l $ to the conic at $ p $ has infinite intersection multiplicity, and is therefore a $3$-fold tangent.
		Furthermore, any other line through $ p $ is tangent only to the component $ l $ of $ Q $ and not to the conic $ C $.
		As the intersection multiplicity of a line with a union of curves is the sum of intersection multiplicities, $ l $ is the only $3$-fold tangent to $ Q $ at $ p $.
		It is also a $4$-fold tangent, so $ (p, l) $ is a tacnode of $ Q $.
		Thus $ \Sing(Q) \in X_1$.

		\begin{minipage}{0.5\textwidth}
			\centering
			\begin{tikzpicture}[scale=0.66]
				\draw[fill] (1,0) circle (1.5pt) node[below]{$p$};
				\draw (-0.5,0) -- (2.5,0);
				\draw (-0.5,-0.5) -- (2,2);
				\draw (1,1) circle (1cm);
			\end{tikzpicture}
		\end{minipage}\hfill
		\begin{minipage}{0.5\textwidth}
			\centering
			\begin{tikzpicture}[scale=0.66]
				\draw[fill] (0,0) circle (1.5pt) node[below]{$p$};
				\draw (-0.5,0) -- (0.5,0);
			\end{tikzpicture}
		\end{minipage}

		\item If the lines intersect the conic in exactly three points, one of which is the intersection point $ p $ of the two lines, then the intersection point is a triple point.
		Then the two lines $ l $ and $ l' $ along with the tangent $ T_pC $ are the directions of the tacnodes of $ Q $, which are the only tacnodes by the same argument as above.
		Thus $ \Sing(Q) \in X_3$.

		\begin{minipage}{0.5\textwidth}
			\centering
			\begin{tikzpicture}[scale=0.66]
				\draw[fill] (0,0) circle (1.5pt) node[below]{$p$};
				\draw (-0.5,-0.5) -- (1.5,1.5);
				\draw (0.5,-0.5) -- (-1.5,1.5);
				\draw (0,1) circle (1cm);
			\end{tikzpicture}
		\end{minipage}\hfill
		\begin{minipage}{0.5\textwidth}
			\centering
			\begin{tikzpicture}[scale=0.66]
				\draw[fill] (0,0) circle (1.5pt) node[below]{$p$};
				\draw (-0.707/2,-0.707/2) -- (0.707/2,0.707/2);
				\draw (-0.5,0) -- (0.5,0);
				\draw (0.707/2,-0.707/2) -- (-0.707/2,0.707/2);
			\end{tikzpicture}
		\end{minipage}

		\item If the lines intersect the conic in exactly two points, none of which is the intersection point of the two lines, then they must be tangent to the conic at two distinct points $ p, q \in C $.
		By the same argument as in Case 6(b), $ \Sing(Q) = \{ (p, l), (q,l') \} $.
		Thus $ \Sing(Q) \in X_{8} $.

		\begin{minipage}{0.5\textwidth}
			\centering
			\begin{tikzpicture}[scale=0.66]
				\draw[fill] (0,1) circle (1.5pt) node[left]{$p$};
				\draw[fill] (1,0) circle (1.5pt) node[below]{$q$};
				\draw (-0.5,0) -- (2,0);
				\draw (0,-0.5) -- (0,2);
				\draw (1,1) circle (1cm);
			\end{tikzpicture}
		\end{minipage}\hfill
		\begin{minipage}{0.5\textwidth}
			\centering
			\begin{tikzpicture}[scale=0.66]
				\draw[fill] (0,1) circle (1.5pt) node[left]{$p$};
				\draw[fill] (1,0) circle (1.5pt) node[below]{$q$};
				\draw (0,0.5) -- (0,1.5);
				\draw (0.5,0) -- (1.5,0);
			\end{tikzpicture}
		\end{minipage}

		\item If the lines intersect the conic in exactly two points, one of which is the intersection point $p$ of the two lines, then the intersection point is a triple point.
		Then the two lines $ l $ and $l' = T_pC $ are the directions of the tacnodes at point $ p \in Q $.
		These two tacnodes are the only tacnodes of $ Q $ by the same argument as above.
		Thus $ \Sing(Q) \in X_{2} $.

		\begin{minipage}{0.5\textwidth}
			\centering
			\begin{tikzpicture}[scale=0.66]
				\draw[fill] (1,0) circle (1.5pt) node[below left]{$p$};
				\draw (-0.5,0) -- (2.5,0);
				\draw (1,-0.5) -- (1,2.5);
				\draw (1,1) circle (1cm);
			\end{tikzpicture}
		\end{minipage}\hfill
		\begin{minipage}{0.5\textwidth}
			\centering
			\begin{tikzpicture}[scale=0.66]
				\draw[fill] (0,0) circle (1.5pt) node[below left]{$p$};
				\draw (-0.5,0) -- (0.5,0);
				\draw (0,-0.5) -- (0,0.5);
			\end{tikzpicture}
		\end{minipage}

		\item If the two lines coincide and intersect the conic in two distinct points, then the double line is a direction at every point on it (as in Case 4).
		Moreover, for each of the two intersection points, the tangent line to the conic at that point is a direction (as in Case 6(d)).
		Thus $ \Sing(Q) \in X_{9}$.

		\begin{minipage}{0.5\textwidth}
			\centering
			\begin{tikzpicture}[scale=0.66]
				\draw[fill] (0,1) circle (1.5pt) node[above left]{$p$};
				\draw[fill] (2,1) circle (1.5pt) node[above right]{$q$};
				\draw (-0.5,1) -- (2.5,1);
				\draw (1,1) circle (1cm);
			\end{tikzpicture}
		\end{minipage}\hfill
		\begin{minipage}{0.5\textwidth}
			\centering
			\begin{tikzpicture}[scale=0.66]
				\draw[fill] (-1,0) circle (1.5pt) node[above left]{$p$};
				\draw[fill] (1,0) circle (1.5pt) node[above right]{$q$};
				\draw (-2,0) -- (2,0);
				\draw (-1,-0.5) -- (-1,0.5);
				\draw (1,-0.5) -- (1,0.5);
			\end{tikzpicture}
		\end{minipage}

		\item If the two lines coincide and intersect the conic in exactly one point, then the line is tangent to $ C $ at that point.
		As in Case 4, the double line is a direction at every point on it.
		The point of tangency is a triple point, but it only has one $4$-fold tangent, namely the line component of $ Q $, and therefore the description of $ \Sing(Q) $ is complete.
		Thus $ \Sing(Q) \in X_4 $.

		\begin{minipage}{0.5\textwidth}
			\centering
			\begin{tikzpicture}[scale=0.66]
				\draw[fill] (1,0) circle (1.5pt) node[below]{$p$};
				\draw (-0.5,0) -- (2.5,0);
				\draw (1,1) circle (1cm);
			\end{tikzpicture}
		\end{minipage}\hfill
		\begin{minipage}{0.5\textwidth}
			\centering
			\begin{tikzpicture}[scale=0.66]
				\draw (-2,0) -- (2,0);
			\end{tikzpicture}
		\end{minipage}

	\end{enumerate}

	\textbf{Case 7. Suppose $ Q = C \cup C' $ is a union of two smooth conics.}

	Note that two conics can intersect tangentially in at most two points, since for a fixed conic and a point on it, tangential intersection at that point imposes two linear conditions on a variable conic.

	\begin{enumerate}[leftmargin=*, label=(\alph*)]

		\item By the local nature of the definition of a tacnode, we can view the points of $ Q $ near the intersection points as a cross -- the intersection of two lines.
		This reduces the argument to Case 1(a).

		An alternative argument can be made.
		If the two conics intersect transversely, then there are no triple points, and at each of the four intersection points there are two $3$-fold tangents: the tangents to each conic.

		\begin{minipage}{0.5\textwidth}
			\centering
			\begin{tikzpicture}[scale=0.66]
				\draw (0,0) ellipse (1.5cm and 1cm);
				\draw (0,0) ellipse (1cm and 1.5cm);
			\end{tikzpicture}
		\end{minipage}
		\begin{minipage}{0.5\textwidth}
			\centering
			$ \emptyset $
		\end{minipage}

		\item If the two conics intersect tangentially in exactly one point $p$, then the set of singularities is \mbox{$\Sing(Q)=\left \{ (p,l=T_pQ) \right \} \in X_1 $}.
		Indeed, the tangents $ T_pC = l = T_pC'$ coincide.
		By additivity of intersection multiplicity for a union of curves, $ l $ is a $4$-fold tangent to $ Q $.
		It is the unique $4$-fold tangent at $ p $.
		Furthermore, every $3$-fold tangent to $ Q $ at $ p $ must be a $4$-fold tangent.
		The other intersection points are as in Case 7(a).

		\begin{minipage}{0.5\textwidth}
			\centering
			\begin{tikzpicture}[scale=0.66]
				\draw[fill] (-1,0) circle (1.5pt) node[left]{$p$};
				\draw (0.5,0) ellipse (1.5cm and 1cm);
				\draw (0,0) ellipse (1cm and 1.5cm);
			\end{tikzpicture}
		\end{minipage}
		\begin{minipage}{0.5\textwidth}
			\centering
			\begin{tikzpicture}[scale=0.66]
				\draw[fill] (0.5,0) circle (1.5pt) node[left]{$p$};
				\draw (0.5,-0.5) -- (0.5,0.5);
			\end{tikzpicture}
		\end{minipage}

		\item If the two conics intersect tangentially in two distinct points $p$ and $q$, then by applying the same argument as in Case 7(b), \mbox{$\Sing(Q)=\left \{ (p,T_pQ),\ (q,T_qQ) \right \} \in X_{8} $}.

		\begin{minipage}{0.5\textwidth}
			\centering
			\begin{tikzpicture}[scale=0.66]
				\draw[fill] (-1,0) circle (1.5pt) node[left]{$p$};
				\draw[fill] (1,0) circle (1.5pt) node[right]{$q$};
				\draw (0,0) ellipse (1cm and 1cm);
				\draw (0,0) ellipse (1cm and 1.5cm);
			\end{tikzpicture}
		\end{minipage}
		\begin{minipage}{0.5\textwidth}
			\centering
			\begin{tikzpicture}[scale=0.66]
				\draw[fill] (-1,0) circle (1.5pt) node[left]{$p$};
				\draw[fill] (1,0) circle (1.5pt) node[right]{$q$};
				\draw (-1,-0.5) -- (-1,0.5);
				\draw (1,-0.5) -- (1,0.5);
			\end{tikzpicture}
		\end{minipage}

		\item If the two conics intersect tangentially in at least three distinct points, then they coincide.
		By applying the same argument as in Case 7(b) to each point, we obtain the set of singularities \mbox{$\Sing(Q)=\left \{ (p,T_pQ)\ |\ p\in Q \right \} \in X_{12} $}.

		\begin{minipage}{0.5\textwidth}
			\centering
			\begin{tikzpicture}[scale=0.66]
				\draw (0,0) circle (1cm);
			\end{tikzpicture}
		\end{minipage}
		\begin{minipage}{0.5\textwidth}
			\centering
			\begin{tikzpicture}[scale=0.66]
				\draw (0,0) circle (1cm);
			\end{tikzpicture}
		\end{minipage}

	\end{enumerate}

	\textbf{Case 8. Suppose $ Q = C \cup l$ is a union of an irreducible cubic curve and a line.}

	The singularities (in the usual sense) of such a quartic $Q$ are precisely the intersection points $l\cap C$ and the singularities of the cubic $C$.

	An irreducible cubic curve has at most one singularity.
	It must be either a node $ y^2z=x^2(x+z) $ or a cusp $ y^2z=x^3 $.

	A smooth irreducible projective cubic curve is equivalent to one of the form $ y^2z = x(x-z)(x-\lambda z)$ (a Legendre cubic) for some $ \lambda \in \C \setminus \{ 0, 1 \} $.

	\begin{enumerate}[leftmargin=*, label=(\alph*)]
		\item If the line $ l $ intersects $ C $ in three distinct points, then the same argument as in Case 6(a) applies.

		\item If $l$ is a tangent to $C$ at a non-singular point $p\in C$ intersecting $C$ in exactly two points, then $\Sing(Q)=\left \{ (p,l) \right \}$.
		Indeed, at the point of tangency the line $ l = T_pC $ is the unique $3$-fold tangent that is a $4$-fold tangent to $ Q $ by the same argument as above.
		The remaining intersection point is a smooth point of $ C $ by counting multiplicity, and the same argument as in Case 8(a) applies.

		\item Suppose $ l $ is a tangent to $ C $ intersecting $ C $ in exactly one non-singular point $ p $.
		In this case the only possible tacnodes are of the form $ (p, l') $, where $  l' \ni p $ are various lines through $ p $, as $ p $ is the only singular point of $ Q $.
		Thus $ \Sing(Q) $ belongs to one of the following: $ X_1, X_2, X_3, X_5 $.

		\item Suppose $ l $ is a tangent to the nodal cubic $ C $ at the singular point $ p $ intersecting $ C $ in exactly one point: \mbox{$ (x+y)(y^2 - x^2(x + 1)) = -x^4-x^3y-(x-y)(x+y)^2 = 0 $}.
		Then by \autoref{thm:tacnode_calc}, $ \Sing(Q) = \{ (p,l), (p,l')\} $, where $ l'= \{x-y=0\}$ is the other tangent to $ C $ at $ p $ intersecting $ C $ in exactly one point.

		\begin{minipage}{0.5\textwidth}
			\centering
			\includegraphics{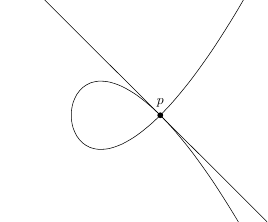}
		\end{minipage}
		\begin{minipage}{0.5\textwidth}
			\centering
			\begin{tikzpicture}[scale=0.66]
				\draw[fill] (0,0) circle (1.5pt) node[above]{$p$};
				\draw (-0.707/2,0.707/2) -- (0.707/2,-0.707/2);
				\draw (-0.707/2,-0.707/2) -- (0.707/2,0.707/2);
			\end{tikzpicture}
		\end{minipage}

		\item Suppose $ l $ is a tangent to the nodal cubic $ C $ at the singular point $ p $ intersecting $ C $ in exactly two points: $ y(y^2 - x^2(x + 1)) = x^3y+y(x-y)(x+y) = 0 $.
		Then by \autoref{thm:tacnode_calc}, \mbox{$ \Sing(Q) = \{ (p,l), (p,l'), (p, l'') \} $}, where $ l'= \{x-y=0\}$ and $ l''=\{x+y=0\} $.

		\begin{minipage}{0.5\textwidth}
			\centering
			\includegraphics{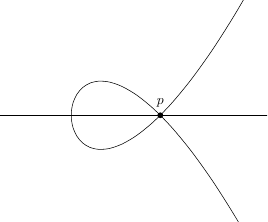}
		\end{minipage}
		\begin{minipage}{0.5\textwidth}
			\centering
			\begin{tikzpicture}[scale=0.66]
				\draw[fill] (0,0) circle (1.5pt) node[above]{$p$};
				\draw (-0.707/2,0.707/2) -- (0.707/2,-0.707/2);
				\draw (-0.707/2,-0.707/2) -- (0.707/2,0.707/2);
				\draw (-0.5,0) -- (0.5,0);
			\end{tikzpicture}
		\end{minipage}

		\item Suppose $ l $ is the $3$-fold tangent to the cuspidal cubic.
		Then in the local chart $ z = 1 $ the quartic $ Q $ is given by $ y(x^3-y^2) = x^3y-y^3 = 0 $ .
		By considering multiplicities, every singular point of $ Q $ lies in this chart.
		By \autoref{thm:tacnode_calc}, we obtain $ \Sing(Q) = \{ (p, l) \} $.

		\begin{minipage}{0.5\textwidth}
			\centering
			\includegraphics{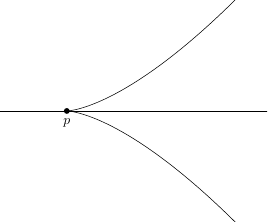}
		\end{minipage}
		\begin{minipage}{0.5\textwidth}
			\centering
			\begin{tikzpicture}[scale=0.66]
				\draw[fill] (0,0) circle (1.5pt) node[below]{$p$};
				\draw (-0.5,0) -- (0.5,0);
			\end{tikzpicture}
		\end{minipage}

		\item Suppose $ l $ is a tangent to the cuspidal cubic $ C $ at the singular point $ p $, intersecting $ C $ in exactly two points: $ (x-y)(x^3-y^2) = x^4-x^3y-y^2(x-y) = 0 $.
		By the same argument as in the previous case, $ \Sing(Q) = \{ (p, l), (p, l') \}$, where $ l'= \{y=0\} $.

		\begin{minipage}{0.5\textwidth}
			\centering
			\includegraphics{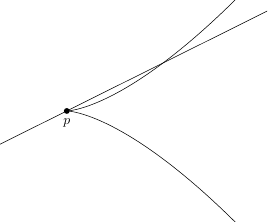}
		\end{minipage}
		\begin{minipage}{0.5\textwidth}
			\centering
			\begin{tikzpicture}[scale=0.66]
				\draw[fill] (0,0) circle (1.5pt) node[below]{$p$};
				\draw (-0.5,0) -- (0.5,0);
				\draw (-0.857493/2,-0.5144958/2) -- (0.857493/2,0.5144958/2);
			\end{tikzpicture}
		\end{minipage}
	\end{enumerate}

	\textbf{Case 9. Suppose that $ Q $ is an irreducible quartic curve.}

	In the following, we are using the classification of singularities of irreducible quartic curves by Namba et al \cite{namba1984geometry}.

	Consider the genus $ g = 3 - \sum \delta_p $ of $ Q $.

	\begin{enumerate}[leftmargin=*, label=(\alph*)]
		\item If $ g = 3 $, then there are no singular points.
		\item If $ g = 2 $, then there is a single singular point.
		It is either a node or a cusp, in which case $ \Sing(Q) = \emptyset $, as in the previous cases.
		\item Let $ g = 1 $.
		If $ Q $ has two distinct singular points, then, as above, they do not satisfy our definition.
		If $ Q $ has one singular point, then, by Namba's classification, it is either a tacnode (in the classical sense) or a double cusp.
		Observe that double cusp differs from a simple cusp in that the unique $3$-fold tangent is also a $4$-fold tangent.
		Therefore, $ \Sing(Q) $ lies in $ X_1 $ in both cases.
		\item Let $ g = 0 $, that is, $ \sum \delta_p = 3 $.
		Suppose first that $ Q $ has a triple point $ p $.
		A triple point has $ \delta_p \geq 3 $, so in this case $ p $ is the only singular point of $ Q $ and $ \delta_p = 3 $.
		Thus, $ \Sing(Q) $ lies in $ X_3 $, if $ p $ is an ordinary triple point (consider trefoil curve).
		$ \Sing(Q) $ is in $ X_{2} $ if $ Q = C_1 \cup C_2 $, where $ p $ is a simple cusp of multiplicity $2$ of one branch $ C_1 $ intersecting another, smooth branch $ C_2 $, and $ C_i $ intersect transversely.
		To see this, consider a local chart around the singular point, where the quartic looks like a transverse intersection of a cuspidal cubic with a line.
		By Case 8(g), $ \Sing(Q) $ is in $ X_{2} $.
		Finally, $ \Sing(Q) $ is in $ X_1 $, if $ p $ is a simple cusp with multiplicity $3$ (consider $ x^4 = y^3 $).

		The only other case left is when every singular point of $ Q $ is a double point: a point of multiplicity $ 4 $ would force $ Q $ to be a union of four lines through it, hence reducible.
		This case includes that of a single singular point, since a ramphoidal cusp and an osnode are double points with $ \delta_p = 3 $.
		$ \Sing(Q) $ is non-empty only if it contains a double cusp, ramphoidal cusp, tacnode or osnode.
		In each of these cases, $ \Sing(Q) $ is in $ X_1 $.
		In case of a ramphoidal cusp, we can use \autoref{thm:tacnode_calc} to check it algebraically.
		As being a singularity of a certain type is a local property, it suffices to check for any such quartic, e.g. the local model $ (y - x^2)^2 = xy^3 $.
		For an osnode $ p \in Q $, as it is the only singularity of the quartic that has it, $ \Sing(Q) $ lies in one of $ X_i $, where $ i = 1, 2, 3, 5 $.
		Note that if $ p $ has more than three directions through it, then $ \Sing(Q) $ is of type $ X_5 $.
	\end{enumerate}
\end{proof}

\begin{proposition}[Condition 2]
	If $ K\in X_i, L \in X_j $ and $K \subsetneq L$, then $i < j $.
\end{proposition}

\begin{proof}
	The arrows in the picture below together with their compositions are all possible inclusions of configurations considered as subsets of $ P = \PP T\C P^2 $.
	Notice the indices of the corresponding configuration spaces.

	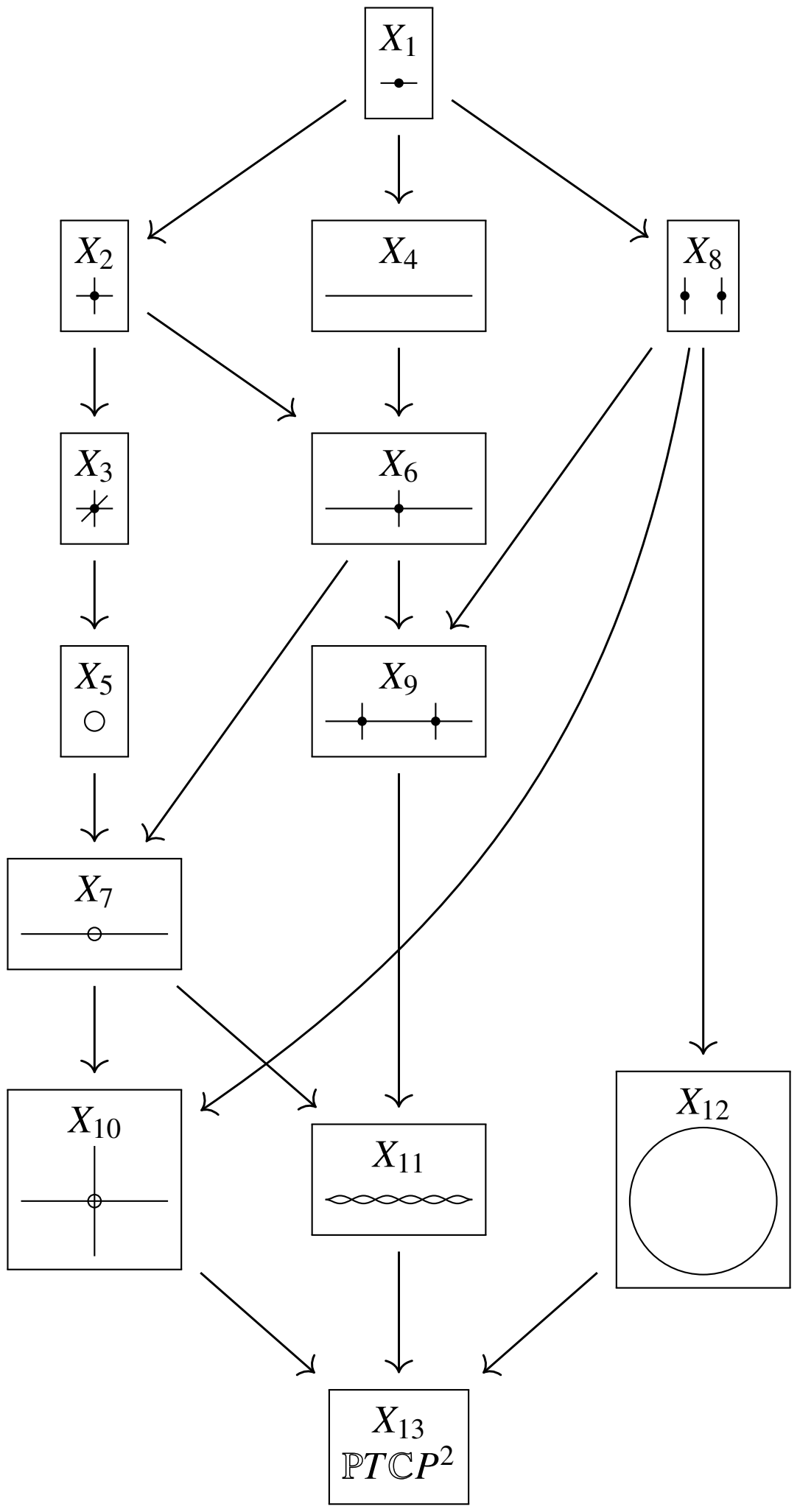
\begin{figure}[htbp]
		\centering
		\begin{tikzcd}
			& \mlnode{$ X_1 $\\$ \configlonetac $} \arrow[ld] \arrow[rd] \arrow[d] &                                \\
			\mlnode{$X_{2}$\\$ \configsmallcross $} \arrow[d] \arrow[rd] & \mlnode{$X_4$\\$ \configline $} \arrow[d]                       & \mlnode{$X_{8}$\\$ \configparallelpair $} \arrow[ldd] \arrow[dddd] \arrow[lldddd, bend left=23] \\
			\mlnode{$X_3$\\$ \configtriple $} \arrow[d]            & \mlnode{$ X_6 $\\$ \configlinetac $} \arrow[d] \arrow[ldd]           &                                \\
			\mlnode{$X_5 $\\$ \configcirc $} \arrow[d]            & \mlnode{$X_{9}$\\$ \configlinepair $} \arrow[dd]                    &                                \\
			\mlnode{$X_7 $\\$ \configlinecirc $} \arrow[d] \arrow[rd] &                                     &                                \\
			\mlnode{$X_{10}$\\$ \configcross{} $} \arrow[rd]        & \mlnode{$X_{11}$\\$ \configfatline $} \arrow[d]                    & \mlnode{$X_{12}$\\$ \configconic $} \arrow[ld]              \\
			& \mlnode{$X_{13}$\\$ \PP T\C P^2 $}                                &
		\end{tikzcd}
		\caption{Relations between configurations}
		\label{fig:config_relations}
	\end{figure}
\end{proof}

\begin{proposition}[Condition 3]
	For any fixed index $ i $, the dimension $ d_i \coloneq \dim_{\C} (L(K)) $ of the complex vector space \mbox{$ L(K) = \left \{ F \in \Pi \ |\ K \subset \Sing(F) \right \} $} is the same for all configurations $ K \in X_i $.
\end{proposition}

\begin{proof}
	\autoref{thm:tacnode_calc} asserts that the conditions for a given pair to be a tacnode are linear with respect to coefficients of the general homogeneous polynomial of degree 4.
	Thus, the polynomials with singularities given by $ K \in X_i $ form a linear subspace in $ \Pi $.

	To see that the dimension of this subspace does not depend on the choice of $ K \in X_i $, we use the action of $ G = \PGL_3(\C) $.
	The group $ \GL_3(\C) $ acts linearly on $ \Pi $ by $ (g\cdot F)(v) = F(g^{-1}v) $.
	Two lifts of an element of $ G $ to $ \GL_3(\C) $ differ by a scalar, so the induced action of $ G $ on the set of linear subspaces of $ \Pi $ is well defined, as are the actions of $ G $ on $ \C P^2 $, on $ \C P^{2\vee} $, and hence on $ P $ and on $ 2^P $.
	As noted after \autoref{thm:tacnode_calc}, tacnodes are preserved by this action, that is, $ \Sing(g\cdot F) = g\Sing(F) $.
	Consequently $ L(gK) = g\cdot L(K) $ for every $ K \subset P $ and every $ g \in G $, and therefore $ \dim_{\C} L(gK) = \dim_{\C} L(K) $.
	Thus $ \dim_{\C} L(K) $ depends only on the $ G $-orbit of $ K $ in $ 2^P $, and it is enough to check that each $ X_i $ is a single $ G $-orbit.
	This is immediate from \autoref{table:configs_def}:
	\begin{itemize}
		\item $ X_1 $ and $ X_7 $ are determined by the incident pair $ (p,l) \in P $, and $ G $ acts transitively on $ P $;
		\item $ X_4 $ and $ X_{11} $ are determined by the line $ l $, and $ X_5 $ by the point $ p $, and $ G $ acts transitively on lines and on points;
		\item $ X_2 $ and $ X_3 $ are determined by $ p $ together with an unordered pair, respectively triple, of distinct lines through $ p $.
		The stabiliser of $ p $ in $ G $ acts on the pencil $ \PP(T_p\C P^2) \cong \C P^1 $ as the whole of $ \PGL_2(\C) $, which is sharply $ 3 $-transitive;
		\item $ X_6 $ is determined by the ordered pair of distinct lines $ (l_1, l_2) $, and $ X_{10} $ by the unordered pair $ \{ l_1, l_2 \} $, and $ G $ acts transitively on both;
		\item $ X_{8} $ and $ X_{9} $ are determined by the points $ p $, $ q $ and $ r \coloneq l_p \cap l_q $ up to interchanging $ p $ and $ q $, since $ l_p = pr $, $ l_q = qr $, and, in the case of $ X_9 $, $ l = pq $.
		These three points are never collinear: $ p, r \in l_p $ and $ q, r \in l_q $, so a line containing all three would have to equal both $ l_p $ and $ l_q $, whereas $ l_p \neq l_q $.
		Since $ G $ acts transitively on ordered triples of points in general position, it acts transitively on $ X_8 $ and on $ X_9 $;
		\item $ X_{12} $ is determined by the smooth conic $ C $, and $ G $ acts transitively on smooth conics;
		\item $ X_{13} = \{ P \} $ consists of a single configuration.
	\end{itemize}
	Hence $ d_i $ is well defined.
	Its value is computed in \autoref{thm:dimensions}.
	The normalisation of $ K $ carried out there is legitimate precisely because of the transitivity just established.
\end{proof}

\begin{proposition}[Condition 4]
	The intersection $ X_i \cap X_j = \emptyset $ is empty if $ i \neq j $.
\end{proposition}

\begin{proof}
	This follows directly from \autoref{table:configs_def}.
\end{proof}

\begin{definition}
	\label{pg:26}
	Let $ K \in \bigcup_i \overline{X}_i $ be a configuration.
	Its \emph{geometrisation} $ \Geom(K) $ is the smallest configuration in $ \bigcup_i X_i $ that contains it.
\end{definition}

For a configuration $ K \in \bigcup X_i $, $ \Geom(K) = K $.

\begin{proposition}[Condition 5]
	\label{prop:condition5}
	For every configuration $ K \in \overline{X}_i \setminus X_i $, $ \Geom(K) \in \bigcup_{j < i} X_j $ and \mbox{$ L(K) = L(\Geom(K)) $}.
\end{proposition}

The proof of \autoref{prop:condition5} computes the closures $ \overline X_i $ of the configuration spaces in $ 2^P $.
For this we use the following properties of the Hausdorff metric $ \rho_H $ on $ 2^P $, which are \cite[Proposition 1]{BerGorKon23}.

\begin{remark}
	\label{rem:hausdorff}
	\begin{enumerate}[leftmargin=*, label=(\arabic*)]
		\item The space $ (2^P, \rho_H) $ is a compact metric space, where $ \rho_H $ is the Hausdorff metric on $ 2^P $.
		\item Let $ \{ K_i \} $ be a convergent sequence in $ (2^P, \rho_H ) $.
		Then $ x \in \lim K_i $, if and only if there exists a sequence $ \{ x_i \} $ that converges to $ x $ and such that every $ x_i \in K_i $.
		\item Consider convergent sequences $ \{ K'_i \} \to K' $ and $ \{ K''_i \} \to K'' $ in $ (2^P, \rho_H ) $.
		We have \mbox{$ \lim (K'_i \cup K''_i) = K' \cup K'' $}.
	\end{enumerate}
\end{remark}

From these properties we derive a criterion for computing such closures, \autoref{lem:closure_scheme} below, and apply it to each $ X_i $ in turn.
The criterion is phrased in terms of the following notion.

\begin{definition}
	Let $ M $ be a metric space and $ M^\circ \subset M $ a subset.
	A map $ \varphi \colon M \to 2^P $ is \emph{continuous along $ M^\circ $} if $ \varphi(m_t) \to \varphi(m) $ for every $ m \in M $ and every sequence $ \{ m_t \} \subset M^\circ $ converging to $ m $.
\end{definition}

\begin{lemma}
	\label{lem:closure_scheme}
	Let $ M $ be a compact metric space and $ M^\circ \subset M $ a dense subset.
	\begin{enumerate}[leftmargin=*, label=(\arabic*)]
		\item If $ \varphi \colon M \to 2^P $ is continuous along $ M^\circ $, then the closure of $ \varphi(M^\circ) $ in $ 2^P $ is $ \varphi(M) $.
		\item A map $ \varphi \colon M \to 2^P $ is continuous along $ M^\circ $ if its "graph" $ \Gamma_\varphi = \{ (m, x) \in M \times P \mid x \in \varphi(m) \} $ is closed in $ M \times P $ and $ \varphi $ satisfies the following \emph{lower bound}:

		For every $ m \in M $ and every sequence $ \{ m_t \} \subset M^\circ $ converging to $ m $, every point of a dense subset of $ \varphi(m) $ is the limit of a sequence $ \{ x_t \} $ with $ x_t \in \varphi(m_t) $.
		\item A finite union of maps $ M \to 2^P $ that are continuous along $ M^\circ $ is continuous along $ M^\circ $.
		\item The maps to $ 2^P $ that send a flag $ (p, m) \in P $ to $ \{ (p, m) \} $, a line $ l $ to the line of tacnodes $ l \times \{ l \} $ and to the thick line $ \{ (q, m) \in P \mid q \in l \} $, and a point $ p $ to the pencil $ \{ p \} \times \PP(T_p\C P^2) $ are continuous.
	\end{enumerate}
\end{lemma}

We use \autoref{lem:closure_scheme} to prove \autoref{prop:condition5} as follows.
For each $ i $ we parametrise the configurations of $ X_i $ together with their degenerations.
That is, we exhibit a compact metric space $ M_i $, a dense subset $ M_i^\circ \subset M_i $ and a map $ \varphi_i \colon M_i \to 2^P $ with $ \varphi_i(M_i^\circ) = X_i $.
Here $ M_i $ is a projective variety of the data that determine a configuration of $ X_i $, such as points, lines and flags, or, for $ i = 12 $, conics.
The subset $ M_i^\circ $ is the locus where these data satisfy the genericity conditions of \autoref{table:configs_def}, and $ \varphi_i $ sends the data to the configuration they determine.
Parts (2)--(4) of the lemma show that $ \varphi_i $ is continuous along $ M_i^\circ $, and part (1) then gives $ \overline X_i = \varphi_i(M_i) $.
Hence every configuration $ K \in \overline X_i \setminus X_i $ is a value of $ \varphi_i $ on the degenerate data $ M_i \setminus M_i^\circ $, and these values can be explicitly listed case by case.
For each of them we determine $ \Geom(K) $ and check that it lies in a configuration space $ X_j $ with $ j < i $.
The equality $ L(K) = L(\Geom(K)) $ is immediate when $ K $ itself lies in some $ X_j $, as then $ \Geom(K) = K $.
For the four families of configurations that lie in no $ X_j $, it is verified in the proof of \autoref{thm:dimensions}.

\begin{proof}[Proof of \autoref{lem:closure_scheme}]
	\begin{enumerate}[leftmargin=*, label=(\arabic*)]
		\item Given $ m \in M $, choose a sequence $ \{ m_t \} \subset M^\circ $ with $ m_t \to m $.
		Then \mbox{$ \varphi(m) = \lim \varphi(m_t) $} lies in the closure of $ \varphi(M^\circ) $.
		Conversely, every configuration $ K $ in the closure of $ \varphi(M^\circ) $ is the limit of a sequence $ \varphi(m_t) $ with $ m_t \in M^\circ $, as $ 2^P $ is a metric space.
		Passing to a subsequence $\{ m'_t \} \subset \{ m_t \}$, $ m'_{t} \to m \in M $ by compactness of $ M $, so that \mbox{$ K = \lim \varphi(m_t) = \lim \varphi(m'_t) = \varphi(m) \in \varphi(M) $}.

		\item Fix $ m \in M $ and a sequence $ \{ m_t \} \subset M^\circ $ with $ m_t \to m $.
		Since the metric space $ 2^P $ is compact by part (1) of \autoref{rem:hausdorff}, it suffices to show that every convergent subsequence of $ \{ \varphi(m_t) \} $ converges to $ \varphi(m) $.
		By part (2) of the remark, the limit $ K $ of such a subsequence consists of the limits of sequences $ \{ x_t \} $ with $ x_t \in \varphi(m_t) $, the index running over the subsequence.
		If $ x_t \to x $ is such a sequence, then $ (m_t, x_t) \in \Gamma_\varphi $ converges to $ (m, x) $, which lies in $ \Gamma_\varphi $ because the graph is closed.
		Hence $ K \subset \varphi(m) $.
		Conversely, $ K $ is closed in $ P $ and contains a dense subset of $ \varphi(m) $ by the lower bound, so $ \varphi(m) \subset K $.

		\item This follows from part (3) of \autoref{rem:hausdorff}.

		\item The first map is an isometry of $ P $ onto its image in $ 2^P $.
		For the other three we use the action of $ G = \PGL_3(\C) $: if $ g_t \to \mathrm{id} $ in $ G $, then \mbox{$ \rho_H(g_t K, K) \leq \sup_{x \in P} d(g_t x, x) \to 0 $} for every $ K \in 2^P $, since the action $ G \times P \to P $ is continuous and $ P $ is compact.
		These three maps commute with the action of $ G $, for instance \mbox{$ g (l \times \{ l \}) = gl \times \{ gl \} $} for every $ g \in G $, because $ G $ acts on $ P $ by bijections that preserve incidence.
		Moreover, $ G $ acts transitively on points and on lines, and the orbit maps \mbox{$ G \to \C P^2 $} and \mbox{$ G \to \C P^{2\vee} $} are locally trivial fibre bundles ($ G $ is a complex Lie group and the stabilisers are closed subgroups), hence admit continuous local sections.
		Therefore, if $ s $ is such a section of $ g \mapsto gl $ near $ l $, any convergent sequence of lines $ l_t \to l $ can be written as $ l_t = g_t l $ with \mbox{$ g_t = s(l_t) s(l)^{-1} \to \mathrm{id} $}.
		By the commutation and the bound above, $ \rho_H(l_t \times \{ l_t \}, l \times \{ l \}) = \rho_H(g_t (l \times \{ l \}), l \times \{ l \}) \to 0 $.
		This shows sequential continuity of $l \mapsto l \times \{ l \}$.
		The same argument applies to thick lines and to pencils.
	\end{enumerate}
\end{proof}

We are now ready to check Condition 5.

\begin{proof}[Proof of \autoref{prop:condition5}]
	We carry out the plan described after \autoref{lem:closure_scheme}, beginning with what all cases share.

	Each $ M_i $ below is a projective variety, in all cases but one cut out by incidence conditions inside a product of copies of $ \C P^2 $ and $ \C P^{2\vee} $ and equipped with the metric of that product, and $ M_i^\circ $ is the locus in $ M_i $ where the genericity conditions of \autoref{table:configs_def} hold, so that $ \varphi_i(M_i^\circ) = X_i $ by construction.
	Every $ M_i $ below is irreducible: it is a point, a projective space, a product of flag varieties, an iterated bundle whose fibres are products of projective lines, or, in the last case, the closure of the image of an irreducible set.
	The genericity conditions are non-degeneracy conditions.
	They require that certain points, or certain lines, be distinct, and, for $ i = 12 $, that certain determinants are non-zero.
	So each of them fails on a Zariski-closed subset, proper because $ X_i \neq \varnothing $.
	Hence $ M_i^\circ $ is a non-empty Zariski-open subset of the irreducible variety $ M_i $, and is therefore dense in $ M_i $ in the complex topology.

	For $ i \neq 12 $, every value of $ \varphi_i $ is a finite union of single flags, lines of tacnodes, pencils and thick lines whose data depend continuously on $ m \in M_i $.
	So $ \varphi_i $ is continuous along $ M_i^\circ $ by parts (3) and (4) of \autoref{lem:closure_scheme}, hence $ \overline X_i = \varphi_i(M_i) $ by part (1), and it remains to list the values of $ \varphi_i $ on $ M_i \setminus M_i^\circ $.
	The map $ \varphi_{12} $, whose values on $ M_{12}^\circ $ are tangency curves of smooth conics, is treated separately in item (12), with part (2) of the lemma.

	\begin{enumerate}[leftmargin=*, label=(\arabic*)]
		\item $ \overline X_1 \setminus X_1 = \emptyset $:
		\nopagebreak

		Take $ M_1 = M_1^\circ = P $ and $ \varphi_1(p, m) = \{ (p, m) \} $.
		Since $ M_1 \setminus M_1^\circ = \emptyset $, the closure is \mbox{$ \overline X_1 = \varphi_1(M_1) = X_1 $}.

		$ \configlonetac $
		converge to $ \configlonetac $.

		\item $ \overline X_{2} \setminus X_{2} = X_1 $:
		\nopagebreak

		Take $ M_2 = \{ (p, m_1, m_2) \in \C P^2 \times (\C P^{2\vee})^2 \mid p \in m_1 \cap m_2 \} $ with $ \varphi_2(p, m_1, m_2) = \{ (p, m_1), (p, m_2) \} $ and $ M_2^\circ = \{ m_1 \neq m_2 \} $.
		If $ m_1 = m_2 $, the two flags coincide and $ \varphi_2(p, m_1, m_1) = \{ (p, m_1) \} \in X_1 $.

		$ \configsmallcross $ converge to $ \configlonetac $, $ \configsmallcross $.

		\item $ \overline X_3 \setminus X_3 = X_1 \cup X_{2} $:
		\nopagebreak

		The same construction with triples of lines through $ p $:
		\begin{equation*}
			\varphi_3(p, m_1, m_2, m_3) = \{ (p, m_1), (p, m_2), (p, m_3) \} .
		\end{equation*}
		If exactly two of the three lines coincide, say $ m_2 = m_3 $, then two of the flags coincide and $ \varphi_3(p, m_1, m_2, m_2) = \{ (p, m_1), (p, m_2) \} \in X_2 $.
		If all three lines coincide, then $ \varphi_3(p, m_1, m_1, m_1) = \{ (p, m_1) \} \in X_1 $.

		$ \configtriple $ converge to $ \configlonetac $, $ \configsmallcross $, $ \configtriple $.

		\item $ \overline X_4 \setminus X_4 = \emptyset $:
		\nopagebreak

		Take $ M_4 = M_4^\circ = \C P^{2\vee} $ and $ \varphi_4(l) = l \times \{ l \} $.

		$ \configline $ converge to $ \configline $.

		\item $ \overline X_5 \setminus X_5 = \emptyset $:
		\nopagebreak

		Take $ M_5 = M_5^\circ = \C P^2 $ and $ \varphi_5(p) = \{ p \} \times \PP(T_p\C P^2) $.

		$ \configcirc $ converge to $ \configcirc $.

		\item $ \overline X_6 \setminus X_6 = X_4 $:
		\nopagebreak

		Take $ M_6 = \{ (l_1, p, l_2) \mid p \in l_1 \cap l_2 \} $ with $ \varphi_6(l_1, p, l_2) = (l_1 \times \{ l_1 \}) \cup \{ (p, l_2) \} $ and $ M_6^\circ = \{ l_2 \neq l_1 \} $.
		If $ l_2 = l_1 $, the extra flag is absorbed into the line: $ \varphi_6(l_1, p, l_1) = l_1 \times \{ l_1 \} \in X_4 $.

		$ \configlinetac $ converge to $ \configline $, $ \configlinetac $.

		\item $ \overline X_7 \setminus X_7 = \emptyset $:
		\nopagebreak

		Take $ M_7 = M_7^\circ = \{ (l, p) \mid p \in l \} $ and $ \varphi_7(l, p) = (l \times \{ l \}) \cup (\{ p \} \times \PP(T_p\C P^2)) $.

		$ \configlinecirc $ converge to $ \configlinecirc $.

		\item $ \overline X_{8} \setminus X_{8} = X_1 \cup X_{2} \cup W_1 \cup W_2 $, where $ W_1 $ and $ W_2 $ are the families of non-geometric configurations described below:
		\nopagebreak

		Take $ M_8 = P \times P $ with $ \varphi_8((p, l_p), (q, l_q)) = \{ (p, l_p), (q, l_q) \} $, and let $ M_8^\circ $ be the locus where $ p $, $ q $ and $ r \coloneq l_p \cap l_q $ are three distinct points.
		On $ M_8 \setminus M_8^\circ $:
		\begin{itemize}
			\item if $ p = q $, then $ \varphi_8((p, l_p), (p, l_q)) \in X_2 $ when $ l_p \neq l_q $, and $ \varphi_8((p, l_p), (p, l_p)) \in X_1 $ when $ l_p = l_q $;
			\item if $ p \neq q $ and $ l_p = l_q $, then the common direction passes through both points, so $ l_p = l_q = pq $, and we obtain the configuration $ \{ (p, pq), (q, pq) \} $.
			The family of these is denoted $ W_1 $;
			\item if $ p \neq q $, $ l_p \neq l_q $ and $ r \in \{ p, q \} $, say $ r = p $, then $ l_q $ passes through both $ p $ and $ q $, so $ l_q = pq $, and we obtain the configuration $ \{ (p, l_p), (q, pq) \} $ with $ l_p \neq pq $.
			The family of these is denoted $ W_2 $.
		\end{itemize}

		$ \configparallelpair $ converge to $ \configlonetac $, $ \configsmallcross $, $ \vcenter{\hbox{ \begin{tikzpicture}[scale=0.66]
					\draw[fill] (-0.75,0) circle (1.5pt);
					\draw[fill] (0.75,0) circle (1.5pt);
					\draw (-1.25,0) -- (-0.25,0);
					\draw (0.25,0) -- (1.25,0);
		\end{tikzpicture} }} \in W_1 $, $ \vcenter{\hbox{ \begin{tikzpicture}[scale=0.66]
					\draw[fill] (-0.25,0) circle (1.5pt);
					\draw[fill] (0.75,0) circle (1.5pt);
					\draw (-0.25,-0.25) -- (-0.25,0.25);
					\draw (0.25,0) -- (1.25,0);
		\end{tikzpicture} }} \in W_2 $, $ \configparallelpair $.

		Neither family contains a configuration from any of $ X_i $.
		Indeed, a configuration consisting of two flags at two distinct points either belongs to $ X_8 $ or to no $ X_i $ at all, and both families fail the genericity condition of $ X_8 $.
		Their geometrisations exist nevertheless.
		Let us find the configuration spaces that contain them.

		Let $ K $ be a configuration from $ W_1 $ or $ W_2 $, and suppose that its geometrisation $\Geom(K)$ exists.
		In both cases $ \Geom(K) $ contains the flag $ (q, pq) $ together with a flag at the distinct point $ p $.
		The configuration spaces $ X_1, X_2, X_3 $ and $ X_5 $ are excluded, their configurations being supported at a single point.
		The space $X_8$ is excluded, just like before.
		The space $ X_{12} $ is also excluded: if $ \PP(TC) $ contained $ (q, pq) $, then $ pq $ would be tangent to $ C $ at $ q $, and a tangent line meets a smooth conic only at its point of tangency, leaving no room for a flag at $ p $.

		The remaining configuration spaces are $ X_4, X_6, X_7, X_9, X_{10}, X_{11} $ and $ X_{13} $.
		For $ K \in W_1$, there is a unique line of tacnodes $ pq \times \{ pq \} = L \in X_4$ such that $K \subset L$.
		It follows from \autoref{fig:config_relations} that $ \Geom(K) = L$, as $L$ is the smallest configuration from $X_i$ containing $K$.
		Similarly, for $K \in W_2 $ we have $ \Geom(K) = (pq \times \{ pq \}) \cup \{ (p, l_p) \} \in X_6$.

		We will see that $ L(\Geom(K)) = L(K) $ for both families in \autoref{thm:dimensions}, items (14) and (15) of its proof.

		\item $ \overline X_{9} \setminus X_{9} = X_4 \cup X_6 \cup W_3 $, where $ W_3 $ is the family of non-geometric configurations described below:
		\nopagebreak

		Take $ M_9 = \{ (l, (p, l_p), (q, l_q)) \in \C P^{2\vee} \times P \times P \mid p, q \in l \} $ with
		\begin{equation*}
			\varphi_9(l, (p, l_p), (q, l_q)) = (l \times \{ l \}) \cup \{ (p, l_p), (q, l_q) \} ,
		\end{equation*}
		and let $ M_9^\circ $ be the locus where $ p $, $ q $ and $ l_p \cap l_q $ are three distinct points.
		This forces $ l_p \neq l \neq l_q $: if, say, $ l_p = l $, then $ l_p \cap l_q = q $.
		On $ M_9 \setminus M_9^\circ $, the value of $ \varphi_9 $ always contains $ l \times \{ l \} $, and the possible degenerations are:
		\begin{itemize}
			\item $ \varphi_9(l, (p, l_p), (q, l_q)) = l \times \{ l \} \in X_4 $ when both extra flags are absorbed into the line: this happens for $ l_p = l_q = l $, and for $ p \neq q $ with $ l_p = l_q $, since a common direction through two distinct points of $ l $ must equal $ pq = l $;
			\item $ \varphi_9(l, (p, l_p), (q, l_q)) \in X_6 $ when exactly one extra flag survives: for $ p = q $ with exactly one of $ l_p, l_q $ equal to $ l $, or with $ l_p = l_q \neq l $.
			This also happens for $ p \neq q $ with $ l_p \neq l_q $ and $ l_p \cap l_q \in \{ p, q \} $.
			If, say, $ l_p \cap l_q = p $, then $ l_q $ contains both $ p $ and $ q $, so $ l_q = pq = l $ is absorbed, while $ l_p \neq l_q = l $ survives;
			\item if $ p = q $ and $ l \neq l_p \neq l_q \neq l $, we obtain the configuration $ (l \times \{ l \}) \cup \{ (p, l_p), (p, l_q) \} $, a line of tacnodes with two extra directions at one of its points.
			The family of these is denoted $ W_3 $.
		\end{itemize}

		$ \configlinepair $ converge to $ \configline $, $ \configlinetac $, $ \vcenter{\hbox{ \begin{tikzpicture}[scale=0.66]
					\draw[fill] (0,0) circle (1.5pt);
					\draw (-1,0) -- (1,0);
					\draw (-0.707/4,-0.707/4) -- (0.707/4,0.707/4);
					\draw (0,-0.25) -- (0,0.25);
		\end{tikzpicture} }} \in W_3 $, $ \configlinepair $.

		A configuration $ K \in W_3 $ lies in no $ X_i $: it contains a line of tacnodes and exactly two extra flags, both at the same point of the line, whereas $ X_6 $ has one extra flag, $ X_9 $ has its two extra flags at distinct points, and $ X_7 $ has the full pencil.
		For the same reason, no configuration of \autoref{table:configs_def} contains $ l \times \{ l \} $ together with exactly two extra directions at $ p $, so any configuration containing $ K $ contains the line and the full pencil at $ p $.
		The smallest such is
		\begin{equation*}
			\Geom(K) = (l \times \{ l \}) \cup (\{ p \} \times \PP(T_p\C P^2)) \in X_7 .
		\end{equation*}
		We will see that $ L(\Geom(K)) = L(K) $ for this family in \autoref{thm:dimensions}, item (16) of its proof.

		\item $ \overline X_{10} \setminus X_{10} = X_7 $:
		\nopagebreak

		Take $ M_{10} = \{ (l_1, l_2, p) \mid p \in l_1 \cap l_2 \} $ with
		\begin{equation*}
			\varphi_{10}(l_1, l_2, p) = (l_1 \times \{ l_1 \}) \cup (l_2 \times \{ l_2 \}) \cup (\{ p \} \times \PP(T_p\C P^2))
		\end{equation*}
		and $ M_{10}^\circ = \{ l_1 \neq l_2 \} $, on which $ p = l_1 \cap l_2 $ is determined.
		If $ l_1 = l_2 $, then $ p $ is an arbitrary point of $ l_1 $ and $ \varphi_{10}(l_1, l_1, p) = (l_1 \times \{ l_1 \}) \cup (\{ p \} \times \PP(T_p\C P^2)) \in X_7 $.

		$ \configcross{} $ converge to $ \configlinecirc $, $ \configcross{} $.

		\item $ \overline X_{11} \setminus X_{11} = \emptyset $:
		\nopagebreak

		Take $ M_{11} = M_{11}^\circ = \C P^{2\vee} $ and $ \varphi_{11}(l) = \{ (q, m) \in P \mid q \in l \} $, the thick line over $ l $.

		$ \configfatline $ converge to $ \configfatline $.

		\item $ \overline X_{12} \setminus X_{12} = X_{10} \cup X_7 \cup W_4 $, where $ W_4 $ is the family of non-geometric configurations described below:
		\nopagebreak
		\setlength{\parskip}{4pt}

		$ \configconic $ converge to $ \configcross{} $, $ \configlinecirc $, $ \vcenter{\hbox{ \begin{tikzpicture}[scale=0.66]
					\draw (-1,0) -- (1,0);
					\draw (-0.5,0) circle (2.5pt);
					\draw (0.5,0) circle (2.5pt);
		\end{tikzpicture} }} \in W_4 $, $ \configconic $.

		A configuration of $ X_{12} $ is the tangency curve $ \PP(TC) $ of a smooth conic $ C $.
		When $ C $ degenerates, one of three things happens.
		\begin{itemize}
			\item $ C $ becomes a pair of distinct lines $ l_1 \cup l_2 $ that meet at a point $ q $.
			The limit consists of the two lines of tacnodes and of the pencil at $ q $, and lies in $ X_{10} $.
			\item $ C $ becomes a double line $ l $, and its tangent lines become the lines through two distinct points $ u, v \in l $.
			The limit consists of the line of tacnodes and of the pencils at $ u $ and at $ v $.
			These configurations form $ W_4 $.
			\item The same with $ u = v $.
			The limit lies in $ X_7 $.
		\end{itemize}
		In the last two cases the double line $ l $ alone does not determine the limit: the points $ u $ and $ v $ depend on how $ C $ degenerates.
		This is why the parameter space $ M_{12} $ below records a conic together with the curve of its tangent lines.

		The proof has three steps.
		In Step 1 we define $ M_{12} $ and $ \varphi_{12} $ and check that $ \varphi_{12}(M_{12}^\circ) = X_{12} $.
		In Step 2 we list the points of $ M_{12} \setminus M_{12}^\circ $ and the values of $ \varphi_{12} $ at them.
		In Step 3 we show that $ \varphi_{12} $ is continuous along $ M_{12}^\circ $.
		Step 3 uses the list of Step 2, which is why it comes last.

		\medskip
		We write points of $ \C P^2 $ as columns $ x $ and lines as rows $ \lambda $, the line being $ \{ x \mid \lambda x = 0 \} $.
		In this item a flag is written $ (x, m) $, where $ x $ is its point and $ m $ is its line.
		A point of the parameter space is written $ ([A],[B]) $.

		A conic is the zero set $ \V(A) = \{ x \mid x^T A x = 0 \} $ of a nonzero symmetric $ 3 \times 3 $ matrix $ A $ determined up to scale.
		It is smooth if and only if $ \det A \neq 0 $.
		A nonzero symmetric matrix $ B $ dually determines the curve
		\begin{equation*}
			\V^\vee(B) = \{ m \in \C P^{2 \vee} \mid m B m^T = 0 \} .
		\end{equation*}
		We use four facts.
		In the first three $ \det A \neq 0 $, and $ \adj A = \det(A) A^{-1} $ is the adjugate matrix.
		\begin{enumerate}[label=(\alph*)]
			\item The tangent line of $ \V(A) $ at $ x \in \V(A) $ is the row $ x^T A $.
			\item The tangent lines of $ \V(A) $ constitute exactly $ \V^\vee(\adj A) $.

			Indeed, $ (x^T A) \adj (A) (x^T A)^T = \det(A) \, x^T A x $, so every tangent line of $ \V(A) $ lies on the conic $ \V^\vee(\adj A) $.
			This conic is smooth, because $ \det \adj A = (\det A)^2 \neq 0 $.
			Hence it is irreducible, and the two curves coincide.
			\item A tangent line of $ \V(A) $ meets $ \V(A) $ only at its point of tangency.

			Indeed, the intersection of a line with a conic has total multiplicity $ 2 $, and tangency uses all of it.
			\item $ \adj (g^T A g) = \det(g)^2 \, g^{-1} (\adj A) g^{-T} $ for every invertible matrix $ g $.
		\end{enumerate}

		\medskip
		\textbf{Step 1.} $ \varphi_{12}(M_{12}^\circ) = X_{12} $.
		\nopagebreak

		Let $ M_{12} \subset \C P^5 \times \C P^5 $ be the variety of \emph{complete conics}: the closure of the set
		\begin{equation*}
			M_{12}^\circ = \left \{ ([A], [\adj A]) \ \middle| \ \det A \neq 0 \right \} ,
		\end{equation*}
		where each $ \C P^5 $ is the projectivised space of symmetric $ 3 \times 3 $ matrices.
		Being the closure of the image of an irreducible set under an algebraic map, $ M_{12} $ is an irreducible projective variety.
		Moreover,
		\begin{equation*}
			M_{12}^\circ = M_{12} \cap \{ \det A \neq 0 \} ,
		\end{equation*}
		so $ M_{12}^\circ $ is open and dense in $ M_{12} $.
		Indeed, let $ ([A],[B]) \in M_{12} $ with $ \det A \neq 0 $ be a limit of $ ([A_t],[\adj A_t]) $ with representatives $ A_t \to A $.
		Then $ \adj A_t \to \adj A \neq 0 $ by continuity of the adjugate, so $ [B] = [\adj A] $.

		Define
		\begin{equation*}
			\varphi_{12}([A],[B]) = \left \{ (x, m) \in P \ \middle| \ x^T A x = 0 , \ m B m^T = 0 \right \} .
		\end{equation*}
		The graph of $ \varphi_{12} $ is cut out by these polynomial conditions together with the incidence condition of $ P $, hence closed.

		On $ M_{12}^\circ $, the value $ \varphi_{12}([A],[\adj A]) $ consists of all flags $ (x, m) $ where $ x $ lies on the smooth conic $ C = \V(A) $ and, by fact (b), $ m $ is a tangent line of $ C $ passing through $ x $.
		By fact (c), $ m $ is the tangent line at $ x $, so $ m = x^T A $ by fact (a), and
		\begin{equation*}
			\varphi_{12}([A], [\adj A]) = \left \{ (x, x^T A) \ \middle| \ x \in C \right \} = \PP(TC) .
		\end{equation*}
		Hence $ \varphi_{12}(M_{12}^\circ) = X_{12} $.

		\medskip
		\textbf{Step 2.} The points of $ M_{12} \setminus M_{12}^\circ $.
		\nopagebreak

		We show that the points $ ([A],[B]) $ of $ M_{12} \setminus M_{12}^\circ $ and the values of $ \varphi_{12} $ at them are exactly the following.
		\begin{center}
			\small
			\setlength{\tabcolsep}{4pt}
			\renewcommand{\arraystretch}{1.6}
			\begin{tabular}{|c|c|c|c|c|c|}
				\hline
				$ A $ & $ B $ & $ \V(A) $ & $ \V^\vee(B) $ & \multicolumn{2}{c|}{$ \varphi_{12}([A],[B]) $} \\
				\hline
				$ \rk A = 2 $ & $ q q^T $ & $ l_1 \cup l_2 $, \ $ l_1 \cap l_2 = q $ & the lines through $ q $ & $ \configcross{} $ & $ \in X_{10} $ \\
				\hline
				$ \rk A = 1 $ & $ u v^T + v u^T $, \ $ u \neq v $ & $ l $, doubled, \ $ u, v \in l $ & the lines through $ u $ or $ v $ & $ \vcenter{\hbox{ \begin{tikzpicture}[scale=0.66]
							\draw (-1,0) -- (1,0);
							\draw (-0.5,0) circle (2.5pt);
							\draw (0.5,0) circle (2.5pt);
				\end{tikzpicture} }} $ & $ \in W_4 $ \\
				\hline
				$ \rk A = 1 $ & $ u u^T $ & $ l $, doubled, \ $ u \in l $ & the lines through $ u $ & $ \configlinecirc $ & $ \in X_7 $ \\
				\hline
			\end{tabular}
		\end{center}
		The second row is the definition of the family $ W_4 $.

		First, let us show that $ A B = 0 $.
		Let $ ([A],[B]) \in M_{12} \setminus M_{12}^\circ $, so that $ \det A = 0 $.
		Choose a sequence in $ M_{12}^\circ $ converging to it,
		\begin{equation*}
			([A_t], [\adj A_t]) \to ([A],[B]) ,
		\end{equation*}
		with representatives chosen (after passing to a subsequence) so that $ A_t \to A $ and $ B_t \to B $, where $ B_t $ is proportional to $ \adj A_t $.
		Since $ A_t \adj A_t = \det (A_t) \, \mathrm{Id} $, the product $ A_t B_t $ is a scalar multiple of the identity:
		\begin{equation*}
			A_t B_t = \delta_t \, \mathrm{Id} .
		\end{equation*}
		In the limit $ A B = \delta \, \mathrm{Id} $.
		Taking determinants, $ \delta^3 = \det A \det B = 0 $.
		Hence
		\begin{equation*}
			A B = 0 .
		\end{equation*}
		Since $ A \neq 0 $, two cases remain: $ \rk A = 2 $ and $ \rk A = 1 $.

		Next, let us show that there are no other points.
		\begin{itemize}
			\item $ \rk A = 2 $.
			The kernel of $ A $ is spanned by a single point $ q $, and $ \V(A) = l_1 \cup l_2 $ is a pair of distinct lines with $ l_1 \cap l_2 = q $.
			Indeed, in coordinates with $ q = (0:0:1) $ the form of $ A $ involves only $ x $ and $ y $ and is a nondegenerate binary quadratic form, hence a product of two nonproportional linear forms.

			The relation $ A B = 0 $ places every column of $ B $ in $ \ker A $, so $ B = q c^T $ for some column $ c $.
			The symmetry of $ B $ forces $ c $ to be proportional to $ q $.
			Thus $ [B] = [q q^T] $.

			\item $ \rk A = 1 $.
			Now $ A = \lambda^T \lambda $ up to scale, and $ \V(A) $ is the line $ l = \{ \lambda x = 0 \} $, doubled.
			The relation $ A B = \lambda^T (\lambda B) = 0 $ forces $ \lambda B = 0 $.

			Diagonalise the symmetric form $ B $, which has rank at most $ 2 $, as $ B = w_1 w_1^T + w_2 w_2^T $ (or $ B = w_1 w_1^T $).
			The condition $ \lambda B = (\lambda w_1) w_1^T + (\lambda w_2) w_2^T = 0 $ gives $ \lambda w_1 = \lambda w_2 = 0 $.
			Setting $ u = w_1 + i w_2 $ and $ v = w_1 - i w_2 $ we obtain, up to scale,
			\begin{equation*}
				B = u v^T + v u^T \qquad \text{with } u, v \in l ,
			\end{equation*}
			where $ u = v $ is allowed and covers the rank-one case.
		\end{itemize}

		Let us show that all these points occur.
		\begin{itemize}
			\item $ \rk A = 2 $.
			The first projection $ M_{12} \to \C P^5 $ has closed image containing the dense set $ \{ \det A \neq 0 \} $, hence is surjective.
			So $ ([A],[B]) \in M_{12} $ for some $ [B] $, and $ [B] = [q q^T] $ by the above.

			\item $ \rk A = 1 $ and $ u \neq v $.
			By fact (d), $ M_{12} $ is invariant under the coordinate changes induced by $ g $.
			Since $ \PGL_3(\C) $ acts transitively on triples (a line, two distinct points on it), we may assume
			\begin{equation*}
				\lambda = (0, 1, 0) , \qquad u = (1,0,0)^T , \qquad v = (0,0,1)^T ,
			\end{equation*}
			i.e. the form of $ A $ is $ y^2 $ and $ l = \{ y = 0 \} $.
			For $ t \neq 0 $, the conic $ A_t $ with the form $ y^2 + 2 t x z $ satisfies $ \det A_t = -t^2 \neq 0 $ and
			\begin{equation*}
				\adj A_t = - t \, ( u v^T + v u^T ) - t^2 \, \lambda^T \lambda .
			\end{equation*}
			So $ [A_t] \to [\lambda^T \lambda] = [A] $ and $ [\adj A_t] = [u v^T + v u^T + t \, \lambda^T \lambda] \to [B] $ as $ t \to 0 $, whence $ ([A],[B]) \in M_{12} $.

			\item $ \rk A = 1 $ and $ u = v $.
			This case follows by taking the limit $ v \to u $ within the fibre of $ M_{12} $ over $ [A] $, which is closed in $ \C P^5 $.
		\end{itemize}

		Finally, let us find the values of $ \varphi_{12} $.
		\begin{itemize}
			\item $ \rk A = 2 $.
			The curve $ \V^\vee(q q^T) = \{ m \mid (m q)^2 = 0 \} $ is the pencil of lines through $ q $, so
			\begin{equation*}
				\varphi_{12}([A],[q q^T]) = \left \{ (x,m) \in P \ \middle| \ x \in l_1 \cup l_2 , \ q \in m \right \} .
			\end{equation*}
			For $ x = q $ every line through $ q $ occurs, while for $ x \in l_i \setminus \{ q \} $ the only line through both $ x $ and $ q $ is $ l_i $.
			Hence
			\begin{equation*}
				\varphi_{12}([A],[q q^T]) = (l_1 \times \{ l_1 \}) \cup (l_2 \times \{ l_2 \}) \cup ( \{ q \} \times \PP(T_q\C P^2)) \in X_{10} .
			\end{equation*}

			\item $ \rk A = 1 $.
			Here $ m B m^T = 2 (m u)(m v) $, so $ \V^\vee(B) $ is the union of the pencils of lines through $ u $ and through $ v $.
			For $ x \in \{ u, v \} $ every line through $ x $ occurs, while for $ x \in l \setminus \{ u, v \} $ a line through $ x $ and through $ u $ or $ v $ must be $ l $ itself.
			Hence
			\begin{equation*}
				\varphi_{12}([A],[B]) = (l \times \{ l \}) \cup (\{ u \} \times \PP(T_u\C P^2)) \cup (\{ v \} \times \PP(T_v\C P^2)) .
			\end{equation*}
			For $ u = v $ this is a configuration in $ X_7 $.
			For $ u \neq v $ it lies in no $ X_i $, and the family of these configurations is denoted $ W_4 $.
		\end{itemize}

		\medskip
		\textbf{Step 3.} $ \varphi_{12} $ is continuous along $ M_{12}^\circ $.
		\nopagebreak

		Fix $ ([A],[B]) \in M_{12} $ and a sequence $ ([A_t],[\adj A_t]) \in M_{12}^\circ $ converging to it, normalised as in Step 2.
		Write $ C_t = \V(A_t) $ and $ K = \varphi_{12}([A],[B]) $.
		As the graph of $ \varphi_{12} $ is closed, by part (2) of \autoref{lem:closure_scheme} it suffices to verify the lower bound:
		\begin{center}
			every flag of a dense subset of $ K $ is the limit of a sequence of flags of $ \PP(TC_t) $.
		\end{center}
		We use two ways to approximate, (P) for points and (L) for lines, and one way to conclude, (U).
		\begin{enumerate}[label=(\Alph*), leftmargin=*]
			\item[(P)] For every $ x \in \V(A) $ there are $ x_t \in C_t $ with $ x_t \to x $.

			Indeed, let $ L $ be a line through $ x $ not contained in $ \V(A) $.
			One exists, as $ \V(A) $ contains at most two lines.
			The restrictions of $ A_t $ to $ L $ are binary quadratic forms converging to the nonzero form $ A|_L $, and the roots in $ L \cong \C P^1 $ of a nonzero binary quadratic form depend continuously on its coefficients.
			As $ x $ is a root of $ A|_L $, there are $ x_t \in \V(A_t) \cap L $ with $ x_t \to x $.

			\item[(L)] If $ ([A],[B]) \notin M_{12}^\circ $, then for every $ m \in \V^\vee(B) $ there are tangent lines $ m_t $ of $ C_t $ with $ m_t \to m $.

			Indeed, dually, let $ \mathcal{L} $ be a pencil containing $ m $ (a line in the dual plane) not contained in $ \V^\vee(B) $.
			One exists, as $ \V^\vee(B) $ is a union of at most two pencils by Step 2.
			The same root-continuity argument yields $ m_t \in \V^\vee(\adj A_t) $ with $ m_t \to m $, and $ m_t $ is a tangent line of $ C_t $ by fact (b).

			\item[(U)] Let $ (x_t, m_t) \in \PP(TC_t) $ with $ x_t \to x $.
			If $ (x, m) $ is the only flag of $ K $ with the point $ x $, then $ (x_t, m_t) \to (x, m) $.
			The same holds with the roles of the point and the line exchanged.

			Indeed, as the graph of $ \varphi_{12} $ is closed, every convergent subsequence of the flags $ (x_t, m_t) $ has its limit in $ K $.
			That limit has the point $ x $, so it equals $ (x, m) $.
			Since $ P $ is compact, the whole sequence converges to $ (x, m) $.
		\end{enumerate}

		Consider first the case $ ([A],[B]) \in M_{12}^\circ $.
		Every flag of $ K = \PP(TC) $ has the form $ (x, x^T A) $ with $ x \in C $.
		Taking $ x_t \to x $ as in (P), the flags $ (x_t, x_t^T A_t) \in \PP(TC_t) $ converge to $ (x, x^T A) $, which proves the lower bound.

		Consider now the case $ ([A],[B]) \notin M_{12}^\circ $.
		We verify the lower bound for two kinds of flags, which together form a dense subset of $ K $.
		\begin{center}
			\small
			\renewcommand{\arraystretch}{1.4}
			\begin{tabular}{|l|c|c|c|}
				\hline
				& $ \rk A = 2 $ & $ \rk A = 1 $ & by \\
				\hline
				the only flag of $ K $ with its point & $ (x, l_i) $, \ $ x \in l_i \setminus \{ q \} $ & $ (x, l) $, \ $ x \in l \setminus \{ u, v \} $ & (P), (U) \\
				\hline
				the only flag of $ K $ with its line & $ (q, m) $, \ $ m \notin \{ l_1, l_2 \} $ & $ (u, m) $ or $ (v, m) $, \ $ m \neq l $ & (L), (U) \\
				\hline
				the flags not covered & $ (q, l_1) $, $ (q, l_2) $ & $ (u, l) $, $ (v, l) $ & \\
				\hline
			\end{tabular}
		\end{center}
		\begin{itemize}
			\item Let $ (x, m) $ be a flag of the first row.
			Take $ x_t \in C_t $ with $ x_t \to x $ by (P) and consider the flags $ (x_t, x_t^T A_t) \in \PP(TC_t) $.
			They converge to $ (x, m) $ by (U).
			\item Let $ (x, m) $ be a flag of the second row.
			It is the only flag of $ K $ with the line $ m $: such a flag has its point in $ (l_1 \cup l_2) \cap m = \{ q \} $ if $ \rk A = 2 $, and in $ l \cap m $, a single point, if $ \rk A = 1 $.
			Take tangent lines $ m_t \to m $ by (L), and let $ x_t \in C_t $ be the point of tangency of $ m_t $.
			The flags $ (x_t, m_t) $ converge to $ (x, m) $ by (U), with the roles of the point and the line exchanged.
		\end{itemize}
		The flags of the third row are limits of flags of the second row.
		So the two kinds are indeed dense in $ K $, and $ \varphi_{12} $ is continuous along $ M_{12}^\circ $.

		\medskip
		\textbf{Conclusion.}
		By part (1) of \autoref{lem:closure_scheme}, $ \overline X_{12} = \varphi_{12}(M_{12}) $.
		By Step 1 and Step 2, this completes the proof that $ \overline X_{12} \setminus X_{12} = X_{10} \cup X_7 \cup W_4 $.

		A configuration $ K = (l \times \{ l \}) \cup (\{ p_+ \} \times \PP(T_{p_+}\C P^2)) \cup (\{ p_- \} \times \PP(T_{p_-}\C P^2)) $ in $ W_4 $ (writing $ p_+ , p_- $ for the two distinct points $ u, v \in l $) lies in no $ X_i $, as it contains full pencils at exactly two points.
		The only configurations of \autoref{table:configs_def} containing full pencils at two distinct points $ p_\pm $ are the thick line over $ p_+ p_- = l $, which contains the full pencil at a point exactly when the point lies on $ l $, and $ P $ itself.
		Both contain $ K $, and the smaller one is the thick line, so
		\begin{equation*}
			\Geom(K) = \{ (q, m) \in P \mid q \in l \} \in X_{11} .
		\end{equation*}
		We will see that $ L(\Geom(K)) = L(K) $ for this family in \autoref{thm:dimensions}, item (17) of its proof.

		\item $ \overline X_{13} \setminus X_{13} = \emptyset $:
		\nopagebreak

		$ X_{13} = \{ P \} $ is a single point of the metric space $ 2^P $, hence closed.
	\end{enumerate}

	Finally, we collect the verification of Condition 5 itself.
	Every geometric limit obtained above lies in a stratum with strictly smaller index.
	The inequalities used are $ 1 < 2 $, $ 1 < 3 $, $ 2 < 3 $, $ 4 < 6 $, $ 1 < 8 $, $ 2 < 8 $, $ 4 < 9 $, $ 6 < 9 $, $ 7 < 10 $, $ 7 < 12 $ and $ 10 < 12 $.
	For such $ K $ we have $ \Geom(K) = K $, so there is nothing more to check.
	The other limits $ K $ lie in the four non-geometric families, and their geometrisations satisfy
	\begin{equation*}
		\Geom(K) \in
		\begin{cases}
			X_4 & \text{for } K \in W_1 \subset \overline X_8 , \\
			X_6 & \text{for } K \in W_2 \subset \overline X_8 , \\
			X_7 & \text{for } K \in W_3 \subset \overline X_9 , \\
			X_{11} & \text{for } K \in W_4 \subset \overline X_{12} ,
		\end{cases}
	\end{equation*}
	with $ 4 < 8 $, $ 6 < 8 $, $ 7 < 9 $ and $ 11 < 12 $ respectively, so $ \Geom(K) \in \bigcup_{j < i} X_j $ in all cases, and the remaining equalities $ L(K) = L(\Geom(K)) $ are verified in items (14)--(17) of the proof of \autoref{thm:dimensions}.
\end{proof}

\begin{proposition}[Condition 7]
	If $ K $ is a finite configuration from $ X_i $, then every nonempty subset $ K' \subsetneq K $ belongs to some $ X_j $ with $ j < i $
\end{proposition}

\begin{proof}
	In our case, the finite configurations can only be found in the following configuration spaces:
	$ X_1, X_{2}, X_3, X_{8} $.
	For each $ i \neq 3 $, every proper nonempty subset is a configuration in $ X_1 $.
	In case of $ X_3 $, such subset is either in $ X_1 $ or $ X_{2} $.
\end{proof}

\section{Computation of the dimensions of the vector spaces \texorpdfstring{$ L(K) $}{}}
\label{sec:v1_dimensions}

To compute the spectral sequence we need to determine the dimensions $ d_i = \dim L(K) $ for $ K \in X_i $.

\begin{theorem}\label{thm:dimensions}
	The dimensions $ d_i $ are as in \autoref{table:configs_dim}.
	In each case we write a local form of a quartic in the affine chart $ U = \{ z = 1 \} $ satisfying the tacnode conditions prescribed by a configuration $K$.

	\begin{table}[!htbp]
		\caption{Dimensions of $ L(K) $}
		\label{table:configs_dim}
		\centering
			\renewcommand{\arraystretch}{1.3}
			\begin{tabular}{|c|c|c|c|}
				\hline
				$i$      & $d_i$ & $F\in L(K)$ (locally, in $ U = \{ z = 1 \} $) & Generic el. of $X_i$                              \\
				\hline
				1        & 9     & $f_4 + ax^2y+bxy^2+cy^3+dy^2$         & $ \configlonetac $                                \\
				\hline
				2        & 7     & $f_4+ax^2y+bxy^2$                     & $ \configsmallcross $                             \\
				\hline
				3        & 6     & $f_4+axy(x-y)$                        & $ \configtriple $                                 \\
				\hline
				4        & 6     & $y^2(ax^2+bxy+cy^2+dx+ey+g)$          & $ \configline $                                   \\
				\hline
				5        & 5     & $f_4$                                 & $ \configcirc $                                   \\
				\hline
				6        & 4     & $y^2(ax^2+bxy+cy^2+dx)$               & $ \configlinetac $                                \\
				\hline

				7        & 3     & $y^2(ax^2+bxy+cy^2)$                  & $ \configlinecirc $                               \\
				\hline
				8        & 3     & $ax^4+bx^2y(y-1)+cy^2(y-1)^2$         & $ \configparallelpair $                           \\
				\hline
				9        & 2     & $x^2(ax^2+by(y-1))$                   & $ \configlinepair $                               \\
				\hline
				10       & 1     & $ax^2y^2$                             & $ \configcross{} $                                \\
				\hline
				11       & 1     & $ay^4$                                & $ \configfatline $                                \\
				\hline
				12       & 1     & $a(x^2+y^2-1)^2$                      & $ \configconic $                                  \\
				\hline
				13       & 0     & 0                                     & $\PP T\C P^2$                                     \\
				\hline
				non-geom & 6     & (see 4)                               & $ \vcenter{\hbox{ \begin{tikzpicture}[scale=0.66]
							\draw[fill] (-0.75,0) circle (1.5pt);
							\draw[fill] (0.75,0) circle (1.5pt);
							\draw (-1.25,0) -- (-0.25,0);
							\draw (0.25,0) -- (1.25,0);
				\end{tikzpicture} }}$ \\
				\hline
				non-geom & 4     & (see 6)                               & $ \vcenter{\hbox{ \begin{tikzpicture}[scale=0.66]
							\draw[fill] (-0.25,0) circle (1.5pt);
							\draw[fill] (0.75,0) circle (1.5pt);
							\draw (-0.25,-0.25) -- (-0.25,0.25);
							\draw (0.25,0) -- (1.25,0);
				\end{tikzpicture} }}$ \\
				\hline
				non-geom & 3     & (see 7)                               & $ \vcenter{\hbox{ \begin{tikzpicture}[scale=0.66]
							\draw[fill] (0,0) circle (1.5pt);
							\draw (-1,0) -- (1,0);
							\draw (-0.707/4,-0.707/4) -- (0.707/4,0.707/4);
							\draw (0,-0.25) -- (0,0.25);
				\end{tikzpicture} }}$ \\
				\hline
				non-geom & 1     & (see 11)                              & $ \vcenter{\hbox{ \begin{tikzpicture}[scale=0.66]
							\draw (-1,0) -- (1,0);
							\draw (-0.5,0) circle (2.5pt);
							\draw (0.5,0) circle (2.5pt);
				\end{tikzpicture} }}$ \\
				\hline
			\end{tabular}
	\end{table}
\end{theorem}

\begin{proof}
	Suppose $ p = (0:0:1)$, $l={y=0}$.
	For a polynomial $f$ in two variables, consider its Taylor expansion around $ (0,0) $, i.e. the expression $ f(x,y) = $
	\begin{gather*}
		= a_{40}x^4+a_{31}x^3y+a_{22}x^2y^2+a_{13}xy^3+a_{04}y^4+\\
		+a_{30}x^3+a_{21}x^2y+a_{12}xy^2+a_{03}y^3+\\
		+a_{20}x^2+a_{11}xy+a_{02}y^2+\\
		+a_{10}x+a_{01}y+\\
		+a_{00} .
	\end{gather*}
	Suppose $ f(x,y) = F(x,y,1) $ is the restriction to the affine chart $ U = \{ z = 1 \} $ of a homogeneous polynomial $ F \in \Pi $.
	\begin{enumerate}[leftmargin=*, label=(\arabic*)]
		\item By \autoref{thm:tacnode_calc}, $ (p, l) $ as above is a tacnode of $ Q = \V(F) $ if and only if
		\[
		0 = a_{00} = a_{10} = a_{01} = a_{20} = a_{11} = a_{30} .
		\]
		Thus, $ d_1 = 9$, and the polynomial $ f $ is of the form $ f(x,y) = $
		\begin{gather*}
			= a_{40}x^4+a_{31}x^3y+a_{22}x^2y^2+a_{13}xy^3+a_{04}y^4+\\
			+a_{21}x^2y+a_{12}xy^2+a_{03}y^3+\\
			+a_{02}y^2 .
		\end{gather*}

		\item Suppose additionally that $ l' = \{ x = 0 \} $.
		Then $ (p,l') $ is also a tacnode of $ Q $ if and only if additionally $ 0 = a_{02} = a_{03} $.
		Thus, $ d_2 = 7 $, and the polynomial $ f $ is of the form $ f(x,y) = $
		\begin{gather*}
			= a_{40}x^4+a_{31}x^3y+a_{22}x^2y^2+a_{13}xy^3+a_{04}y^4+\\
			+a_{21}x^2y+a_{12}xy^2 .
		\end{gather*}

		\item Suppose additionally that $ l'' = \{ x - y = 0 \}$.
		Then $ (p, l'') $ is also a tacnode of $ Q $ if and only if additionally $ a_{21} = - a_{12}$.
		Thus, $ d_3 = 6 $, and the polynomial $ f $ is of the form $ f(x,y) = $
		\begin{gather*}
			= a_{40}x^4+a_{31}x^3y+a_{22}x^2y^2+a_{13}xy^3+a_{04}y^4+\\
			+a_{21}xy(x-y) .
		\end{gather*}

		\item Suppose $ f $ is as in (1).
		Assume that every point on the line $ l $ is a tacnode with a direction $ l $.
		Then for any translation $ x \mapsto x + t$, the polynomial $ g(x,y) = f(x+t,y) $ must satisfy conditions in \autoref{thm:tacnode_calc}.
		$ 0 = g(0,0) = f(t,0) = a_{40}t^4$ for each $t$ which implies $a_{40} = 0$ (take $t = 1$).
		Then $ g'_x(0,0) = f'_x(t,0) = 0$ as a polynomial in $ t $, and $ g'_y(0,0) = f'_y(t,0) = a_{31}t^3 + a_{21}t^2$ implies $0 = a_{31} = a_{21}$.
		Thus,
		\begin{gather*}
			g_t(x,y) = g(x,y) = f(x+t,y) =  \\[0.4em]
			= a_{22}(x+t)^2y^2+a_{13}(x+t)y^3+a_{04}y^4+\\
			+a_{12}(x+t)y^2+a_{03}y^3+\\
			+a_{02}y^2=\\[0.4em]
			= a_{22}x^2y^2+a_{13}xy^3+a_{04}y^4+\\
			+ (2a_{22}t+a_{12})xy^2+(a_{13}t+a_{03})y^3+\\
			+ (a_{22}t^2+a_{02})y^2 .
		\end{gather*}
		This clearly satisfies conditions in \autoref{thm:tacnode_calc} for each $ t \in \C $.
		Hence $ d_4 = 6 $.

		\item Suppose $ f $ as in (1).
		If every line through $ p $ is a tacnode direction at $ p $, then clearly $ f_2(x,y) = 0 = f_3(x,y)$.
		Therefore, $ d_5 = 5 $, and
		\begin{gather*}
			f(x,y) = f_4(x,y) = a_{40}x^4+a_{31}x^3y+a_{22}x^2y^2+a_{13}xy^3+a_{04}y^4 .
		\end{gather*}

		\item Suppose $ f $ as in (4), $ l' = \{ x = 0 \} $.
		As in (2), $(p, l')$ is a tacnode if additionally $ a_{02} = 0 = a_{03}$.
		Thus, $ d_6 = 4 $, and $ f(x,y) = $
		\begin{gather*}
			= a_{22}x^2y^2+a_{13}xy^3+a_{04}y^4+\\
			+ a_{12}xy^2 .
		\end{gather*}

		\item Suppose $ f $ as in (6).
		If every line through $ p $ is a tacnode direction at $ p $, then $ a_{12} = 0 $.
		Hence $ d_7 = 3 $, and
		\begin{gather*}
			f(x,y) = a_{22}x^2y^2+a_{13}xy^3+a_{04}y^4 .
		\end{gather*}

		\item Suppose $ f $ as in (1), $ q = (0:1:1) $, $ l' = \{ y-z = 0 \} $.
		If $ (q,l') $ is also a tacnode of $ Q $, then by \autoref{thm:tacnode_calc}, $ f(0,1) = a_{04} + a_{03} + a_{02} = 0$, $ f'_x(0,1) = a_{13} + a_{12} = 0 $, and $ f'_y(0,1) = 4a_{04} + 3a_{03} + 2a_{02} = 0 $.
		Therefore, $ a_{12} = - a_{13} $, $ a_{02} = a_{04}$, and $ a_{03} = -2a_{04}$.
		Consider $ g(x,y) = f(x,y+1) = $
		\begin{gather*}
			= a_{40}x^4+a_{31}x^3(y+1)+a_{22}x^2(y+1)^2+a_{13}x(y+1)^3+a_{04}(y+1)^4+\\
			+a_{21}x^2(y+1)-a_{13}x(y+1)^2-2a_{04}(y+1)^3+\\
			+a_{04}(y+1)^2 = \\[0.4em]
			=a_{40}x^4+a_{31}x^3y+a_{22}x^2y^2+a_{13}xy^3+a_{04}y^4+\\
			+a_{31}x^3+(2a_{22}+a_{21})x^2y+2a_{13}xy^2+2a_{04}y^3+\\
			+(a_{22}+a_{21})x^2+a_{13}xy+a_{04}y^2 .
		\end{gather*}
		Again, by \autoref{thm:tacnode_calc}, $a_{22} + a_{21} = 0 = a_{13} = a_{31}$, and therefore $ g(x,y) = $
		\begin{gather*}
			=a_{40}x^4+a_{22}x^2y^2+a_{04}y^4+\\
			+a_{22}x^2y+2a_{04}y^3+\\
			+a_{04}y^2 = \\[0.4em]
			= a_{40}x^4 + a_{22}x^2y(y+1)+a_{04}y^2(y+1)^2 .
		\end{gather*}
		Hence, $ d_{8} = 3 $.

		\item Suppose $ f $ is as in (8).
		If the line $ l''=\{ x = 0 \} $ is a direction for each point on it, then for any $ t \in \C $, $f(0,t) = a_{04}t^2(t+1)^2 \implies a_{04} = 0$.
		Consider $ g(x,y) = f(x,y+t) = $
		\begin{gather*}
			=a_{40}x^4+a_{22}x^2(y+t)^2+\\
			+a_{22}x^2(y+t)= \\[0.4em]
			= a_{40}x^4 + a_{22}x^2(y+t)(y+t+1) = \\[0.4em]
			= a_{40}x^4 + a_{22}x^2y^2 + \\
			a_{22}(2t+1)x^2y + \\
			a_{22}(t+1)x^2 .
		\end{gather*}
		Clearly, the conditions in \autoref{thm:tacnode_calc} hold for each point $ (0,t) \in l'' $.
		Thus, $d_{9} = 2 $.

		\item Suppose $ f $ as in (7).
		Applying the same arguments, but switching the roles of $ x $ and $ y $, we get $ f = a_{22}x^2y^2 $.
		Hence $ d_{10} = 1 $.

		\item Suppose $ f $ as in (7).
		As each point $ (t,0) $ on $ l = \{ y = 0 \}$ is a tacnode with every line through it a direction, from \autoref{thm:tacnode_calc} applied to $ 	g(x,y) = f(x+t,y) = $
		\begin{gather*}
			= a_{22}(x+t)^2y^2+a_{13}(x+t)y^3+a_{04}y^4 = \\[0.4em]
			= a_{22}x^2y^2 + a_{13}xy^3 + a_{04}y^4 + \\
			+ 2a_{22}txy^2 + a_{13}ty^3+\\
			+ a_{22}t^2y^2 ,
		\end{gather*}
		we see that $ a_{22} = 0 = a_{13} $, and therefore $ f(x,y) = a_{04}y^4 $.
		Hence $ d_{11} = 1$.

		\item Consider a configuration in $ X_{12} $ with every tacnode being a pair $(p,T_pC)$ where $C$ is given by $x^2+y^2-z^2$.
		Consider a polynomial $ f(x,y) = F(x,y,1) $ in the affine chart $ U = \{ z = 1 \} $ such that the configuration is equal to $ \Sing(F) $.
		Since every point of $ C $ is on $ Q $, $ f(x,y) $ is divisible by $ c(x,y) = x^2 + y^2 - 1 $.
		Furthermore, $ f(x,y) / c(x,y) $ is also a conic.
		Therefore, it is uniquely defined by any five points on it, and as every point of $ Q $ that lies on $ C $ is singular, it must be a double point, i.e. the intersection of the two conics - the components of $ Q $.
		Hence $ f(x,y) = \lambda \cdot c(x,y)^2 $, so $ d_{12} = 1 $.

		\item Obviously, $ L(\{ \PP T\C P^2 \}) = \{ f(x,y) = 0 \} $ is a one-point space.
		Thus, $ d_{13} = 0 $.

		\item Suppose $ f $ as in (1).
		Let $ q = (1:0:1) $.
		Then $ (q, l) $ is also a tacnode of $ Q $ if and only if \mbox{$ g(x,y) = f(x,y+1) $} satisfies conditions in \autoref{thm:tacnode_calc}.
		In particular, \mbox{$ 0 = g(0,0) = f(1,0) = a_{40} $}, and hence $ 0 = g'_x(0,0) = f'_x(1,0) $ as a polynomial, and \mbox{$ 0 = g'_y(0,0) = f'_y(1,0) = a_{31} + a_{21} $}.
		Thus, $ g(x,y) = f(x+1,y) = $
		\begin{gather*}
			= a_{31}(x+1)^3y+a_{22}(x+1)^2y^2+a_{13}(x+1)y^3+a_{04}y^4-\\
			-a_{31}(x+1)^2y+a_{12}(x+1)y^2+a_{03}y^3+\\
			+a_{02}y^2=\\[0.4em]
			= a_{31}x^3y+a_{22}x^2y^2+a_{13}xy^3+a_{04}y^4+\\
			+2a_{31}x^2y+(a_{12}+2a_{22})xy^2+(a_{13}+a_{03})y^3+\\
			+a_{31}xy+(a_{12}+a_{22}+a_{02})y^2 .
		\end{gather*}
		By \autoref{thm:tacnode_calc}, $ 0 = a_{31} = a_{21} $, and therefore $ f(x,y) = $
		\begin{gather*}
			= a_{22}x^2y^2+a_{13}xy^3+a_{04}y^4+\\
			+a_{12}xy^2+a_{03}y^3+\\
			+a_{02}y^2 .
		\end{gather*}
		Notice also that this is exactly the polynomial we got in (4).
		Hence $ L(\Geom(K)) = L(K) $.

		\item Suppose $ f $ as in (1), $ l' = ( x - z = 0)$, $q = (1:0:1)$.
		Then $ (q,l')$ is a tacnode implies \mbox{$ f(1,0) = a_{40} = 0$}, $ f'_x(1,0) = 4a_{40} = 0$, $ f'_y(1,0) = a_{31} + a_{21} = 0 $.
		Consider \mbox{$ g(x,y) = f(x+1,y) = $}
		\begin{gather*}
			= a_{31}(x+1)^3y+a_{22}(x+1)^2y^2+a_{13}(x+1)y^3+a_{04}y^4-\\
			-a_{31}(x+1)^2y+a_{12}(x+1)y^2+a_{03}y^3+\\
			+a_{02}y^2=\\[0.4em]
			= a_{31}x^3y+a_{22}x^2y^2+a_{13}xy^3+a_{04}y^4+\\
			+2a_{31}x^2y+(a_{12}+2a_{22})xy^2+(a_{13}+a_{03})y^3+\\
			+a_{31}xy+(a_{12}+a_{22}+a_{02})y^2 .
		\end{gather*}
		By \autoref{thm:tacnode_calc}, $ l' $ is a direction implies $ 0 = a_{31} = a_{12} + a_{22} + a_{02} = a_{13} + a_{03} $.
		Hence $ g(x,y) = $
		\begin{gather*}
			= a_{22}x^2y^2+a_{13}xy^3+a_{04}y^4+\\
			+(a_{12}+2a_{22})xy^2 = \\[0.4em]
			= a_{22}x^2y^2+a_{13}xy^3+a_{04}y^4+\\
			+a'_{12}xy^2 ,
		\end{gather*}
		where $ a'_{12} = a_{12}+2a_{22}=a_{22}-a_{02}$.
		Notice that this is the same equation as in (6).
		Hence $ L(\Geom(K)) = L(K) $.

		\item Suppose $ f $ as in (6).
		If the line $ l'' = \{ x-y = 0 \} $ is a tacnode direction at $ p $, then $ a_{12} = 0 $.
		Hence
		\begin{gather*}
			f(x,y) = a_{22}x^2y^2+a_{13}xy^3+a_{04}y^4 .
		\end{gather*}
		Notice that this is the same equation as in (7).
		Hence $ L(\Geom(K)) = L(K) $.

		\item Consider the non-geometric limit $ K = (l \times \{ l \}) \cup (\{ p_+ \} \times \PP(T_{p_+}\C P^2)) \cup (\{ p_- \} \times \PP(T_{p_-}\C P^2)) $ of $ X_{12} $ obtained in the proof of Condition 5, and take $ l = \{ y = 0 \} $, $ p_+ = (0:0:1) $, $ p_- = (1:0:1) $.
		The subconfiguration $ (l \times \{ l \}) \cup (\{ p_+ \} \times \PP(T_{p_+}\C P^2)) $ belongs to $ X_7 $, so $ f $ is as in (7).
		As every line through $ p_- $ is a tacnode direction at $ p_- $, \autoref{thm:tacnode_calc} applied to $ g(x,y) = f(x+1,y) = $
		\begin{gather*}
			= a_{22}(x+1)^2y^2+a_{13}(x+1)y^3+a_{04}y^4 = \\[0.4em]
			= a_{22}x^2y^2 + a_{13}xy^3 + a_{04}y^4 + \\
			+ 2a_{22}xy^2 + a_{13}y^3 + \\
			+ a_{22}y^2
		\end{gather*}
		gives $ 0 = a_{22} = a_{13} $, and therefore $ f(x,y) = a_{04}y^4 $.
		Notice that this is the same equation as in (11).
		Hence $ L(\Geom(K)) = L(K) $ for $ \Geom(K) = \configfatline \in X_{11} $, and the corresponding dimension is $ 1 = d_{11} $.

	\end{enumerate}
\end{proof}

\section{Computation of the first page of the main spectral sequence}
\label{sec:v1_first_page}

Before we proceed, we need to establish some notions.
Here we give a short summary of the tools we use.
More details are given in the following section, but it may also be skipped and consulted when needed while reading the proofs.
Finally, this introduction and \autoref{sec:defs_and_core_tools} essentially summarise \cite[Section 1]{BerGorKon23}, included here for the ease of reading.
We recommend the reader to consult the source for detailed explanations, proofs, constructions, and further useful details.

We will use the \emph{conical resolution} $ \sigma $ of $ \Sigma $.
It is a particular subset of $ \Sigma \times \widetilde{\Lambda} $, where $ \widetilde{\Lambda} $ is the union of coherent simplices.

There is a filtration $ \emptyset = \Phi_0 \subset \Phi_1 \subset \ldots \subset \Phi_{\mathscr{N}} = \widetilde{\Lambda} $ on $ \widetilde{\Lambda} $, defined in \autoref{sec:defs_and_core_tools} from the family $ Y_1, \ldots, Y_{\mathscr{N}} $ that refines the configuration spaces $ X_1, \ldots, X_N $.
It induces a filtration \mbox{$ \emptyset = F_0 \subset F_1 \subset \ldots \subset F_{\mathscr{N}} = \sigma $} via $ F_i = p^{-1}(\Phi_i) $, where $ p $ is the projection as in the diagram below:
\begin{equation*}
	\begin{tikzcd}
		F_i \arrow[r, hook] \arrow[d] & \sigma \arrow[d, "p"] \arrow[r, "\pi"] \arrow[rd, hook] & \Sigma                                            \\
		\Phi_i \arrow[r, hook]        & \widetilde{\Lambda}                                   & \Sigma \times \widetilde{\Lambda} \arrow[u] \arrow[l]
	\end{tikzcd}
\end{equation*}
Moreover, the projection $ \pi \colon \sigma \to \Sigma $ is proper, and induces an isomorphism of Borel--Moore homology groups $ \overline{H}_*(\sigma) \cong \overline{H}_*(\Sigma) $.

The filtrations on $ \sigma $ and $ \widetilde{\Lambda} $ allow us to compute the BM homology of $ \sigma $, which in turn yields the BM homology of $ \Sigma $.
To exploit these filtrations we use a collection of core tools, which we present in the form of lemmas in the next section.

\subsection{Definitions and main tools}\label{sec:defs_and_core_tools}

First, there is a way to go from the configuration spaces $ X_1, \ldots, X_N $ satisfying conditions $ 1$ -- $4 $ and $ 5- $ from \cite[Condition list 1]{BerGorKon23} to some configuration spaces $ Y_1, \ldots, Y_{\mathscr{N}} $ satisfying the following useful properties (\cite[Lemma 5]{BerGorKon23}):
\begin{enumerate}
	\item There is $ \alpha : \{1, \ldots, N\} \to \{1, \ldots, \mathscr{N}\} $ a strictly increasing function;
	\item $ Y_{\alpha(i)} = X_i $ for $ i = 1, \ldots, N $;
	\item $ \bigcup_{j=1}^{\mathscr{N}} Y_j = \bigcup_{i=1}^N \overline{X}_i $;
	\item If $ \alpha(i - 1) < j \leq \alpha(i) $, then the geometrisation of every $ K $ in $ Y_j $ belongs to $ X_i $;
	\item The spaces $ Y_1, \ldots, Y_{\mathscr{N}} $ satisfy conditions $ 1$ -- $4 $ and $ 5+ $ from \cite[Condition list 1]{BerGorKon23}.
\end{enumerate}

\label{pg:35}
Next, we consider $ \mathbb{Y} = \bigcup_{j=1}^{\mathscr N} Y_j $ and its \emph{$k$-th symmetric join} $ \mathbb{Y}^{*k} $ (not to be confused with join $ \mathbb{Y} * \ldots * \mathbb{Y} $ of $ \mathbb{Y} $ with itself $ k $ times).
The definition and details of construction of a symmetric join can be found in \cite[Section 1.2.1]{BerGorKon23}.
It can be intuitively viewed as the union of all $ (k-1) $-simplices $ \Delta \subset \R^{\Omega} $ such that vertices of $ \Delta $ belong to distinct copies of $ \iota(\mathbb{Y}) $, where $ \iota $ is a continuous embedding $ \iota \colon \mathbb{Y} \to \R^{\Omega} $.
In what follows, we consider points of $ \mathbb{Y} $ as both configurations (i.e. compact subspaces of $ P $ ) and vertices of simplices in $ \mathbb{Y}^{*k} $.

We will call a simplex $ \Delta \subset \mathbb{Y}^{*k} $ \emph{coherent} if for any two of its vertices, one is a subset of the other.
Among the vertices of a coherent simplex $ \Delta $ there is a configuration that contains all the other: the \emph{principal vertex} of $ \Delta $.
Let $ \widetilde{\Lambda} $ be the union of all coherent simplices.
\label{pg:36}

For every $ K \in \mathbb{Y} $ we denote by $ \widetilde{\Lambda}(K) $ the union of all coherent simplices with principal vertex $ K $.

For every $ K \in \mathbb{X} \coloneq \bigcup_{i=1}^N \overline{X}_i $ let $ \Lambda(K) $ be the union of all coherent simplices $ \Delta $ such that the principal vertex of $ \Delta $ is $ K $, and for every vertex $ K' \neq K $ of $ \Delta $ we have $ \Geom(K') \neq K $.

For every $ K \in X_i = Y_{\alpha(i)} $, we set
\begin{equation*}
	\partial \Lambda(K) \coloneq \bigcup_{j < i} \bigcup_{\substack{K' \in X_j = Y_{\alpha(j)}, \\  K' \subset K}} \widetilde{\Lambda}(K') .
\end{equation*}
and for $K \in Y_i$, we set
\begin{equation*}
	\partial \widetilde \Lambda(K) \coloneq \bigcup_{j < i} \bigcup_{\substack{K' \in Y_j, \\ K' \subset K}} \widetilde{\Lambda}(K') .
\end{equation*}
We set $ \Phi_0 \coloneq \emptyset $ and $ \Phi_i \coloneq \bigcup_{j \leq i} \bigcup_{K \in Y_j} \widetilde{\Lambda}(K) $ for $ i = 1, \ldots, \mathscr{N} $.
These sets define a filtration on $ \widetilde{\Lambda} $:
\[
	\emptyset = \Phi_0 \subset \Phi_1 \subset \ldots \subset \Phi_{\mathscr{N}} = \widetilde{\Lambda} .
\]
The \emph{conical resolution} $ \sigma $ of $ \Sigma $ is the subset of $ \Sigma \times \widetilde{\Lambda} $ formed by all pairs $ (f, x) $ such that $ f \in \Sigma $, $ x \in \widetilde \Lambda (\Sing(f)) $.

There is a filtration $ \emptyset = F_0 \subset F_1 \subset \ldots \subset F_{\mathscr{N}} = \sigma $, where $ F_i = p^{-1}(\Phi_i) $ and $ p $ is the projection as in the diagram below:
\begin{equation*}
	\begin{tikzcd}
		F_i \arrow[r, hook] \arrow[d] & \sigma \arrow[d, "p"] \arrow[r, "\pi"] \arrow[rd, hook] & \Sigma                                            \\
		\Phi_i \arrow[r, hook]        & \widetilde{\Lambda}                                   & \Sigma \times \widetilde{\Lambda} \arrow[u] \arrow[l]
	\end{tikzcd}
\end{equation*}
Finally, we set $ \Psi_i \coloneq (\Phi_{\alpha(i)} \setminus \Phi_{\alpha(i-1)}) \cap (\bigcup_{K \in X_i} \Lambda(K)) $.

These objects satisfy several important properties, which can be found in \cite{BerGorKon23}.
We focus on a collection of core tools that we will use repeatedly throughout the section.

The results below are stated in \cite{BerGorKon23} under one further standing assumption (\cite[p.~15]{BerGorKon23}): the space $ \mathbb{Y} $ must be homeomorphic to a finite CW-complex.
We verify this for our family.
By property (3) above, $ \mathbb{Y} = \bigcup_{i=1}^{13} \overline X_i $.
The proof of Condition 5 exhibits each $ \overline X_i $ as the image of a continuous map $ \varphi_i \colon M_i \to 2^P $ from a compact complex projective variety: an incidence variety for $ i \neq 12 $, and the variety of complete conics for $ i = 12 $.
Let $ M $ be the disjoint union of the $ M_i $, realised as a compact semialgebraic set by placing real algebraic models of the $ M_i $ (\cite[Section 3.4]{BCR98}) in disjoint affine subspaces, and let $ \varphi \colon M \to 2^P $ be the combined map.
As a continuous surjection from a compact space to a metric space, $ \varphi $ maps $ M $ onto its image $ \mathbb{Y} $ as a topological quotient, so $ \mathbb{Y} \cong M / R $, where the equivalence relation $ R = (\varphi \times \varphi)^{-1}(\Delta_{2^P}) $ is closed in $ M \times M $.
Moreover, $ R $ is semialgebraic: each incidence relation $ \Gamma_i = \{ (m, x) \in M_i \times P \mid x \in \varphi_i(m) \} $ is cut out by finitely many bihomogeneous polynomial equations, hence is real algebraic, and membership of $ (m, m') $ in $ R $ is the first-order condition $ \forall x \in P \colon \left( (m, x) \in \Gamma \Leftrightarrow (m', x) \in \Gamma \right) $ on the disjoint union $ \Gamma \subset M \times P $ of the $ \Gamma_i $, so the Tarski--Seidenberg theorem (\cite[Section 2.2]{BCR98}) applies.
By \cite[Corollary 1.5]{Brumfiel87}, a closed semialgebraic equivalence relation on a closed and bounded affine semialgebraic set admits an affine semialgebraic quotient: there is a semialgebraic map $ f \colon M \to Y $ onto an affine semialgebraic set $ Y $ with $ R = (f \times f)^{-1}(\Delta_Y) $.
(This is the form in which we need \cite[Theorem 1.4]{Brumfiel87}, whose hypothesis that the projection $ R \to M $ be proper is automatic here because $ M $ is compact.)
By \cite[Remark 1.3]{Brumfiel87}, such an $ f $ induces the quotient topology on $ Y $ in the ordinary sense.
An arbitrary, not merely semialgebraic, subset of $ Y $ is closed if and only if its preimage is.
So $ Y \cong M / R \cong \mathbb{Y} $.
Being the continuous image of a compact set, $ Y $ is closed and bounded, and every compact semialgebraic set is homeomorphic to a finite simplicial complex by \cite[Theorem 9.2.1]{BCR98}.
Hence $ \mathbb{Y} $ is homeomorphic to a finite CW-complex, as required.

First, the following result is obtained from \cite[Theorem 4]{BerGorKon23}.
\begin{lemma}\label{lem:tool_vb}
	\begin{enumerate}[(1)]
		\item For $ i = 1, \ldots, N $, there is a $ \C^{d_i} $-vector bundle
		\begin{equation*}
			p \mid_{F_{\alpha (i)} \setminus F_{\alpha (i-1)}} \colon F_{\alpha (i)} \setminus F_{\alpha (i-1)} \to \Phi_{\alpha (i)} \setminus \Phi_{\alpha (i-1)}.
		\end{equation*}
		\item For $ j = 1, \ldots, \mathscr{N} $, there is a $ \C^{d'_j} $-vector bundle
		\begin{equation*}
			p \mid_{F_{j} \setminus F_{j-1}} \colon F_{j} \setminus F_{j-1} \to \Phi_{j} \setminus \Phi_{j-1}.
		\end{equation*}
	\end{enumerate}

\end{lemma}

We will use it in combination with the following fact:
\begin{lemma}\label{lem:tool_vb_homology}
	Let $ E \to B $ be an algebraic vector bundle from the previous lemma.
	Let $ d $ be its complex rank.
	Then
	\begin{equation*}
		\overline H_*(E,\Q) = \Q(d) \otimes \overline H_*(B,\Q),
	\end{equation*}
	where $ \Q(d) $ is in degree $ 2d $.
\end{lemma}

\begin{proof}
	See \cite[Theorem 8]{BerGorKon23} and \cite[Theorem 9]{BerGorKon23}.
\end{proof}

Next, from \cite[Theorem 5]{BerGorKon23}, we get the following two lemmas:
\begin{lemma}\label{lem:tool_fibre_bundle}
	\begin{enumerate}[(1)]
		\item For $ i = 1, \ldots , N $, the space $ \Phi_{\alpha(i)} \setminus \Phi_{\alpha(i-1)} $ is fibred over $ X_i $ with generic fibres $ \widetilde \Lambda(K) \setminus  \partial \widetilde \Lambda(K) $ over the point $ K \in X_i $.
		\item For $ i = 1, \ldots , N $, the space $ \Psi_{i} $ is fibred over $ X_{i} $ with generic fibres $ \Lambda(K) \setminus \partial \Lambda(K) $ over the point $ K \in X_{i} $.
		\item If $ M $ is finitely generated, then the inclusion $ \Psi_i \subset \Phi_{\alpha(i)} \setminus \Phi_{\alpha(i-1)} $ induces an isomorphism on the BM homology groups $ \overline{H}_*(-; M) $ for $ i = 1, \ldots, N $.
	\end{enumerate}
\end{lemma}

\begin{lemma}\label{lem:tool_cone}
	Let $ X $ be a compact CW-complex.
	Then
	\begin{equation*}
		\overline H_*(\Cone(X) \setminus X) = H_*(\Cone(X), X) = \widetilde H_{*-1}(X) .
	\end{equation*}
\end{lemma}
This is a classical fact, and we are going to use it repeatedly for $ X = \partial \Lambda(K) $, where $ K $ is a configuration.
By \cite[Lemma 8]{BerGorKon23} the spaces $ \Lambda(K), \partial \Lambda(K), \widetilde \Lambda(K),  \partial \widetilde \Lambda(K) $ are closed in the compact space $ \mathbb{Y}^{*\mathscr{N}} $, hence compact.
That $ \Lambda(K) $ and $ \widetilde \Lambda(K) $ are cones over $ \partial \Lambda(K) $ and $  \partial \widetilde \Lambda(K) $ is established in the discussion preceding that lemma in \cite[Section 1.2.2]{BerGorKon23}.

\begin{lemma}\label{lem:tool_contr}
	If CW complexes $ A, B $ and $ A \cap B $ are all contractible, then the space $ A \cup B $ is also contractible.
\end{lemma}

Finally, there is the following basic result from \cite[Section 22.1]{fomenko2016homotopical} that we are going to use.
\begin{lemma}\label{lem:tool_ss}
	For the spectral sequence of a homologically simple fibration (e.g. a fibration with simply connected base), the second page is given in terms of base $ B $ and fibre $ F $ homologies by
	\begin{equation*}
		E^2_{pq} = H_p(B; H_q(F;R)) .
	\end{equation*}
\end{lemma}
An analogous statement holds for cohomology and BM homology spectral sequences.

\subsection{First page of the spectral sequence}

\begin{theorem}\label{thm:spectral_groups}
	The first page $E^{1}$ of the spectral sequence is given in \autoref{table:main_sp_seq_p1}.
\end{theorem}

\begin{table}[!htbp]
	\centering
	\caption{$ E^{1} $}
	\label{table:main_sp_seq_p1}
		\small
		\renewcommand{\arraystretch}{1.1}
		\setlength{\tabcolsep}{4pt}
		\begin{tabular}{|r|>{\centering\arraybackslash}m{0.85cm}|*{12}{>{\centering\arraybackslash}m{0.7cm}|}}\hline
			& 1                                           & 2                            & 3 & 4       & 5 & 6 & 7 & 8                           & 9                           & 10                          & 11 & 12 & 13                          \\
			\hline
			23 & $ \Q(12)$                                   &                              &   &         &   &   &   &                             &                             &                             &    &    &                             \\
			\hline
			22 &                                             &                              &   &         &   &   &   &                             &                             &                             &    &    &                             \\
			\hline
			21 & \footnotesize $2\Q(11)$                     &                              &   &         &   &   &   &                             &                             &                             &    &    &                             \\
			\hline
			20 &                                             &                              &   &         &   &   &   &                             &                             &                             &    &    &                             \\
			\hline
			19 & \footnotesize $2\Q(10)$ & $\Q(10)$ &   &         &   &   &   &                             &                             &                             &    &    &                             \\
			\hline
			18 &                                             &                              &   &         &   &   &   &                             &                             &                             &    &    &                             \\
			\hline
			17 & $\Q(9)$                 & $\Q(9)$  &   &         &   &   &   &                             &                             &                             &    &    &                             \\
			\hline
			16 &                                             &                              &   &         &   &   &   &                             &                             &                             &    &    &                             \\
			\hline
			15 &                                             & $\Q(8)$                      &   & $\Q(9)$ &   &   &   &                             &                             &                             &    &    &                             \\
			\hline
			14 &                                             &                              &   &         &   &   &   &                             &                             &                             &    &    &                             \\
			\hline
			13 &                                             &                              &   & $\Q(8)$ &   &   &   &                             &                             &                             &    &    &                             \\
			\hline
			12 &                                             &                              &   &         &   &   &   &                             &                             &                             &    &    &                             \\
			\hline
			11 &                                             &                              &   & $\Q(7)$ &   &   &   &                             &                             &                             &    &    &                             \\
			\hline
			10 &                                             &                              &   &         &   &   &   &                             &                             &                             &    &    &                             \\
			\hline
			9  &                                             &                              &   &         &   &   &   & $\Q(8)$                     &                             &                             &    &    &                             \\
			\hline
			8  &                                             &                              &   &         &   &   &   &                             &                             &                             &    &    &                             \\
			\hline
			7  &                                             &                              &   &         &   &   &   & $\Q(7)$ & $\Q(7)$ &                             &    &    &                             \\
			\hline
			6  &                                             &                              &   &         &   &   &   &                             &                             &                             &    &    &                             \\
			\hline
			5  &                                             &                              &   &         &   &   &   & $\Q(6)$ & $\Q(6)$ &                             &    &    &                             \\
			\hline
			4  &                                             &                              &   &         &   &   &   &                             &                             & $\Q(6)$                     &    &    &                             \\
			\hline
			3  &                                             &                              &   &         &   &   &   &                             & $\Q(5)$                     &                             &    &    &                             \\
			\hline
			2  &                                             &                              &   &         &   &   &   &                             &                             & $\Q(5)$ &    &    &                             \\
			\hline
			1  &                                             &                              &   &         &   &   &   &                             &                             &                             &    &    &                             \\
			\hline
			0  &                                             &                              &   &         &   &   &   &                             &                             & $\Q(4)$ &    &    & $\Q(5)$ \\
			\hline
			-1 &                                             &                              &   &         &   &   &   &                             &                             &                             &    &    &                             \\
			\hline
			-2 &                                             &                              &   &         &   &   &   &                             &                             &                             &    &    & $\Q(4)$ \\
			\hline
			& 1                                           & 2                            & 3 & 4       & 5 & 6 & 7 & 8                           & 9                           & 10                          & 11 & 12 & 13                          \\
			\hline
		\end{tabular}
\end{table}

While reading the proof, the reader may find it helpful to refer to \autoref{table:configs_def} and \autoref{fig:config_relations} to keep track of the configurations and their relations.

\begin{proof}

	{\boldmath $ E^1_{1,*} = \overline H_{*+1}(F_{\alpha(1)}) $.}
	\nopagebreak

	By \autoref{lem:tool_vb}, there is a complex vector bundle of rank $ d_1 = 9 $
	\begin{equation*}
		F_{\alpha(1)} \to \Phi_{\alpha(1)},
	\end{equation*}
	from which, by \autoref{lem:tool_vb_homology}, it follows that
	\begin{equation*}
		\overline H_*(F_{\alpha(1)}) = \overset{18}{\Q(9)} \otimes \overline H_{*}(\Phi_{\alpha(1)}),
	\end{equation*}
	where $ \Q(9) $ is in degree $ 18 $.

	We consider the fibre bundle provided by \autoref{lem:tool_fibre_bundle}:
	\begin{equation*}
		\widetilde \Lambda (K) \setminus  \partial \widetilde \Lambda (K) \to \Phi_{\alpha(1)} \to X_1 \ni K \coloneq \configlonetac .
	\end{equation*}
	Since for $ K \in X_1 $ the set $ \widetilde \Lambda (K) = pt $ consists of a single vertex, its boundary is empty, so \mbox{$ \widetilde \Lambda (K) \setminus  \partial \widetilde \Lambda (K) = \widetilde \Lambda (K) = pt $}.
	Since $ X_{1} = \PP T\C P^2 $, from the Borel--Moore homology spectral sequence associated with the above fibre bundle we get
	\begin{equation*}
		\overline H_*(\Phi_{\alpha(1)}) = \overline H_*(pt) \otimes \overline H_*(\PP T\C P^2) = H_*(\PP T\C P^2).
	\end{equation*}
	The last equality holds because $ \PP T\C P^2 $ is compact, being the projectivisation of a vector bundle over the compact base $ \C P^2 $.

	There is a fibre bundle with a simply connected base
	\begin{equation*}
		\C P^1 \to \PP T\C P^2 \to \C P^2 ,
	\end{equation*}
	which gives rise to a spectral sequence $ e $ with \mbox{$ e^{\infty}_{p,q} = e^2_{p,q} = H_p(\C P^2, H_q(\C P^1, \Q)) $}:
	\begin{equation*}
		{
			\begin{array}{|c|c|c|c|c|c|}
				\hline
				2 & \Q(1) &   & \Q(2) &   & \Q(3) \\ \hline
				1 &       &   &       &   &       \\ \hline
				0 & \Q(0) &   & \Q(1) &   & \Q(2) \\ \hline
				& 0     & 1 & 2     & 3 & 4     \\ \hline
			\end{array}
		}
	\end{equation*}
	Here we used \autoref{lem:tool_ss} and identity
	\begin{equation*}
		H_p(\C P^2, H_{2k}(\C P^1, \Q)) = H_p(\C P^2,\Q(k)) = [
		\overset{0}{\Q(k)},
		\overset{1}{0},
		\overset{2}{\Q(k+1)},
		\overset{3}{0},
		\overset{4}{\Q(k+2)}
		],
	\end{equation*}
	where the numbers on top indicate the degrees of the corresponding vector spaces.

	Since the coefficient ring is a field, the groups $ e^\infty $ are torsion-free, and therefore
	\begin{equation*}
		H_n(\PP T\C P^2) = \bigoplus_{p+q=n} e^\infty_{p,q} .
	\end{equation*}
	Consequently
	\begin{equation*}
		\overline H_n(\Phi_{\alpha(1)}) = H_*(\PP T\C P^2) = [
		\overset{0}{\Q(0)},
		\overset{1}{0},
		\overset{2}{2\Q(1)},
		\overset{3}{0},
		\overset{4}{2\Q(2)},
		\overset{5}{0},
		\overset{6}{\Q(3)}
		] ,
	\end{equation*}
	and by \autoref{lem:tool_vb_homology}
	\begin{equation*}
		\overline H_*(F_{\alpha(1)}) = [
		\overset{18}{\Q(9)},
		\overset{19}{0},
		\overset{20}{2\Q(10)},
		\overset{21}{0},
		\overset{22}{2\Q(11)},
		\overset{23}{0},
		\overset{24}{\Q(12)}
		] ,
	\end{equation*}
	and therefore
	\begin{equation*}
		E^1_{1,*} = [
		\overset{17}{\Q(9)},
		\overset{18}{0},
		\overset{19}{2\Q(10)},
		\overset{20}{0},
		\overset{21}{2\Q(11)},
		\overset{22}{0},
		\overset{23}{\Q(12)}
		].
	\end{equation*}

	{\boldmath $ E^1_{2,*} = \overline H_{*+2}(F_{\alpha(2)}\setminus F_{\alpha(1)}) $.}
	\nopagebreak

	By \autoref{lem:tool_vb}, there is a complex vector bundle of rank $ d_2 = 7 $
	\begin{equation*}
		\C^7 \longrightarrow F_{\alpha(2)} \setminus F_{\alpha(1)} \longrightarrow \Phi_{\alpha(2)} \setminus \Phi_{\alpha(1)} ,
	\end{equation*}
	from which by \autoref{lem:tool_vb_homology} we obtain
	\begin{equation*}
		\overline H_*(F_{\alpha(2)} \setminus F_{\alpha(1)}) = \overset{14}{\Q(7)} \otimes \overline H_{*}(\Phi_{\alpha(2)} \setminus \Phi_{\alpha(1)}).
	\end{equation*}
	By \autoref{lem:tool_fibre_bundle}, there is a fibre bundle
	\begin{equation*}
		\Lambda(K) \setminus \partial \Lambda(K) \to \Psi_2 \to X_{2} = B(\C P^{2\vee}, 2) \ni K \coloneq \configsmallcross ,
	\end{equation*}
	whose generic fibre is an open interval.
	Indeed, for a configuration $ K = \{ (p, l_1) , (p, l_2) \} \in X_{2} $, $ \Lambda(K) $ is a union of two $1$-simplices with vertices $ [ \{ (p, l_1) \} ] $, $ [ \{ (p, l_1) , (p, l_2) \} ] $, and $ [ \{(p, l_2) \} ] $.
	Accordingly, $ \partial \Lambda(K) $ is the union of two disjoint $0$-simplices $ [ \{ (p, l_1) \} ] $ and $ [ \{ (p, l_2) \} ] $.

	If $ E \to B $ is a locally trivial fibre bundle with fibre $ \R^n $, then
	\begin{equation*}
		\overline H_*(E) = \overline H_{*-n}(B; \mathscr{L}),
	\end{equation*}
	where $ \mathscr{L} $ is the local system corresponding to orientation on $ E $.

	In our case we have an $ \R^1 $-bundle $ \Psi_2 \to B(\C P^{2\vee}, 2) $.
	This gives rise to the alternating local system $ \pm \Q $.

	By \cite[Lemma 2B]{Vass99}, $ \overline H_*(B(\C P^n, k), \pm \R) $ is isomorphic to $ H_{*-k(k-1)}(G_k(\C^{n+1}), \R) $, where $ G_k(\C^{n+1}) $ denotes the Grassmannian of $ k $-dimensional subspaces of $ \C^{n+1} $.
	Vassiliev states his homological lemmas with real coefficients, but $ B(\C P^n, k) $ is a smooth complex variety, hence an oriented manifold of the homotopy type of a finite CW-complex.
	Poincar\'e duality therefore identifies $ \overline H_*(B(\C P^n, k), \pm F) $ with $ H^{2nk - *}(B(\C P^n, k), \pm F) $ for any field $ F $, and singular cohomology of a finite CW-complex with coefficients in a local system of $ \Q $-vector spaces commutes with the extension of scalars $ \Q \subset \R $, so the lemma holds with $ \pm \Q $ coefficients as well.
	Since $ \C P^{2\vee} \cong \C P^2 $, it applies here with $ n = k = 2 $, and $ G_2(\C^3) \cong \C P^{2\vee} $ has one-dimensional homology in the degrees $ 0 $, $ 2 $ and $ 4 $.
	The Tate twists are those of the strata used in the proof of that lemma: relative to a fixed flag $ \C P^0 \subset \C P^1 \subset \C P^{2\vee} $, the space $ B(\C P^{2\vee}, 2) $ is the union of three locally closed strata that are complex affine spaces of dimensions $ 1 $, $ 2 $ and $ 3 $, and on each of them the flag distinguishes the two points of a configuration, so the local system $ \pm \Q $ is trivial there.
	The strata therefore contribute $ \Q(1) $, $ \Q(2) $ and $ \Q(3) $ in the degrees $ 2 $, $ 4 $ and $ 6 $, and, these degrees being distinct, the spectral sequence of the stratification degenerates.
	Hence
	\begin{equation*}
		\overline H_k(B(\C P^{2\vee}, 2), \pm \Q) =
		\begin{cases}
			\Q \left( \frac{k}{2} \right) & \text{if $ k = 2, 4, 6 $} \\
			0                             & \text{otherwise}
		\end{cases}
		.
	\end{equation*}
	Combining this with \autoref{lem:tool_fibre_bundle},
	\begin{equation*}
		\overline H_{*}(\Psi_2) = \overline H_{*}(\Phi_{\alpha(2)} \setminus \Phi_{\alpha(1)}) = \overset{1}{\Q(0)} \otimes \overline H_*(B(\C P^{2\vee}, 2), \pm \Q) = [
		\overset{3}{\Q(1)},
		\overset{4}{0},
		\overset{5}{\Q(2)},
		\overset{6}{0},
		\overset{7}{\Q(3)}
		] ,
	\end{equation*}
	where $ \overset{1}{\Q(0)} $ corresponds to the fibre $ \R^1 $.

	Hence by \autoref{lem:tool_vb_homology}
	\begin{equation*}
		\overline H_*(F_{\alpha(2)} \setminus F_{\alpha(1)}) = [
		\overset{17}{\Q(8)},
		\overset{18}{0},
		\overset{19}{\Q(9)},
		\overset{20}{0},
		\overset{21}{\Q(10)}
		],
	\end{equation*}
	and therefore
	\begin{equation*}
		E^1_{2,*} = [
		\overset{15}{\Q(8)},
		\overset{16}{0},
		\overset{17}{\Q(9)},
		\overset{18}{0},
		\overset{19}{\Q(10)}
		].
	\end{equation*}

	{\boldmath $ E^1_{3,*} = \overline H_{*+3}(F_{\alpha(3)}\setminus F_{\alpha(2)}) $.}
	\nopagebreak

	By \autoref{lem:tool_vb}, there is a complex vector bundle of rank $ d_3 = 6 $
	\begin{equation*}
		\C^6 \longrightarrow F_{\alpha(3)} \setminus F_{\alpha(2)} \longrightarrow \Phi_3 \setminus \Phi_{\alpha(2)} ,
	\end{equation*}
	from which by \autoref{lem:tool_vb_homology} we obtain
	\begin{equation*}
		\overline H_*(F_{\alpha(3)} \setminus F_{\alpha(2)}) = \Q(6) \otimes \overline H_{*}(\Phi_{\alpha(3)} \setminus \Phi_{\alpha(2)}).
	\end{equation*}

	For each point $ x \in \Phi_{\alpha(3)} \setminus \Phi_{\alpha(2)} $ there exists a unique point $ p \in \C P^2 $ such that $ x $ lies in the connected component of the principal vertex $ [\{ (p, l_1), (p, l_2), (p, l_3) \}] = [\configtriple] $, where \mbox{$ K = \configtriple \in X_{3} $} is a configuration.

	We consider the map $ \Phi_{\alpha(3)} \setminus \Phi_{\alpha(2)} \to \C P^2 $ which sends $ x $ to the unique $ p $ corresponding to $ [\{ (p, l_1), (p, l_2), (p, l_3) \}] $.

	The fibre over $ p $ is then the union of all connected components for principal vertices $[K] = [\configtriple] $ with varying $ l_i $.
	That is, a union of simplices $ \Delta(\configtriple) $ with the boundary $ \partial \Delta(\configtriple) $ removed (see \mbox{\cite[Theorem 6]{BerGorKon23}}) parameterised by $ B(\C P^1, 3) $:

	\begin{center}
		\begin{tikzpicture}[scale=0.66]
			\draw[fill] (0,1) circle (1.5pt) node[left]{$p$};
			\draw[fill] (-1,0) circle (1.5pt) node[above left]{$p_1$};
			\draw[fill] (0,0) circle (1.5pt) node[below left]{$p_2$};
			\draw[fill] (1,0) circle (1.5pt) node[below left]{$p_3$};
			\draw (-1.5,-0.5) -- (0.5, 1.5);
			\draw (1.5,-0.5) -- (-0.5, 1.5);
			\draw (0,-0.5) -- (0, 1.5);
			\draw (-1.5,0) -- (1.5,0);
		\end{tikzpicture}
	\end{center}

	This defines a fibre bundle.
	The corresponding Borel--Moore homology spectral sequence gives
	\begin{equation*}
		\overline H_* (\Phi_{\alpha(3)} \setminus \Phi_{\alpha(2)}) = \overline H_* (B(\C P^1, 3), \pm \Q) \otimes \overline H_* (\C P^2),
	\end{equation*}
	where the local system comes from choices of orientation on $ \Delta(\configtriple) \setminus \partial \Delta(\configtriple) $: indeed, switching $ l_1 $ and $ l_2 $ in the definition of a connected component preserves the fibre but reverses the orientation.

	By \cite[Lemma 2B]{Vass99} again, the group $ \overline H_* (B(\C P^n, k), \pm \R) $ vanishes whenever $ k > n + 1 $, and the vanishing passes to $ \pm \Q $ coefficients as explained in the computation of $ E^1_{2,*} $ above.
	In particular $ \overline H_* (B(\C P^1, 3), \pm \Q) = 0 $, and therefore $ \overline H_* (\Phi_{\alpha(3)} \setminus \Phi_{\alpha(2)}) = 0 $.
	Hence
	\begin{equation*}
		E^1_{3,*} = 0 .
	\end{equation*}

	{\boldmath $ E^1_{4,*} = \overline H_{4 + *}(F_{\alpha(4)}\setminus F_{\alpha(3)}) $.}
	\nopagebreak

	By \autoref{lem:tool_vb}, there is a complex vector bundle of rank $ d_4 = 6 $
	\begin{equation*}
		\C^6 \longrightarrow F_{\alpha(4)} \setminus F_{\alpha(3)} \longrightarrow \Phi_{\alpha(4)} \setminus \Phi_{\alpha(3)} ,
	\end{equation*}
	from which, by \autoref{lem:tool_vb_homology},it follows that
	\begin{equation*}
		\overline H_*(F_{\alpha(4)} \setminus F_{\alpha(3)}) = \overset{12}{\Q(6)} \otimes \overline H_{*}(\Phi_{\alpha(4)} \setminus \Phi_{\alpha(3)}) .
	\end{equation*}
	By \autoref{lem:tool_fibre_bundle}, we have
	\begin{equation*}
		\Lambda(K) \setminus \partial\Lambda(K) \to \Psi_4 \to X_4 \ni K \coloneq \configline .
	\end{equation*}
	This bundle is isomorphic to the following fibre bundle over a simply connected base:
	\begin{equation*}
		\mathring{\Cone}(\C P^1) = \Cone (\C P^1) \setminus \C P^1 \to \Psi_4 \to \C P^{2\vee} .
	\end{equation*}
	Hence there is a spectral sequence such that
	\begin{equation*}
		e^2_{p,q} = e^\infty_{p,q} = \overline H_p(\C P^{2\vee}, \overline H_q(\mathring{\Cone}(\C P^1))) .
	\end{equation*}
	By \autoref{lem:tool_cone},
	\begin{equation*}
		\overline H_q(\mathring{\Cone}(\C P^1)) = H_q(\Cone (\C P^1), \C P^1) =  \widetilde{H}_{q-1}(\C P^1) =
		\begin{cases}
			\Q(1) & \text{if $ q = 3$} \\
			0     & \text{otherwise}
		\end{cases} .
	\end{equation*}
	Since $\C P^{2\vee} = \C P^2$ is compact, $ \overline H_*(\C P^{2\vee},\Q) = \overline H_*(\C P^2,\Q) = H_*(\C P^2,\Q) $.

	The page $ e^2 = e^\infty $:
	\begin{equation*}
		{
			\begin{array}{|c|c|c|c|c|c|}
				\hline
				3 & \Q(1)   &   & \Q(2)   &   &   \Q(3) \\ \hline
				& 0       & 1 & 2       & 3 & 4       \\ \hline
			\end{array}
		}
	\end{equation*}
	Consequently
	\begin{equation*}
		\overline H_*(\Phi_{\alpha(4)} \setminus \Phi_{\alpha(3)}) = \overline H_*(\Psi_4) = [
		\overset{3}{\Q(1)},
		\overset{4}{0},
		\overset{5}{\Q(2)},
		\overset{6}{0},
		\overset{7}{\Q(3)}
		].
	\end{equation*}
	Hence we obtain
	\begin{equation*}
		\overline H_*(F_{\alpha(4)} \setminus F_{\alpha(3)}) = [
		\overset{15}{\Q(7)},
		\overset{16}{0},
		\overset{17}{\Q(8)},
		\overset{18}{0},
		\overset{19}{\Q(9)}
		],
	\end{equation*}
	and therefore
	\begin{equation*}
		E^1_{4,*} = [
		\overset{11}{\Q(7)},
		\overset{12}{0},
		\overset{13}{\Q(8)},
		\overset{14}{0},
		\overset{15}{\Q(9)}
		].
	\end{equation*}

	{\boldmath $ E^1_{5,*} = \overline H_{5 + *}(F_{\alpha(5)}\setminus F_{\alpha(4)}) $.}
	\nopagebreak

	By \autoref{lem:tool_vb}, there is a complex vector bundle of rank $ d_5 = 5 $
	\begin{equation*}
		\C^5 \longrightarrow F_{\alpha(5)} \setminus F_{\alpha(4)} \longrightarrow \Phi_{\alpha(5)} \setminus \Phi_{\alpha(4)} ,
	\end{equation*}
	which implies that
	\begin{equation*}
		\overline H_*(F_{\alpha(5)} \setminus F_{\alpha(4)}) = \Q(5) \otimes \overline H_{*}(\Phi_{\alpha(5)} \setminus \Phi_{\alpha(4)}) .
	\end{equation*}
	By \autoref{lem:tool_fibre_bundle},
	\begin{equation*}
		\overline H_{*}(\Phi_{\alpha(5)} \setminus \Phi_{\alpha(4)}) = \overline H_*(\widetilde \Lambda(K) \setminus  \partial \widetilde \Lambda(K)) \otimes \overline H_*(X_5) .
	\end{equation*}
	By \autoref{lem:tool_cone},
	\begin{equation*}
		\overline H_*(\widetilde \Lambda(K) \setminus  \partial \widetilde \Lambda(K)) = H_*(\Cone( \partial \widetilde \Lambda(K)),  \partial \widetilde \Lambda(K)) = \widetilde H_{*-1}( \partial \widetilde \Lambda (K)) .
	\end{equation*}
	For $ K \coloneq \configcirc \in X_{5} $, we have $  \partial \widetilde \Lambda (K) = (\C P^1)^{*3} $, as it is a union of all possible $ \configtriple $ going through a fixed point, defined by $ K $.
	By \cite[Lemma 21]{BerGorKon23}, the self-join $ (\C P^1)^{*k} $ has the rational homology of a point for every $ k \geq 2 $.
	See also \cite{Vass99}, where the notation $ S^2 = \C P^1 $ is used.

	Therefore, \begin{equation*}
		E^1_{5,*} = 0 .
	\end{equation*}

	{\boldmath $ E^1_{6,*} = \overline H_{6 + *}(F_{\alpha(6)} \setminus F_{\alpha(5)}) $.}
	\nopagebreak

	By \autoref{lem:tool_vb}, there is a complex vector bundle of rank $ d_6 = 4 $
	\begin{equation*}
		\C^4 \longrightarrow F_{\alpha(6)} \setminus F_{\alpha(5)} \longrightarrow \Phi_{\alpha(6)} \setminus \Phi_{\alpha(5)} ,
	\end{equation*}
	which implies that
	\begin{equation*}
		\overline H_*(F_{\alpha(6)} \setminus F_{\alpha(5)}) = \Q(4) \otimes \overline H_{*}(\Phi_{\alpha(6)} \setminus \Phi_{\alpha(5)}) .
	\end{equation*}
	By \autoref{lem:tool_fibre_bundle},
	\begin{equation*}
		\overline H_{*}(\Phi_{\alpha(6)} \setminus \Phi_{\alpha(5)}) = \overline H_{*}(\Psi_6) =  \overline H_*(\Lambda(K) \setminus \partial \Lambda(K)) \otimes \overline H_*(X_6) .
	\end{equation*}
	By \autoref{lem:tool_cone},
	\begin{equation*}
		\overline H_*(\Lambda(K) \setminus \partial \Lambda(K)) = H_*(\Cone(\partial \Lambda(K)), \partial \Lambda(K)) = \widetilde H_{*-1}(\partial \Lambda (K)) .\end{equation*}
	By definition,
	\begin{equation*}
		\partial \Lambda(K) = \cup \widetilde \Lambda(K') ,
	\end{equation*}
	where the union is taken over all maximal geometric $ K' \subsetneq K $.
	Hence, for $ K \coloneq \configlinetac \in X_{6} $,
	\begin{equation*}
		\partial \Lambda\left( \configlinetac \right) = \widetilde \Lambda\left( \configline \right) \cup \widetilde \Lambda\left( \configsmallcross \right) .
	\end{equation*}
	Let $ A \coloneq \widetilde \Lambda\left( \configline \right) $ and $ B \coloneq  \widetilde \Lambda\left( \configsmallcross \right) $.

	$ A $ is contractible as a cone over the corresponding $  \partial \widetilde \Lambda(K') $.
	The same holds for $ B $.
	By \cite[Lemma 6]{BerGorKon23}, their intersection is $ A \cap B = \widetilde \Lambda\left( \configlonetac \right) $, which is also contractible.
	By \autoref{lem:tool_contr}, their union $ A \cup B $ is then also contractible.
	Therefore, $ \partial \Lambda\left( \configlinetac \right) $ has the homology of a point, and hence
	\begin{equation*}
		E^1_{6,*} = 0.
	\end{equation*}

	{\boldmath $ E^1_{7,*} = \overline H_{7 + *}(F_{\alpha(7)}\setminus F_{\alpha(6)}) $.}
	\nopagebreak

	By \autoref{lem:tool_vb}, there is a complex vector bundle of rank $ d_7 = 3 $
	\begin{equation*}
		\C^3 \longrightarrow F_{\alpha(7)} \setminus F_{\alpha(6)} \longrightarrow \Phi_{\alpha(7)} \setminus \Phi_{\alpha(6)} ,
	\end{equation*}
	which implies that
	\begin{equation*}
		\overline H_*(F_{\alpha(7)}\setminus F_{\alpha(6)}) = \Q(3) \otimes \overline H_{*}(\Phi_{\alpha(7)} \setminus \Phi_{\alpha(6)}) .
	\end{equation*}
	By \autoref{lem:tool_fibre_bundle},
	\begin{equation*}
		\overline H_{*}(\Phi_{\alpha(7)} \setminus \Phi_{\alpha(6)}) = \overline H_{*}(\Psi_7) = \overline H_*(\Lambda(K) \setminus \partial \Lambda(K)) \otimes \overline H_*(X_7) .
	\end{equation*}
	By \autoref{lem:tool_cone},
	\begin{equation*}
		\overline H_*(\Lambda(K) \setminus \partial \Lambda(K)) =  H_*(\Cone(\partial \Lambda(K)), \partial \Lambda(K)) = \widetilde H_{*-1}(\partial \Lambda(K)) .
	\end{equation*}
	By definition,
	\begin{equation*}
		\partial \Lambda(K) = \bigcup \widetilde \Lambda(K') ,
	\end{equation*}
	where the union is taken over all maximal geometric $ K' \subsetneq K $.
	Hence, for $ K \coloneq \configlinecirc \in X_{7} $,
	\begin{equation*}
		\partial \Lambda\left( \configlinecirc \right) = \widetilde \Lambda\left( \configcirc \right) \cup \bigcup \widetilde \Lambda\left( \configlinetac \right) .
	\end{equation*}
	To compute the homology $ H_*(\partial \Lambda\left( \configlinecirc \right)) = \overline H_*(\partial \Lambda\left( \configlinecirc \right)) $, consider the filtration
	\[
	\emptyset \subset f_1 \subset f_1' \subset f_2 = \partial \Lambda\left( \configlinecirc \right) ,
	\]
	where
	\[
	f_1 = \widetilde \Lambda\left( \configcirc \right),\quad
	f_1' = f_1 \cup \widetilde \Lambda\left( \configline \right),\quad
	f_2 = f_1' \cup \bigcup \widetilde \Lambda\left( \configlinetac \right) .
	\]
	As in the computation of $ E^1_{6,*} $ above, we use throughout that $ \widetilde \Lambda(K') \cap \widetilde \Lambda(K'') = \widetilde \Lambda(K' \cap K'') $ by \cite[Lemma 6]{BerGorKon23}.
	Let $ e^* $ denote the associated Borel--Moore homological spectral sequence.
	Its first page $ e^1 $ has columns given by
	\[
	e^1_{1,*} = \overline H_{*+1}(f_1) = \overset{-1}{\Q}=\Q[-1],\quad
	e^1_{2,*} = \overline H_{*+2}(f_1' \setminus f_1),\quad
	e^1_{3,*} = \overline H_{*+3}(f_2 \setminus f_1') .
	\]
	The second column vanishes.
	Indeed, $ \configline \cap \configcirc = \configlonetac $, so that $ \widetilde \Lambda\left( \configline \right) \cap f_1 = \widetilde \Lambda\left( \configlonetac \right) $ is a single point and $ f_1' \setminus f_1 = \widetilde \Lambda\left( \configline \right) \setminus \widetilde \Lambda\left( \configlonetac \right) $.
	From the long exact sequence of the pair $ ( \widetilde \Lambda\left( \configline \right) , \widetilde \Lambda\left( \configlonetac \right) ) $, both terms of which are contractible, we obtain
	\[
	\overline H_*(f_1' \setminus f_1) = H_*( \widetilde \Lambda\left( \configline \right) , \widetilde \Lambda\left( \configlonetac \right) ) = 0 .
	\]
	For the third column, note first that
	\[
	\widetilde \Lambda\left( \configlinetac \right) \cap f_1' = \widetilde \Lambda\left( \configsmallcross \right) \cup \widetilde \Lambda\left( \configline \right) = \partial \Lambda\left( \configlinetac \right) ,
	\]
	and that two distinct configurations $ \configlinetac \subsetneq K $ intersect precisely in $ \configline $.
	The sets $ \widetilde \Lambda\left( \configlinetac \right) \setminus f_1' $ are therefore pairwise disjoint, and they exhaust $ f_2 \setminus f_1' $.
	They are indexed by the choice of the second direction at the point of $ K $, which ranges over $ \C P^1 $ with the direction of the line of $ K $ excluded, and they are permuted transitively by the stabiliser of $ K $.
	Consider therefore the fibre bundle
	\[
	\widetilde \Lambda\left( \configlinetac \right) \setminus \partial \Lambda \left( \configlinetac \right) \longrightarrow f_2 \setminus f_1' \longrightarrow \C^1 = \C P^1 \setminus \C P^0 .
	\]
	From the long exact sequence of the pair $( \widetilde \Lambda\left( \configlinetac \right) , \partial \Lambda \left( \configlinetac \right) )$ we obtain
	\[
	\overline H_*(\widetilde \Lambda\left( \configlinetac \right) \setminus \partial \Lambda \left( \configlinetac \right)) = H_*(\widetilde \Lambda\left( \configlinetac \right) , \partial \Lambda \left( \configlinetac \right)) = 0 ,
	\]
	since $ \widetilde \Lambda\left( \configlinetac \right) $ is a cone, and $ \partial \Lambda\left( \configlinetac \right) $ was shown to be contractible in the computation of $ E^1_{6,*} $ above.
	Therefore, from the Leray spectral sequence associated to the last fibre bundle, $ \overline H_*(f_2 \setminus f_1') = 0 $.
	Hence $ \widetilde H_{*-1}(\partial \Lambda(K)) = 0 $, and consequently
	\begin{equation*}
		E^1_{7,*} = 0 .
	\end{equation*}

	{\boldmath $ E^1_{8,*} = \overline H_{8 + *}(F_{\alpha(8)}\setminus F_{\alpha(7)}) $.}
	\nopagebreak

	By \autoref{lem:tool_vb}, there is a complex vector bundle of rank $ d_{8} = 3 $
	\begin{equation*}
		\C^3 \longrightarrow F_{\alpha(8)} \setminus F_{\alpha(7)} \longrightarrow \Phi_{\alpha(8)} \setminus \Phi_{\alpha(7)} ,
	\end{equation*}
	which implies that
	\begin{equation*}
		\overline H_*(F_{\alpha(8)}\setminus F_{\alpha(7)}) = \overset{6}{\Q(3)} \otimes \overline H_{*}(\Phi_{\alpha(8)} \setminus \Phi_{\alpha(7)}) .
	\end{equation*}
	By \autoref{lem:tool_fibre_bundle} applied to
	\begin{equation*}
		\R = \widetilde \Lambda(K) \setminus  \partial \widetilde \Lambda(K) \to \Phi_{\alpha(8)} \setminus \Phi_{\alpha(7)} \to X_{8} \ni K \coloneq \configparallelpair
	\end{equation*}
	we get
	\begin{equation*}
		\overline H_{*}(\Phi_{\alpha(8)} \setminus \Phi_{\alpha(7)}) = \overline H_{*-1}(X_{8}; \pm \Q) .
	\end{equation*}
	The fibre over a configuration $ K = \{ (p_1, l_1), (p_2, l_2) \} $ is the interval
	\begin{equation*}
		\begin{tikzpicture}
			\draw[fill] (0,0) circle (1.5pt) node[above]{$[ (p_1, l_1), (p_2, l_2) ]$};
			\draw (-2.5,0) -- (2.5,0);
			\draw (-2.5,0) circle (1.5pt) node[above]{$[ (p_1, l_1) ]$};
			\draw (2.5,0) circle (1.5pt) node[above]{$[ (p_2, l_2) ]$};
		\end{tikzpicture}
	\end{equation*}
	with end points $ [(p_i, l_i)] $ removed, and is therefore homotopy equivalent to $ \R $.
	The local system $ \pm \Q $ arises from the fundamental group of $ X_{8} $ reversing the orientation of the interval.

	Consider the fibre bundle
	\begin{equation*}
		\C^2=\C P^2 \setminus \C P^1 \to X_{8} \to B(\C P^2,2) .
	\end{equation*}
	Using the corresponding Leray spectral sequence for Borel--Moore homology, together with the groups $ \overline H_{*}(B(\C P^2, 2); \pm \Q) $ found in the computation of $ E^1_{2,*} $, we obtain
	\begin{gather*}
		\overline H_{*}(X_{8}; \pm \Q) = \overset{4}{\Q(2)} \otimes \overline H_{*}(B(\C P^2,2); \pm \Q) = \overset{4}{\Q(2)} \otimes [
		\overset{2}{\Q(1)},
		\overset{3}{0},
		\overset{4}{\Q(2)},
		\overset{5}{0},
		\overset{6}{\Q(3)}
		] =\\
		= [
		\overset{6}{\Q(3)},
		\overset{7}{0},
		\overset{8}{\Q(4)},
		\overset{9}{0},
		\overset{10}{\Q(5)}
		].
	\end{gather*}
	Hence
	\begin{equation*}
		\overline H_*(\Phi_{\alpha(8)} \setminus \Phi_{\alpha(7)}) = [
		\overset{7}{\Q(3)},
		\overset{8}{0},
		\overset{9}{\Q(4)},
		\overset{10}{0},
		\overset{11}{\Q(5)}
		].
	\end{equation*}
	Consequently
	\begin{equation*}
		\overline H_*(F_{\alpha(8)} \setminus F_{\alpha(7)}) = \overset{6}{\Q(3)} \otimes \overline H_{*}(\Phi_{\alpha(8)} \setminus \Phi_{\alpha(7)}) = [
		\overset{13}{\Q(6)},
		\overset{14}{0},
		\overset{15}{\Q(7)},
		\overset{16}{0},
		\overset{17}{\Q(8)}
		],
	\end{equation*}
	and therefore
	\begin{equation*}
		E^1_{8,*} = [
		\overset{5}{\Q(6)},
		\overset{6}{0},
		\overset{7}{\Q(7)},
		\overset{8}{0},
		\overset{9}{\Q(8)}
		].
	\end{equation*}

	{\boldmath $ E^1_{9,*} = \overline H_{9 + *}(F_{\alpha(9)} \setminus F_{\alpha(8)}) $.}
	\nopagebreak

	By \autoref{lem:tool_vb}, there is a complex vector bundle of rank $ d_{9} = 2 $
	\begin{equation*}
		\C^2 \longrightarrow F_{\alpha(9)} \setminus F_{\alpha(8)} \longrightarrow \Phi_{\alpha(9)} \setminus \Phi_{\alpha(8)} ,
	\end{equation*}
	which implies that
	\begin{equation*}
		\overline H_*(F_{\alpha(9)}\setminus F_{\alpha(8)}) = \Q(2) \otimes \overline H_{*}(\Phi_{\alpha(9)} \setminus \Phi_{\alpha(8)}).
	\end{equation*}
	By \autoref{lem:tool_fibre_bundle},
	\begin{equation*}
		\overline H_{*}(\Phi_{\alpha(9)} \setminus \Phi_{\alpha(8)}) = \overline H_{*}(\Psi_{9}) .
	\end{equation*}
	By the same lemma, there is a fibre bundle
	\begin{equation*}
		\Lambda(K) \setminus \partial\Lambda(K) \to \Psi_{9} \to X_{9} .
	\end{equation*}
	The generic fibre can be described as follows.
	By \autoref{lem:tool_cone},
	\begin{equation*}
		\overline H_*(\Lambda(K) \setminus \partial \Lambda(K)) =  H_*(\Cone(\partial \Lambda(K)), \partial \Lambda(K)) = \widetilde H_{*-1}(\partial \Lambda(K)) .
	\end{equation*}
	By definition,
	\begin{equation*}
		\partial \Lambda(K) = \cup \widetilde \Lambda(K') ,
	\end{equation*}
	where the union is taken over all maximal geometric $ K' \subsetneq K $.
	Hence, for $ K \coloneq \configlinepair \in X_{9} $,
	\begin{equation*}
		\partial \Lambda\left( \configlinepair \right) = \widetilde \Lambda\left( \configparallelpair \right) \cup \widetilde \Lambda\left(\vcenter{\hbox{ \begin{tikzpicture}[scale=0.66]
					\draw[fill] (-0.5,0) circle (1.5pt);
					\draw (-1,0) -- (1,0);
					\draw (-0.5,-0.25) -- (-0.5,0.25);
		\end{tikzpicture} }}\right) \cup \widetilde \Lambda\left(\vcenter{\hbox{ \begin{tikzpicture}[scale=0.66]
					\draw[fill] (0.5,0) circle (1.5pt);
					\draw (-1,0) -- (1,0);
					\draw (0.5,-0.25) -- (0.5,0.25);
		\end{tikzpicture} }}\right) .
	\end{equation*}
	Let $ A \coloneq \widetilde \Lambda\left(\vcenter{\hbox{ \begin{tikzpicture}[scale=0.66]
				\draw[fill] (-0.5,0) circle (1.5pt);
				\draw (-1,0) -- (1,0);
				\draw (-0.5,-0.25) -- (-0.5,0.25);
	\end{tikzpicture} }}\right) $ and $ B \coloneq \widetilde \Lambda\left(\vcenter{\hbox{ \begin{tikzpicture}[scale=0.66]
				\draw[fill] (0.5,0) circle (1.5pt);
				\draw (-1,0) -- (1,0);
				\draw (0.5,-0.25) -- (0.5,0.25);
	\end{tikzpicture} }}\right) $.

	Since both are cones over corresponding $  \partial \widetilde \Lambda (K') $, the spaces $ A $ and $ B $ are contractible.
	By \cite[Lemma 6]{BerGorKon23}, \mbox{$ A \cap B =  \widetilde \Lambda\left( \configline \right) $}, which is contractible for the same reason.
	By \autoref{lem:tool_contr}, their union $ A \cup B $ is then also contractible.

	Contracting the union within $ \partial \Lambda (K) $ yields a space homeomorphic to the closed interval $ \widetilde \Lambda\left( \configparallelpair \right) $ with the boundary points glued together, that is, a $1$-sphere $ S^1 $.
	Thus
	\begin{equation*}
		\overline H_*(\Lambda(K) \setminus \partial \Lambda(K)) = \widetilde H_{*-1}(\partial \Lambda(K)) = \widetilde H_{*-1}(S^1) = [
		\overset{2}{\Q(0)}
		].
	\end{equation*}
	The fundamental group of the base inverts the orientation on the sphere, as it interchanges the two boundary points, and therefore we obtain the equality
	\begin{equation*}
		\overline H_{*}(\Psi_{9}) = \overline H_*(\Lambda(K) \setminus \partial \Lambda(K)) \otimes \overline H_*(X_{9}; \pm \Q) ,
	\end{equation*}
	with non-trivial monodromy giving the local system $ \pm \Q $.

	To compute  $\overline H_*(X_{9}; \pm \Q)$, observe that there is a fibre bundle
	\begin{equation*}
		\C\times \C \to X_{9} \to B(\C P^2, 2) ,
	\end{equation*}
	where each $ \C $ parametrises the choices of directions through one of the points in the two-set \mbox{$ \{ \vcenter{\hbox{ \begin{tikzpicture}[scale=0.66]
				\draw[fill] (-0.25,0) circle (1.5pt);
				\draw[fill] (0.25,0) circle (1.5pt);
	\end{tikzpicture} }} \} = S \in B(\C P^2, 2) $} corresponding to a configuration $ \configlinepair $ from $ X_{9} $.

	Thus, using the Leray BM homology spectral sequence and again the groups of $ B(\C P^2, 2) $ from the computation of $ E^1_{2,*} $,
	\begin{gather*}
		\overline H_*(X_{9}; \pm \Q) = \overline H_*(\C^2) \otimes \overline H_*(B(\C P^2,2); \pm \Q) = \overset{4}{\Q(2)} \otimes [
		\overset{2}{\Q(1)},
		\overset{3}{0},
		\overset{4}{\Q(2)},
		\overset{5}{0},
		\overset{6}{\Q(3)}
		] =\\
		= [
		\overset{6}{\Q(3)},
		\overset{7}{0},
		\overset{8}{\Q(4)},
		\overset{9}{0},
		\overset{10}{\Q(5)}
		].
	\end{gather*}
	Hence
	\begin{equation*}
		\overline H_*(\Phi_{\alpha(9)}\setminus \Phi_{\alpha(8)}) = \overline H_{*}(\Psi_{9}) =  \overset{2}{\Q(0)} \otimes [
		\overset{6}{\Q(3)},
		\overset{7}{0},
		\overset{8}{\Q(4)},
		\overset{9}{0},
		\overset{10}{\Q(5)}
		] = [
		\overset{8}{\Q(3)},
		\overset{9}{0},
		\overset{10}{\Q(4)},
		\overset{11}{0},
		\overset{12}{\Q(5)}
		] .
	\end{equation*}
	Consequently
	\begin{equation*}
		\overline H_*(F_{\alpha(9)}\setminus F_{\alpha(8)}) = \overset{4}{\Q(2)} \otimes [
		\overset{8}{\Q(3)},
		\overset{9}{0},
		\overset{10}{\Q(4)},
		\overset{11}{0},
		\overset{12}{\Q(5)}
		] = [
		\overset{12}{\Q(5)},
		\overset{13}{0},
		\overset{14}{\Q(6)},
		\overset{15}{0},
		\overset{16}{\Q(7)}
		],
	\end{equation*}
	and therefore
	\begin{equation*}
		E^1_{9,*} = [
		\overset{3}{\Q(5)},
		\overset{4}{0},
		\overset{5}{\Q(6)},
		\overset{6}{0},
		\overset{7}{\Q(7)}
		].
	\end{equation*}

	{\boldmath $ E^1_{10,*} = \overline H_{10 + *}(F_{\alpha(10)} \setminus F_{\alpha(9)}) $.}
	\nopagebreak

	By \autoref{lem:tool_vb}, there is a complex vector bundle of rank $ d_{10} = 1 $
	\begin{equation*}
		\C^1 \longrightarrow F_{\alpha(10)} \setminus F_{\alpha(9)} \longrightarrow \Phi_{\alpha(10)} \setminus \Phi_{\alpha(9)} ,
	\end{equation*}
	which implies that
	\begin{equation*}
		\overline H_*(F_{\alpha(10)}\setminus F_{\alpha(9)}) = \overset{2}{\Q(1)} \otimes \overline H_{*}(\Phi_{\alpha(10)} \setminus \Phi_{\alpha(9)}) .
	\end{equation*}
	Similarly to the previous case,
	\begin{equation*}
		\overline H_{*}(\Phi_{\alpha(10)} \setminus \Phi_{\alpha(9)}) = \overline H_{*}(\Psi_{10}) ,
	\end{equation*}
	and there is a fibre bundle
	\begin{equation*}
		\Lambda(K) \setminus \partial\Lambda(K) \to \Psi_{10} \to X_{10} .
	\end{equation*}
	The generic fibre can be described as follows.
	By definition,
	\begin{equation*}
		\partial \Lambda(K) = \cup \widetilde \Lambda(K') ,
	\end{equation*}
	where the union is taken over all maximal geometric $ K' \subsetneq K $.
	Hence, for $ K \coloneq \configcross{x} \in X_{10} $,
	\begin{equation*}
		\partial \Lambda\left( \configcross{x} \right)  = \widetilde \Lambda\left( \configlinecirc \right) \bigcup \widetilde \Lambda\left(\vcenter{\hbox{ \begin{tikzpicture}[scale=0.66]
					\draw (0,0) circle (2.5pt);
					\draw (0,-0.75) -- (0,0.75);
		\end{tikzpicture} }}\right) \bigcup \left( \bigcup_{p,q \neq x} \widetilde \Lambda\left( \configpair[0.6] \right) \right) ,
	\end{equation*}
	where $ p $ and $ q $ run over two different lines defining the configuration $ K $.

	Let $ X \coloneq \partial \Lambda\left( \configcross{x} \right) $, and set $ A \coloneq \widetilde \Lambda\left( \configlinecirc \right) $, $ B \coloneq \widetilde \Lambda\left(\vcenter{\hbox{ \begin{tikzpicture}[scale=0.66]
				\draw (0,0) circle (2.5pt);
				\draw (0,-0.75) -- (0,0.75);
	\end{tikzpicture} }}\right) $ and $ C \coloneq \bigcup\limits_{p,q \neq x} \widetilde \Lambda\left( \configpair[0.6] \right)$.
	Let $ Y \coloneq A \cup B $.

	The spaces $ X $ and $ Y $ are compact, since they are closed \cite[Lemma 8]{BerGorKon23} in the compact space $ \mathbb{Y}^{*\mathscr{N}} $.
	Here we use two facts.
	First, $ \mathbb{Y} $ is compact, as a finite union of compact sets.
	Second, by definition, for any finite $ n $, the $ n^\text{th} $ self-join of a compact space is compact.
	Since $ Y \subset X $ is closed, there is a long exact sequence of BM homology groups
	\begin{equation*}
		\dots \to \overline H_i(Y) \to \overline H_i(X) \to \overline H_i(X\setminus Y) \to \overline H_{i-1}(Y) \to \dots
	\end{equation*}
	that can be rewritten as
	\begin{equation*}
		\dots \to H_i(Y) \to H_i(X) \to \overline H_i(X\setminus Y) \to H_{i-1}(Y) \to \dots \ .
	\end{equation*}
	We have
	$ X \setminus Y = C \setminus Y = \bigcup\limits_{p,q \neq x} \left(\vcenter{\hbox{ 	\begin{tikzpicture}
				\draw[fill] (0,0) circle (1.5pt) node[above]{$[ (p, l_1), (q, l_2) ]$};
				\draw (-2.5,0) -- (2.5,0);
				\draw (-2.5,0) circle (1.5pt) node[above]{$[ (p, l_1) ]$};
				\draw (2.5,0) circle (1.5pt) node[above]{$[ (q, l_2) ]$};
	\end{tikzpicture} }}\right) = \C \times \C \times \R $.
	Therefore
	\begin{equation*}
		\overline H_*(X\setminus Y) = \overline H_*(\C^2 \times \R) = \overset{5}{\Q(2)}
	\end{equation*}
	is concentrated in degree $5$.

	By \cite[Lemma 6]{BerGorKon23}, $ A \cap B = \widetilde \Lambda\left( \configcirc \right) $.
	As $ A $, $ B $ and $ A \cap B $ are contractible, $ Y = A \cup B $ is also contractible.
	Consequently
	\begin{equation*}
		\overline H_*(\Lambda(K) \setminus \partial \Lambda(K)) = \widetilde H_{*-1}(\partial \Lambda(K)) = \widetilde H_{*-1}(X) = \overline H_{*-1}(X\setminus Y) = \overset{6}{\Q(2)}
	\end{equation*}
	is concentrated in degree $6$.

	As before, the fundamental group of the base acts by switching the lines $ l_1 $ and $ l_2 $ with each other, therefore reversing the orientation of each interval indexed by $ p, q \neq x $ as above.
	This yields the local system $ \pm \Q $ in the equation below.
	\begin{equation*}
		\overline H_{*}(\Psi_{10}) = \overline H_*(\Lambda(K) \setminus \partial \Lambda(K)) \otimes \overline H_*(X_{10}; \pm \Q) .
	\end{equation*}
	Since $ X_{10} = B(\C P^{2\vee}, 2) $, we obtain
	\begin{equation*}
		\overline H_{*}(\Phi_{\alpha(10)} \setminus \Phi_{\alpha(9)}) = \overline H_{*}(\Psi_{10}) = \overset{6}{\Q(2)} \otimes [
		\overset{2}{\Q(1)},
		\overset{3}{0},
		\overset{4}{\Q(2)},
		\overset{5}{0},
		\overset{6}{\Q(3)}
		] = [
		\overset{8}{\Q(3)},
		\overset{9}{0},
		\overset{10}{\Q(4)},
		\overset{11}{0},
		\overset{12}{\Q(5)}
		] .
	\end{equation*}
	Consequently
	\begin{equation*}
		\overline H_*(F_{\alpha(10)}\setminus F_{\alpha(9)}) = \overset{2}{\Q(1)} \otimes \overline H_{*}(\Phi_{\alpha(10)} \setminus \Phi_{\alpha(9)}) = [
		\overset{10}{\Q(4)},
		\overset{11}{0},
		\overset{12}{\Q(5)},
		\overset{13}{0},
		\overset{14}{\Q(6)}
		],
	\end{equation*}
	and therefore
	\begin{equation*}
		E^1_{10,*} = [
		\overset{0}{\Q(4)},
		\overset{1}{0},
		\overset{2}{\Q(5)},
		\overset{3}{0},
		\overset{4}{\Q(6)}
		] .
	\end{equation*}

	{\boldmath $ E^1_{11,*} = \overline H_{11 + *}(F_{\alpha(11)} \setminus F_{\alpha(10)}) $.}
	\nopagebreak

	By \autoref{lem:tool_vb}, there is a complex vector bundle of rank $ d_{11} = 1 $
	\begin{equation*}
		\C^1 \longrightarrow F_{\alpha(11)} \setminus F_{\alpha(10)} \longrightarrow \Phi_{\alpha(11)} \setminus \Phi_{\alpha(10)} ,
	\end{equation*}
	which implies that
	\begin{equation*}
		\overline H_*(F_{\alpha(11)}\setminus F_{\alpha(10)}) = \Q(1) \otimes \overline H_{*}(\Phi_{\alpha(11)} \setminus \Phi_{\alpha(10)}) .
	\end{equation*}
	As before,
	\begin{equation*}
		\overline H_{*}(\Phi_{\alpha(11)} \setminus \Phi_{\alpha(10)}) = \widetilde H_{*-1}(\partial \Lambda(K)) \otimes \overline H_*(X_{11}) .
	\end{equation*}
	Consider $ K \coloneq \configfatline \in X_{11} $.
	To compute the homology of $\partial \Lambda(K)$ we introduce a filtration reflecting the different types of maximal subconfigurations in $ K $.

	Set $ f_0 \coloneq \widetilde \Lambda ( \configline ) $, where $ \configline \subset K $ is the unique line corresponding to the configuration $ K $.
	Let $ f_1 \coloneq \bigcup\limits_{p} A_p $, where $ A_p \coloneq \widetilde \Lambda\left(\vcenter{\hbox{ \begin{tikzpicture}[scale=0.66]
				\draw (0,0) circle (2.5pt) node[above]{$p$} ;
				\draw (-1,0) -- (1,0);
	\end{tikzpicture} }}\right) $ and $ p $ runs over all the points on the line defining $ K $.
	Denote $ f_2' \coloneq \bigcup \left( \widetilde \Lambda\left( \configlinepair \right) \setminus \left( \widetilde \Lambda\left(\vcenter{\hbox{ \begin{tikzpicture}[scale=0.66]
				\draw[fill] (-0.5,0) circle (1.5pt);
				\draw (-1,0) -- (1,0);
				\draw (-0.5,-0.25) -- (-0.5,0.25);
	\end{tikzpicture} }}\right) \cup \widetilde \Lambda\left(\vcenter{\hbox{ \begin{tikzpicture}[scale=0.66]
				\draw[fill] (0.5,0) circle (1.5pt);
				\draw (-1,0) -- (1,0);
				\draw (0.5,-0.25) -- (0.5,0.25);
	\end{tikzpicture} }}\right) \right) \right) $, where the union is over $ \configpair[0.45] \subset K $.
	For such pair of tacnodes, let us denote $ L_{p,q} \coloneq \configlinepair , L_{x} \coloneq \configlinetac $, where $ x \in \{ p, q \} $.

	Consider the filtration
	\begin{equation*}
		f_0 \subset f_1 \subset f_2 \coloneq f_1 \sqcup f_2' = \partial \Lambda(K) .
	\end{equation*}
	We use the corresponding Borel--Moore homology spectral sequence converging to the BM homology of $ \partial \Lambda(K) $.
	We compute its first page.

	\begin{enumerate}[leftmargin=*, label=(\arabic*)]
		\item The first column contains the Borel--Moore homology $ \overline H_* ( f_0 ) = H_* ( f_0 ) = H_* ( pt ) $.

		\item The second column contains the Borel--Moore homology of $ f_1 \setminus f_0 = \bigcup\limits_{p} \left( A_p \setminus f_0 \right) $.

		Using the long exact sequence of BM homology groups
		\begin{equation*}
			\dots \to \overline H_i(f_0) \to \overline H_i(A_p) \to \overline H_i(A_p \setminus f_0) \to \overline H_{i-1}(f_0) \to \dots \ ,
		\end{equation*}
		one checks that $ \overline H_*(A_p \setminus f_0) = 0 $.
		Considering the Borel--Moore homology spectral sequence of the fibre bundle
		\begin{equation*}
			A_p \setminus f_0 \to \bigcup\limits_{p} \left( A_p \setminus f_0 \right) \to \C P^1
		\end{equation*}
		we see that the second column
		\begin{equation*}
			\overline H_*(f_1 \setminus f_0 ) = \overline H_*\left( \bigcup\limits_{p} \left( A_p \setminus f_0 \right) \right) = \overline H_*\left( A_p \setminus f_0 \right) \otimes \overline H_*( \C P^1 ) = 0
		\end{equation*}
		is zero.

		\item The third column contains the Borel--Moore homology of $ f_2 \setminus f_1 = f_2' $.

		From the computation of $ E^1_{9,*} $ we already know that $ \widetilde \Lambda\left( L_{p} \right) \cup \widetilde \Lambda\left( L_{q} \right) $ is homotopy equivalent to a point.
		Exactly as before, the long exact sequence of the pair implies that
		\begin{equation*}
			\overline H_*\left( \widetilde \Lambda\left( L_{p,q} \right) \setminus \left( \widetilde \Lambda\left( L_{p} \right) \cup \widetilde \Lambda\left( L_{q} \right) \right) \right) = 0 ,
		\end{equation*}
		and, as before, the spectral sequence of the fibre bundle
		\begin{equation*}
			\widetilde \Lambda\left( L_{p,q} \right) \setminus \left( \widetilde \Lambda\left( L_{p} \right) \cup \widetilde \Lambda\left( L_{q} \right) \right) \to f_2' \to B(\C P^1, 2) \times \C^2 ,
		\end{equation*}
		where $ f_2 \setminus f_1 = f_2' = \bigcup \left( \widetilde \Lambda\left( L_{p,q} \right) \setminus \left( \widetilde \Lambda\left( L_{p} \right) \cup \widetilde \Lambda\left( L_{q} \right) \right) \right) $, implies that the second column
		\begin{equation*}
			\overline H_*(f_2 \setminus f_1 ) = \overline H_*\left( \widetilde \Lambda\left( L_{p,q} \right) \setminus \left( \widetilde \Lambda\left( L_{p} \right) \cup \widetilde \Lambda\left( L_{q} \right) \right) \right) \otimes \overline H_*\left( B(\C P^1, 2) \times \C^2 \right) = 0
		\end{equation*}
		is zero.
	\end{enumerate}

	Thus the spectral sequence has only one non-zero group on the first page, located at the intersection of the $ \nth{0} $ row and $ \nth{0} $ column.
	Therefore, its limit Borel--Moore homology \mbox{$ \overline H_{*}(\partial \Lambda(K)) = H_{*}(\partial \Lambda(K)) = \Q $} is concentrated in degree $ 0 $, and consequently
	\begin{equation*}
		\widetilde H_{*}(\partial \Lambda(K)) = 0 \implies \overline H_{*}(\Phi_{\alpha(11)} \setminus \Phi_{\alpha(10)}) = 0 \implies E^1_{11,*} = 0 .
	\end{equation*}

	{\boldmath $ E^1_{12,*} = \overline H_{12 + *}(F_{\alpha(12)} \setminus F_{\alpha(11)}) $.}
	\nopagebreak

	By \autoref{lem:tool_vb}, there is a complex vector bundle of rank $ d_{12} = 1 $
	\begin{equation*}
		\C^1 \longrightarrow F_{\alpha(12)} \setminus F_{\alpha(11)} \longrightarrow \Phi_{\alpha(12)} \setminus \Phi_{\alpha(11)} ,
	\end{equation*}
	which implies that
	\begin{equation*}
		\overline H_*(F_{\alpha(12)}\setminus F_{\alpha(11)}) = \Q(1) \otimes \overline H_{*}(\Phi_{\alpha(12)} \setminus \Phi_{\alpha(11)}) .
	\end{equation*}
	By \autoref{lem:tool_fibre_bundle},
	\begin{equation*}
		\overline H_{*}(\Phi_{\alpha(12)} \setminus \Phi_{\alpha(11)}) = \overline H_{*}(\Psi_{12}) = \overline H_*(\Lambda(K) \setminus \partial \Lambda(K)) \otimes \overline H_*(X_{12}) .
	\end{equation*}
	By \autoref{lem:tool_cone},
	\begin{equation*}
		\overline H_*(\Lambda(K) \setminus \partial \Lambda(K)) =  H_*(\Cone(\partial \Lambda(K)), \partial \Lambda(K)) = \widetilde H_{*-1}(\partial \Lambda(K)) .
	\end{equation*}
	By definition,
	\begin{equation*}
		\partial \Lambda(K) = \cup \widetilde \Lambda(K') ,
	\end{equation*}
	where the union is taken over all maximal geometric $ K' \subsetneq K $.

	Recall that the configuration $ K \coloneq \configconic[0.45] \in X_{12} $ consists of all $ (p,l) $ such that $ l $ is tangent to a fixed non-singular conic $ C $ at the point $ p $.
	Therefore, the maximal geometric configurations are pairs of such tacnodes.
	We claim that $ \partial \Lambda(K) = (C)^{*2} \cong (\C P^1)^{*2} $.
	For maximal subconfigurations $ K' = \configpair[0.45] \subset K $, the space $ \cup \Lambda (K') $ is a union of closed intervals parametrised by their ends $ [(p,l_p)] $ and $ [(q, l_q)] $.
	Choosing a tacnode $ (p, l_p) $, with $ p \in C $ and $ l_p $ a tangent to $ C $ at $ p $, is the same as choosing a point $ p \in C \cong \C P^1 $.
	Consequently, the union $ \cup \Lambda (K') $ is exactly $ C^{*2} = \bigcup\limits_{p,q \in C}[p,q] \supset \bigcup\limits_{p \in C}[p,p] = C^{*1} $.
	By \cite[Lemma 21]{BerGorKon23}, the self-join $ (\C P^1)^{*k} $ has the rational homology of a point for every $ k \geq 2 $.
	See also \cite{Vass99}.
	Hence
	\begin{equation*}
		E^1_{12,*} = 0.
	\end{equation*}

	{\boldmath $ E^1_{13,*} = \overline H_{13 + *}(F_{\alpha(13)} \setminus F_{\alpha(12)}) $.}
	\nopagebreak

	By \autoref{lem:tool_vb}, there is a complex vector bundle of rank $ d_{13} = 0 $
	\begin{equation*}
		\C^0 \longrightarrow F_{\alpha(13)} \setminus F_{\alpha(12)} \longrightarrow \Phi_{\alpha(13)} \setminus \Phi_{\alpha(12)} ,
	\end{equation*}
	which implies that
	\begin{equation*}
		\overline H_*(F_{\alpha(13)}\setminus F_{\alpha(12)}) = \Q(0) \otimes \overline H_{*}(\Phi_{\alpha(13)} \setminus \Phi_{\alpha(12)}) = \overline H_{*}(\Phi_{\alpha(13)} \setminus \Phi_{\alpha(12)}).
	\end{equation*}
	As before,
	\begin{equation*}
		\overline H_{*}(\Phi_{\alpha(13)} \setminus \Phi_{\alpha(12)}) = \overline H_{*}(\Psi_{13}) = \widetilde H_{*-1}(\partial \Lambda(K)) \otimes \overline H_*(X_{13}).
	\end{equation*}
	Since $ X_{13} = pt = \left\{ \PP T\C P^2 \right\} $, it suffices to compute $ \widetilde H_{*-1}(\partial \Lambda(K)) $.
	By definition,
	\begin{equation*}
		\partial \Lambda(K) = \cup \widetilde \Lambda(K') ,
	\end{equation*}
	where the union is taken over all maximal geometric $ K' \subsetneq K $.
	More explicitly,
	\begin{align*}
		\partial \Lambda(K) & = \bigcup \widetilde \Lambda \left( \configcross[0.55]{} \right) \cup \bigcup \widetilde \Lambda(\configfatline) \cup \bigcup \widetilde \Lambda \left( \configconic[0.45] \right) =                                              \\
		& = \left( \bigcup_{K \in X_{10}} \widetilde \Lambda \left( K \right) \right) \cup \left(\bigcup_{K' \in X_{11}} \widetilde \Lambda \left( K' \right) \right) \cup \left(\bigcup_{K'' \in X_{12}} \widetilde \Lambda(K'') \right) .
	\end{align*}
	We are going to find the homology of $ \partial \Lambda(K) $ via a BM homology spectral sequence associated with the following filtration:
	\begin{align*}
		\cup \widetilde \Lambda(\configcirc)
		&\subset \cup \widetilde \Lambda(\configfatline) \subset \\
		&\subset \cup \widetilde \Lambda(\configfatline) \bigcup \cup \widetilde \Lambda\left(\configcross[0.55]{}\right) \subset \\
		&\subset \cup \widetilde \Lambda(\configfatline) \bigcup \cup \widetilde \Lambda\left(\configcross[0.55]{}\right) \bigcup \cup \widetilde \Lambda\left(\configconic[0.45]\right)
		= \partial \Lambda(K).
	\end{align*}
	As the space $ \partial \Lambda(K) $ is compact, it suffices to find its BM homology $ \overline H_{*}(\partial \Lambda(K)) = H_{*}(\partial \Lambda(K))  $.
	The above filtration yields a spectral sequence $ e^* $.
	Let us compute its first page, column by column.

	\begin{enumerate}[leftmargin=*, label=(\arabic*), start=0]
		\item $ e^1_{0,*} = \overline H_*(\sqcup \widetilde \Lambda(\configcirc)) $

		Consider the fibre bundle
		\begin{equation*}
			\widetilde \Lambda(\configcirc ) \to \sqcup \widetilde \Lambda(\configcirc) \to  \C P^2 .
		\end{equation*}
		By \autoref{lem:tool_ss}, the second page of the corresponding spectral sequence consists of groups $ \overline H_p(\C P^2; \overline H_q(\widetilde \Lambda(\configcirc))) $.
		Observe that $ \overline H_q(\widetilde \Lambda(\configcirc)) = H_q(\widetilde \Lambda(\configcirc)) = H_q(pt) $, so
		\begin{equation*}
			\overline H_{*}(\cup \widetilde \Lambda(\configcirc)) = \overline H_{*}(\C P^2) = H_{*}(\C P^2) = [
			\overset{0}{\Q(0)},
			\overset{1}{0},
			\overset{2}{\Q(1)},
			\overset{3}{0},
			\overset{4}{\Q(2)}
			] .
		\end{equation*}
		\item $ e^1_{1,*} = \overline H_{1+*}(\cup \widetilde \Lambda(\configfatline) \setminus \cup \widetilde \Lambda( \configcirc )) $

		Fix some line $ \configfatline $.
		Let $ X \coloneq \widetilde \Lambda (\configfatline) $ and $ Y \coloneq \cup \widetilde \Lambda( \configcirc ) $, where the union is over the points $ p $ of the chosen line $ l $.
		Two configurations $ \configcirc $ belonging to distinct points are disjoint, and one belonging to a point outside $ l $ is disjoint from $ \configfatline $, so the union is disjoint and equals $ X \cap \cup \widetilde \Lambda(\configcirc) $.
		Sending a point of $ \widetilde \Lambda(\configcirc) $ to the corresponding $ p $ therefore gives a fibre bundle
		\begin{equation*}
			\widetilde \Lambda(\configcirc) \to Y \to l = \C P^1 ,
		\end{equation*}
		whose fibre is a cone and hence contractible.
		So, exactly as in the case of $ e^1_{0,*} $ above, \autoref{lem:tool_ss} gives $ \overline H_*(Y) = \overline H_*(\C P^1) = H_*(\C P^1) $.
		From the long exact sequence of BM homology of the pair $ (X, Y) $
		\begin{equation*}
			\dots \to \overline H_i(Y) \to \overline H_i(X) \to \overline H_i(X\setminus Y) \to \overline H_{i-1}(Y) \to \dots
		\end{equation*}
		we find that $ \overline H_*( \widetilde \Lambda(\configfatline) \setminus \cup \widetilde \Lambda(\configcirc) ) = \overline H_*(X\setminus Y) = [ \overset{3}{\Q(1)} ] $.

		Consider the fibre bundle
		\begin{equation*}
			\widetilde \Lambda(\configfatline) \setminus \cup \widetilde \Lambda(\configcirc) \to \sqcup ( \widetilde \Lambda(\configfatline) \setminus \cup \widetilde \Lambda(\configcirc) ) \to  \C P^{2\vee} .
		\end{equation*}
		Following the same steps as in the case of $ e^1_{0,*} $, we get
		\begin{equation*}
			\overline H_{*}(\sqcup (\widetilde \Lambda (\configfatline) \setminus \cup \widetilde \Lambda(\configcirc))) = [
			\overset{3}{\Q(1)}
			] \otimes \overline H_{*}(\C P^{2\vee}) = [
			\overset{3}{\Q(1)},
			\overset{4}{0},
			\overset{5}{\Q(2)},
			\overset{6}{0},
			\overset{7}{\Q(3)}
			].
		\end{equation*}
		Hence \begin{equation*}
			e^1_{1,*} = \overline H_{1+*}(\sqcup (\widetilde \Lambda (\configfatline) \setminus \cup \widetilde \Lambda(\configcirc))) = [
			\overset{2}{\Q(1)},
			\overset{3}{0},
			\overset{4}{\Q(2)},
			\overset{5}{0},
			\overset{6}{\Q(3)}
			].
		\end{equation*}

		\item $ e^1_{2,*} = \overline H_{2+*} \left( \cup \widetilde \Lambda \left( \configcross[0.55]{} \right) \setminus \cup \widetilde \Lambda(\configfatline) \right) $

		As all maximal configurations $ M \subset \configcross[0.55]{} $ are contained in some $ \configfatline $, we see that
		\begin{equation*}
			\overline H_* \left( \cup \widetilde \Lambda \left( \configcross[0.55]{} \right) \setminus \cup \widetilde \Lambda(\configfatline) \right) = \overline H_* \left( \cup \left(\widetilde \Lambda \left( \configcross[0.55]{} \right) \setminus \cup \widetilde \Lambda(M)\right) \right) ,
		\end{equation*}
		where the union $ \cup \widetilde \Lambda(M) $ is taken over all maximal configurations $ M $ lying inside of $ \configcross[0.55]{} $.

		Let $ X \coloneq \widetilde \Lambda \left( \configcross[0.55]{} \right)$, $ Y \coloneq \cup \widetilde \Lambda(M) = \partial \Lambda \left( \configcross[0.55]{} \right) $ for some fixed configuration $ \configcross[0.55]{} $.
		From computations in the case of $ E^1_{10,*} $, we know that
		\begin{equation*}
			H_*(Y) = \overline H_*(Y) = [
			\overset{0}{\Q(0)},
			\overset{5}{\Q(2)}
			].
		\end{equation*}
		Consider the BM homology long exact sequence of the pair $ (X,Y) = \left( \widetilde \Lambda \left( \configcross[0.55]{} \right) , \cup \widetilde \Lambda(M) \right) $:
		\begin{equation*}
			\dots \to H_i(Y) \to H_i(X) \to \overline H_i(X\setminus Y) \to H_{i-1}(Y) \to \dots \ \ .
		\end{equation*}
		As $ \overline H_*(X) = H_*(X) = H_*(pt) $, from the long exact sequence we get
		\begin{equation*}
			\overline H_* \left( X \setminus Y \right) = [
			\overset{6}{\Q(2)}
			].
		\end{equation*}
		Consider the fibre bundle
		\begin{equation*}
			\widetilde \Lambda \left( \configcross[0.55]{} \right) \setminus \cup \widetilde \Lambda (M) \to \cup \left(\widetilde \Lambda \left( \configcross[0.55]{} \right) \setminus \cup \widetilde \Lambda (M) \right) \to B(\C P^{2\vee}, 2) .
		\end{equation*}
		Following the same arguments as in the case of $ e^1_{0,*} $, and using once more the groups of $ B(\C P^{2\vee}, 2) $ from the computation of $ E^1_{2,*} $, we get
		\begin{gather*}
			\overline H_* \left( \cup \left(\widetilde \Lambda \left( \configcross[0.55]{} \right) \setminus \cup \widetilde \Lambda (M) \right) \right) = \overline H_*\left(\widetilde \Lambda \left( \configcross[0.55]{} \right) \setminus \cup \widetilde \Lambda (M)\right) \otimes \overline H_*(B(\C P^{2\vee}, 2); \pm \Q ) =\\
			= [
			\overset{6}{\Q(2)}
			] \otimes [
			\overset{2}{\Q(1)},
			\overset{3}{0},
			\overset{4}{\Q(2)},
			\overset{5}{0},
			\overset{6}{\Q(3)}
			] = [
			\overset{8}{\Q(3)},
			\overset{10}{\Q(4)},
			\overset{12}{\Q(5)}
			],
		\end{gather*}
		and hence
		\begin{equation*}
			e^1_{2,*} = [
			\overset{6}{\Q(3)},
			\overset{8}{\Q(4)},
			\overset{10}{\Q(5)}
			].
		\end{equation*}

		\item $ e^1_{3,*} = \overline H_{3+*} \left( \partial \Lambda(K) \setminus \left(\cup \widetilde \Lambda \left( \configcross[0.55]{} \right) \bigcup \cup \widetilde \Lambda(\configfatline) \right) \right) $

		We have
		\begin{equation*}
			\overline H_{*}\left(\partial \Lambda(K) \setminus \left(\cup \widetilde \Lambda \left( \configcross[0.55]{} \right) \bigcup \cup \widetilde \Lambda(\configfatline)\right)\right) = \overline H_* \left( \bigsqcup \left(\widetilde \Lambda \left( \configconic[0.45] \right) \setminus \cup \widetilde \Lambda \left( \configpair[0.55] \right) \right) \right) .
		\end{equation*}
		For a fixed conic, let $ X \coloneq \widetilde \Lambda \left( \configconic[0.45] \right) $, and let $ Y \coloneq \cup \widetilde \Lambda \left( \configpair[0.55] \right) = \partial \Lambda \left( \configconic[0.45] \right) $.
		By previous computations, $ \overline H_*(X) = H_*(X) = H_*(pt) = H_*(Y) = \overline H_*(Y) $ (here we use the fact that the $ \nth{2} $ self-join $ Y = (\C P^1)^{*2} $ of a compact space is compact).
		Therefore, from the long exact sequence of BM homology for the pair $ (X, Y) $, we know that $ \overline H_*(X \setminus Y) = 0 $.

		The fibre bundle
		\begin{equation*}
			\widetilde \Lambda \left( \configconic[0.45] \right) \setminus \cup \widetilde \Lambda \left( \configpair[0.55] \right) \to \bigsqcup \left(\widetilde \Lambda \left( \configconic[0.45] \right) \setminus \cup \widetilde \Lambda \left( \configpair[0.55] \right) \right) \to \C P^1 =  \configconic[0.45]
		\end{equation*}
		implies that the BM homology in question vanishes, as the BM homology of the generic fibre is zero:
		\begin{equation*}
			\overline H_* \left( \bigsqcup \left( \widetilde \Lambda \left( \configconic[0.45] \right) \setminus \cup \widetilde \Lambda \left( \configpair[0.55] \right) \right) \right) = \overline H_*(\C P^1) \otimes \overline H_* \left( X \setminus Y \right) = 0 .
		\end{equation*}
		But then also $  e^1_{3,*} = 0 $.
	\end{enumerate}
	The spectral sequence $ e^* $ now computes the cohomology $ E^1_{13,*} $, but there are 3 potentially non-zero differentials.
	We now show that they are all in fact non-trivial.

	First, let us record what $ E^1_{13,*} $ would look like if all three differentials were trivial, so that the entry each of them removes can be read off.
	We would have
	\begin{equation*}
		H_{*}(\partial \Lambda(K)) = \overline H_{*}(\partial \Lambda(K)) = [
		\overset{0}{\Q(0)},
		\overset{2}{\Q(1)?},
		\overset{3}{\Q(1)?},
		\overset{4}{\Q(2)?},
		\overset{5}{\Q(2)?},
		\overset{7}{\Q(3)?},
		\overset{8}{\Q(3)?},
		\overset{10}{\Q(4)},
		\overset{12}{\Q(5)}
		].
	\end{equation*}
	The question marks indicate entries that would vanish if the differentials were non-trivial.
	So,
	\begin{gather*}
		\overline H_{*}(F_{\alpha(13)} \setminus F_{\alpha(12)}) = \overline H_{*}(\Phi_{\alpha(13)} \setminus \Phi_{\alpha(12)}) =\\
		= [
		\overset{3}{\Q(1)?},
		\overset{4}{\Q(1)?},
		\overset{5}{\Q(2)?},
		\overset{6}{\Q(2)?},
		\overset{8}{\Q(3)?},
		\overset{9}{\Q(3)?},
		\overset{11}{\Q(4)},
		\overset{13}{\Q(5)}
		],
	\end{gather*}
	and therefore
	\begin{equation*}
		E^1_{13,*} = \overline H_{13+*}(F_{\alpha(13)} \setminus F_{\alpha(12)}) = [
		\overset{-10}{\Q(1)?},
		\overset{-9}{\Q(1)?},
		\overset{-8}{\Q(2)?},
		\overset{-7}{\Q(2)?},
		\overset{-5}{\Q(3)?},
		\overset{-4}{\Q(3)?},
		\overset{-2}{\Q(4)},
		\overset{0}{\Q(5)}
		].
	\end{equation*}
	The question marks come in three pairs, one for each of the three differentials: the differential out of $ e^1_{1,2} $ accounts for the entries in the degrees $ -10 $ and $ -9 $, the one out of $ e^1_{1,4} $ for those in the degrees $ -8 $ and $ -7 $, and the one out of $ e^1_{2,6} $ for those in the degrees $ -5 $ and $ -4 $.
	Each of the three is a map between two copies of $ \Q $ of the same weight, so it is either zero or an isomorphism, and in the latter case both entries of its pair vanish.
	We therefore rule out the three separately.

	The six entries in question lie on the diagonals $ p + q = 3, 4, 5, 6, 8, 9 $ of the main spectral sequence, one on each.
	By the computations of the columns $ 1, \ldots, 12 $ above, every non-zero entry of $ E^1 $ outside the $ \nth{13} $ column lies on a diagonal $ p+q \geq 10 $.
	Hence none of the six is the target of a differential, since the $ \nth{13} $ column is the last one, and none of them is the source of a non-zero one, since a differential out of $ E^r_{13,q} $ lands in a column $ p < 13 $ on the diagonal $ p + q - 1 \leq 8 $, which is empty there.
	Each of the six would therefore survive to $ E^{\infty} $ and, being alone on its diagonal, would be the whole of $ \overline H_{p+q}(\Sigma) $.
	So a vanishing differential out of $ e^1_{1,2} $ would give $ \overline H_3(\Sigma) \neq 0 $, one out of $ e^1_{1,4} $ would give $ \overline H_5(\Sigma) \neq 0 $, and one out of $ e^1_{2,6} $ would give $ \overline H_8(\Sigma) \neq 0 $.

	Recall that
	\begin{equation*}
		\widetilde H^*(\Pi \setminus \Sigma) = \overline H_{\dim_{\R} \Pi - * - 1}(\Sigma) = \overline H_{29 - *}(\Sigma) ,
	\end{equation*}
	where $ \Pi = \C_4[x_0,x_1,x_2] $ is of (complex) dimension $ \binom{4 + 2}{4} = 15 $.

	From this, the three cases would give respectively $ H^{26}(\Pi \setminus \Sigma) \neq 0 $, $ H^{24}(\Pi \setminus \Sigma) \neq 0 $ and $ H^{21}(\Pi \setminus \Sigma) \neq 0 $.

	We now show that $ H^{\geq 21}(\Pi \setminus \Sigma) = 0 $.
	This contradicts all three at once, and therefore shows that each of the three differentials is non-trivial and that only the groups $ \Q(4) $ and $ \Q(5) $ survive in the $ \nth{13} $ column.

	Let $ \UU \coloneq \Pi \setminus \Sigma $.
	Since $ \UU $ is an open subset of the metric space $ \Pi $, its dimension satisfies $ \dim_{\C} \UU = \dim \UU = \dim \Pi = 15 $: any open ball $ B \subset \Pi $ has dimension $ \dim B = \dim \Pi $.
	But dimension is monotone with respect to inclusions of subsets of a metric space, so
	\begin{equation*}
		15 = \dim B \leq \dim \UU \leq \dim \Pi = 15 \implies \dim \UU = 15 .
	\end{equation*}
	Let $ G \coloneq \GL_3(\C) $.
	Then $ \dim G = 9 $, and therefore $ \dim \UU / G = \dim \UU - \dim G = 6 $.
	From this, $ H^{\geq 12}(\UU / G) = 0 $ (the top cohomology is zero, as $ \UU / G $ is non-compact).
	Let $ {}^{G}E^{}_{*} $ be the spectral sequence associated with the quotient map $ \UU \to \UU / G $.
	By \mbox{\cite[Corollary A.1]{BerGorKon23}}, its second page is given by $ {}^{G}E^{pq}_{*} = H^{p}(\UU/G; H^{q}(G)) $.
	If $ p \geq 12 $ or $ q > 9 $, then \mbox{$ H^{p}(\UU/G; H^{q}(G)) = 0 $} by earlier observations.
	Since it converges to the cohomology of $ \UU $, this remark shows that \mbox{$ 0 = H^{\geq 21}(\UU) = \overline H_{\leq 8}(\Sigma) $}.
	Contradiction.

	It follows that
	\begin{equation*}
		E^1_{13,*} = [
		\overset{-2}{\Q(4)},
		\overset{0}{\Q(5)}
		]. \qedhere
	\end{equation*}
\end{proof}

\section{Differentials in the main spectral sequence: the result and the tools}
\label{sec:main_spse_diffs}

In \autoref{sec:v1_first_page} we proved \autoref{thm:spectral_groups} which yields \autoref{table:main_sp_seq_p1} -- the first page of our main Borel--Moore homology spectral sequence $ E^* $.
We now determine the limiting term $ E^{\infty} $ by computing the differentials in $ E^ * $, and consequently compute the cohomology groups of the space $ \UU = \Pi \setminus \Sigma $ of polynomials corresponding to GIT-stable quartics.

\begin{theorem}\label{thm:ch2_einfty}
	The $E^{\infty}$ page of the spectral sequence is given in \autoref{table:main_sp_seq_pinf}.
\end{theorem}

\begin{table}[!htbp]
	\centering
	\caption{$ E^{\infty} $}
	\label{table:main_sp_seq_pinf}
		\small
		\renewcommand{\arraystretch}{1.1}
		\setlength{\tabcolsep}{4pt}
		\begin{tabular}{|r|>{\centering\arraybackslash}m{0.85cm}|*{12}{>{\centering\arraybackslash}m{0.7cm}|}}\hline
			& 1                       & 2       & 3 & 4       & 5 & 6 & 7 & 8       & 9       & 10      & 11 & 12 & 13 \\
			\hline
			23 & $ \Q(12)$               &         &   &         &   &   &   &         &         &         &    &    &    \\
			\hline
			22 &                         &         &   &         &   &   &   &         &         &         &    &    &    \\
			\hline
			21 & \footnotesize $2\Q(11)$ &         &   &         &   &   &   &         &         &         &    &    &    \\
			\hline
			20 &                         &         &   &         &   &   &   &         &         &         &    &    &    \\
			\hline
			19 & $ \Q(10) $              &         &   &         &   &   &   &         &         &         &    &    &    \\
			\hline
			18 &                         &         &   &         &   &   &   &         &         &         &    &    &    \\
			\hline
			17 &                         &         &   &         &   &   &   &         &         &         &    &    &    \\
			\hline
			16 &                         &         &   &         &   &   &   &         &         &         &    &    &    \\
			\hline
			15 &                         & $\Q(8)$ &   & $\Q(9)$ &   &   &   &         &         &         &    &    &    \\
			\hline
			14 &                         &         &   &         &   &   &   &         &         &         &    &    &    \\
			\hline
			13 &                         &         &   & $\Q(8)$ &   &   &   &         &         &         &    &    &    \\
			\hline
			12 &                         &         &   &         &   &   &   &         &         &         &    &    &    \\
			\hline
			11 &                         &         &   & $\Q(7)$ &   &   &   &         &         &         &    &    &    \\
			\hline
			10 &                         &         &   &         &   &   &   &         &         &         &    &    &    \\
			\hline
			9  &                         &         &   &         &   &   &   & $\Q(8)$ &         &         &    &    &    \\
			\hline
			8  &                         &         &   &         &   &   &   &         &         &         &    &    &    \\
			\hline
			7  &                         &         &   &         &   &   &   &         &         &         &    &    &    \\
			\hline
			6  &                         &         &   &         &   &   &   &         &         &         &    &    &    \\
			\hline
			5  &                         &         &   &         &   &   &   &         &         &         &    &    &    \\
			\hline
			4  &                         &         &   &         &   &   &   &         &         & $\Q(6)$ &    &    &    \\
			\hline
			3  &                         &         &   &         &   &   &   &         & $\Q(5)$ &         &    &    &    \\
			\hline
			2  &                         &         &   &         &   &   &   &         &         &         &    &    &    \\
			\hline
			1  &                         &         &   &         &   &   &   &         &         &         &    &    &    \\
			\hline
			0  &                         &         &   &         &   &   &   &         &         &         &    &    &    \\
			\hline
			-1 &                         &         &   &         &   &   &   &         &         &         &    &    &    \\
			\hline
			-2 &                         &         &   &         &   &   &   &         &         &         &    &    &    \\
			\hline
			& 1                       & 2       & 3 & 4       & 5 & 6 & 7 & 8       & 9       & 10      & 11 & 12 & 13 \\
			\hline
		\end{tabular}
\end{table}

\begin{proof}
	We claim that every differential in $ E^r $ that can be non-trivial is indeed non-trivial.
	These are the differentials from the $\nth{2}$ column to the $\nth{1}$ and from the $\nth{9}$ column to the $\nth{8}$ on $ E^1 $, and from the $\nth{13}$ column to the $\nth{10}$ on $ E^3 $.
	To prove this, we develop a method for computing such differentials in \autoref{sec:hvb} and \autoref{sec:diff_vs_leq}.
	The three claims are established in \autoref{sec:diff_2_to_1}, \autoref{sec:diff_9_to_8} and \autoref{sec:diff_13_to_10} respectively, and the passage from them to \autoref{table:main_sp_seq_pinf} is carried out in \autoref{sec:assembly}.
\end{proof}

\begin{theorem}
	The spectral sequence converges to the Borel--Moore homology of $\Sigma$, whose non-zero groups are listed below.
	\begin{equation*}
		{
			\begin{array}{c|cccccccc}
				n & 12 & 14 & 15 & 17                & 19 & 20 & 22                & 24 \\ \hline
				\overline H_n(\Sigma) & \Q(5) & \Q(6) & \Q(7) & \Q(8)^{\oplus 3} & \Q(9) & \Q(10) & \Q(11)^{\oplus 2} & \Q(12)
			\end{array}
		}
	\end{equation*}
\end{theorem}

\begin{theorem}\label{thm:cohom_of_equations}
	Applying Alexander duality (for $i > 0$) $ H^i(\Pi \setminus \Sigma) = \overline{H}_{2 \dim_{\C}(\Pi) - 1 -i}(\Sigma) $ we obtain the cohomology groups of the space of polynomials corresponding to GIT-stable plane quartics:
	\begin{equation*}
		\resizebox{0.95\textwidth}{!}{$
			\begin{array}{c|ccccccccc}
				n      					  & 0     & 5      & 7                 & 9      & 10     & 12                & 14     & 15     & 17      \\ \hline
				H^n(\Pi \setminus \Sigma) & \Q(0) & \Q(-3) & \Q(-4)^{\oplus 2} & \Q(-5) & \Q(-6) & \Q(-7)^{\oplus 3} & \Q(-8) & \Q(-9) & \Q(-10)
			\end{array}
		$}
	\end{equation*}
\end{theorem}

\subsection{Preliminary definitions and constructions}

We will soon need to compute some Chern classes.
We begin by establishing some preliminary notation.
First, set $ P \coloneq \PP T\C P^2 $ and $ G \coloneq \SL_3 (\C) $.
Let $ B \coloneq  \left\{
\begin{bmatrix}
	a_{11} & a_{12} & a_{13} \\
	0      & a_{22} & a_{23} \\
	0      & 0      & a_{33}
\end{bmatrix}
\in G \right\} $ be the standard Borel subgroup of $ G $.

The standard Torus subgroup $ T \subset G $ can be described as
\begin{equation*}
	T = G \cap \{\text{diagonal matrices}\} = \left\{ T(a,b) \coloneq
	\begin{bmatrix}
		(ab)^{-1} & 0 & 0 \\
		0         & a & 0 \\
		0         & 0 & b
	\end{bmatrix}
	\in G \right\} .
\end{equation*}
The corresponding Lie algebra $ t $ can be described as
\begin{equation*}
	\mathfrak{t}  \coloneq \Lie(T) = \left\{
	\begin{bmatrix}
		-a-b & 0 & 0 \\
		0    & a & 0 \\
		0    & 0 & b
	\end{bmatrix}
	\right\} \subset \mathfrak{gl}_3(\C) = \Mat_{3\times 3}(\C) .
\end{equation*}
The dual space $\mathfrak{t}^*$ can be written as
\begin{equation*}
	\mathfrak{t}^{*} = \langle \varepsilon_0, \varepsilon_1, \varepsilon_2 \rangle  / \langle \varepsilon_0 + \varepsilon_1 + \varepsilon_2  \rangle ,
\end{equation*}
where
\begin{equation*}
	\varepsilon_i \colon \mathfrak{t} \to \C \colon \varepsilon_i\left(
	\begin{bmatrix}
		a_{11} & 0      & 0      \\
		0      & a_{22} & 0      \\
		0      & 0      & a_{33}
	\end{bmatrix} \right) \coloneq a_{ii}.
\end{equation*}
We denote by $c(E) = \sum c_i(E)$ the total Chern class of a complex vector bundle $E$, where $c_i(E)$ is the $i^\text{th}$ Chern class of $E$.

We will need the following two results.

\begin{proposition}\label{thm:ptcp2_cohom_structure}
	Let $ I \coloneq (s_1, s_2, s_3) \subset \Q[\varepsilon_0, \varepsilon_1, \varepsilon_2] $ be the ideal generated by the elementary symmetric polynomials $ s_i $ in the polynomial ring in three variables over $\Q$.
	Then
	\begin{equation*}
		H^*(G / B; \Q) = \Q[\varepsilon_0, \varepsilon_1, \varepsilon_2] / I .
	\end{equation*}
	Moreover, an analogous representation exists for coefficients in $ \Z $ instead of $ \Q $.
\end{proposition}

\begin{proof}
	This is \cite[Lemma 17]{BerGorKon23}.
\end{proof}
We work with this explicit form of the cohomology ring to compute the Chern classes.
Let \mbox{$ Q \coloneq \left\{
\begin{bmatrix}
	a_{11} & a_{12} & a_{13} \\
	0      & a_{22} & a_{23} \\
	0      & a_{32} & a_{33}
\end{bmatrix}
\in G \right\} $}.

\begin{proposition}
	Let $ p \coloneq (1:0:0) $, $ q \coloneq (0:1:0) \in \C P^2 $.
	Let ${ \ l \coloneq \overline{pq} }$ be the projective line containing them.
	Then the following maps are homeomorphisms:
	\begin{enumerate}[(1)]
		\item $ G / Q \xrightarrow{\sim} \C P^2 \colon g\cdot Q \mapsto g \cdot p $;
		\item $ G / B \xrightarrow{\sim} P \colon g\cdot B \mapsto g \cdot (p,l) $.
	\end{enumerate}

\end{proposition}

\begin{proof}
	\autoref{rmk:orbit-stab} below + the Orbit--Stabiliser Theorem.
\end{proof}

\begin{lemma}\label{thm:transitivity}
	There is a transitive action of $ G $ on $ P $.
\end{lemma}

\begin{proof}
	We view $P$ as the flag variety
	\begin{equation*}
		P = \PP T\C P^2 = F(3; 1,2) = \left\{ (\langle v_1 \rangle , \langle v_1, v_2 \rangle ) \in \Gr(1, \C^3) \times \Gr(2, \C^3) = \C P^2 \times \C P^{2\vee} \right\} .
	\end{equation*}
	Define $ g \cdot (\langle v_1 \rangle , \langle v_1, v_2 \rangle ) \coloneq (\langle g v_1 \rangle  , \langle g v_1, g v_2 \rangle ) $.
	Let $ v_3 \in \C^3 $ be such that $ \langle v_1, v_2, v_3 \rangle  = \C^3 $.

	Then there is an element $ g \coloneq [v_1, v_2, v_3] \in \GL_3(\C) $ such that $ g e_i = v_i $.
	Since $ \det(g) \neq 0 $, there exists an element \mbox{$ h \coloneq \left[v_1, v_2, \frac{v_3}{\det(g)}\right] \in \SL_3(\C) $} such that $ h e_i = v_i $ for $ i = 1, 2 $, and $ h e_3 = \frac{1}{\det(g)} v_3 $.
	Therefore,
	\begin{equation*}
		h \cdot (\langle e_1 \rangle  , \langle e_1, e_2 \rangle ) = (\langle v_1 \rangle , \langle v_1, v_2 \rangle ) .
	\end{equation*}
	As the flag spanned by the standard basis can be mapped to any other flag, we are done.
\end{proof}

Observe that in the proof above the points of $ \C P^2 $ determined by $ h e_3 $ and $ v_3 $ coincide.

\begin{corollary}
	The action of $ G $ on $ \C P^2 $ given by $ g \cdot \PP(\langle v \rangle ) = \PP(\langle gv \rangle ) $ is transitive.
\end{corollary}

\begin{remark}\label{rmk:orbit-stab}
	Let $ p \coloneq (1:0:0) $, $ q \coloneq (0:1:0) \in \C P^2 $.
	Let ${ \ l \coloneq \overline{pq} }$ be the projective line containing them.
	Then
	\begin{enumerate}[(i)]
		\item The orbit of $ (p,l) $ is $ G \cdot (p, l) = P $, and the stabiliser is $ G_{(p, l)} = B =  \left\{ \begin{bmatrix}
			* & * & * \\
			0         & * & * \\
			0         & 0 & *
		\end{bmatrix} \in G \right \} $;
		\item $ G \cdot p = \C P^2 $, and $ G_p = Q = \left\{ \begin{bmatrix}
			* & * & * \\
			0         & * & * \\
			0         & * & *
		\end{bmatrix} \in G \right \} $.
	\end{enumerate}
\end{remark}

\begin{remark}
	From the inclusions $ B \subset Q \subset G $ follows that $ Q / B \subset G / B  $.
	Consequently, there is a well-defined surjective map $ G / B \to G / Q \colon gB \mapsto gQ $ such that the diagram below is commutative:
	\begin{equation*}
		\begin{tikzcd}
			G / B \arrow[rrr] \arrow[ddd, "\sim" {rotate=90, anchor=north}] &                                          &                       & G / Q \arrow[ddd, "\sim" {rotate=90, anchor=north}] \\
			& gB \arrow[r, maps to] \arrow[d, maps to] & gQ \arrow[d, maps to] &                           \\
			& {g\cdot (p,l)} \arrow[r, maps to]        & g \cdot p             &                           \\
			P \arrow[rrr]                         &                                          &                       & \C P^2
		\end{tikzcd}
	\end{equation*}

\end{remark}

\begin{remark}\label{rmk:homot-TB}
	Let $ U_3 = \C^3 $ be the subgroup of all unipotent matrices in $ G $.
	Then the map
	\[
	B \to T \times U_3 \colon X \mapsto (\diag X, (\diag X)^{-1} \cdot X)
	\]
	is a homeomorphism, where $ \diag X $ is the diagonal part of the matrix $ X $.
	But then from the spectral sequence of the fibre bundle
	\[
	pt \sim  B / T \to G / T \to G / B ,
	\]
	it follows that $ H^*(G / B) = H^*(G / T) $.
\end{remark}

\subsection{Homogeneous Vector Bundles}\label{sec:hvb}

We will make extensive use of homogeneous vector bundles (HVBs) in what follows.
The reason is that for HVBs one can easily compute their Chern classes in terms of the variables $ \varepsilon_0, \varepsilon_1, \varepsilon_2 $ introduced earlier.

From the homeomorphism $ B = T \times U_3 $, we obtain that every character $ T \to \C^\times $ extends uniquely to a character $ B \to \C^\times $.
Conversely, every character on $ B $ restricts to a character on $ T $.
These two operations are inverse to each other.
Thus we can talk about characters of $ T $ and characters of $ B $ interchangeably.

Let $ R \colon B \to \GL(V) $ be a representation on a finite-dimensional vector space $ V $.
Define
\mbox{$ G \times_B V \coloneq G \times_R V \coloneq (G \times V) / B $},
where $ B $ acts on the product via $ b \cdot (g,v) \coloneq (gb^{-1}, R(b)v) $.
With the obvious projection map, $ G \times_R V $ is a vector bundle over $ G / B $ with generic fibre $ V $.
Vector bundles constructed in this way are \emph{homogeneous vector bundles}.

Let $ G $ act transitively on a space $ X $, so that $ X = G / H $ for a subgroup $ H \subset G $.
The following two definitions are equivalent.

\begin{definition}
	\label{pg:61}
	\begin{enumerate}[(1)]
		\item The vector bundle $ p \colon E \to X $ is a \emph{homogeneous vector bundle (HVB)} if it is isomorphic $ E \simeq G \times_H V $ to a vector bundle constructed as above, where $ V $ is an $ H $-representation.
		\item The vector bundle $ E \to X $ is a \emph{homogeneous vector bundle} if $ G $ acts on $ E $ in a way that is linear on each fibre and agrees with $ E \to X $ (i.e. $ p(g \cdot t) = g \cdot p(t) $ for $ g \in G, t \in E $).
	\end{enumerate}
\end{definition}

We will need the following result later:

\begin{proposition}\label{thm:hvb}
	Let $ E_i \to G / B $ be a HVB, $ i = 1, 2 $.
	Fix $ x_0 \in G / B $, for example $ x_0 = e \cdot B $.

	If the fibres $ E_{1, x_0} \simeq E_{2, x_0} $ are isomorphic as $ G_{x_0} $-representations (in the case $ x_0 = e \cdot B $, we have $ G_{x_0} = B $), then $ E_1 \simeq E_2 $ as HVBs.
\end{proposition}

We now show that if $ \chi : T \to \C^\times $ is a character, then $ \dif \chi $ lies in the $ \Z $-module $ \mathfrak{t}^{*} (\Z) \subset \mathfrak{t}^{*} $.
This allows us to relate linear HVBs to the variables $ \varepsilon_0, \varepsilon_1, \varepsilon_2 $.

\begin{proposition}
	If $ f \colon \C^\times \to \C^\times $ is a holomorphic group homomorphism, then it is of the form $ z \mapsto z^n $ for some $ n\in \Z $.
\end{proposition}

The hypothesis cannot be weakened to continuity: complex conjugation $ z \mapsto \overline z $ is a continuous group homomorphism $ \C^\times \to \C^\times $ and is not of the form $ z \mapsto z^n $.

\begin{proof}
	Since $ f $ is a continuous homomorphism of Lie groups, it is smooth.
	Hence its differential $ f_* \colon \C \to \C $ is a well-defined real-linear map: $ f_*(z) = az + b \overline z $, for some $ a,b \in \C $.
	As $ f $ is holomorphic, $ f_* $ is in fact complex-linear, so $ b = 0 $.
	From the commutativity of the following general diagram
	\begin{equation*}
		\begin{tikzcd}
			G \arrow[r, "f"]                               & G                    \\
			\Lie(G) \arrow[u, "\exp"'] \arrow[r, "f_*"]    & \Lie(G) \arrow[u, "\exp"']
		\end{tikzcd}
	\end{equation*}
	we then obtain $ f(\eu^{z}) = \eu^{f_*(z)} = \eu^{az + b \overline z} $.

	Hence $ 1 = f(1) = f(\eu^{2 \pi i}) = \eu^{2 \pi i (a - b)} $ implies that $ a - b = n \in \Z $.

	Write $ z = r\eu^{i\theta} $, $ r \in \R_{>0} $.
	Then
	\begin{equation*}
		f(z) = f(r\eu^{i\theta}) = f(r)\eu^{i\theta n} = f(\eu^{\ln{r}})\eu^{i\theta n} = \eu^{ n \ln{r}} \eu^{i\theta n} = r^n \eu^{i\theta  n} = z^n .
	\end{equation*}
\end{proof}

\begin{corollary}
	If $ R \colon T \to \C^\times $ is a holomorphic group homomorphism, then it is of the form $ T(a,b) \mapsto a^nb^k $ for some $ n, k \in \Z $, where $ a, b \in \C^\times $.
\end{corollary}

\begin{proof}
	As $ R $ is a homomorphism, from the previous proposition and the commutativity of
	\begin{equation*}
		\begin{tikzcd}
			T \arrow[r]                                                              & \C^{\times} \\
			\C^{\times} \times \C^{\times} \arrow[u, "\sim" {rotate=90, anchor=north, inner sep=.4mm}] \arrow[ru, "proj"'] &
		\end{tikzcd}
	\end{equation*}
	by the universal property of the product, we obtain $ R(T(a,b)) = a^nb^k $.
\end{proof}

\begin{corollary}
	For any representation $ R \colon T \to \C^\times $, its differential $ \dif R \colon \mathfrak{t} \to \C $ is in $ \mathfrak{t}^{*}(\Z) \subset \mathfrak{t}^{*} $.
\end{corollary}

\begin{proof}

	Recall the following diagram:
	\begin{equation*}
		\begin{tikzcd}
			T \arrow[r, "R"]                    & \C^{\times}          \\
			\mathfrak{t} \arrow[u, "\exp"] \arrow[r, "\dif R"] & \C \arrow[u, "\exp"']
		\end{tikzcd}
	\end{equation*}
	Let $ X = \begin{bmatrix}
		-a-b & 0      & 0     \\
		0      & a & 0     \\
		0      & 0      & b
	\end{bmatrix} \in \mathfrak{t}.
	$ Then for a real-valued variable $ s $,
	\begin{equation*}
		\eu^{sX} = T(\eu^{sa}, \eu^{sb}) = \begin{bmatrix}
			\eu^{-s(a+b)} & 0      & 0     \\
			0      & \eu^{sa} & 0     \\
			0      & 0      & \eu^{sb}
		\end{bmatrix} ,
	\end{equation*} and therefore, by above, $ R(\eu^{sX}) = \eu^{s(na + kb)} $.
	But then
	\begin{equation*}
		\dif R(X) = \frac{\dif}{\dif s}R(\eu^{s X})\mid_{s=0} = \frac{\dif}{\dif s}(\eu^{s(na+kb)})\mid_{s=0} = na + kb \implies \dif R = n \varepsilon_1 + k \varepsilon_2 .
	\end{equation*}
	Thus
	\begin{equation*}
		\dif R \in \mathfrak{t}^{*}(\Z) \subset \mathfrak{t}^{*} = \hom_{\C}(t, \C) = \langle  \varepsilon_0, \varepsilon_1, \varepsilon_2 \rangle _{\C-mod} / \langle  \varepsilon_0 + \varepsilon_1 + \varepsilon_2 \rangle _{\C-mod} ,
	\end{equation*}
	where $ \varepsilon_i \colon \begin{bmatrix}
		x_0 & 0      & 0     \\
		0      & x_1 & 0     \\
		0      & 0      & x_2
	\end{bmatrix} = X \mapsto \varepsilon_i(X) = x_i $.
\end{proof}

\begin{corollary}\label{cor:char_diff}
	For a character $ R \colon T \to \C^{\times} \colon T(a,b) \mapsto a^nb^k $, the corresponding differential is \mbox{$ \dif R = n \varepsilon_1 + k \varepsilon_2 $}.
\end{corollary}

\begin{corollary}
	$ \mathfrak{t}^{*} (\Z) = \{ \dif R\ |\ R \in \mathfrak{X}(T) \} $, where $ \mathfrak{X}(T) $ is the set of all characters on $ T $.
\end{corollary}

\begin{proof}
	$ \supset $ : See above.

	$ \subset $ :
	\begin{equation*}
		f \in \mathfrak{t}^{*}(\Z) \implies f = n_0 \varepsilon_0 + n_1 \varepsilon_1 + n_2 \varepsilon_2 = (n_1 - n_0) \varepsilon_1 + (n_2 - n_0) \varepsilon_2 .
	\end{equation*}
	Let $ k_1 \coloneq n_1 - n_0 , k_2 \coloneq n_2 - n_0 $.
	Then $ f = \dif( T(a,b) \mapsto a^{k_1} b^{k_2} ) $.
\end{proof}

Let $ L \to G / B $ be a line bundle obtained from a character $ \chi \colon B \to \C^\times $.
Then the isomorphism in \autoref{thm:ptcp2_cohom_structure} takes
\begin{equation*}
	H^2(G / B; \Q) \ni c_1(L) \mapsto \dif (\chi_T) = (\dif \chi)_{\mathfrak{t}} \in \mathfrak{t}^{*} (\Z) .
\end{equation*}

\begin{proposition}
	Every representation $ T \to \GL(V) $ of the torus decomposes into 1-dimensional representations.
\end{proposition}

\begin{proof}
	\cite[Chapter $3, \text{\textsection} 2$, Theorem $ 3 $]{OniVin11}.
\end{proof}

\begin{theorem}
	Let $ R \colon B \to \GL(V) $ be a representation, and let $ E \to G / B $ be the induced vector bundle.
	Then the restricted representation $R|_T$ decomposes into a direct sum of one-dimensional representations.
\end{theorem}

\begin{proof}
	Follows from the previous proposition and remarks on $ G/B \sim G/T $.
\end{proof}

\begin{corollary}\label{cor:chern_prod}
	If $ \chi_i \colon T \to \C^\times $ are the characters occurring in the decomposition of $ R|_T $ above, then under the isomorphism in \autoref{thm:ptcp2_cohom_structure} ${ c(E) \mapsto \Pi (1 + \dif \chi_i) }$.
\end{corollary}

\begin{proof}
	\mbox{\cite[Lemma 17, part 2]{BerGorKon23}}.
\end{proof}

The last corollary allows us to compute a Chern class of a vector bundle $ E $ by finding a group action, corresponding representation $ R $ and its characters, such that $ E $ is isomorphic to the homogeneous vector bundle constructed from $ R $.
This observation will be used repeatedly in what follows.

\subsection{Differentials via the long exact sequence of a pair } \label{sec:diff_vs_leq}

For details on this section, see \cite[Lemma 18]{BerGorKon23}.

Using long exact sequences of pairs in Borel--Moore homology, we are going to construct a diagram, from which the differentials can be identified.
More precisely, they appear as the connecting morphisms of the long exact sequence of a pair.

Let $ (X, Y) $ be a CW pair.
Let $ E_1 \to X $ and $ E_2 \to Y $ be vector bundles of ranks $ r_1 \leq r_2 $.
\begin{equation*}
	\begin{tikzcd}
		E_1 \arrow[rr, "\C^{r_1}", two heads] &  & X                  \\
		&  &                    \\
		E_2 \arrow[rr, "\C^{r_2}", two heads] &  & Y \arrow[uu, hook]
	\end{tikzcd}
\end{equation*}
Suppose there is an isomorphism of vector bundles over $ Y $:
\[
\iota \colon E_1|_Y \xrightarrow{\sim} E'_2 \subset E_2 .
\]
Then the quotient $E_2/E'_2$ is a vector bundle over $Y$ of rank $r_2-r_1$.
Moreover, one can define $ Z \coloneq E_1 \cup_{\iota} E_2 \supset E_2 $, and $ (Z, E_2) $ is a CW pair.

Consider the following commutative diagram.
\begin{equation*}
	\begin{tikzcd}
		\overline H_*(Z\setminus E_2) \arrow[dd, equal] \arrow[rrrrdd, dashrightarrow, "{\partial_{(Z,E_2)}}"]             &  &                                                                             &  &                                                              \\
		&  &                                                                             &  &                                                              \\
		\overline H_*(E_1\setminus (E_1|_Y)) \arrow[dd, "Th" {inner sep=2mm}, "\sim" {rotate=90, anchor=north, inner sep=.5mm}] \arrow[rr, "{\partial_{(E_1,E_1|_Y)}}"{yshift=1.2ex}] &  & \overline H_{*-1}(E_1 |_Y) \arrow[dd, "Th" {inner sep=2mm}, "\sim" {rotate=90, anchor=north, inner sep=.5mm}] \arrow[rr] &  & \overline H_{*-1}(E_2) \arrow[dd, "Th" {inner sep=2mm}, "\sim" {rotate=90, anchor=north, inner sep=.5mm}] \\
		&  &                                                                             &  &                                                              \\
		\overline H_{*-2r_1}(X\setminus Y) \arrow[rr, "{\partial_{(X,Y)}}"{yshift=1.2ex}]                                               &  & \overline H_{* - 2r_1 - 1}(Y) \arrow[rr, "-\frown e(E_2/E'_2)"{yshift=1.2ex}]                &  & \overline H_{* - 2r_2 - 1}(Y)
	\end{tikzcd}
\end{equation*}
From this diagram, it follows that the connecting morphism $\partial_{(Z,E_2)}$ is non-trivial whenever both the connecting morphism $\partial_{(X,Y)}$ and the multiplication by the Euler class are non-trivial.

From here on all vector bundles are complex, and $ r_1 $ and $ r_2 $ always denote \emph{complex} ranks, as the labels $ \C^{r_1} $ and $ \C^{r_2} $ above indicate.
Accordingly, the Thom isomorphism of a bundle of complex rank $ r $ shifts the degree by $ 2r $, which is why the bottom row of the diagram is indexed by $ * - 2r_1 $, $ * - 2r_1 - 1 $ and $ * - 2r_2 - 1 $, and the Euler class $ e(E_2/E'_2) $ lies in $ H^{2(r_2-r_1)}(Y;\Q) $, so that capping with it lowers the degree by $ 2(r_2 - r_1) $.
The same doubling is in force at every occurrence of these shifts in the two applications of the diagram below.
Note finally that \cite[Lemma 18]{BerGorKon23}, the source of the diagram, is stated for \emph{real} oriented bundles and carries the signs $ (-1)^{r_1} $ and $ (-1)^{r_1(r_2-r_1)} $ on its two squares.
For complex bundles those real ranks are even, both signs are $ +1 $, and we suppress them.

The homology sequence of the pair $ (Z, E_2) $ can be reinterpreted as a spectral sequence associated with the two-term filtration $ \emptyset \subset E_2 \subset Z $.
The connecting morphisms then coincide with the differentials on the first page of this spectral sequence.
We will use this interpretation to compute the differentials in the spectral sequence under consideration.

\section{Differentials from the second to the first column}\label{sec:diff_2_to_1}

This section computes the differentials from the second column to the first.
Only two of them have a non-zero source and a non-zero target, and \autoref{thm:2to1_nontrivial} states that one is an isomorphism and the other has rank $ 1 $.

By \autoref{sec:two_to_one_goal}, these differentials are connecting morphisms of the pair $ (F_{\alpha(2)}, F_{\alpha(1)}) $.
The construction of \autoref{sec:diff_vs_leq} computes such a morphism when the vector bundle over $ \Phi_{\alpha(2)} \setminus \Phi_{\alpha(1)} $ extends to all of $ \Phi_{\alpha(2)} $.
Here it does not extend (\autoref{sec:closure_failure}).
A point of $ \Phi_{\alpha(1)} $ stands for one flag, a nearby point stands for a pair of flags, and the fibre there depends on both flags.
So the fibres have different limits along different ways of approaching the point.
We therefore replace $ \Phi_{\alpha(1)} $ by the space of ordered pairs of flags with a common point, which remembers the second flag.
The result is the \emph{resolved pair} of \autoref{sec:resolved_pair}.
It agrees with the original pair away from $ \Phi_{\alpha(1)} $, and over it the bundle does extend.
The extension uses the limit of the fibre as the two flags merge (\autoref{sec:config_lim}).

For the resolved pair, the construction of \autoref{sec:diff_vs_leq} writes the differential as a composition of three maps: the boundary map of the resolved pair, the cap product with the Euler class of a vector bundle of rank $ 2 $, and the pushforward that forgets the second flag.
We find the image of the boundary map in \autoref{sec:resolved_cohomology}.
In \autoref{sec:map_alpha} we find the Euler class, and a short calculation in a cohomology ring gives the other two maps on that image.
In \autoref{sec:two_to_one_diff} we show that the differential is this composition, and the theorem follows.

\subsection{The differentials to be computed}\label{sec:two_to_one_goal}

This subsection states the result and reduces it to a connecting morphism.
We then describe $ \Phi_{\alpha(2)} $: every $ K \in X_2 $ contributes an arc $ \widetilde \Lambda(K) $ that joins the vertices of its two flags, and a point outside $ \Phi_{\alpha(1)} $ lies on exactly one arc.
Finally, we identify $ d^1 $ with the connecting morphism of the pair $ (F_{\alpha(2)}, F_{\alpha(1)}) $.

Let $ P \coloneq \PP T\C P^2 $, and let $\Pi_{4,2} \coloneq \Pi = \C[x_0,x_1,x_2]_4$.

The indexing function $ \alpha $ can be chosen so that $ \alpha(1) = 1 $ and $ \alpha(2) = 2 $.
We use this throughout this section, although we keep the notation $ \alpha(i) $ for consistency.

By the proof of \autoref{thm:spectral_groups}, the first two columns of the spectral sequence are given by \mbox{$ E^1_{1,q} = \overline H_{q+1}(F_{\alpha(1)}) $} and $ E^1_{2,q} = \overline H_{q+2}(F_{\alpha(2)} \setminus F_{\alpha(1)}) $.
Reading \autoref{table:main_sp_seq_p1}, the only rows $ q $ in which both $ E^1_{2,q} $ and $ E^1_{1,q} $ are non-zero are $ q = 17 $ and $ q = 19 $.
The aim of this section is to prove that the differential $ d^1 \colon E^1_{2,17} \to E^1_{1,17} $ is an isomorphism and that the differential $ d^1 \colon E^1_{2,19} \to E^1_{1,19} $ has rank $ 1 $ (\autoref{thm:2to1_nontrivial}).

To compute these differentials we need an explicit description of the spaces $ \Phi_{\alpha(1)} \subset \Phi_{\alpha(2)} $ and of the fibres of $ F_{\alpha(2)} $ over their points.
We give it now, in the notation of \autoref{sec:defs_and_core_tools}.

Recall that the self-join $ \mathbb{Y}^{*\mathscr{N}} $ is a union of simplices in a Euclidean space $ \R^\Omega $, whose vertices are the images of the configurations under a continuous embedding $ \iota \colon \mathbb{Y} \to \R^\Omega $.
For a configuration $ K \in \mathbb{Y} $ denote by $ \delta_K \coloneq \iota(K) $ the corresponding vertex, and for $ v \in P $ write $ \delta_v \coloneq \delta_{ \left\{ v \right\} } $.
The embedding is generic, so every point of $ \mathbb{Y}^{*\mathscr{N}} $ lies in the interior of exactly one simplex (see \cite[Sections 1.2.1--1.2.2]{BerGorKon23}).

The only coherent simplex with principal vertex $ \delta_v $ is $ \delta_v $ itself.
So $ \widetilde \Lambda(\left\{ v \right\}) = \left\{ \delta_v \right\} $ and \mbox{$ \Phi_{\alpha(1)} = \left\{ \delta_v \mid v \in P \right\} \cong X_1 \cong P $}.
Let $ U \coloneq F_{\alpha(1)} $.
By \autoref{lem:tool_vb}, $ U $ is a complex vector bundle of rank $ 9 $ over $ \Phi_{\alpha(1)} $, and its fibre over $ \delta_v $ is $ L(\left\{ v \right\}) $.

Now let $ K = \left\{ w_1, w_2 \right\} \in X_2 $.
The configurations contained in $ K $ are $ K $ itself and the singletons $ \left\{ w_1 \right\} $ and $ \left\{ w_2 \right\} $, and neither singleton contains the other.
So the coherent simplices with principal vertex $ \delta_K $ are $ \delta_K $ itself and the two edges that join it to $ \delta_{w_1} $ and to $ \delta_{w_2} $, and
\begin{equation*}
	\widetilde \Lambda(K) = \left[ \delta_{w_1}, \delta_K \right] \cup \left[ \delta_K, \delta_{w_2} \right]
\end{equation*}
(compare the description of $ \Lambda(K) $ in the proof of \autoref{thm:spectral_groups}).
This is an arc from $ \delta_{w_1} $ to $ \delta_{w_2} $, made of two edges that meet at $ \delta_K $.

Since $ \alpha(1) = 1 $ and $ \alpha(2) = 2 $, the space $ \Phi_{\alpha(2)} $ is the union of these arcs, and $ \Phi_{\alpha(1)} $ is the set of their endpoints.
So every point $ b \in \Phi_{\alpha(2)} \setminus \Phi_{\alpha(1)} $ lies on an edge $ \left[ \delta_v, \delta_K \right] $ of some arc and differs from its end $ \delta_v $:
\begin{equation}\label{eq:edge_coords}
	b = s \, \delta_v + ( 1 - s ) \, \delta_K , \qquad 0 \leq s < 1 , \quad v \in K \in X_2 .
\end{equation}
For $ s = 0 $ the point $ b $ is the vertex $ \delta_K $, whichever of the two $ v \in K $ is taken.
The configuration $ K $ is determined by $ b $.
Indeed, $ b $ lies in the interior of exactly one simplex, which is the edge $ \left[ \delta_v, \delta_K \right] $ if $ s > 0 $ and the vertex $ \delta_K $ if $ s = 0 $.
In both cases the principal vertex of this simplex is $ \delta_K $.
We denote the configuration $ K $ by $ K_b \in X_2 $.

Again by \autoref{lem:tool_vb}, $ F_{\alpha(2)} \setminus F_{\alpha(1)} $ is a complex vector bundle of rank $ 7 $ over $ \Phi_{\alpha(2)} \setminus \Phi_{\alpha(1)} $, and its fibre over $ b $ is $ L(K_b) $.
Thus over an arc $ \widetilde \Lambda(K) $ the fibre of $ F_{\alpha(2)} $ is the space $ L(K) $ of dimension $ 7 $, except at the two endpoints $ \delta_v $, where it is the larger space $ L(\left\{ v \right\}) $ of dimension $ 9 $.

Consider the fibre bundle $ (\C P^1)^{*2} \to \Phi_{\alpha(2)} \to \C P^2 $.
Here, $ \Phi_{\alpha(2)} = \Phi_{2}$ is the union of all coherent simplices that have a principal vertex in $ X_{2} $.
It decomposes as a disjoint union of subspaces indexed by the points $ p \in \C P^2 $ of intersection of the two directions $ l_1, l_2 $ of $ \left\{ (p,l_1), (p,l_2) \right\} \in X_{2} $.
The projection map sends everything in a connected component to the corresponding point in $ \C P^2 $.
Correspondingly, the preimage of a point $ p \in \C P^2 $ is the said connected component, the union of the arcs $ \widetilde \Lambda(K) $ over all pairs $ K $ of directions through $ p $.
Straightening each two-edge arc onto the segment joining its endpoints identifies this union with the second self-join $ (\C P^1)^{*2} $ of the space of directions through $ p $.

Recall that $ (\C P^1)^{*2} $ has the rational homology of a point (see \mbox{\cite[Lemma 21]{BerGorKon23}}, whose configuration-space input goes back to \cite{Vass99}).
From the spectral sequence of this fibre bundle (noting that the base is simply connected), we obtain that $ \Phi_{\alpha(2)} $ has the same rational homology as $ \C P^2 $.

The differential $ d^1 \colon E^1_{2,q} \to E^1_{1,q} $ coincides with the connecting morphism
\begin{equation*}
	\partial \colon \overline H_{q+2} ( F_{\alpha(2)} \setminus F_{\alpha(1)} ) \longrightarrow \overline H_{q+1} ( F_{\alpha(1)} )
\end{equation*}
of the pair $ (F_{\alpha(2)}, F_{\alpha(1)}) $ in Borel--Moore homology.
Indeed, the spectral sequence of the two-term filtration $ \emptyset \subset F_{\alpha(1)} \subset F_{\alpha(2)} $ has these connecting morphisms as its first-page differentials (see the end of \autoref{sec:diff_vs_leq}), and the inclusion of this filtration into the full filtration $ \left\{ F_{\alpha(i)} \right\} $ induces the identity on the entries $ (1, q) $ and $ (2, q) $ of the first pages, hence intertwines the differentials between them.

\subsection{The limit of the fibre as two flags merge}\label{sec:config_lim}

To apply the construction of \autoref{sec:diff_vs_leq} to the pair $ (\Phi_{\alpha(2)}, \Phi_{\alpha(1)}) $, the vector bundle $ F_{\alpha(2)} \setminus F_{\alpha(1)} $ would have to extend to a vector bundle over all of $ \Phi_{\alpha(2)} $, that is, over the vertices $ \delta_v $ as well.
So we ask what its fibre $ L(K_b) $ does as $ b $ approaches a vertex.

For $ K \in X_2 $ the space $ L(K) $ has dimension $ 7 $.
This subsection shows that $ L(K) $ has a limit in the Grassmannian when the two flags of $ K $ merge into one, and computes the limit (\autoref{lem:config_lim}, \autoref{cor:config_lim_ext}).
The limit spaces form a vector bundle $ W $ of rank $ 7 $ over $ X_1 \cong P $.

\begin{lemma}\label{lem:config_lim}

	Let $ \left\{ K_n \right\} \subset X_{2} \subset 2^P $ be a sequence converging to $ K \in X_1 \subset 2^P $ within $ 2^P $.
	Then the following hold:
	\begin{enumerate}[(1)]
		\item The sequence $ \left\{ L(K_n) \right\} \subset \Gr(7, \Pi_{4,2}) $ converges within $ \Gr(7, \Pi_{4,2}) $ to a subspace $L$ of $ L(K) $.
		\item If $ \left\{ K'_n \right\} \subset X_{2} $ is another sequence converging to $ K $, then $ \lim_n L(K'_n) = L = \lim_n L(K_n) $.
	\end{enumerate}
\end{lemma}

To better understand the lemma, let us recall some basic definitions.

\begin{definition}
	\label{pg:65}
	The \emph{Grassmann manifold}
	\begin{equation*}
		G_n(\C^{n+k}) = \Gr(n, \C^{n+k}) \coloneq \left\{ V \subset \C^{n+k} \mid V \text{ is an } n\text{-dimensional vector subspace of } \C^{n+k} \right\}
	\end{equation*}
	An \emph{$ n $-frame} in $ \C^{n+k} $ is a tuple $ (v_1, \ldots, v_n) \in (\C^{n+k})^{\times n} $ such that the vectors $ v_i $ are linearly independent.

	The \emph{Stiefel manifold} $ V_n(\C^{n+k}) \coloneq \left\{ (v_1, \ldots, v_n) \in (\C^{n+k})^{\times n } \mid (v_1, \ldots, v_n) \text{ is an $ n $-frame} \right\} $ is an open set in $ (\C^{n+k})^{\times n} $.

	Let $ q \colon V_n(\C^{n+k}) \to \Gr(n,\C^{n+k}) $ be defined by $ (v_1, \ldots, v_n) \mapsto \langle v_1, \ldots, v_n \rangle _{\C} \simeq \C^{n} $.

	The topology on $ \Gr(n,\C^{n+k}) $ is defined as the quotient topology: $ U \subset \Gr(n,\C^{n+k}) $ is open if and only if $ q^{-1}(U) $ is open in $ V_n(\C^{n+k}) $.

	Equivalently, denote $ V^0_n(\C^{n+k}) \coloneq \left\{ (v_1, \ldots, v_n) \in V_n(\C^{n+k}) \mid v_i \text{ are orthonormal} \right\} $.
	Consider the commutative diagram
	\begin{equation*}
		\begin{tikzcd}
			V^0_n(\C^{n+k}) \arrow[r, "i", hook] \arrow[rd, "q_0", two heads] & V_n(\C^{n+k}) \arrow[r, "j", two heads] \arrow[d, "q", two heads] & V^0_n(\C^{n+k}) \arrow[ld, "q_0", two heads] \\
			& \Gr(n,\C^{n+k})                                                     &
		\end{tikzcd}
	\end{equation*}
	where $ i $ is the inclusion, and $ j(v_1, \ldots, v_n) = (v'_1, \ldots, v'_n) $ is obtained via the Gram--Schmidt process.

	Then $ \Gr(n,\C^{n+k}) = q_0(V^0_n(\C^{n+k})) \simeq V^0_n(\C^{n+k}) / U(n) $.

\end{definition}

\begin{remark}
	\begin{itemize}
		\item $ \Gr(n,\C^{n+k}) $ is a compact topological manifold of real dimension $ 2nk $.
		\item $ \Gr(n,\C^{n+k}) \to \Gr(k,\C^{n+k}) \colon X \mapsto X^{\bot} $ is a homeomorphism.
	\end{itemize}
\end{remark}

\begin{definition}
	A canonical vector bundle over $ \Gr(n,\C^{n+k}) $:
	\begin{equation*}
		\gamma^n = \gamma^n(\C^{n+k}) \colon F_{X_0} \hookrightarrow E \twoheadrightarrow \Gr(n,\C^{n+k}) \ni X_0,
	\end{equation*}
	where $ \C^{n} \simeq X_0 = F_{X_0} $, $ E \coloneq \left\{ (X,v) \in \Gr(n,\C^{n+k}) \times \C^{n+k} \colon v\in X \right\} $.
\end{definition}

\begin{proof}{\it (\autoref{lem:config_lim})}

	Without loss of generality, we may assume that $ K = \{ (p, l) \} $, where $ p = (1:0:0) $ and \mbox{$ l = \overline{(1:0:0), (0:1:0)} $}.
	Consider the local affine chart $ U = \{ X = 1 \} = \C^2 $.
	In local coordinates, $ p = (0,0) $ and $ l = \{ z=0 \} $.

	We claim that for any sequence $ \left\{ K_n \right\} $ converging to $ K $, the sequence $ \left\{ L(K_n) \right\} $ converges to
	\begin{equation*}
		L \coloneq {\left \langle \begin{array}{c} y^4 , y^3z , y^2z^2 , yz^3 , z^4 , \\ xyz^2 , \\ xz^3 \end{array} \right \rangle } .
	\end{equation*}
	From this point onwards, we will identify all objects by their local descriptions in $ U $, going back to the global description of $ L $ at the end.

	Without loss of generality, we may assume that $ p_n \in U $ for all $ n $ (otherwise, shift the sequence by $ N $, where $ \forall n > N $, $ p_n \in U $ by the definition of a limit -- we can do so as topology on $ P $ is induced from $ \C P^2 \times \C P^{2\vee} $, and therefore we can consider open neighbourhoods of $ p_n $ and $ l^i_n $ separately; see \autoref{rem:top_on_ptcp2}).

	Let $ p_n = (a_n, b_n) $, $ l^i_n = \{ \alpha^i_n y + \beta^i_n z - \alpha^i_n a_n - \beta^i_n b_n = 0 \} $, $ i = 1, 2 $, $ K_n = \{ (p_n, l^1_n), (p_n, l^2_n) \} $.

	Let $ F(1,y,z) = f(y, z) $ be a polynomial.
	Then $ K_n \subset \Sing(f) $ if and only if
	\begin{equation*}
		f(y + a_n, z + b_n) = g_n(y,z) = g_{n,4}(y, z) + (\alpha^1_n y + \beta^1_n z)(\alpha^2_n y + \beta^2_n z)(\alpha y + \beta z),
	\end{equation*}
	(here $ \alpha $ and $ \beta $ define an arbitrary line) i.e. \begin{equation*}
		f(y + a_n, z + b_n) = g_n(y,z) \in  {\left \langle \begin{array}{c} y^4 , y^3z , y^2z^2 , yz^3 , z^4 , \\ \alpha^1_n \alpha^2_n y^3 + (\alpha^1_n \beta^2_n + \beta^1_n \alpha^2_n) y^2z + \beta^1_n \beta^2_n yz^2 , \\ \alpha^1_n \alpha^2_n y^2z + (\alpha^1_n \beta^2_n + \beta^1_n \alpha^2_n) yz^2 + \beta^1_n \beta^2_n z^3 \end{array} \right \rangle }.
	\end{equation*}
	Equivalently, the polynomial $f(y, z)$ lies in the space
	\begin{equation*}
		L(K_n) = {\left \langle \begin{array}{c} (y-a_n)^4 , (y-a_n)^3(z-b_n) , (y-a_n)^2(z-b_n)^2 , (y-a_n)(z-b_n)^3 , (z-b_n)^4 , \\ \alpha^1_n \alpha^2_n (y-a_n)^3 + (\alpha^1_n \beta^2_n + \beta^1_n \alpha^2_n) (y-a_n)^2(z-b_n) + \beta^1_n \beta^2_n (y-a_n)(z-b_n)^2 , \\ \alpha^1_n \alpha^2_n (y-a_n)^2(z-b_n) + (\alpha^1_n \beta^2_n + \beta^1_n \alpha^2_n) (y-a_n)(z-b_n)^2 + \beta^1_n \beta^2_n (z-b_n)^3 \end{array} \right \rangle }.
	\end{equation*}
	Without loss of generality, we may assume that $ \beta^i_n = 1 $ (otherwise, change to the affine neighbourhood where $ \beta^i_n \neq 0 $, i.e. $ l^i_n \neq \{ y = 0 \} $, and then replace $ \alpha^i_n $ with $ \frac{\alpha^i_n}{\beta^i_n} $; see \autoref{rem:top_on_ptcp2}).

	Now, since $ K_n \to K $, we have $ a_n, b_n, \alpha^1_n, \alpha^2_n \to 0 $.
	Therefore, returning to global coordinates, and by the definition of the topology on $ \Gr(7, \Pi_{4,2}) $, we obtain
	\begin{equation*}
		L(K_n) \to L = {\left \langle \begin{array}{c} y^4 , y^3z , y^2z^2 , yz^3 , z^4 , \\ xyz^2 , \\ xz^3 \end{array} \right \rangle } .
	\end{equation*}
	Indeed, consider the fibre bundle
	\begin{equation*}
		\begin{tikzcd}
			\GL_7(\C) \arrow[r, hook] & V_7(\C^{15}) \arrow[r, "q", two heads]                            & G_7(\C^{15}) \\
			& q^{-1}(W) \arrow[r, "q", two heads] \arrow[d, equal] & W             \\
			& \GL_7(\C) \times W \arrow[ru, "proj_2"]                                 &
		\end{tikzcd}
	\end{equation*}
	together with a trivialising neighbourhood $ W $ around $ L $.
	Here the fibre is $ \GL_7(\C) $ because $ V_7(\C^{15}) $ is the manifold of all $ 7 $-frames.
	Had we used the orthonormal frames $ V^0_7(\C^{15}) $ of the definition above, it would be $ U(7) $.
	Fix the identity $ I \in \GL_7(\C) $.
	Now each space $ L(K_n) $ can be lifted to a frame in $ V_7(\C^{15}) $.
	It follows, using the fact that $ V_7(\C^{15}) $ has the topology induced from product topology, that convergence of basis vectors in $ L(K_n) $ to a basis of $ L $ implies convergence of $ L(K_n) $ to $ L $.
\end{proof}

\begin{remark}\label{rem:top_on_ptcp2}
	\begin{enumerate}[(1)]
		\item The topology on $ P $ is induced from $ \C P^2 \times \C P^{2\vee} $.
		\item $ K_n \to K $ implies $ a_n, b_n, \alpha^1_n, \alpha^2_n \to 0 $ (see $ 1 $).
	\end{enumerate}
\end{remark}

The following description of the spaces $ L(K) $ for $ K \in X_2 $ will be used repeatedly below.

\begin{lemma}\label{lem:two_flag_fibre}
	Let $ K = \left\{ (p, l_1), (p, l_2) \right\} \in X_2 $.
	Work in an affine chart centred at $ p $ as in \autoref{thm:tacnode_calc}, so that every polynomial has a Taylor expansion $ f = f_0 + f_1 + f_2 + f_3 + f_4 $ around $ p $, and let $ L_1 $, $ L_2 $ be linear equations of $ l_1 $, $ l_2 $ in this chart.
	Then
	\begin{equation*}
		L(K) = \left\{ f \mid f_0 = f_1 = f_2 = 0 \text{ and } L_1 L_2 \text{ divides } f_3 \right\} .
	\end{equation*}
\end{lemma}

\begin{proof}
	By definition $ L(K) = L(\left\{ (p, l_1) \right\}) \cap L(\left\{ (p, l_2) \right\}) $, and \autoref{thm:tacnode_calc} describes both terms of the intersection.
	The conditions $ f_0 = f_1 = 0 $ are common to both.
	The two conditions on $ f_2 $ read $ f_2 = \lambda_1 L_1^2 = \lambda_2 L_2^2 $.
	Since $ l_1 \neq l_2 $, the forms $ L_1 $ and $ L_2 $ are non-proportional, and by unique factorisation a non-zero multiple of $ L_1^2 $ is never a multiple of $ L_2^2 $, so $ \lambda_1 = \lambda_2 = 0 $, i.e. $ f_2 = 0 $.
	Finally, $ f_3 $ vanishes on both lines, so the coprime linear forms $ L_1 $ and $ L_2 $ both divide $ f_3 $, i.e. $ L_1 L_2 \mid f_3 $.
	The converse implications are clear.
	As a consistency check, the dimension of the right-hand side is $ 2 + 5 = 7 = d_2 $ (the multiples of $ L_1 L_2 $ among cubics $f_3$, plus all quartics $f_4$), in agreement with the row $ i = 2 $ of \autoref{table:configs_dim}.
\end{proof}

\autoref{lem:config_lim} extends the assignment $ K \mapsto L(K) $ continuously to singleton configurations.
We record the extension once and for all.

\begin{corollary}\label{cor:config_lim_ext}
	Define
	\begin{equation*}
		\widehat L \colon X_2 \cup X_1 \to \Gr(7, \Pi_{4,2}) , \qquad
		\widehat L(K) \coloneq
		\begin{cases}
			L(K)            & K \in X_2 , \\
			\lim_n L(K_n)   & K \in X_1 ,
		\end{cases}
	\end{equation*}
	where in the second case $ \left\{ K_n \right\} \subset X_2 $ is any sequence converging to $ K $ (the limit exists and does not depend on the choice by \autoref{lem:config_lim}).
	Then:
	\begin{enumerate}[(1)]
		\item $ \widehat L $ is continuous;
		\item in the local notation of \autoref{thm:tacnode_calc}, for $ K = \left\{ (p, l) \right\} \in X_1 $ with local equation $ L $ of $ l $,
		\begin{equation*}
			\widehat L \left( \left\{ (p, l) \right\} \right) = \left\{ f \mid f_0 = f_1 = f_2 = 0 \text{ and } L^2 \text{ divides } f_3 \right\} .
		\end{equation*}
	\end{enumerate}
\end{corollary}

\begin{proof}
	First observe that, with $ \beta^i_n = 1 $ as in the proof of \autoref{lem:config_lim}, the sixth and seventh spanning polynomials of $ L(K_n) $ displayed there factor:
	\begin{equation*}
		\begin{aligned}
			v_6 &= (y - a_n) \left( \alpha^1_n (y - a_n) + (z - b_n) \right) \left( \alpha^2_n (y - a_n) + (z - b_n) \right) = (y - a_n) \widetilde L^1_n \widetilde L^2_n , \\
			v_7 &= (z - b_n) \widetilde L^1_n \widetilde L^2_n ,
		\end{aligned}
	\end{equation*}
	where $ \widetilde L^i_n \coloneq \alpha^i_n (y - a_n) + (z - b_n) $ is the local equation of $ l^i_n $.
	Together with the five quartics in $ (y - a_n, z - b_n) $, these seven polynomials depend continuously on the parameters $ (a_n, b_n, \alpha^1_n, \alpha^2_n) $, and they form a frame whenever the spanned space is $ 7 $-dimensional.

	We prove $ (2) $ first.
	Fix $ (p, l) $ with $ p = (a, b) $ in the chart and $ l = \{ \widetilde L = 0 \} $, $ \widetilde L = \alpha (y - a) + (z - b) $, and take $ K_n = \left\{ (p, l), (p, l'_n) \right\} $ with $ l'_n \to l $, $ l'_n \neq l $.
	By the factorisation above, the frames of $ L(K_n) $ converge to the seven polynomials
	\begin{equation*}
		(y - a) \widetilde L^2 , \quad (z - b) \widetilde L^2 , \quad (y - a)^4 , \ (y - a)^3 (z - b) , \ \ldots , \ (z - b)^4 ,
	\end{equation*}
	which are linearly independent and span exactly $ \left\{ f \mid f_0 = f_1 = f_2 = 0 \text{ and } \widetilde L^2 \text{ divides } f_3 \right\} $: the cubic parts $ (y-a) \widetilde L^2 $ and $ (z-b) \widetilde L^2 $ span the multiples of $ \widetilde L^2 $ among the cubics in $ (y - a, z - b) $, and the quartic parts are unrestricted.
	By the Stiefel-manifold argument concluding the proof of \autoref{lem:config_lim}, the spaces $ L(K_n) $ converge to this span, which by \autoref{lem:config_lim} is $ \widehat L(\left\{ (p,l) \right\}) $.

	For $ (1) $, let $ \left\{ K_n \right\} \subset X_2 \cup X_1 $ converge to $ K $.
	It suffices to show that every subsequence admits a further subsequence along which $ \widehat L(K_n) \to \widehat L(K) $, and after passing to a subsequence we may assume that all $ K_n $ lie in the same stratum.
	Note that a sequence of singletons cannot converge to a two-point configuration: the Hausdorff distance from a one-point set to $ \left\{ v, w \right\} $ is at least half the distance between $ v $ and $ w $.
	Three cases remain.
	If $ K_n \in X_2 $ and $ K \in X_1 $, this is \autoref{lem:config_lim}.
	If $ K_n \in X_2 $ and $ K \in X_2 $, then, after the chart normalisations of \autoref{rem:top_on_ptcp2}, the parameters of $ K_n $ converge to those of $ K $, so the frames above converge to the corresponding frame of $ L(K) $, which consists of $ 7 $ linearly independent vectors because $ \dim L(K) = 7 $ by \autoref{thm:dimensions}.
	The Stiefel-manifold argument applies again.
	If $ K_n \in X_1 $ and $ K \in X_1 $, then by $ (2) $ the spaces $ \widehat L(K_n) $ are spanned by the frames
	\begin{equation*}
		(y - a_n) \widetilde L_n^2 , \quad (z - b_n) \widetilde L_n^2 , \quad \text{the five quartics in } (y - a_n, z - b_n) ,
	\end{equation*}
	which again depend continuously on the parameters $ (a_n, b_n, \alpha_n) $ and converge to the corresponding frame for $ K $.
\end{proof}

We single out the resulting family of limit spaces over $ X_1 \cong P $.
By \autoref{cor:config_lim_ext} and \cite[Lemma 10]{BerGorKon23}, the set
\begin{equation*}
	W \coloneq \left\{ (f, (p, l)) \in \Pi_{4,2} \times P \mid f \in \widehat L \left( \left\{ (p, l) \right\} \right) \right\}
\end{equation*}
is the total space of a rank-$ 7 $ complex vector bundle over $ P $, namely the pullback of the tautological bundle $ \gamma^7 $ under $ (p,l) \mapsto \widehat L(\left\{ (p,l) \right\}) $.

\subsection{Why the bundle does not extend over \texorpdfstring{$ \Phi_{\alpha(2)} $}{Phi}}\label{sec:closure_failure}

A vertex $ \delta_{(p,l)} \in \Phi_{\alpha(1)} $ is an end of many arcs.
There is one arc $ \widetilde \Lambda(K) $ for every configuration
\begin{equation*}
	K = \left\{ (p,l), (p,l') \right\} \in X_2 , \qquad l' \neq l \text{ a line through } p .
\end{equation*}
Along this arc the fibre of $ F_{\alpha(2)} \setminus F_{\alpha(1)} $ is the constant space $ L(K) $ (\autoref{sec:two_to_one_goal}), and this space depends on $ l' $.
So the fibre has different limits at the vertex along different arcs, and $ F_{\alpha(2)} \setminus F_{\alpha(1)} $ does not extend to a subbundle of the trivial bundle over $ \Phi_{\alpha(2)} $ (\autoref{prop:no_bundle_ext}).
The closure of $ F_{\alpha(2)} \setminus F_{\alpha(1)} $ contains all these limits.
Its fibre over the vertex is the union of the spaces $ \widehat L ( \left\{ (p,l), (p,l') \right\} ) $ over all lines $ l' $ through $ p $, and it has dimension $ 8 $ instead of $ 7 $ (\autoref{lem:closure_jump}).

Throughout this subsection $ (p,l) \in P $ is a flag and
\begin{equation*}
	x_0 \coloneq \delta_{(p,l)} \in \Phi_{\alpha(1)}
\end{equation*}
is its vertex.
We work in an affine chart centred at $ p $ as in \autoref{thm:tacnode_calc}.
Every polynomial has a Taylor expansion $ f = f_0 + f_1 + f_2 + f_3 + f_4 $ around $ p $, and $ L $, $ L' $ are linear equations of lines $ l $, $ l' $ through $ p $ in this chart.
The letter $ L $ followed by a configuration, as in $ L(K) $, still denotes the space of polynomials of that configuration.

We first describe the simplest way of approaching the vertex, along one arc.
Let $ l' \neq l $ be a line through $ p $ and $ K \coloneq \left\{ (p,l), (p,l') \right\} \in X_2 $.
Consider the points
\begin{equation}\label{eq:approach_arc}
	b_n \coloneq \left( 1 - \tfrac{1}{n} \right) \delta_{(p,l)} + \tfrac{1}{n} \, \delta_{K} , \qquad n \geq 1 . \tag{A}
\end{equation}
They converge to $ x_0 $.
They lie on the edge $ [ \delta_{(p,l)}, \delta_K ] $ of the arc $ \widetilde \Lambda(K) $ and differ from its end $ \delta_{(p,l)} $.
So $ b_n \in \Phi_{\alpha(2)} \setminus \Phi_{\alpha(1)} $ and $ K_{b_n} = K $ by \autoref{eq:edge_coords}.
By \autoref{lem:two_flag_fibre}, the fibre of $ F_{\alpha(2)} \setminus F_{\alpha(1)} $ over every $ b_n $ is the same space
\begin{equation*}
	L(K) = \left\{ f \mid f_0 = f_1 = f_2 = 0 \text{ and } L L' \text{ divides } f_3 \right\} .
\end{equation*}

\begin{proposition}\label{prop:no_bundle_ext}
	No topological subbundle of rank $ 7 $ of the trivial bundle $ \Pi_{4,2} \times \Phi_{\alpha(2)} $ restricts to $ F_{\alpha(2)} \setminus F_{\alpha(1)} $ over $ \Phi_{\alpha(2)} \setminus \Phi_{\alpha(1)} $.
\end{proposition}

\begin{proof}
	Such a subbundle would have a continuous classifying map $ \Phi_{\alpha(2)} \to \Gr(7, \Pi_{4,2}) $ that extends the map $ b \mapsto L(K_b) $.
	We show that the map $ b \mapsto L(K_b) $ has no limit at one vertex.

	Take $ p = (1:0:0) $ and $ l = \{ z = 0 \} $, with the chart of the proof of \autoref{lem:config_lim}, so that $ L = z $.
	Take the sequences \eqref{eq:approach_arc} for the two lines
	\begin{equation*}
		l' = \{ y = 0 \} \qquad \text{and} \qquad l'' = \{ y = z \} .
	\end{equation*}
	Both converge to $ x_0 $.
	Along them the map $ b \mapsto L(K_b) $ has the constant values
	\begin{equation*}
		\begin{aligned}
			L \left( \left\{ (p,l), (p,l') \right\} \right) & = \left\{ f \mid f_0 = f_1 = f_2 = 0 \text{ and } y z \text{ divides } f_3 \right\} , \\
			L \left( \left\{ (p,l), (p,l'') \right\} \right) & = \left\{ f \mid f_0 = f_1 = f_2 = 0 \text{ and } z ( y - z ) \text{ divides } f_3 \right\} .
		\end{aligned}
	\end{equation*}
	These two spaces are different.
	Indeed, the polynomial $ x y^2 z $ lies in the first and not in the second, because its cubic part $ y^2 z $ at $ p $ is divisible by $ y z $ but not by $ z ( y - z ) $.
	So no continuous extension to $ x_0 $ exists.
\end{proof}

The same two sequences show that the map $ b \mapsto K_b $ does not extend continuously to the vertex.
Along the first its value is constantly $ \left\{ (p,l), (p,l') \right\} $, and along the second it is constantly $ \left\{ (p,l), (p,l'') \right\} $.

By \autoref{prop:no_bundle_ext}, \cite[Lemma 18]{BerGorKon23} cannot be applied to the pair $ (\Phi_{\alpha(2)}, \Phi_{\alpha(1)}) $ at all.
The lemma requires its bundle $ E_1 $ to be defined over the whole ambient space of the pair.

We now turn to the closure of the bundle.
Let
\begin{equation*}
	\overline F \coloneq \overline{ F_{\alpha(2)} \setminus F_{\alpha(1)} }
\end{equation*}
denote the closure taken inside $ F_{\alpha(2)} $.
Over $ \Phi_{\alpha(2)} \setminus \Phi_{\alpha(1)} $ this is the bundle $ F_{\alpha(2)} \setminus F_{\alpha(1)} $ of rank $ 7 $.
Over the vertices, however, the fibres jump.
We compute the fibre over a vertex, because these fibres return in \autoref{sec:map_alpha} as the fibres of a vector bundle over $ P $.

Three spaces of polynomials are attached to the flag $ (p,l) $, where $ l' \neq l $ is again a line through $ p $:
\begin{equation*}
	\begin{aligned}
		L \left( \left\{ (p,l), (p,l') \right\} \right) & = \left\{ f \mid f_0 = f_1 = f_2 = 0 \text{ and } L L' \text{ divides } f_3 \right\} , \\
		\widehat L \left( \left\{ (p,l) \right\} \right) & = \left\{ f \mid f_0 = f_1 = f_2 = 0 \text{ and } L^2 \text{ divides } f_3 \right\} , \\
		W'_{(p,l)} & \coloneq \left\{ f \mid f_0 = f_1 = f_2 = 0 \text{ and } L \text{ divides } f_3 \right\} .
	\end{aligned}
\end{equation*}
The first two are described in \autoref{lem:two_flag_fibre} and \autoref{cor:config_lim_ext}(2), and they have dimension $ 7 $.
The third contains the first two, and \autoref{lem:closure_jump} below shows that it is the fibre of $ \overline F $ over $ x_0 $.

In the second way of approaching the vertex, the second flag merges into the first.
Let $ l'_n \neq l $ be lines through $ p $ with $ l'_n \to l $, and
\begin{equation}\label{eq:approach_merge}
	K_n \coloneq \left\{ (p,l), (p,l'_n) \right\} \in X_2 , \qquad b_n \coloneq \delta_{K_n} . \tag{B}
\end{equation}
Then $ b_n \in \widetilde \Lambda(K_n) \setminus \Phi_{\alpha(1)} $, and $ b_n \to x_0 $.
Indeed, $ K_n \to \left\{ (p,l) \right\} $ in $ 2^P $, and the vertex embedding $ \iota $ is continuous.
The fibre over $ b_n $ is $ L(K_n) $, and $ L(K_n) \to \widehat L ( \left\{ (p,l) \right\} ) $ with converging frames (proof of \autoref{cor:config_lim_ext}).

\begin{lemma}\label{lem:closure_jump}
	Let $ x_0 = \delta_{(p,l)} \in \Phi_{\alpha(1)} $ be a vertex, and let $ L $ be a local equation of $ l $ as in \autoref{thm:tacnode_calc}.
	Then the fibre of $ \overline F $ over $ x_0 $, that is $ \left\{ f \mid (f, x_0) \in \overline F \right\} $, equals
	\begin{equation*}
		W'_{(p,l)} \coloneq \left\{ f \mid f_0 = f_1 = f_2 = 0 \text{ and } L \text{ divides } f_3 \right\} ,
	\end{equation*}
	a linear space of dimension $ 8 $.
	In particular $ \overline F $ is not a vector bundle over $ \Phi_{\alpha(2)} $: its fibres have dimension $ 7 $ over $ \Phi_{\alpha(2)} \setminus \Phi_{\alpha(1)} $ and dimension $ 8 $ over $ \Phi_{\alpha(1)} $.
\end{lemma}

\begin{proof}
	The dimension count is as in \autoref{lem:two_flag_fibre}: $ 3 $ multiples of $ L $ among the cubics plus $ 5 $ quartics.

	First, let us show that the fibre is contained in $ W'_{(p,l)} $.
	Let $ (f, x_0) \in \overline F $.
	Then
	\begin{equation*}
		(f, x_0) = \lim_n (f_n, b_n) , \qquad b_n \in \Phi_{\alpha(2)} \setminus \Phi_{\alpha(1)} , \quad f_n \in L(K_{b_n}) .
	\end{equation*}
	Write $ K_n \coloneq K_{b_n} = \left\{ v_n, w_n \right\} \in X_2 $, labelled so that the edge coordinates of \autoref{eq:edge_coords} read
	\begin{equation*}
		b_n = s_n \delta_{v_n} + (1 - s_n) \delta_{K_n} , \qquad s_n \in [0, 1) .
	\end{equation*}
	Passing to a subsequence, we may assume $ s_n \to s $, $ v_n \to v $ and $ w_n \to w $ in the compact spaces involved.
	Both flags of $ K_n $ lie over a common point of $ \C P^2 $, so their limits do as well, and therefore
	\begin{equation*}
		K_\infty \coloneq \left\{ v, w \right\} \in X_2 \cup X_1 , \qquad K_n \to K_\infty \text{ in } 2^P .
	\end{equation*}

	Let us show that the flag $ (p,l) $ belongs to $ K_\infty $.
	By the continuity of $ \iota $ and of convex combinations in $ \R^\Omega $, the points $ b_n $ converge to $ s \delta_v + (1 - s) \delta_{K_\infty} $.
	They also converge to $ x_0 = \delta_{(p,l)} $, so
	\begin{equation*}
		\delta_{(p,l)} = s \delta_v + (1 - s) \delta_{K_\infty} .
	\end{equation*}
	If $ K_\infty = \left\{ v \right\} $, the right-hand side equals $ \delta_v $, so $ v = (p, l) $ and $ K_\infty = \left\{ (p,l) \right\} $.
	Otherwise $ v \in K_\infty $ makes $ \left[ \delta_v, \delta_{K_\infty} \right] $ a $ 1 $-simplex of the self-join.
	Its points other than the two endpoints lie in the interior of a $ 1 $-simplex and are therefore not vertices.
	Hence either $ s = 1 $ and $ v = (p, l) $, or $ s = 0 $ and $ K_\infty = \left\{ (p,l) \right\} $.
	In every case $ (p, l) \in K_\infty $.

	It remains to show that $ f $ lies in $ \widehat L(K_\infty) $ and that $ \widehat L(K_\infty) \subseteq W'_{(p,l)} $.
	By \autoref{cor:config_lim_ext}, $ L(K_n) \to \widehat L(K_\infty) $ in the Grassmannian.
	The total space of the tautological bundle $ \gamma^7 $ is closed in $ \Pi_{4,2} \times \Gr(7, \Pi_{4,2}) $.
	So the limit $ f $ of $ f_n \in L(K_n) $ lies in $ \widehat L(K_\infty) $.
	Finally, $ K_\infty $ contains $ (p,l) $, so $ \widehat L(K_\infty) $ is one of the first two of the three spaces displayed before the lemma (\autoref{lem:two_flag_fibre} and \autoref{cor:config_lim_ext}(2)).
	It consists of polynomials with $ f_0 = f_1 = f_2 = 0 $ and $ L \mid f_3 $, that is, $ \widehat L(K_\infty) \subseteq W'_{(p,l)} $.

	Next, let us show that $ W'_{(p,l)} $ is contained in the fibre.
	Let $ f \in W'_{(p,l)} $, and write $ f_3 = L q $ with $ q $ a quadratic form in the local coordinates.
	There are two cases, and each uses one of the two ways of approaching the vertex.

	Suppose first that $ q $ has a linear factor $ L' $ non-proportional to $ L $.
	This happens whenever $ q \neq 0 $ is not a multiple of $ L^2 $, and also when $ q = 0 $, in which case any $ L' $ will do.
	Let $ l' $ be the line through $ p $ with equation $ L' $, and $ K \coloneq \left\{ (p, l), (p, l') \right\} \in X_2 $.
	Then $ f \in L(K) $ by \autoref{lem:two_flag_fibre}.
	Let $ b_n $ be the points \eqref{eq:approach_arc}.
	The fibre over $ b_n $ is $ L(K) $, so $ (f, b_n) \in F_{\alpha(2)} \setminus F_{\alpha(1)} $, and $ (f, b_n) \to (f, x_0) $.
	Hence $ (f, x_0) \in \overline F $.

	In the remaining case $ q $ is a non-zero multiple of $ L^2 $, so $ f \in \widehat L(\left\{ (p,l) \right\}) $ by \autoref{cor:config_lim_ext}(2).
	Let $ K_n $ and $ b_n $ be as in \eqref{eq:approach_merge}.
	Since $ L(K_n) \to \widehat L(\left\{ (p,l) \right\}) $ with converging frames, there are $ f_n \in L(K_n) $ with $ f_n \to f $.
	Then $ (f_n, b_n) \in F_{\alpha(2)} \setminus F_{\alpha(1)} $ and $ (f_n, b_n) \to (f, x_0) $, so $ (f, x_0) \in \overline F $.
\end{proof}

\begin{example}
	Take $ p = (1:0:0) $ and $ l = \{ z = 0 \} $ as in the proof of \autoref{prop:no_bundle_ext}, so that $ L = z $, and consider
	\begin{equation*}
		f = x y^2 z + y^4 + z^4 .
	\end{equation*}
	Its cubic part at $ p $ is $ y^2 z $, which is divisible by $ z $ but not by $ z^2 $.
	So $ f $ lies in $ W'_{(p,l)} $ but not in the fibre $ \widehat L( \left\{ (p,l) \right\} ) $ of $ W $.
	It is covered by the first case of the proof above, with $ l' = \{ y = 0 \} $.
	The points \eqref{eq:approach_arc} approach the vertex along one edge of the fixed arc $ \widetilde \Lambda(K) $, $ K = \left\{ (p,l), (p,l') \right\} $, and the fibre over each of them is the constant space $ L(K) \ni f $.
\end{example}

The next subsection produces a substitute pair, over which a bundle as in \cite[Lemma 18]{BerGorKon23} does exist.

\subsection{The resolved pair}\label{sec:resolved_pair}

We construct the resolved pair in three steps.
First, we describe the space that replaces $ \Phi_{\alpha(1)} $, namely the space of ordered pairs of flags through a common point, and its quotient by the exchange of the two flags (\autoref{lem:xbar_quotient}).
Second, we define the space that replaces $ \Phi_{\alpha(2)} $ and show that it differs from $ \Phi_{\alpha(2)} $ only over $ \Phi_{\alpha(1)} $ (\autoref{lem:phihat}).
Third, we extend the vector bundle to this space (\autoref{lem:res_bundles}).

Let
\begin{equation*}
	\widehat Y \coloneq P \times_{\C P^2} P = \left\{ \left( (p, l), (p, l') \right) \in P \times P \right\}
\end{equation*}
be the space of ordered pairs of flags with a common point.
The following maps are used throughout the rest of this section.
\begin{equation*}
	\begin{aligned}
		\psi & \colon P \to \C P^2 , & (p, l) & \mapsto p , \\
		\rho & \colon \widehat Y \to P , & ((p, l), (p, l')) & \mapsto (p, l) , \\
		\tau & \colon \widehat Y \to \widehat Y , & ((p, l), (p, l')) & \mapsto ((p, l'), (p, l)) , \\
		\Delta & \colon P \to \widehat Y , & (p, l) & \mapsto ((p, l), (p, l)) , \\
		\eta & \colon \widehat Y \to 2^P , & ((p, l), (p, l')) & \mapsto \left\{ (p, l), (p, l') \right\} .
	\end{aligned}
\end{equation*}
So $ \rho $ is the first projection and $ \rho \circ \tau $ the second, $ \tau $ is the involution exchanging the two flags, and $ \Delta $ is the diagonal embedding, whose image is the fixed locus of $ \tau $.
The map $ \eta $ sends an ordered pair of flags to the set of the two flags.

For $ p \in \C P^2 $ and $ l \in \C P^{2 \vee} $ let $ \overline p \subset \overline l \subset \C^3 $ denote the corresponding line and plane.
Let $ \mathcal{O}(-1) \subset \C^3 \times \C P^2 $ be the tautological line bundle and
\begin{equation*}
	\mathcal{Q} \coloneq \left( \C^3 \times \C P^2 \right) / \mathcal{O}(-1)
\end{equation*}
the universal quotient bundle of rank $ 2 $, whose fibre over $ p $ is $ \C^3 / \overline p $.
Sending $ (p, l) $ to the line $ \overline l / \overline p \subset \mathcal{Q}_p $ identifies $ P $ with the projectivisation $ \PP(\mathcal{Q}) $ over $ \C P^2 $, and likewise identifies $ \widehat Y $ with $ \PP( \psi^* \mathcal{Q} ) $ over $ P $.
In particular $ \widehat Y $ is a smooth projective variety of complex dimension $ 4 $, and $ \rho $ is a $ \C P^1 $-bundle.

Recall that $ \overline X_2 $ denotes the closure of $ X_2 $ in $ 2^P $, and that $ \overline X_2 = X_2 \sqcup X_1 $ (see item (2) of the proof of \autoref{prop:condition5}).
The next lemma identifies the quotient of $ \widehat Y $ by $ \tau $ with $ \overline X_2 $, and also with a bundle of projective planes over $ \C P^2 $.
The first description relates the quotient to the configurations, and the second makes it a smooth projective variety.

To state the lemma we need the following notation.
For a complex vector space $ V $ denote by $ \Sym^2 V $ its symmetric square, that is, the quotient of $ V \otimes V $ by the subspace spanned by the elements $ u \otimes u' - u' \otimes u $.
For $ u, u' \in V $ denote by $ u \odot u' \in \Sym^2 V $ the image of $ u \otimes u' $.
Denote by $ \Sym^2 \mathcal{Q} $ the vector bundle of rank $ 3 $ over $ \C P^2 $ whose fibre over $ p $ is $ \Sym^2 \mathcal{Q}_p $, and by $ \PP( \Sym^2 \mathcal{Q} ) $ its projectivisation.
A point of $ \PP( \Sym^2 \mathcal{Q} ) $ is a pair $ (p, [\, q \,]) $, where $ p \in \C P^2 $ and $ [\, q \,] $ is the line spanned by a non-zero element $ q \in \Sym^2 \mathcal{Q}_p $.

\begin{lemma}\label{lem:xbar_quotient}
	\begin{enumerate}[(1)]
		\item The map $ \eta $ is continuous, and its image is $ \overline{X}_2 $.
		Moreover $ \eta \circ \tau = \eta $, and $ \eta(\hat y) = \eta(\hat y\,') $ if and only if $ \hat y\,' \in \left\{ \hat y, \tau(\hat y) \right\} $.
		\item The map $ m \colon \widehat Y \to \PP( \Sym^2 \mathcal{Q} ) $, $ ((p, l), (p, l')) \mapsto (p, [\, u \odot u' \,]) $, where $ u $ and $ u' $ span the lines $ \overline l / \overline p $ and $ \overline{l'} / \overline p $ in $ \mathcal{Q}_p $, is a surjective morphism of varieties, and $ m(\hat y) = m(\hat y\,') $ if and only if $ \hat y\,' \in \left\{ \hat y, \tau(\hat y) \right\} $.
		\item There is a unique homeomorphism $ \overline{X}_2 \cong \PP( \Sym^2 \mathcal{Q} ) $ under which $ \eta $ corresponds to $ m $.
	\end{enumerate}
\end{lemma}

\begin{proof}
	$ (1) $ Continuity is immediate from the definition of the Hausdorff metric:
	\begin{equation*}
		\rho_H \left( \left\{ a, b \right\}, \left\{ a', b' \right\} \right) \leq \max \left( d(a, a'), d(b, b') \right) .
	\end{equation*}
	The image consists of the sets $ \left\{ (p, l), (p, l') \right\} $, which form $ X_2 $ when $ l \neq l' $ and $ X_1 $ when $ l = l' $.
	So the image is $ X_2 \sqcup X_1 = \overline X_2 $.
	The last two statements are clear (a two-element set determines its elements up to order).

	$ (2) $ The map is well defined ($ u $ and $ u' $ are unique up to scalars) and is a morphism (choose local frames).
	Every point of a fibre $ \PP( \Sym^2 \mathcal{Q}_p ) $ is a class $ [\, q \,] $ with $ q \neq 0 $, and every element of $ \Sym^2 $ of a $ 2 $-dimensional complex vector space factors as $ q = u \odot u' $, uniquely up to order and scalars.
	This gives surjectivity and the description of the fibres.

	$ (3) $ Both $ \eta $ (onto $ \overline X_2 $) and $ m $ are continuous surjections from a compact space onto Hausdorff spaces, hence closed maps, hence quotient maps.
	By $ (1) $ and $ (2) $ they identify exactly the same pairs of points, namely the $ \tau $-orbits.
	Therefore the bijection $ \PP(\Sym^2 \mathcal{Q}) \to \overline X_2 $ sending $ m(\hat y) \mapsto \eta(\hat y) $ is well defined, and it is continuous by the universal property of quotient maps.
	A continuous bijection from a compact space to a Hausdorff space is a homeomorphism.
\end{proof}

From now on we identify $ \overline X_2 $ with $ \PP( \Sym^2 \mathcal{Q} ) $ by the homeomorphism of \autoref{lem:xbar_quotient}(3).
So $ \overline X_2 $ carries the structure of a smooth projective variety of complex dimension $ 4 $, and $ \eta = m $ is a morphism.
By \autoref{lem:xbar_quotient}(1) the fibres of $ \eta $ are the orbits of $ \tau $, so $ \eta $ is two-to-one exactly outside the diagonal.

\begin{corollary}\label{cor:xbar_cohomology}
	$ \dim_\Q H^{2k} ( \overline X_2 ; \Q ) = 1, 2, 3, 2, 1 $ for $ k = 0, \ldots, 4 $, and the cohomology vanishes in odd degrees.
	The same holds for homology.
\end{corollary}

\begin{proof}
	$ \overline X_2 \cong \PP(\Sym^2 \mathcal{Q}) $ is a $ \C P^2 $-bundle over $ \C P^2 $.
	Let $ \xi $ be the first Chern class of its tautological line subbundle.
	The restrictions of $ 1, \xi, \xi^2 $ to each fibre are (up to sign) the powers of the hyperplane class, hence a basis of the cohomology of the fibre.
	By the Leray--Hirsch theorem \cite[Theorem 4D.1]{hatcher2002}, $ H^*(\overline X_2; \Q) $ is a free module over $ H^*(\C P^2; \Q) $ with basis $ 1, \xi, \xi^2 $, which gives the stated dimensions.
	Alternatively, we could have used the spectral sequence of the bundle.
\end{proof}

We can now define the resolved pair.
Let
\begin{equation*}
	\widehat \Phi \coloneq \left\{ (b, K) \in \Phi_{\alpha(2)} \times \overline X_2 \mid b \in \widetilde \Lambda(K) \right\} ,
\end{equation*}
and let $ \widehat \pi \colon \widehat \Phi \to \overline X_2 $ and $ \widehat \rho \colon \widehat \Phi \to \Phi_{\alpha(2)} $ be the two projections.

\begin{lemma}\label{lem:phihat}
	\begin{enumerate}[(1)]
		\item $ \widehat \Phi $ is a closed, hence compact, subspace of $ \Phi_{\alpha(2)} \times \overline X_2 $.
		\item $ \widehat \pi $ is proper, and its fibre over $ K \in \overline X_2 $ is $ \widetilde \Lambda(K) \times \{ K \} $.
		\item $ \widehat \rho $ is surjective, and $ \widehat \rho^{\,-1} ( \Phi_{\alpha(1)} ) = \left\{ ( \delta_v , K ) \mid v \in K \right\} $.
		\item The map $ \widehat Y \to \widehat \rho^{\,-1} ( \Phi_{\alpha(1)} ) $, $ ((p, l), (p, l')) \mapsto \left( \delta_{(p,l)}, \left\{ (p,l), (p,l') \right\} \right) $, is a homeomorphism.
		\item $ \widehat \rho $ restricts to a homeomorphism $ \widehat \Phi \setminus \widehat \rho^{\,-1} ( \Phi_{\alpha(1)} ) \to \Phi_{\alpha(2)} \setminus \Phi_{\alpha(1)} $ with inverse $ b \mapsto (b, K_b) $.
	\end{enumerate}
\end{lemma}

\begin{proof}
	$ (1) $ The self-join $ \mathbb{Y}^{*\mathscr{N}} $ is compact because $ \mathbb{Y} $ is, being the image of the compact space $ \mathbb{Y}^{\mathscr{N}} \times \left\{ t \in \R^{\mathscr{N}} \mid t_j \geq 0, \sum_j t_j = 1 \right\} $ under $ ( K_1, \ldots, K_{\mathscr{N}}, t ) \mapsto \sum_j t_j \iota(K_j) $, and $ \Phi_{\alpha(2)} $ is closed in it by \cite[Lemma 8]{BerGorKon23} (cf.\ the discussion following \autoref{lem:tool_cone}), hence compact.
	The space $ \overline X_2 $ is compact as well, so it suffices to show that $ \widehat \Phi $ is closed in $ \Phi_{\alpha(2)} \times \overline X_2 $.
	Let $ (b_n, K_n) \in \widehat \Phi $ converge to $ (b, K) $.
	Every point of $ \widetilde \Lambda(K_n) $ lies on an edge $ [ \delta_{v_n}, \delta_{K_n} ] $ with $ v_n \in K_n $ (a degenerate edge, reduced to the point $ \delta_{K_n} $, when $ K_n \in X_1 $), so we may write $ b_n = s_n \delta_{v_n} + (1 - s_n) \delta_{K_n} $ with $ s_n \in [0, 1] $ and $ v_n \in K_n $.
	Passing to a subsequence, $ s_n \to s $ and $ v_n \to v \in P $.
	From $ v_n \in K_n $ and $ K_n \to K $ in the Hausdorff metric we get $ \operatorname{dist}(v_n, K) \to 0 $, so $ v \in K $.
	By the continuity of $ \iota $ and of convex combinations in $ \R^\Omega $, $ b = s \delta_v + (1 - s) \delta_K $.
	Since $ v \in K $, this point lies on the edge $ [ \delta_v, \delta_K ] \subseteq \widetilde \Lambda(K) $ (equal to the point $ \delta_K $ if $ K = \left\{ v \right\} $), so $ (b, K) \in \widehat \Phi $.

	$ (2) $ $ \widehat \Phi $ is compact by $ (1) $, so $ \widehat \pi $ is proper.
	The fibre over $ K $ is $ \widetilde \Lambda(K) \times \{ K \} $ by definition.

	$ (3) $ For $ b = \delta_v \in \Phi_{\alpha(1)} $ we have $ (b, \left\{ v \right\}) \in \widehat \Phi $, and for $ b \in \Phi_{\alpha(2)} \setminus \Phi_{\alpha(1)} $ we have $ (b, K_b) \in \widehat \Phi $.
	This proves surjectivity.
	A pair $ (\delta_v, K) $ lies in $ \widehat \Phi $ if and only if $ \delta_v $ is a point of $ \widetilde \Lambda(K) $, i.e. if and only if $ v \in K $.

	$ (4) $ The map $ \widehat Y \to \widehat \rho^{\,-1}(\Phi_{\alpha(1)}) $ is continuous, and it is bijective: given $ (\delta_v, K) $, the second flag is recovered as the element of $ K \setminus \{ v \} $ if $ K \in X_2 $ and as $ v $ itself if $ K = \{ v \} $.
	A continuous bijection from the compact $ \widehat Y $ to a Hausdorff space is a homeomorphism.

	$ (5) $ For $ b \in \Phi_{\alpha(2)} \setminus \Phi_{\alpha(1)} $, membership $ b \in \widetilde \Lambda(K) $ with $ K \in \overline X_2 $ forces $ K = K_b $: for $ K \in X_1 $ the space $ \widetilde \Lambda(K) $ is a single vertex, which lies in $ \Phi_{\alpha(1)} $ and so cannot contain $ b $, and for $ K \in X_2 $ this is the uniqueness of the principal vertex established in \autoref{sec:two_to_one_goal}.
	So $ \widehat \rho $ is bijective over $ \Phi_{\alpha(2)} \setminus \Phi_{\alpha(1)} $.
	It is also proper there, being the restriction of the proper map $ \widehat \rho $ over an open subset of the target.
	A proper continuous bijection of locally compact Hausdorff spaces is a homeomorphism.
\end{proof}

\begin{corollary}\label{cor:phihat_homology}
	The map $ \widehat \pi $ induces an isomorphism $ H_* ( \widehat \Phi ; \Q ) \cong H_* ( \overline X_2 ; \Q ) $.
	In particular the odd homology of $ \widehat \Phi $ vanishes.
\end{corollary}

\begin{proof}
	By \autoref{lem:phihat}(2) the map $ \widehat \pi $ is proper, and its fibres are the spaces $ \widetilde \Lambda(K) \times \{ K \} $.
	These are an arc (the union of two $ 1 $-simplices) for $ K \in X_2 $, a point for $ K \in X_1 $, hence contractible.
	Both $ \widehat \Phi $ and $ \overline X_2 $ are compact, locally compact Hausdorff spaces, and all fibres are contractible, so \cite[Lemma 11]{BerGorKon23} gives an isomorphism on Borel--Moore homology, which for compact spaces coincides with singular homology.
	The vanishing of the odd part is \autoref{cor:xbar_cohomology}.
\end{proof}

From now on we identify $ \widehat Y $ with the subspace $ \widehat \rho^{\,-1} ( \Phi_{\alpha(1)} ) $ of $ \widehat \Phi $ by the homeomorphism of \autoref{lem:phihat}(4).
By the definition of this homeomorphism, $ \widehat \rho |_{\widehat Y} \colon \widehat Y \to \Phi_{\alpha(1)} \cong P $ is then the projection $ \rho $, and $ \widehat \pi |_{\widehat Y} \colon \widehat Y \to \overline X_2 $ is the quotient map $ \eta $.
In these terms \autoref{lem:phihat}(5) says that $ \widehat \rho $ restricts to a homeomorphism $ \widehat \Phi \setminus \widehat Y \to \Phi_{\alpha(2)} \setminus \Phi_{\alpha(1)} $.
So the pair $ ( \widehat \Phi, \widehat Y ) $ differs from the pair $ ( \Phi_{\alpha(2)}, \Phi_{\alpha(1)} ) $ only over $ \Phi_{\alpha(1)} $.

The next lemma exhibits $ \widehat \Phi $ as a mapping cylinder of $ \eta $, collapsed along the diagonal, and equips the pair $ ( \widehat \Phi, \widehat Y ) $ with a CW structure.
The latter is a hypothesis of \cite[Lemma 18]{BerGorKon23}, to be applied in \autoref{sec:two_to_one_diff}.

\begin{lemma}\label{lem:phihat_cw}
	Define
	\begin{align*}
		\Xi \colon \left( \widehat Y \times [0, 1] \right) \sqcup \overline X_2 & \longrightarrow \widehat \Phi , \\
		\Xi ( \hat y, s ) & = \left( s \, \delta_{\rho(\hat y)} + ( 1 - s ) \, \delta_{\eta(\hat y)} , \, \eta(\hat y) \right) , \qquad
		\Xi ( K ) = ( \delta_K , K ) .
	\end{align*}
	Then $ \Xi $ is a quotient map identifying $ \widehat \Phi $ with the adjunction space
	\begin{equation*}
		\overline X_2 \cup_g \left( \widehat Y \times [0, 1] \right) , \qquad
		g \colon \left( \widehat Y \times \{ 0 \} \right) \cup \left( \Delta(P) \times [0, 1] \right) \to \overline X_2 , \quad
		g ( \hat y, s ) = \eta ( \hat y ) ,
	\end{equation*}
	under which $ \widehat Y \subset \widehat \Phi $ corresponds to the image of $ \widehat Y \times \{ 1 \} $.
	Consequently $ ( \widehat \Phi, \widehat Y ) $ is homeomorphic to a finite CW pair.
\end{lemma}

\begin{proof}
	The map $ \Xi $ is continuous ($ \rho $, $ \eta $, $ \iota $ and convex combination in $ \R^\Omega $ are continuous), and its values lie in $ \widehat \Phi $: the point $ s \delta_{\rho(\hat y)} + (1-s) \delta_{\eta(\hat y)} $ lies on the edge $ [ \delta_{\rho(\hat y)}, \delta_{\eta(\hat y)} ] \subseteq \widetilde \Lambda ( \eta(\hat y) ) $, and $ \delta_K \in \widetilde \Lambda(K) $ (the common vertex of the two edges for $ K \in X_2 $, the unique point of $ \widetilde \Lambda(K) $ for $ K \in X_1 $).

	It is surjective: a point $ (b, K) \in \widehat \Phi $ with $ b = \delta_K $ is $ \Xi(K) $, and otherwise $ K \in X_2 $ and, in the edge coordinates of \autoref{eq:edge_coords}, $ b = s \delta_v + (1 - s) \delta_K $ with $ s \in (0, 1] $ and $ v \in K $, so $ (b, K) = \Xi ( \hat y, s ) $ for the ordered pair $ \hat y \in \widehat Y $ with $ \rho(\hat y) = v $ and $ \eta(\hat y) = K $.
	As a continuous surjection from a compact space onto a Hausdorff space, $ \Xi $ is a closed map, hence a quotient map.

	The fibres of $ \Xi $ are exactly the classes of the equivalence relation generated by $ a \sim g(a) $ for $ a $ in the domain of $ g $, which identifies $ \widehat \Phi $ with $ \overline X_2 \cup_g ( \widehat Y \times [0,1] ) $.

	Indeed, $ \Xi(\hat y, 0) = ( \delta_{\eta(\hat y)}, \eta(\hat y) ) = \Xi ( \eta(\hat y) ) $ and $ \Xi ( \Delta(u), s ) = ( \delta_u, \left\{ u \right\} ) = \Xi ( \left\{ u \right\} ) $ for all $ s $, so $ \Xi $ factors through the adjunction space.

	Conversely, suppose $ \Xi(\hat y, s) = \Xi(\hat y\,', s') $ with $ \hat y \notin \Delta(P) $ and $ s > 0 $.
	Then $ K \coloneq \eta(\hat y) = \eta(\hat y\,') \in X_2 $, and the common image point $ b = s \delta_{\rho(\hat y)} + (1-s) \delta_K $ determines $ s $ and $ \delta_{\rho(\hat y)} $: for $ s < 1 $ it lies in the interior of the $ 1 $-simplex $ [ \delta_{\rho(\hat y)}, \delta_K ] $, whose vertices and affine coordinates it determines (the singleton vertex being distinguished from $ \delta_K $ by the cardinality of the underlying configuration), and for $ s = 1 $ it is the vertex $ \delta_{\rho(\hat y)} $ itself, while no point with $ s' < 1 $ is that vertex.
	Hence $ s' = s $ and $ \rho(\hat y\,') = \rho(\hat y) $, which together with $ \eta(\hat y\,') = \eta(\hat y) $ gives $ \hat y\,' = \hat y $.

	The remaining cases ($ s = 0 $, or $ \hat y \in \Delta(P) $, or points taken in different parts of the disjoint union) reduce to the two identities displayed above and to the injectivity of $ K \mapsto ( \delta_K, K ) $.

	For the CW structure we realise $ \widehat \Phi $, together with the subspace corresponding to $ \widehat Y $, as a compact semialgebraic set.
	Both $ \widehat Y $ and $ \overline X_2 $ are complex projective varieties, hence Zariski closed in some $ \C P^N $, and $ \C P^N $ is an affine nonsingular real algebraic variety by \cite[Proposition 3.4.6]{BCR98}.
	This realises them as affine real algebraic sets $ \widehat Y \subset \R^m $ and $ \overline X_2 \subset \R^{m'} $, the device already used in \autoref{sec:defs_and_core_tools}.
	The maps $ \tau $ and $ \eta $ are morphisms of complex projective varieties (\autoref{lem:xbar_quotient}), so their graphs are Zariski closed and both maps are semialgebraic in these models.

	Consider the continuous semialgebraic map
	\begin{equation*}
		\Theta \colon \widehat Y \times [0,1] \longrightarrow \overline X_2 \times \R^m , \qquad
		\Theta ( \hat y, s ) \coloneq \left( \eta ( \hat y ) , \ s ( \hat y - \tau \hat y ) \right) ,
	\end{equation*}
	the difference being taken in $ \R^m $.
	Its fibres are exactly those of the restriction of $ \Xi $ to $ \widehat Y \times [0,1] $, which, as shown above, are the classes of $ a \sim g(a) $ that meet $ \widehat Y \times [0,1] $.

	Indeed, $ \Theta ( \hat y, s ) = \Theta ( \hat y\,', s' ) $ forces $ \eta ( \hat y ) = \eta ( \hat y\,' ) $, hence $ \hat y\,' \in \left\{ \hat y, \tau \hat y \right\} $ by \autoref{lem:xbar_quotient}(1).
	If $ \hat y \in \Delta(P) $, then $ \hat y = \tau \hat y $ and both sides equal $ ( \eta ( \hat y ), 0 ) $ for all $ s, s' $, so the segment $ \left\{ \hat y \right\} \times [0,1] $ is collapsed to a point.
	And if $ \hat y \notin \Delta(P) $, then $ u \coloneq \hat y - \tau \hat y \neq 0 $, so $ \hat y\,' = \hat y $ gives $ s u = s' u $ and hence $ s = s' $, while $ \hat y\,' = \tau \hat y $ gives $ s u = - s' u $ and hence $ s = s' = 0 $.

	Since $ \eta $ is surjective, $ \Xi $ is surjective already on $ \widehat Y \times [0,1] $.
	So $ \Theta $ and $ \Xi |_{\widehat Y \times [0,1]} $ are continuous surjections from the same compact space onto Hausdorff spaces and have the same fibres, and both are therefore quotient maps.
	Hence $ \widehat \Phi \cong S \coloneq \Theta ( \widehat Y \times [0,1] ) $, and under this homeomorphism $ \widehat Y $ corresponds to $ S_1 \coloneq \Theta ( \widehat Y \times \left\{ 1 \right\} ) $: the restriction of $ \Xi $ to $ \widehat Y \times \left\{ 1 \right\} $ is the embedding $ \hat y \mapsto ( \delta_{\rho(\hat y)}, \eta(\hat y) ) $ of \autoref{lem:phihat}(4), whose image is $ \widehat \rho^{\,-1} ( \Phi_{\alpha(1)} ) = \widehat Y $.

	The sets $ S $ and $ S_1 $ are semialgebraic by \cite[Proposition 2.2.7]{BCR98} and compact, hence closed and bounded, so by \cite[Theorem 9.2.1]{BCR98} there are a finite simplicial complex $ K $ and a semialgebraic homeomorphism $ |K| \to S $ under which $ S_1 $ is a union of open simplices.
	Being closed, $ S_1 $ is the realisation of a subcomplex.
	So $ ( \widehat \Phi, \widehat Y ) $ is homeomorphic to a finite simplicial pair, in particular to a finite CW pair.
\end{proof}

Recall that $ U \coloneq F_{\alpha(1)} $.

\begin{lemma}\label{lem:res_bundles}
	Let
	\begin{equation*}
		\mathcal{T} \coloneq \left\{ ( f, (b, K) ) \in \Pi_{4,2} \times \widehat \Phi \mid f \in \widehat L (K) \right\} .
	\end{equation*}
	\begin{enumerate}[(1)]
		\item $ \mathcal{T} $ is a rank-$ 7 $ complex vector bundle over $ \widehat \Phi $, the pullback of the tautological bundle $ \gamma^7 $ under $ \widehat L \circ \widehat \pi $.
		\item Under the homeomorphism of \autoref{lem:phihat}(5), the restriction $ \mathcal{T} |_{\widehat \Phi \setminus \widehat Y} $ is identified with the vector bundle $ F_{\alpha(2)} \setminus F_{\alpha(1)} \to \Phi_{\alpha(2)} \setminus \Phi_{\alpha(1)} $ by the map $ (f, (b, K_b)) \mapsto (f, b) $, a homeomorphism of total spaces.
		\item The restriction $ \widehat W \coloneq \mathcal{T} |_{\widehat Y} $ has fibre
		\begin{equation*}
			\widehat W_{((p,l),(p,l'))} = \left\{ f \mid f_0 = f_1 = f_2 = 0 \text{ and } L L' \text{ divides } f_3 \right\}
		\end{equation*}
		(in local notation as before), and $ \widehat W $ is a subbundle of $ \rho^* U $.
		On the diagonal the condition reads $ L^2 \mid f_3 $.
		The quotient
		\begin{equation*}
			\mathcal{N} \coloneq \rho^* U / \widehat W
		\end{equation*}
		is a rank-$ 2 $ complex vector bundle over $ \widehat Y $.
	\end{enumerate}
\end{lemma}

\begin{proof}
	$ (1) $ The map $ \widehat L \circ \widehat \pi \colon \widehat \Phi \to \Gr(7, \Pi_{4,2}) $ is continuous by \autoref{cor:config_lim_ext} and \autoref{lem:phihat}, and $ \mathcal{T} $ is exactly the pullback of $ \gamma^7 $ described in \cite[Lemma 10]{BerGorKon23}.

	$ (2) $ Over $ b \in \Phi_{\alpha(2)} \setminus \Phi_{\alpha(1)} $ the fibre of $ F_{\alpha(2)} \setminus F_{\alpha(1)} $ is $ L(K_b) = \widehat L(K_b) $, so the displayed map is a bijection commuting with the projections and linear on the fibres.
	It is continuous with continuous inverse because the homeomorphism of \autoref{lem:phihat}(5) is, and the $ f $-coordinate is untouched.
	Moreover, by \autoref{lem:tool_vb} (or \cite[Theorem 4]{BerGorKon23}) the right-hand side is a vector bundle, so the identification is one of vector bundles.

	$ (3) $ The fibre of $ \widehat W $ over $ \hat y $ is $ \widehat L ( \eta( \hat y ) ) $, which is the displayed space by \autoref{lem:two_flag_fibre} (off the diagonal) and by \autoref{cor:config_lim_ext}(2) (on the diagonal).
	Each such fibre is contained in the corresponding fibre $ L ( \left\{ (p, l) \right\} ) $ of $ \rho^* U $: the conditions $ f_0 = f_1 = 0 $, $ f_2 = 0 = 0 \cdot L^2 $ and $ L \mid f_3 $ of \autoref{thm:tacnode_calc} hold.
	Thus the inclusion of total spaces $ \widehat W \subset \rho^* U $ is a continuous fibrewise-linear injection of vector bundles of constant ranks $ 7 \leq 9 $.
	Locally, a frame of $ \widehat W $ extends to a frame of $ \rho^* U $ (complete a frame at a point and shrink the neighbourhood), so $ \widehat W $ is a subbundle and the quotient $ \mathcal{N} $ is a vector bundle of rank $ 2 $.
\end{proof}

\subsection{The boundary morphism of the resolved pair}\label{sec:resolved_cohomology}

We now compute the cohomology of $ \widehat Y $, the action of $ \tau $ on it, and the kernel of $ \eta_* $.
The kernel is the image of the boundary map of the resolved pair (\autoref{prop:boundary_resolved}).
It consists of the classes that are anti-invariant under $ \tau $, and it has dimension $ 1 $ in each of the degrees $ 2 $, $ 4 $ and $ 6 $ (\autoref{lem:ker_eta}).
We also compute the integration along the fibre of $ \rho $, which gives the pushforward to $ P $ (\autoref{lem:gysin_rho}).
Throughout, $ H^*(P; \Q) = \Q[\varepsilon_0, \varepsilon_1, \varepsilon_2] / I $ as in \autoref{thm:ptcp2_cohom_structure}, where $ I $ is generated by the elementary symmetric polynomials $ s_1, s_2, s_3 $ in $ \varepsilon_0, \varepsilon_1, \varepsilon_2 $.
Recall from \autoref{sec:resolved_pair} the projections $ \psi \colon P \to \C P^2 $ and $ \rho \colon \widehat Y \to P $, the involution $ \tau $ of $ \widehat Y $, the quotient map $ \eta \colon \widehat Y \to \overline X_2 $, and the bundle $ \mathcal{Q} $ over $ \C P^2 $ with fibre $ \C^3 / \overline p $ over $ p $.

\begin{lemma}\label{lem:chern_dictionary}
	The three tautological line bundles over $ P $,
	\begin{equation*}
		\mathcal{S}_0 \text{ with fibre } \overline p , \qquad
		\mathcal{S}_1 \text{ with fibre } \overline l / \overline p , \qquad
		\mathcal{S}_2 \text{ with fibre } \C^3 / \overline l
	\end{equation*}
	(sub- and quotient bundles of the trivial bundle $ \C^3 \times P $), are homogeneous vector bundles, and
	\begin{equation*}
		c_1 ( \mathcal{S}_0 ) = \varepsilon_0 , \qquad
		c_1 ( \mathcal{S}_1 ) = \varepsilon_1 , \qquad
		c_1 ( \mathcal{S}_2 ) = \varepsilon_2 .
	\end{equation*}
\end{lemma}

\begin{proof}
	The group $ G = \SL_3(\C) $ acts on $ \C^3 \times P $ by $ g \cdot (v, (p,l)) = (g v, g (p,l)) $, and this action preserves each of the three families, linearly on the fibres.
	So each is a homogeneous vector bundle, and by \autoref{thm:hvb} it is isomorphic to $ G \times_B ( \text{fibre at } x_0 ) $, where $ x_0 $ is the standard flag $ p = (1:0:0) $, $ l = \overline{(1:0:0), (0:1:0)} $.
	At $ x_0 $ we have $ \overline p = \langle e_0 \rangle $ and $ \overline l = \langle e_0, e_1 \rangle $, and the torus element $ T(a, b) = \diag \left( (ab)^{-1}, a, b \right) $ acts on the three fibres $ \langle e_0 \rangle $, $ \overline l / \overline p = \langle \overline e_1 \rangle $, $ \C^3 / \overline l = \langle \overline e_2 \rangle $ by the characters
	\begin{equation*}
		\chi_0 (a, b) = a^{-1} b^{-1} , \qquad \chi_1 (a, b) = a , \qquad \chi_2 (a, b) = b .
	\end{equation*}
	By \autoref{cor:char_diff} and the identification $ c_1 = \dif \chi $ stated after it, $ c_1(\mathcal{S}_0) = - \varepsilon_1 - \varepsilon_2 = \varepsilon_0 $ (using $ \varepsilon_0 + \varepsilon_1 + \varepsilon_2 = 0 $), $ c_1(\mathcal{S}_1) = \varepsilon_1 $ and $ c_1(\mathcal{S}_2) = \varepsilon_2 $.
	The description of $ H^*(P; \Q) $ by the classes $ \varepsilon_i $ is that of \cite[Lemma 17]{BerGorKon23} (\autoref{thm:ptcp2_cohom_structure}), so the signs obtained here can be compared with those of \cite{BerGorKon23}.
	The fibre of $ \mathcal{S}_0 $ over $ (p, l) $ is the line $ \overline p $, so $ \mathcal{S}_0 $ is the pullback $ \psi^* \mathcal{O}(-1) $ of the tautological line bundle of $ \C P^2 $.
	The proof of \cite[Lemma 19]{BerGorKon23} gives $ c_1 ( \psi^* \mathcal{O}(-1) ) = \varepsilon_0 $, in agreement with the value of $ c_1 ( \mathcal{S}_0 ) $ found above.
\end{proof}

We will also need the multiplicative relations in $ H^*(P; \Q) $ in explicit form.
From $ s_1 = 0 $ we get $ \varepsilon_2 = - \varepsilon_0 - \varepsilon_1 $.
Substituting this into $ s_2 = 0 $ gives
\begin{equation}\label{eq:relation_eps1sq}
	\varepsilon_1^2 = - \varepsilon_0^2 - \varepsilon_0 \varepsilon_1 ,
\end{equation}
and substituting into $ s_3 = 0 $ gives $ \varepsilon_0 \varepsilon_1 ( \varepsilon_0 + \varepsilon_1 ) = 0 $, i.e. $ \varepsilon_0 \varepsilon_1^2 = - \varepsilon_0^2 \varepsilon_1 $.
Multiplying \autoref{eq:relation_eps1sq} by $ \varepsilon_0 $ and comparing, we obtain $ \varepsilon_0^3 = 0 $.
Modulo the ideal $ I $, the graded vector space $ H^*(P) $ is then
\begin{equation*}
	H^*(P) = \underbrace{\langle 1 \rangle}_{H^0(P)} \oplus \underbrace{\langle \varepsilon_0, \varepsilon_1 \rangle}_{H^2(P)} \oplus \underbrace{\langle \varepsilon_0^2, \varepsilon_0 \varepsilon_1 \rangle}_{H^4(P)} \oplus \underbrace{\langle \varepsilon_0^2 \varepsilon_1 \rangle}_{H^6(P)} ,
\end{equation*}
in agreement with the Betti numbers of $ P $ computed in the proof of \autoref{thm:spectral_groups} (the case $ E^1_{1,*} $).

\begin{lemma}\label{lem:cohom_yhat}
	Let $ \mathcal{S}' $ be the tautological line subbundle of $ \rho^* \psi^* \mathcal{Q} $ over $ \widehat Y = \PP(\psi^* \mathcal{Q}) $, i.e. the line bundle with fibre $ \overline{l'} / \overline p $ over $ ((p,l),(p,l')) $, and set
	\begin{equation*}
		\varepsilon_1' \coloneq c_1 ( \mathcal{S}' ) \in H^2 ( \widehat Y ; \Q ) .
	\end{equation*}
	In what follows we suppress the pullback $ \rho^* $ from the notation and write $ \varepsilon_i $ also for $ \rho^* \varepsilon_i $.
	\begin{enumerate}[(1)]
		\item $ H^* ( \widehat Y ; \Q ) $ is a free module over $ H^*(P; \Q) $ with basis $ 1, \varepsilon_1' $, and its ring structure is determined by the relation
		\begin{equation}\label{eq:relation_eps1prime}
			( \varepsilon_1' )^2 = - \varepsilon_0 \varepsilon_1' - \varepsilon_0^2 .
		\end{equation}
		In particular $ \dim_\Q H^{2k}(\widehat Y) = 1, 3, 4, 3, 1 $ for $ k = 0, \ldots, 4 $, the odd cohomology vanishes, and a $ \Q $-basis is given by the products of the basis monomials of $ H^*(P) $ above with $ 1 $ and $ \varepsilon_1' $.
		\item $ \tau^* $ fixes $ \varepsilon_0 $ and exchanges $ \varepsilon_1 \leftrightarrow \varepsilon_1' $.
		\item The anti-invariant part $ H^{2k} ( \widehat Y )^- \coloneq \left\{ u \mid \tau^* u = - u \right\} $ is one-dimensional for $ 2k = 2, 4, 6 $, spanned respectively by
		\begin{equation*}
			\varepsilon_1 - \varepsilon_1' , \qquad
			\varepsilon_0 ( \varepsilon_1 - \varepsilon_1' ) , \qquad
			\varepsilon_0^2 ( \varepsilon_1 - \varepsilon_1' ) ,
		\end{equation*}
		and zero for $ 2k = 0, 8 $.
	\end{enumerate}
\end{lemma}

\begin{proof}
	$ (1) $ The restriction of $ \mathcal{S}' $ to a fibre $ \rho^{-1}((p,l)) \cong \PP(\mathcal{Q}_p) \cong \C P^1 $ is the tautological line bundle of that projective line, so the restrictions of $ 1 $ and $ \varepsilon_1' $ form a basis of the cohomology of every fibre, and the Leray--Hirsch theorem \cite[Theorem 4D.1]{hatcher2002} applies.
	For the relation, consider the exact sequence of vector bundles over $ \widehat Y $
	\begin{equation*}
		0 \to \mathcal{S}' \to \psi^* \mathcal{Q} \to \psi^* \mathcal{Q} / \mathcal{S}' \to 0
	\end{equation*}
	(pullbacks suppressed).
	Over $ P $, the sequence $ 0 \to \mathcal{S}_0 \to \C^3 \to \psi^* \mathcal{Q} \to 0 $ and \autoref{lem:chern_dictionary} give $ c ( \psi^* \mathcal{Q} ) = ( 1 + \varepsilon_0 )^{-1} = 1 - \varepsilon_0 + \varepsilon_0^2 $ (using $ \varepsilon_0^3 = 0 $).
	Whitney's formula for the sequence over $ \widehat Y $ reads $ ( 1 + \varepsilon_1' ) \left( 1 + c_1 ( \psi^* \mathcal{Q} / \mathcal{S}' ) \right) = 1 - \varepsilon_0 + \varepsilon_0^2 $.
	In degree $ 1 $ it gives $ c_1 ( \psi^* \mathcal{Q} / \mathcal{S}' ) = - \varepsilon_0 - \varepsilon_1' $, and in degree $ 2 $ it gives $ \varepsilon_1' \cdot c_1 ( \psi^* \mathcal{Q} / \mathcal{S}' ) = \varepsilon_0^2 $.
	Together, $ \varepsilon_1' \cdot ( - \varepsilon_0 - \varepsilon_1' ) = \varepsilon_0^2 $, which is \autoref{eq:relation_eps1prime}.
	The dimension count follows from freeness.

	$ (2) $ As subbundles of the trivial bundle, $ \tau^* ( \rho^* \mathcal{S}_0 ) = \rho^* \mathcal{S}_0 $ (the fibre $ \overline p $ depends only on the common point), while $ \tau^* \mathcal{S}' $ has fibre $ \overline l / \overline p $ at $ ((p,l),(p,l')) $, i.e. $ \tau^* \mathcal{S}' = \rho^* \mathcal{S}_1 $, and symmetrically \mbox{$ \tau^* ( \rho^* \mathcal{S}_1 ) = \mathcal{S}' $}.
	Taking first Chern classes and using \autoref{lem:chern_dictionary} gives the claim.

	$ (3) $ On the $ \Q $-basis of $ (1) $, the involution $ \tau^* $ acts by the permutation exchanging $ \varepsilon_1 \leftrightarrow \varepsilon_1' $ (and expanding via \autoref{eq:relation_eps1sq} and \autoref{eq:relation_eps1prime} where needed).
	In degree $ 2 $, the basis $ \varepsilon_0, \varepsilon_1, \varepsilon_1' $ is permuted with $ \varepsilon_0 $ fixed, so the anti-invariants are spanned by $ \varepsilon_1 - \varepsilon_1' $.
	In degree $ 4 $, the basis $ \varepsilon_0^2, \varepsilon_0 \varepsilon_1, \varepsilon_0 \varepsilon_1', \varepsilon_1 \varepsilon_1' $ is permuted with $ \varepsilon_0^2 $ and $ \varepsilon_1 \varepsilon_1' $ fixed, so the anti-invariants are spanned by $ \varepsilon_0 ( \varepsilon_1 - \varepsilon_1' ) $.
	In degree $ 6 $, $ \tau^* $ exchanges $ \varepsilon_0^2 \varepsilon_1 \leftrightarrow \varepsilon_0^2 \varepsilon_1' $ and fixes $ \varepsilon_0 \varepsilon_1 \varepsilon_1' $, so the anti-invariants are spanned by $ \varepsilon_0^2 ( \varepsilon_1 - \varepsilon_1' ) $.
	In degrees $ 0 $ and $ 8 $ the action is trivial ($ \varepsilon_0^2 \varepsilon_1 \varepsilon_1' $ is symmetric).
\end{proof}

The next lemma describes the integration along the fibre of the $ \C P^1 $-bundle $ \rho \colon \widehat Y \to P $.

\begin{lemma}\label{lem:gysin_rho}
	Let $ PD $ denote the Poincar\'e duality isomorphisms of the compact complex manifolds $ \widehat Y $ and $ P $, and let
	\begin{equation*}
		\rho_! \coloneq PD^{-1} \circ \rho_* \circ PD \colon H^{k} ( \widehat Y ; \Q ) \to H^{k - 2} ( P ; \Q ) .
	\end{equation*}
	Then:
	\begin{enumerate}[(1)]
		\item $ \rho_! $ is $ H^*(P) $-linear: $ \rho_! ( \rho^* (u) \smile v ) = u \smile \rho_! (v) $;
		\item $ \rho_! (1) = 0 $, and $ c \coloneq \rho_! ( \varepsilon_1' ) \in H^0(P; \Q) = \Q $ is non-zero;
		\item $ \rho_! \left( ( \varepsilon_1' )^2 \right) = - c \, \varepsilon_0 $.
	\end{enumerate}
\end{lemma}

\begin{proof}
	$ (1) $ Using the projection formula $ \rho_* ( \rho^* (u) \frown z ) = u \frown \rho_* (z) $ and the equality \mbox{$ PD (u \smile v) = u \frown PD(v) $}:
	\begin{equation*}
		\rho_! ( \rho^* u \smile v ) = PD^{-1} \rho_* ( \rho^* u \frown PD (v) ) = PD^{-1} ( u \frown \rho_* PD (v) ) = u \smile \rho_! (v) .
	\end{equation*}

	$ (2) $ $ \rho_!(1) $ lies in $ H^{-2}(P) = 0 $.
	If $ c $ were zero then, by $ (1) $ and the freeness of \autoref{lem:cohom_yhat}(1), $ \rho_! $ would vanish identically, hence so would $ \rho_* = PD \circ \rho_! \circ PD^{-1} $ on homology.
	But $ \rho_* \colon H_0(\widehat Y) \to H_0(P) $ is an isomorphism, as both spaces are non-empty and path-connected.
	(With the standard orientation conventions one can check that $ c = -1 $, the degree of the tautological bundle on each fibre.
	Only $ c \neq 0 $ will be used below.)

	$ (3) $ By \autoref{eq:relation_eps1prime} and $ (1) $--$ (2) $:
	$ \rho_! ( ( \varepsilon_1' )^2 ) = \rho_! ( - \varepsilon_0 \varepsilon_1' - \varepsilon_0^2 ) = - \varepsilon_0 \, c - 0 = - c \, \varepsilon_0 $.
\end{proof}

\begin{lemma}\label{lem:ker_eta}
	For every $ k $ the map $ \eta_* \colon H_k ( \widehat Y ; \Q ) \to H_k ( \overline X_2 ; \Q ) $ is surjective, and
	\begin{equation*}
		\ker \eta_* = H_k ( \widehat Y ; \Q )^- \coloneq \left\{ z \mid \tau_* z = - z \right\} ,
	\end{equation*}
	which corresponds under Poincar\'e duality to the anti-invariant cohomology of \autoref{lem:cohom_yhat}(3).
	In particular $ \ker \eta_* $ is one-dimensional in degrees $ 2, 4, 6 $ and zero otherwise.
\end{lemma}

\begin{proof}
	Since $ \tau $ is a biholomorphic involution, it preserves the orientation, so $ \tau_* [\widehat Y] = [\widehat Y] $.
	By the projection formula, $ \tau_* PD ( \tau^* u ) = \tau_* \left( \tau^* u \frown [\widehat Y] \right) = u \frown \tau_* [\widehat Y] = PD (u) $.
	Replacing $ u $ by $ \tau^* v $ and using $ ( \tau^* )^2 = \mathrm{id} $ gives $ \tau_* PD (v) = PD ( \tau^* v ) $.
	Hence $ PD $ exchanges the $ (\pm 1) $-eigenspaces of $ \tau^* $ on cohomology with those of $ \tau_* $ on homology.
	That is, if $ \tau^* v = \pm v $ then $ \tau_* PD (v) = \pm PD (v) $ with the same sign, so $ PD $ maps $ H^k ( \widehat Y )^- $ isomorphically onto $ H_{8 - k} ( \widehat Y )^- $.
	The dimension statement now follows from \autoref{lem:cohom_yhat}(3).

	The inclusion $ H_k ( \widehat Y )^- \subseteq \ker \eta_* $ is immediate from $ \eta \circ \tau = \eta $:
	for $ \tau_* z = - z $ we get $ \eta_* z = \eta_* \tau_* z = - \eta_* z $, so $ \eta_* z = 0 $.

	For the reverse inclusion and the surjectivity we use the smooth structure of $ \overline X_2 = \PP( \Sym^2 \mathcal{Q} ) $ and the map $ \eta_! \coloneq PD^{-1} \circ \eta_* \circ PD $, defined by the same formula as $ \rho_! $ in \autoref{lem:gysin_rho}.
	First, $ \eta_* [ \widehat Y ] = 2 [ \overline X_2 ] $.
	To see this, let $ D \coloneq \eta ( \Delta(P) ) \subset \overline X_2 $ be the branch locus, a closed subvariety of real dimension $ 6 $, and $ V \coloneq \overline X_2 \setminus D $.
	In the long exact sequence of Borel--Moore homology of the pair, $ \overline H_8 (D) = \overline H_7 (D) = 0 $ for dimension reasons, so the restriction $ H_8 ( \overline X_2 ) \to \overline H_8 ( V ) $ is an isomorphism.
	Over $ V $ the map $ \eta $ restricts to a two-sheeted covering $ \widehat Y \setminus \Delta(P) \to V $ of connected oriented $ 8 $-manifolds, and the pushforward of the Borel--Moore fundamental class of a two-sheeted covering is twice the fundamental class of the base.
	Since $ [\widehat Y] $ and $ [\overline X_2] $ restrict to the fundamental classes of $ \widehat Y \setminus \Delta(P) $ and $ V $, the claim follows.

	Consequently, for every $ u \in H^* ( \overline X_2 ) $,
	\begin{equation*}
		\eta_! \eta^* (u) = PD^{-1} \eta_* ( \eta^* u \frown [ \widehat Y ] ) = PD^{-1} ( u \frown \eta_* [ \widehat Y ] ) = PD^{-1} ( u \frown 2 [ \overline X_2 ] ) = 2 u .
	\end{equation*}
	Hence $ \eta_! $ is surjective (equivalently, $ \eta_* $ is surjective).

	For the reverse inclusion, the surjectivity of $ \eta_! $ gives
	\begin{equation*}
		\dim \ker \eta_! = \dim H^* ( \widehat Y ) - \dim H^* ( \overline X_2 ) = 12 - 9 = 3
	\end{equation*}
	by \autoref{lem:cohom_yhat}(1) and \autoref{cor:xbar_cohomology}.
	The anti-invariant part $ H^* ( \widehat Y )^- $, of the same total dimension $ 3 $, is contained in $ \ker \eta_! $: this is the $ PD $-conjugate of the inclusion proved above.
	Therefore \mbox{$ \ker \eta_! = H^* ( \widehat Y )^- $}, and conjugating by $ PD $ once more gives $ \ker \eta_* = H_* ( \widehat Y )^- $.
\end{proof}

We can now describe the boundary map of the resolved pair.

\begin{proposition}\label{prop:boundary_resolved}
	Let $ \partial_{( \widehat \Phi, \widehat Y )} \colon \overline H_{j+1} ( \widehat \Phi \setminus \widehat Y ) \to H_j ( \widehat Y ) $ be the connecting morphism of the pair $ ( \widehat \Phi, \widehat Y ) $ in Borel--Moore homology.
	\begin{enumerate}[(1)]
		\item $ \overline H_{2k} ( \widehat \Phi \setminus \widehat Y ) = 0 $ for every $ k $.
		\item For every $ k $ the map $ \partial_{( \widehat \Phi, \widehat Y )} $ is an isomorphism of $ \overline H_{2k+1} ( \widehat \Phi \setminus \widehat Y ) $ onto $ \ker \eta_* = H_{2k} ( \widehat Y )^- $.
		\item $ \overline H_* ( \widehat \Phi \setminus \widehat Y ) $ is $ \Q $ in the degrees $ 3 $, $ 5 $, $ 7 $ and zero otherwise.
		The image of $ \overline H_5 ( \widehat \Phi \setminus \widehat Y ) $ is spanned by $ PD \left( \varepsilon_0 ( \varepsilon_1 - \varepsilon_1' ) \right) $, and the image of $ \overline H_7 ( \widehat \Phi \setminus \widehat Y ) $ is spanned by $ PD ( \varepsilon_1 - \varepsilon_1' ) $.
	\end{enumerate}
\end{proposition}

\begin{proof}
	Both $ \widehat \Phi $ and $ \widehat Y $ are compact, so their Borel--Moore homology is ordinary homology, and the long exact sequence of the pair reads
	\begin{equation*}
		H_{j+1} ( \widehat \Phi ) \to \overline H_{j+1} ( \widehat \Phi \setminus \widehat Y ) \xrightarrow{\ \partial\ } H_j ( \widehat Y ) \xrightarrow{\ \iota_* \ } H_j ( \widehat \Phi ) ,
	\end{equation*}
	where $ \iota \colon \widehat Y \hookrightarrow \widehat \Phi $ is the inclusion.
	The odd homology of $ \widehat \Phi $ and of $ \widehat Y $ vanishes (\autoref{cor:phihat_homology}, \autoref{lem:cohom_yhat}(1)).
	Since $ \widehat \pi_* $ is an isomorphism (\autoref{cor:phihat_homology}) and $ \widehat \pi \circ \iota = \eta $ (by the identification of $ \widehat Y $ with a subspace of $ \widehat \Phi $), the map $ \iota_* $ has the kernel of $ \eta_* $, and it is surjective because $ \eta_* $ is (\autoref{lem:ker_eta}).

	$ (1) $ The group $ H_{2k-1} ( \widehat Y ) $ is zero, so $ \overline H_{2k} ( \widehat \Phi \setminus \widehat Y ) $ is the cokernel of $ \iota_* \colon H_{2k} ( \widehat Y ) \to H_{2k} ( \widehat \Phi ) $, which is zero.

	$ (2) $ For $ j = 2k $ the group on the left of the sequence vanishes, so $ \partial $ is injective.
	By exactness its image is $ \ker \iota_* = \ker \eta_* $, which is $ H_{2k} ( \widehat Y )^- $ by \autoref{lem:ker_eta}.

	$ (3) $ By \autoref{lem:ker_eta}, the space $ H_{2k} ( \widehat Y )^- $ is one-dimensional for $ 2k = 2, 4, 6 $ and zero otherwise, and it is the image of $ H^{8 - 2k} ( \widehat Y )^- $ under $ PD $.
	The classes that span $ H^4 ( \widehat Y )^- $ and $ H^2 ( \widehat Y )^- $ are given in \autoref{lem:cohom_yhat}(3).
\end{proof}

\begin{remark}\label{rem:2to1_consistency}
	By \autoref{prop:boundary_resolved}(3) and \autoref{lem:phihat}(5), $ \overline H_* ( \Phi_{\alpha(2)} \setminus \Phi_{\alpha(1)} ) $ is $ \Q $ in the degrees $ 3 $, $ 5 $, $ 7 $ and zero otherwise, i.e. $ E^1_{2,q} = \Q $ exactly for $ q = 15, 17, 19 $.
	So the resolved pair reproduces the second column of \autoref{table:main_sp_seq_p1} independently of the computation of that column in the proof of \autoref{thm:spectral_groups}.
\end{remark}

\subsection{The cap product with the Euler class}\label{sec:map_alpha}

The Euler class enters the computation through the following diagram.
By \cite[Corollary 3]{BerGorKon23}, for the subbundle $ \widehat W \subset \rho^* U $ of \autoref{lem:res_bundles} over the smooth compact base $ \widehat Y $, the map induced on Borel--Moore homology by the inclusion of total spaces satisfies
\begin{equation*}
	\begin{tikzcd}[column sep=small]
		\rho^* U \arrow[rrd, "\C^9"] &  &            & \overline H_*(\rho^* U) \arrow[rr, "Th^{-1}"] &  & H_{*-18}(\widehat Y) \arrow[rr, "PD^{-1}"] &  & H^{26-*}(\widehat Y) \\
		&  & \widehat Y &                                               &  &                                            &  &                      \\
		\widehat W \arrow[uu, hook] \arrow[rru, "\C^7"] &  &            & \overline H_*(\widehat W) \arrow[uu] \arrow[rr, "Th^{-1}"] &  & H_{*-14}(\widehat Y) \arrow[rr, "PD^{-1}"] \arrow[uu, "-\frown e(\mathcal{N})"] &  & H^{22-*}(\widehat Y) \arrow[uu, "-\smile e(\mathcal{N})"]
	\end{tikzcd}
\end{equation*}
That is, under the Thom and Poincar\'e duality isomorphisms the inclusion becomes the cup product with the (equivariant, in the homogeneous situations below) Euler class $ e ( \mathcal{N} ) $ of the quotient bundle $ \mathcal{N} = \rho^* U / \widehat W $.
This subsection computes the Euler classes of the two quotient bundles that appear in this way: the class $ e ( U / W ) $ over $ P $, which is the diagonal restriction, and the class $ e ( \mathcal{N} ) $ over $ \widehat Y $.
The first is a homogeneous vector bundle, and its Euler class is read off from the characters of the torus (\autoref{thm:euler_class_via_hvb}).
The second is an extension of two line bundles.
One of them is pulled back from $ P $, and the first Chern class of the other is determined by its restrictions to the diagonal and to a fibre of $ \rho $ (\autoref{lem:euler_resolved}).
At the end we cap the anti-invariant classes with $ e ( \mathcal{N} ) $ and push the result forward to $ P $ (\autoref{prop:cap_pushforward}).

Recall that $ U = F_{\alpha(1)} $ is a vector bundle of rank $ 9 $ over $ P \cong \Phi_{\alpha(1)} $, and recall the rank-$ 7 $ bundle $ W $ over $ P $ defined after \autoref{cor:config_lim_ext}.
Its fibres $ \widehat L ( \left\{ (p,l) \right\} ) $ are subspaces of the fibres $ L ( \left\{ (p,l) \right\} ) $ of $ U $ (the case $ \lambda = 0 $ of \autoref{thm:tacnode_calc}), so $ W \subset U $.
Let
\begin{equation*}
	W' \coloneq \left\{ ( f, (p, l) ) \in \Pi_{4,2} \times P \mid f_0 = f_1 = f_2 = 0 \text{ and } L \text{ divides } f_3 \right\}
\end{equation*}
be the family of the jump fibres of \autoref{lem:closure_jump}.

\begin{lemma}\label{lem:hvb_chain}
	$ W' $ is the total space of a rank-$ 8 $ subbundle of $ U $ containing $ W $, and the bundles $ U \supset W' \supset W $ and the quotients $ U / W' $, $ W' / W $ and $ U / W $ are homogeneous vector bundles with respect to the action
	\begin{equation*}
		g \cdot ( f, (p, l) ) \coloneq \left( f \circ g^{-1} , \ g \cdot (p, l) \right) , \qquad g \in G = \SL_3(\C) .
	\end{equation*}
\end{lemma}

\begin{proof}
	The action is well defined on $ U $: looking back at \autoref{thm:tacnode_calc}, all conditions are tautologically satisfied after a projective change of coordinates, i.e. $ \Sing ( f \circ g^{-1} ) = g \cdot \Sing (f) $ (cf. the proof of \autoref{thm:dimensions}).
	The same equivariance shows that the action preserves $ W' $: its fibre consists of the polynomials for which $ p $ is a point of multiplicity at least $ 3 $ whose tangent cone contains $ l $, an intrinsic condition.
	It preserves $ W $ as well: on $ X_2 $ we have $ L ( g K ) = g \cdot L(K) $, and by the continuity of $ \widehat L $ (\autoref{cor:config_lim_ext}) and of the $ G $-action this passes to the limit spaces, $ \widehat L ( g K ) = g \cdot \widehat L (K) $ for $ K \in X_1 $.

	The fibre dimensions are constant ($ 9 $, $ 8 $ and $ 7 $): they are computed at the standard flag in \autoref{thm:dimensions}, \autoref{lem:closure_jump} and \autoref{cor:config_lim_ext} respectively, and $ G $ acts transitively on $ P $.

	So far each of $ U $, $ W' $ and $ W $ is a $ G $-invariant subset of $ \Pi_{4,2} \times P $ whose fibres over $ P $ are vector subspaces of $ \Pi_{4,2} $.
	We now show that such a subset is a vector bundle, and a homogeneous one.

	Let $ E $ be one of the three, and let $ x_0 \coloneq (p, l) \in P $ be the flag with $ p = (1:0:0) $ and $ l = \overline{(1:0:0), (0:1:0)} $.
	The stabiliser of $ x_0 $ in $ G $ is $ B $ (this is how $ P $ is identified with $ G / B $), so $ B $ preserves the fibre $ E_{x_0} $, and $ E_{x_0} $ is a representation of $ B $.

	Consider the map
	\begin{equation*}
		G \times_B E_{x_0} \to E , \qquad [g, f] \mapsto \left( f \circ g^{-1} , \ g \cdot x_0 \right) .
	\end{equation*}
	It is well defined, because $ ( g b^{-1} , f \circ b^{-1} ) $ and $ (g, f) $ have the same image for every $ b \in B $, and it takes values in $ E $, because $ E $ is $ G $-invariant.

	It is continuous and commutes with the projections to $ P = G / B $.
	On the fibre over $ g \cdot x_0 $ it is the linear map $ f \mapsto f \circ g^{-1} $ from $ E_{x_0} $ to $ E_{g \cdot x_0} $, which is a bijection with inverse $ f \mapsto f \circ g $, again by the $ G $-invariance of $ E $.
	So the map is a bijection.

	Its inverse is continuous: over every open subset of $ P $ on which the projection $ G \to G / B $ has a section $ s $, the inverse is $ (f, x) \mapsto [ s(x) , f \circ s(x) ] $.
	The map is also $ G $-equivariant for the action $ h \cdot [g, f] = [hg, f] $ on the left-hand side.
	Hence $ E $ is a vector bundle isomorphic to $ G \times_B E_{x_0} $, that is, a homogeneous vector bundle.

	Quotients of homogeneous bundles by homogeneous subbundles are again homogeneous (induced by the quotient representation).
\end{proof}

Recall from \autoref{thm:ptcp2_cohom_structure} that $ H^*(P; \Q) = \Q[\varepsilon_0, \varepsilon_1, \varepsilon_2] / I $.
The next statement is the diagonal case of the Euler-class computation.

\begin{proposition}\label{thm:euler_class_via_hvb}
	The Euler class of $ U / W $ is $ e ( U / W ) = - 2 \varepsilon_1^2 $.
\end{proposition}

\begin{proof}
	Let $ p = (1:0:0), l = \overline{(1:0:0), (0:1:0)} $.
	We now show that the fibre of $ U / W $ over $ (p, l) $ is generated by monomials $ xy^2z, x^2z^2 $.
	To do so, it suffices to consider the diagram
	\begin{equation*}
		\small
		\begin{tikzcd}[column sep=small]
			W_{(p,l)} \arrow[r, equal] \arrow[d]                & {\left \langle \begin{array}{c} xyz^2 , xz^3 , \\ y^4 , y^3z, y^2z^2 , yz^3 , \\ z^4  \end{array} \right \rangle } \arrow[r, equal] \arrow[d]                 & {\widehat L ( \left\{ (p,l) \right\} )} \arrow[r] \arrow[d] & W \arrow[rd] \arrow[d] &               \\
			{U_{(p,l)}} \arrow[r, equal] \arrow[d] & {\left \langle \begin{array}{c} x^2z^2 , xy^2z , xyz^2 , xz^3 , \\ y^4 , y^3z , y^2z^2 , yz^3, \\ z^4 \end{array} \right \rangle } \arrow[r, equal] \arrow[d] & {L(\left \{ (p,l) \right\})} \arrow[r] \arrow[d]          & U \arrow[r] \arrow[d]         & {P \ni (p,l)} \\
			{U_{(p,l)} / W_{(p,l)}} \arrow[r, equal]       & {\langle  x^2z^2 , xy^2z  \rangle } \arrow[r, equal]                                                                                                        & {L(\left \{ (p,l) \right\}) / \widehat L ( \left\{ (p,l) \right\} )} \arrow[r]     & U / W \arrow[ru]     &
		\end{tikzcd}
	\end{equation*}
	where the fibre $ W_{(p,l)} = \widehat L ( \left\{ (p,l) \right\} ) $ is the limit space computed in the proof of \autoref{lem:config_lim} (see \autoref{cor:config_lim_ext}).

	The group $ G = \SL_3(\C) $ acts on the space of polynomials via $ (g . f)(v) = f(g^{-1}v) $.
	For \mbox{$ g = \diag(\lambda_0, \lambda_1, \lambda_2) $} the monomial $ x^n y^m z^k $ spans a character of weight $ \lambda_0^{-n} \lambda_1^{-m} \lambda_2^{-k} $.
	Its differential is
	\mbox{$ \dif \chi = -(n \varepsilon_0 + m \varepsilon_1 + k \varepsilon_2) $}.

	By \autoref{lem:hvb_chain}, $ U / W $ is a homogeneous vector bundle.
	This means that we can apply \autoref{thm:hvb}.
	Recall that taking $ (p, l) $ as above, it is easy to check that $ B $ stabilises $ (p, l) $, and the isomorphism $ P = G / B $ maps $ (p, l) $ to $ x_0 \coloneq e \cdot B $.
	From this, $ U / W $ is determined by its fibre at $ x_0 $ as a $ B $-representation, that is, by the representation $ R \colon B \to \GL(F) $ on \mbox{$ F \coloneq {\langle  x^2z^2 , xy^2z  \rangle } $}.

	Now, we can find the Euler class of HVB $ U / W $ via techniques described in \autoref{sec:hvb}.
	Namely, it has an obvious decomposition into line bundles.
	By finding their corresponding characters and their differentials, we will find the complete Chern class of the vector bundle.

	Let us find the characters of the representation corresponding to this action.

	From the definition of the action, for the generating polynomial $ f = xy^2z $ and a group element $ g = T(s,t) \in G $,
	\begin{equation*}
		T(s,t)(xy^2z) = (stx)(s^{-1}y)^2(t^{-1}z) = s^{-1}xy^2z = \chi_1(s,t) xy^2z .
	\end{equation*}
	Thus $ \chi_1(s,t) = s^{-1} $.
	Consequently, by \autoref{cor:char_diff} $ \dif \chi_1 = -\varepsilon_1 $.

	Similarly,
	\begin{equation*}
		T(s,t)(x^2z^2) = s^2x^2z^2 = \chi_2(s,t) x^2z^2 \implies \dif \chi_2 = 2\varepsilon_1 .
	\end{equation*}
	Hence, by \autoref{cor:chern_prod}
	\begin{equation*}
		c(U / W) = (1 - \varepsilon_1)(1 + 2\varepsilon_1) = 1 + \varepsilon_1 - 2\varepsilon_1^2 \in \Q[\varepsilon_0, \varepsilon_1, \varepsilon_2] / I .
	\end{equation*}
	Since $ U / W $ has rank $2$, its Euler class coincides with the top Chern class:
	\[
	e(U / W) = c_2(U / W) = - 2\varepsilon_1^2 .
	\]
\end{proof}

Two by-products of this proof will be used below: the sub-line $ \langle x y^2 z \rangle $ of $ F $ maps isomorphically onto the fibre of the line bundle $ W' / W $ at $ x_0 $, and the quotient line $ \langle x^2 z^2 \rangle $ onto the fibre of $ U / W' $.
Since the characters computed above are those of the corresponding $ T $-eigenlines, we get, by \autoref{cor:char_diff} applied to these homogeneous line bundles (\autoref{lem:hvb_chain}),
\begin{equation}\label{eq:c1_diag_lines}
	c_1 ( W' / W ) = - \varepsilon_1 , \qquad
	c_1 ( U / W' ) = 2 \varepsilon_1 .
\end{equation}

\begin{lemma}\label{lem:euler_resolved}
	Over $ \widehat Y $ there is a short exact sequence of vector bundles
	\begin{equation*}
		0 \to \rho^* W' / \widehat W \to \mathcal{N} \to \rho^* ( U / W' ) \to 0 ,
	\end{equation*}
	with
	\begin{equation*}
		c_1 \left( \rho^* ( U / W' ) \right) = 2 \varepsilon_1 , \qquad
		c_1 \left( \rho^* W' / \widehat W \right) = \varepsilon_1 - 2 \varepsilon_1' ,
	\end{equation*}
	and consequently
	\begin{equation*}
		e ( \mathcal{N} ) = c_2 ( \mathcal{N} ) = 2 \varepsilon_1 ( \varepsilon_1 - 2 \varepsilon_1' ) .
	\end{equation*}
\end{lemma}

\begin{proof}
	Fibrewise, $ \widehat W \subseteq \rho^* W' $: the conditions $ f_0 = f_1 = f_2 = 0 $ and $ L L' \mid f_3 $ (or $ L^2 \mid f_3 $ on the diagonal) imply $ L \mid f_3 $.
	Both inclusions $ \widehat W \subseteq \rho^* W' \subseteq \rho^* U $ are subbundles as in the proof of \autoref{lem:res_bundles}(3), so the displayed sequence exists, with quotients of ranks $ 8 - 7 = 1 $ and $ 9 - 8 = 1 $.
	The first Chern class of $ \rho^* (U / W') $ is $ 2 \varepsilon_1 $ by \autoref{eq:c1_diag_lines}.

	Write $ \beta \coloneq c_1 ( \rho^* W' / \widehat W ) = a \varepsilon_0 + b \varepsilon_1 + c' \varepsilon_1' $ in the basis of $ H^2 ( \widehat Y ) $ from \autoref{lem:cohom_yhat}(1).
	We determine the coefficients by two restrictions.

	\emph{Restriction to the diagonal.}
	The diagonal embedding $ \Delta \colon P \to \widehat Y $ satisfies $ \rho \circ \Delta = \mathrm{id}_P $, and $ \Delta^* \widehat W = W $: both sides are the pullbacks of $ \gamma^7 $ under the same map $ (p, l) \mapsto \widehat L ( \left\{ (p,l) \right\} ) $, by \autoref{lem:res_bundles}(3).
	Hence $ \Delta^* ( \rho^* W' / \widehat W ) = W' / W $, whose first Chern class is $ - \varepsilon_1 $ by \autoref{eq:c1_diag_lines}.
	On classes, $ \Delta^* \varepsilon_0 = \varepsilon_0 $, $ \Delta^* \varepsilon_1 = \varepsilon_1 $, and $ \Delta^* \varepsilon_1' = \varepsilon_1 $ (the restriction of $ \mathcal{S}' $ to the diagonal is $ \mathcal{S}_1 $).
	So $ a \varepsilon_0 + ( b + c' ) \varepsilon_1 = - \varepsilon_1 $ in $ H^2(P) $, giving $ a = 0 $ and $ b + c' = -1 $.

	\emph{Restriction to a fibre.}
	Fix the standard flag $ (p, l) $ as in the proof of \autoref{thm:euler_class_via_hvb} and its chart $ \{ X = 1 \} $ with the coordinates $ (y, z) $, so that $ p = (0, 0) $ and $ L = z $.
	The map $ (1 : y : z) \mapsto (0, y, z) + \overline p $ identifies this chart with the vector space $ \C^3 / \overline p $.
	It sends $ p $ to the origin and a line $ l' $ through $ p $ to the line $ \overline{l'} / \overline p $.
	So the parts $ f_i $ of the Taylor expansion of a polynomial $ f $ around $ p $ are forms of degree $ i $ on $ \C^3 / \overline p $, and an equation $ L' $ of $ l' $ in the chart is a linear form on $ \C^3 / \overline p $ with kernel $ \overline{l'} / \overline p $.

	Consider the fibre $ \rho^{-1} ( (p,l) ) \cong \C P^1 $, over which $ \rho^* W' $ is the trivial bundle with fibre $ W'_{(p,l)} $.
	For $ f \in W'_{(p,l)} $ write $ f_3 = L \, q_f $ with $ q_f $ a quadratic form on $ \C^3 / \overline p $.
	Define, for each point $ ((p,l), (p,l')) $ of the fibre,
	\begin{equation*}
		\Theta ( f ) \coloneq q_f |_{\overline{l'} / \overline p} \in \Sym^2 \left( ( \overline{l'} / \overline p )^{\vee} \right) = \left( ( \mathcal{S}' )^{\vee} \right)^{\otimes 2} \Big|_{((p,l),(p,l'))} .
	\end{equation*}
	As the point varies, $ \Theta $ is a continuous map of vector bundles over $ \rho^{-1} ( (p,l) ) $ from the trivial bundle with fibre $ W'_{(p,l)} $ to $ ( ( \mathcal{S}' )^{\vee} )^{\otimes 2} $.

	The kernel of $ \Theta $ at the point $ ((p,l),(p,l')) $ is
	\begin{equation*}
		\begin{aligned}
			\left\{ f \in W'_{(p,l)} \mid q_f \text{ vanishes on } \overline{l'} / \overline p \right\}
			& = \left\{ f \in W'_{(p,l)} \mid L' \text{ divides } q_f \right\} \\
			& = \left\{ f \in W'_{(p,l)} \mid L L' \text{ divides } f_3 \right\}
			= \widehat W_{((p,l),(p,l'))} .
		\end{aligned}
	\end{equation*}
	The first equality holds because a binary quadratic form that vanishes on a line through the origin is divisible by the equation of this line.
	Moreover, $ \Theta $ is surjective over every point of the fibre, because for every line through the origin there is a quadratic form that does not vanish on it.
	Hence $ \Theta $ induces an isomorphism
	\begin{equation*}
		\left( \rho^* W' / \widehat W \right) \Big|_{\rho^{-1}((p,l))} \cong \left( ( \mathcal{S}' )^{\vee} \right)^{\otimes 2} \Big|_{\rho^{-1}((p,l))} .
	\end{equation*}

	We now restrict $ \beta = a \varepsilon_0 + b \varepsilon_1 + c' \varepsilon_1' $ to the fibre $ \rho^{-1} ( (p,l) ) $.
	The classes $ \varepsilon_0 $ and $ \varepsilon_1 $ are pulled back from $ P $ along $ \rho $, and $ \rho $ is constant on the fibre, so they restrict to zero.
	Hence $ \beta $ restricts to $ c' \cdot ( \varepsilon_1' |_{\text{fibre}} ) $.

	On the other hand, by the isomorphism above, $ \beta $ restricts to the first Chern class of $ ( ( \mathcal{S}' )^{\vee} )^{\otimes 2} $ on the fibre.
	The first Chern class of a line bundle changes sign under dualisation and is additive under tensor products, so this class is the restriction of $ - 2 c_1 ( \mathcal{S}' ) = - 2 \varepsilon_1' $.
	Therefore $ c' \cdot ( \varepsilon_1' |_{\text{fibre}} ) = - 2 \, ( \varepsilon_1' |_{\text{fibre}} ) $.

	The fibre $ \rho^{-1} ( (p,l) ) $ is the projective line $ \PP ( \C^3 / \overline p ) $, whose points are the lines $ \overline{l'} / \overline p $, and the restriction of $ \mathcal{S}' $ to it is the tautological line bundle.
	So $ \varepsilon_1' |_{\text{fibre}} $ is the first Chern class of the tautological line bundle of $ \C P^1 $, which generates $ H^2 ( \C P^1 ) $ and in particular is not zero.
	This gives $ c' = -2 $, and then $ b = 1 $ from $ b + c' = -1 $.

	Therefore $ \beta = \varepsilon_1 - 2 \varepsilon_1' $, and by the Whitney formula
	$ c ( \mathcal{N} ) = ( 1 + 2 \varepsilon_1 ) ( 1 + \varepsilon_1 - 2 \varepsilon_1' ) $, whose degree-$ 2 $ part is $ c_2 ( \mathcal{N} ) = 2 \varepsilon_1 ( \varepsilon_1 - 2 \varepsilon_1' ) $.
	For a complex bundle the Euler class is the top Chern class.
\end{proof}

\begin{remark}\label{rem:euler_consistency}
	As a coherence check: $ \Delta^* \mathcal{N} = U / W $ (diagonal restriction of \autoref{lem:res_bundles}(3)), and correspondingly
	\begin{equation*}
		\Delta^* e ( \mathcal{N} ) = 2 \varepsilon_1 ( \varepsilon_1 - 2 \varepsilon_1 ) = - 2 \varepsilon_1^2 = e ( U / W ) ,
	\end{equation*}
	in agreement with \autoref{thm:euler_class_via_hvb}.
	This is not evidence independent of \autoref{lem:euler_resolved}.
	The diagonal restriction was already used there to pin down $ \beta $.
	But the check confirms that the two computations combine consistently.
\end{remark}

It remains to cap the anti-invariant classes with $ e ( \mathcal{N} ) $ and to push the result forward to $ P $.

\begin{proposition}\label{prop:cap_pushforward}
	For every $ u \in H^* ( \widehat Y ) $,
	\begin{equation*}
		\rho_* \left( PD ( u ) \frown e ( \mathcal{N} ) \right) = PD \left( \rho_! ( u \smile e ( \mathcal{N} ) ) \right) .
	\end{equation*}
	Moreover,
	\begin{align*}
		\rho_! \left( ( \varepsilon_1 - \varepsilon_1' ) \smile e ( \mathcal{N} ) \right) & = c \left( 6 \varepsilon_0^2 + 2 \varepsilon_0 \varepsilon_1 \right) , \\
		\rho_! \left( \varepsilon_0 ( \varepsilon_1 - \varepsilon_1' ) \smile e ( \mathcal{N} ) \right) & = 2 c \, \varepsilon_0^2 \varepsilon_1 ,
	\end{align*}
	where $ c \neq 0 $ is the number of \autoref{lem:gysin_rho}(2), and both classes are non-zero.
\end{proposition}

\begin{proof}
	For $ u \in H^* ( \widehat Y ) $ we have $ PD ( u ) \frown e ( \mathcal{N} ) = PD ( u \smile e ( \mathcal{N} ) ) $ and $ \rho_* PD = PD \, \rho_! $ (\autoref{lem:gysin_rho}).
	This gives the first equality.

	By \autoref{lem:euler_resolved}, $ e ( \mathcal{N} ) = 2 \varepsilon_1 ( \varepsilon_1 - 2 \varepsilon_1' ) $.
	For $ u = \varepsilon_1 - \varepsilon_1' $ we expand the product in powers of $ \varepsilon_1' $:
	\begin{equation*}
		( \varepsilon_1 - \varepsilon_1' ) \smile 2 \varepsilon_1 ( \varepsilon_1 - 2 \varepsilon_1' ) = 2 \varepsilon_1^3 - 6 \varepsilon_1^2 \varepsilon_1' + 4 \varepsilon_1 ( \varepsilon_1' )^2 .
	\end{equation*}
	The map $ \rho_! $ is $ H^*(P) $-linear, and $ \rho_! (1) = 0 $, $ \rho_! ( \varepsilon_1' ) = c $, $ \rho_! \left( ( \varepsilon_1' )^2 \right) = - c \, \varepsilon_0 $ with $ c \neq 0 $ (\autoref{lem:gysin_rho}).
	Applying it term by term, we get
	\begin{equation*}
		\rho_! \left( ( \varepsilon_1 - \varepsilon_1' ) \smile e ( \mathcal{N} ) \right) = - 6 c \, \varepsilon_1^2 - 4 c \, \varepsilon_0 \varepsilon_1 .
	\end{equation*}
	By \autoref{eq:relation_eps1sq}, the right-hand side equals
	\begin{equation*}
		- 6 c \left( - \varepsilon_0^2 - \varepsilon_0 \varepsilon_1 \right) - 4 c \, \varepsilon_0 \varepsilon_1 = c \left( 6 \varepsilon_0^2 + 2 \varepsilon_0 \varepsilon_1 \right) ,
	\end{equation*}
	which is non-zero because $ \varepsilon_0^2 $ and $ \varepsilon_0 \varepsilon_1 $ form a basis of $ H^4 ( P ) $.

	For $ u = \varepsilon_0 ( \varepsilon_1 - \varepsilon_1' ) $, by $ H^*(P) $-linearity and $ \varepsilon_0^3 = 0 $,
	\begin{equation*}
		\rho_! \left( \varepsilon_0 ( \varepsilon_1 - \varepsilon_1' ) \smile e ( \mathcal{N} ) \right) = c \, \varepsilon_0 \left( 6 \varepsilon_0^2 + 2 \varepsilon_0 \varepsilon_1 \right) = 2 c \, \varepsilon_0^2 \varepsilon_1 ,
	\end{equation*}
	which is non-zero because $ \varepsilon_0^2 \varepsilon_1 $ spans $ H^6 ( P ) $.
\end{proof}

\subsection{The differentials as a composition of three maps}\label{sec:two_to_one_diff}

This subsection proves \autoref{thm:2to1_nontrivial}.
\autoref{lem:transport} shows that the connecting morphism built from the resolved pair is the differential $ d^1 $.
\autoref{prop:d1_composition} then writes $ d^1 $ as the composition of three maps, and the theorem follows from \autoref{prop:boundary_resolved} and \autoref{prop:cap_pushforward}.
\autoref{cor:2to1_einfty} lists what survives in the first two columns.

We now assemble the pieces, applying the construction of \autoref{sec:diff_vs_leq} to the resolved pair $ ( \widehat \Phi, \widehat Y ) $ of \autoref{sec:resolved_pair} and to the bundles $ \mathcal{T} $ and $ \widehat W $ of \autoref{lem:res_bundles}.
Take
\begin{equation*}
	(X, Y) = ( \widehat \Phi, \widehat Y ) , \qquad
	E_1 = \mathcal{T} , \qquad
	E_2 = \rho^* U , \qquad
	E_2' = \widehat W = \mathcal{T} |_{\widehat Y} \subset \rho^* U ,
\end{equation*}
with the identity as the glueing isomorphism $ E_1 |_Y \cong E_2' $ (all bundles are complex, so canonically oriented), and let
\begin{equation*}
	Z \coloneq \mathcal{T} \cup_{\widehat W} \rho^* U
\end{equation*}
be the glued space, as in \autoref{sec:diff_vs_leq}.

\begin{lemma}\label{lem:transport}
	There is a proper map of pairs
	\begin{equation*}
		h \colon ( Z, \rho^* U ) \longrightarrow ( F_{\alpha(2)}, F_{\alpha(1)} )
	\end{equation*}
	with the following properties.
	\begin{enumerate}[(1)]
		\item $ h^{-1} ( F_{\alpha(1)} ) = \rho^* U $, and $ h $ restricts to a homeomorphism $ Z \setminus \rho^* U \to F_{\alpha(2)} \setminus F_{\alpha(1)} $, namely the bundle identification of \autoref{lem:res_bundles}(2).
		\item Over $ F_{\alpha(1)} = U $, the restriction $ h |_{\rho^* U} \colon \rho^* U \to U $ is the canonical map of the pullback bundle, covering $ \rho \colon \widehat Y \to P \cong \Phi_{\alpha(1)} $.
		\item The connecting morphisms of the two pairs in Borel--Moore homology satisfy
		\begin{equation*}
			\partial_{( F_{\alpha(2)}, F_{\alpha(1)} )} = ( h |_{\rho^* U} )_* \circ \partial_{( Z, \rho^* U )} \circ \left( ( h |_{Z \setminus \rho^* U} )_* \right)^{-1} ,
		\end{equation*}
		and the Thom isomorphisms of $ \rho^* U $ and $ U $ intertwine $ ( h |_{\rho^* U} )_* $ with $ \rho_* $:
		\begin{equation*}
			Th \circ ( h |_{\rho^* U} )_* = \rho_* \circ Th .
		\end{equation*}
	\end{enumerate}
\end{lemma}

\begin{proof}
	Points of $ \mathcal{T} $ are triples $ (f, (b, K)) $ with $ f \in \widehat L(K) $, and points of $ \rho^* U $ are pairs $ (f, \hat y) $ with $ f \in L ( \left\{ \rho( \hat y ) \right\} ) $.
	Define
	\begin{equation*}
		h ( f, (b, K) ) \coloneq ( f, b ) , \qquad
		h ( f, \hat y ) \coloneq ( f, \delta_{\rho ( \hat y )} ) .
	\end{equation*}
	The two prescriptions agree on the glued subspace: a point of $ \mathcal{T} |_{\widehat Y} $ is $ (f, ( \delta_v, K )) $ with $ v \in K $, identified with $ (f, \hat y) \in \widehat W \subset \rho^* U $ where $ \rho ( \hat y ) = v $, and both prescriptions send it to $ (f, \delta_v) $.
	Since $ Z $ carries the quotient topology of $ \mathcal{T} \sqcup \rho^* U $, the map $ h $ is well defined and continuous.
	Its image lies in $ F_{\alpha(2)} $ because $ f \in \widehat L (K) \subseteq L ( \left\{ v \right\} ) $ for every $ v \in K $, and the fibres of $ F_{\alpha(2)} $ are the spaces $ L ( K_b ) $, respectively $ L ( \left\{ v \right\} ) $ over the vertex $ \delta_v $.
	($ Z $ is Hausdorff, being the adjunction of two metrisable spaces along a closed subspace.)

	Properness: let $ C \subset F_{\alpha(2)} $ be compact.
	Then the polynomial coordinates of the points of $ C $ range over a compact subset $ C_\Pi \subset \Pi_{4,2} $.
	Note that both $ \widehat \Phi $ and $ \widehat Y $ are compact.
	The preimage of $ C $ in $ \mathcal{T} $ is a closed subset of $ C_\Pi \times \widehat \Phi $, hence compact, and similarly for $ \rho^* U $ (with $ \widehat Y $).
	The preimage of $ C $ in $ Z $ is the image of the union of these two compacta under the quotient map, hence compact.

	$ (1) $ $ h ( f, (b, K) ) \in F_{\alpha(1)} $ if and only if $ b \in \Phi_{\alpha(1)} $, i.e. $ b = \delta_v $.
	But the points of $ \mathcal{T} $ with $ b \in \Phi_{\alpha(1)} $ are exactly those glued into $ \rho^* U $.
	So $ h^{-1} ( F_{\alpha(1)} ) = \rho^* U $, and on $ Z \setminus \rho^* U = \mathcal{T} |_{\widehat \Phi \setminus \widehat Y} $ the map $ h $ is the identification of \autoref{lem:res_bundles}(2).

	$ (2) $ Immediate from the definition.

	$ (3) $ $ h $ is a proper map of pairs, so it induces a morphism of the long exact sequences of the pairs in Borel--Moore homology.
	The square containing the connecting morphisms reads
	\mbox{$ \partial_{( F_{\alpha(2)}, F_{\alpha(1)} )} \circ ( h |_{Z \setminus \rho^* U} )_* = ( h |_{\rho^* U} )_* \circ \partial_{( Z, \rho^* U )} $},
	and the factor $ ( h |_{Z \setminus \rho^* U} )_* $ on the left is invertible by $ (1) $.

	For the Thom compatibility, recall the chain-level construction of the Thom isomorphism in the proof of \cite[Proposition 10]{BerGorKon23}.
	Denote by $ \pi_U \colon U \to P $ and $ \pi_{\rho^* U} \colon \rho^* U \to \widehat Y $ the bundle projections.
	Let $ u $ be a cochain that represents the Thom class of $ U $.
	The Thom isomorphism of $ U $ is induced by the chain map
	\begin{equation*}
		s \mapsto ( \pi_U )_\# ( s \frown u )
	\end{equation*}
	on locally finite chains.
	The chain $ s \frown u $ is finite: $ u $ vanishes on the simplices that miss the zero section, and the zero section is compact because $ P $ is.
	So the restriction to the chains subordinate to a cover, which is made in \cite{BerGorKon23} for a general base, is not needed here.

	Denote by $ \widetilde h $ the restriction $ h |_{\rho^* U} $, so that $ \pi_U \circ \widetilde h = \rho \circ \pi_{\rho^* U} $.
	The cochain $ \widetilde h^\# u $ represents the Thom class of $ \rho^* U $, and $ \widehat Y $ is compact, so the Thom isomorphism of $ \rho^* U $ is induced in the same way by $ \widetilde h^\# u $.
	By the projection formula for chains,
	\begin{equation*}
		( \pi_U )_\# \left( \widetilde h_\# s \frown u \right)
		= ( \pi_U )_\# \widetilde h_\# \left( s \frown \widetilde h^\# u \right)
		= \rho_\# ( \pi_{\rho^* U} )_\# \left( s \frown \widetilde h^\# u \right) .
	\end{equation*}
	Passing to homology, we get $ Th \circ \widetilde h_* = \rho_* \circ Th $.
\end{proof}

\begin{proposition}\label{prop:d1_composition}
	For every $ q $ the following diagram commutes.
	\begin{equation*}
		\begin{tikzcd}
			\overline H_{q+2} ( F_{\alpha(2)} \setminus F_{\alpha(1)} ) \arrow[rr, "d^1"]                                          &                                                                 & \overline H_{q+1} ( U ) \arrow[rd, "Th" {inner sep=2mm}, "\sim" {sloped, inner sep=.5mm}]                        &                \\
			\overline H_{q+2} ( Z \setminus \rho^* U ) \arrow[u, "h_*"' {inner sep=2mm}, "\sim" {rotate=90, anchor=north, inner sep=.5mm}] \arrow[d, "Th" {inner sep=2mm}, "\sim" {rotate=90, anchor=north, inner sep=.5mm}] \arrow[rr, "{\partial_{( Z, \rho^* U )}}"] &                                                                 & \overline H_{q+1} ( \rho^* U ) \arrow[u, "h_*"' {inner sep=2mm}] \arrow[d, "Th" {inner sep=2mm}, "\sim" {rotate=90, anchor=north, inner sep=.5mm}] & H_{q-17} ( P ) \\
			\overline H_{q-12} ( \widehat \Phi \setminus \widehat Y ) \arrow[r, "{\partial_{( \widehat \Phi, \widehat Y )}}"]      & H_{q-13} ( \widehat Y ) \arrow[r, "- \frown e ( \mathcal{N} )"] & H_{q-17} ( \widehat Y ) \arrow[ru, "\rho_*"']                   &
		\end{tikzcd}
	\end{equation*}
	The map $ h_* $ in the left column is induced by the homeomorphism $ h |_{Z \setminus \rho^* U} $.
	The maps $ Th $ are the Thom isomorphisms of $ \mathcal{T} |_{\widehat \Phi \setminus \widehat Y} $ (rank $ 7 $), of $ \rho^* U $ and of $ U $ (rank $ 9 $).
	Consequently, after we identify the source of $ d^1 $ with $ \overline H_{q-12} ( \widehat \Phi \setminus \widehat Y ) $ along the left column and its target with $ H_{q-17} ( P ) $ by the Thom isomorphism of $ U $, the differential $ d^1 $ becomes the composition
	\begin{equation*}
		\rho_* \circ \left( - \frown e ( \mathcal{N} ) \right) \circ \partial_{( \widehat \Phi, \widehat Y )} \colon \overline H_{q-12} ( \widehat \Phi \setminus \widehat Y ) \to H_{q-17} ( P ) .
	\end{equation*}
\end{proposition}

\begin{proof}
	By \autoref{sec:two_to_one_goal}, $ d^1 \colon E^1_{2,q} \to E^1_{1,q} $ is the connecting morphism of the pair $ ( F_{\alpha(2)}, F_{\alpha(1)} ) $.
	So the upper rectangle and the triangle on the right commute by \autoref{lem:transport}(3).

	The lower rectangle is the diagram of \autoref{sec:diff_vs_leq} for the pair $ ( \widehat \Phi, \widehat Y ) $, which is compact (\autoref{lem:phihat}(1)) and homeomorphic to a finite CW pair (\autoref{lem:phihat_cw}).
	The bundle data $ E_1 = \mathcal{T} $, $ E_2 = \rho^* U $, $ E_2' = \widehat W $ is supplied by \autoref{lem:res_bundles}, with $ E_1 $ defined over \emph{all} of $ \widehat \Phi $.
	The signs are $ +1 $ since all bundles are complex (see the remark at the end of \autoref{sec:diff_vs_leq}).
\end{proof}

\begin{theorem}\label{thm:2to1_nontrivial}
	The differential $ d^1 \colon E^1_{2,17} \to E^1_{1,17} $ is an isomorphism, and the differential $ d^1 \colon E^1_{2,19} \to E^1_{1,19} $ has rank $ 1 $.
\end{theorem}

\begin{proof}
	By \autoref{prop:d1_composition}, it is enough to prove both claims for the composition of that proposition on $ \overline H_{q-12} ( \widehat \Phi \setminus \widehat Y ) $ for $ q = 17 $ and $ q = 19 $, that is, in the degrees $ 5 $ and $ 7 $.
	By \autoref{prop:boundary_resolved}, this group is one-dimensional, and $ \partial_{( \widehat \Phi, \widehat Y )} $ maps it onto the line spanned by $ PD ( u ) $, where $ u = \varepsilon_0 ( \varepsilon_1 - \varepsilon_1' ) $ for $ q = 17 $ and $ u = \varepsilon_1 - \varepsilon_1' $ for $ q = 19 $.
	By \autoref{prop:cap_pushforward}, the class $ \rho_* \left( PD ( u ) \frown e ( \mathcal{N} ) \right) $ is non-zero.
	So in both cases $ d^1 $ is non-zero and its source is one-dimensional, hence $ d^1 $ has rank $ 1 $.

	For $ q = 17 $ the target is $ E^1_{1,17} \cong H_0 ( P ) \cong \Q $, so $ d^1 $ is an isomorphism.
\end{proof}

\begin{corollary}\label{cor:2to1_einfty}
	The entries of the columns $ 1 $ and $ 2 $ of $ E^1 $ that survive the differentials \mbox{$ d^1 \colon E^1_{2,q} \to E^1_{1,q} $} are those of the columns $ 1 $ and $ 2 $ of \autoref{table:main_sp_seq_pinf}.
\end{corollary}

\begin{proof}
	Only the columns $ 1 $ and $ 2 $ take part, so the sole differentials between them are those of the first page, $ d^1 \colon E^1_{2,q} \to E^1_{1,q} $, and $ E^2 = E^{\infty} $ there.
	Reading \autoref{table:main_sp_seq_p1} against \autoref{thm:2to1_nontrivial} gives
	\begin{equation*}
		\begin{array}{c|cc|c|cc}
			q  & E^1_{1,q} & E^1_{2,q} & \rk d^1 & E^{\infty}_{1,q} & E^{\infty}_{2,q} \\ \hline
			15 & 0         & \Q(8)     & 0       & 0                & \Q(8)            \\
			17 & \Q(9)     & \Q(9)     & 1       & 0                & 0                \\
			19 & 2\Q(10)   & \Q(10)    & 1       & \Q(10)           & 0                \\
			21 & 2\Q(11)   & 0         & 0       & 2\Q(11)          & 0                \\
			23 & \Q(12)    & 0         & 0       & \Q(12)           & 0
		\end{array}
	\end{equation*}
	where the rows $ q = 15, 21, 23 $ have a zero target or a zero source and so impose nothing.
	In the row $ q = 17 $ the differential is an isomorphism by \autoref{thm:2to1_nontrivial}, so both entries die.
	In the row $ q = 19 $ it has rank $ 1 $, so $ E^{\infty}_{2,19} = 0 $ and $ E^{\infty}_{1,19} = \Q(10) $.
	The surviving entries are precisely the columns $ 1 $ and $ 2 $ of \autoref{table:main_sp_seq_pinf}.
\end{proof}

\section{Differentials from the ninth to the eighth column}\label{sec:diff_9_to_8}

Let us use the construction in \autoref{sec:diff_vs_leq}.
The CW pair is
\begin{align*}
	X & \coloneq \Phi_{\alpha(9)} \setminus \Phi_{\alpha(7)} , \\
	Y & \coloneq \Phi_{\alpha(8)} \setminus \Phi_{\alpha(7)} .
\end{align*}
Hence $ X \setminus Y = \Phi_{\alpha(9)} \setminus \Phi_{\alpha(8)} $.

There is a commutative diagram
\begin{equation*}
	\begin{tikzcd}
		Y \arrow[rdd] \arrow[rr, hook] &                   & X \arrow[ldd] \arrow[rdd] &        & X \setminus Y \arrow[ll, hook'] \arrow[ldd] \\
		&                   &                           &        &                                                      \\
		& X_{8} \arrow[rr, "\sim" {anchor=south, inner sep=.4mm}, "\phi" {above, yshift=0.3em}] &                           & X_{9} &
	\end{tikzcd}
\end{equation*}
The isomorphism $ \phi $ arises from the fact that a pair $ \{ (p,l_p), (q,l_q) \} = \configparallelpair \in X_{8} $ corresponds uniquely to an element $ \{ (p,l_p), (q,l_q) \} \cup \{ (r,\overline{pq}) \text{ st } r \in \overline{pq}  \} = \configlinepair \in X_{9} $.
The maps from $ X $ are defined via this isomorphism.

The diagram above allows us to define a $ \C^{d_9} = \C^2 $-vector bundle over $ X = \Phi_{\alpha(9)} \setminus \Phi_{\alpha(7)}$ (recall $ d_9 = 2 $ from \autoref{table:configs_dim}):
\begin{equation*}
	E_1 \coloneq \left\{ (f,b) \in \Sigma \times \Phi_{\alpha(9)} \setminus \Phi_{\alpha(7)} \ \mid\ f \in L(K), \text{ where $ K $ is the image of $ b $ in } X_{9} \right\} .
\end{equation*}
Let
\begin{equation*}
	E_2 \coloneq F_{\alpha(8)} \setminus F_{\alpha(7)} \to \Phi_{\alpha(8)} \setminus \Phi_{\alpha(7)} = Y.
\end{equation*}
be the $ \C^3 $-vector bundle defined previously.

As in \autoref{sec:diff_vs_leq}, we obtain the diagram
\begin{equation*}
	\begin{tikzcd}
		E_1 \arrow[rr, "\C^{2}", two heads] &  & \Phi_{\alpha(9)} \setminus \Phi_{\alpha(7)}                  \\
		&  &                    \\
		 F_{\alpha(8)} \setminus F_{\alpha(7)} \arrow[rr, "\C^{3}", two heads] &  & \Phi_{\alpha(8)} \setminus \Phi_{\alpha(7)} \arrow[uu, hook]
	\end{tikzcd}
\end{equation*}

Let
\[\resizebox{\textwidth}{!}{$
E'_2 \coloneq E_1 |_{\Phi_{\alpha(8)} \setminus \Phi_{\alpha(7)}} = \left\{ (f,b) \in \Sigma \times \Phi_{\alpha(8)} \setminus \Phi_{\alpha(7)} \ \mid\ f \in L(K), \text{ where $ K $ is the image of $ b $ in } X_{9} \right\} . $}
\]
There is a natural inclusion $ \iota \colon E'_2 \to E_2 $.
Locally over a chart $ U \subset \Phi_{\alpha(8)} \setminus \Phi_{\alpha(7)} $, the bundle $ E'_2 $ is isomorphic to $ L( \configlinepair ) \times U = \C^2 \times U $, while $ E_2 $ is isomorphic to $ L( \configparallelpair ) \times U = \C^3 \times U $, where $ \phi( \configparallelpair ) = \configlinepair $.

Consider $ Z \coloneq E_1 \cup_{\iota} E_2 $.
One checks that $ Z \setminus E_2 = E_1 \setminus (E_1 |_{\Phi_{\alpha(8)} \setminus \Phi_{\alpha(7)}}) = F_{\alpha(9)} \setminus F_{\alpha(8)} $ as topological spaces.
But then the connecting morphism $ \partial_{(Z, E_2)} = \partial_{(Z, F_{\alpha(8))} \setminus F_{\alpha(7)})} $ in the long exact sequence of the pair $ (Z, E_2) $ is precisely the differential on the first page of the spectral sequence of filtration
\begin{equation*}
	\varnothing \subset F_{\alpha(8)} \setminus F_{\alpha(7)} \subset F_{\alpha(9)} \setminus F_{\alpha(7)} ,
\end{equation*}
which coincides with the differential from the ninth to the eighth column in the original spectral sequence.
Here we use the interpretation of the long exact sequence of a pair as the spectral sequence associated to a two-step filtration.
\begin{equation*}
	\begin{tikzcd}[column sep=small]
		\overline H_*(Z\setminus (F_{\alpha(8)} \setminus F_{\alpha(7)})) \arrow[dd, equal] \arrow[rrrrdd, "{\partial_{(Z,F_{\alpha(8)} \setminus F_{\alpha(7)})}}"] &  &                                                                                          &  &                                                                                           \\
		&  &                                                                                          &  &                                                                                           \\
		\overline H_*\left( E_1\setminus (E_1|_{\Phi_{\alpha(8)} \setminus \Phi_{\alpha(7)}}) \right) \arrow[dd, "Th" {inner sep=2mm}, "\sim" {rotate=90, anchor=north, inner sep=.5mm}] \arrow[rr, "{\partial_{(E_1,E_1|_{\Phi_{\alpha(8)} \setminus \Phi_{\alpha(7)}})}}"{yshift=1.2ex}]                               &  & \overline H_{*-1}(E_1 |_{\Phi_{\alpha(8)} \setminus \Phi_{\alpha(7)}}) \arrow[dd, "Th" {inner sep=2mm}, "\sim" {rotate=90, anchor=north, inner sep=.5mm}] \arrow[rr]     &  & \overline H_{*-1}(F_{\alpha(8)} \setminus F_{\alpha(7)}) \arrow[dd, "Th" {inner sep=2mm}, "\sim" {rotate=90, anchor=north, inner sep=.5mm}] \\
		&  &                                                                                          &  &                                                                                           \\
		\overline H_{*-2r_1}(\Phi_{\alpha(9)} \setminus \Phi_{\alpha(8)}) \arrow[rr, "{\partial_{(\Phi_{\alpha(9)} \setminus \Phi_{\alpha(7)},\Phi_{\alpha(8)} \setminus \Phi_{\alpha(7)})}}"{yshift=1.5ex}] &  & \overline H_{* - 2r_1 - 1}(\Phi_{\alpha(8)} \setminus \Phi_{\alpha(7)}) \arrow[rr, "-\frown e((F_{\alpha(8)} \setminus F_{\alpha(7)})/E'_2)"{yshift=1.5ex}] &  & \overline H_{* - 2r_2 - 1}(\Phi_{\alpha(8)} \setminus \Phi_{\alpha(7)})
	\end{tikzcd}
\end{equation*}
To show that the differentials from the ninth to the eighth column are all non-trivial, it suffices to show that both maps
\[
\partial_{(\Phi_{\alpha(9)} \setminus \Phi_{\alpha(7)},\Phi_{\alpha(8)} \setminus \Phi_{\alpha(7)})} = \partial_{(X,Y)} \text{ and } -\frown e((F_{\alpha(8)} \setminus F_{\alpha(7)})/E'_2) = -\frown e(E_2/E'_2)
\]
are non-trivial.

\subsection{Boundary morphism \texorpdfstring{$\partial_{(X,Y)}$}{}}

Consider $ K \coloneq \configlinepair \in X_{9} $.
Let $ M_1 \coloneq \vcenter{\hbox{\begin{tikzpicture}[scale=0.66]
			\draw[fill] (-0.5,0) circle (1.5pt);
			\draw (-0.5,-0.25) -- (-0.5,0.25);
			\draw (-1,0) -- (1,0);
\end{tikzpicture}}} $ and $ M_2 \coloneq \vcenter{\hbox{\begin{tikzpicture}[scale=0.66]
			\draw[fill] (0.5,0) circle (1.5pt);
			\draw (0.5,-0.25) -- (0.5,0.25);
			\draw (-1,0) -- (1,0);
\end{tikzpicture}}} $.
The configurations $ M_1 $, $ M_2 $ and $ \phi^{-1}(K) = \configparallelpair $ are precisely the maximal subconfigurations of $ K $.
Consider the following fibre bundle:
\begin{equation*}
	\begin{tikzcd}
		\widetilde{\Lambda}(K) \setminus \bigl( \widetilde{\Lambda}(M_1) \cup \widetilde{\Lambda}(M_2) \bigr) \arrow[r, hook] & \Phi_{\alpha(9)} \setminus \Phi_{\alpha(7)} \arrow[r, two heads, "\pi"] & X_{9} \ni K
	\end{tikzcd} .
\end{equation*}
The projection map on elements of $ \Phi_{\alpha(9)} \setminus \Phi_{\alpha(8)} $ is given as in the following fibre bundle:
\begin{equation*}
	\begin{tikzcd}
		\widetilde{\Lambda}(K) \setminus  \partial \widetilde \Lambda(K) \arrow[r, hook] & \Phi_{\alpha(9)} \setminus \Phi_{\alpha(8)} \arrow[r, two heads] & X_{9} \ni K
	\end{tikzcd} .
\end{equation*}
The projection map on elements of $ \Phi_{\alpha(8)} \setminus \Phi_{\alpha(7)} $ is the composition of the projection map of
\begin{equation*}
	\begin{tikzcd}
		\widetilde{\Lambda}(\phi^{-1}(K)) \setminus  \partial \widetilde \Lambda(\phi^{-1}(K)) \arrow[r, hook] & \Phi_{\alpha(8)} \setminus \Phi_{\alpha(7)} \arrow[r, two heads] & X_{8} \ni \phi^{-1}(K)
	\end{tikzcd}
\end{equation*}
and the isomorphism $ \phi \colon X_{8} \to X_{9} $ as above.

By definition, $  \partial \widetilde \Lambda(K) $ is the union of $ \widetilde{\Lambda}(M) $ indexed by $ M \subset K $ maximal subconfigurations.
So in our case $  \partial \widetilde \Lambda(K) = \widetilde{\Lambda}(M_1) \cup \widetilde{\Lambda}(M_2) \cup \widetilde{\Lambda}(\phi^{-1}(K)) $.
Therefore, the fibre is
\begin{equation*}
	\pi^{-1}(K) = \widetilde{\Lambda}(K) \setminus  \partial \widetilde \Lambda(K) \bigsqcup \widetilde{\Lambda}(\phi^{-1}(K)) \setminus  \partial \widetilde \Lambda(\phi^{-1}(K)) = \widetilde{\Lambda}(K) \setminus \bigl( \widetilde{\Lambda}(M_1) \cup \widetilde{\Lambda}(M_2) \bigr) .
\end{equation*}
The union $ \widetilde{\Lambda}(M_1) \cup \widetilde{\Lambda}(M_2) $ has the homology of a point: both terms are cones, and by \cite[Lemma 6]{BerGorKon23} their intersection is $ \widetilde{\Lambda}(M_1 \cap M_2) = \widetilde{\Lambda}(\configline) $, which is a cone as well, so \autoref{lem:tool_contr} applies.
Consequently, by \autoref{lem:tool_cone} $ \overline{H}_*(\pi^{-1}(K)) = 0 $.
Hence by \autoref{lem:tool_fibre_bundle} we obtain \mbox{$ \overline H_*(X) = \overline H_*(\Phi_{\alpha(9)} \setminus \Phi_{\alpha(7)}) = 0 $}, which implies that the connecting morphism \mbox{$ \partial_{(X,Y)} \colon \overline{H}_*(X \setminus Y) \to \overline{H}_{*-1}(Y) $} is an isomorphism in every possible degree.

\subsection{Cap product with Euler class \texorpdfstring{$ e(E_2 / E'_2) $}{}}

We redefine the problem via the following diagram:
\begin{equation*}
	\begin{tikzcd}
		E'_2 \arrow[rr, hook] \arrow[dd] &  & E_2 \arrow[rr, two heads] \arrow[dd] &  & E \arrow[rr] \arrow[dd]             &  & \Phi_{\alpha(8)} \setminus \Phi_{\alpha(7)} \arrow[dd, "\pi"]                 \\
		&  &                                      &  &                                     &  &                                    \\
		\overline{E}'_2 \arrow[rr, hook] &  & \overline{E}_2 \arrow[rr, two heads] &  & \overline{E} \arrow[rr]             &  & X_{8}                             \\
		&  &                                      &  &                                     &  &                                    \\
		&  &                                      &  & \widetilde{E} \arrow[rr] \arrow[uu] &  & \widetilde{X}_{8} \arrow[uu, "h"']
	\end{tikzcd}
\end{equation*}
We will explain what each part of the diagram means as we go.
The reader may find it useful to refer back to this diagram throughout the following lemmas and theorems.
These results form a sequence of reductions that reduces the problem of non-triviality of the cap product in question to an equivalent problem of non-triviality of a cup product.
We then explicitly calculate the cup product in terms of a known ring structure.

From now on, by ``non-trivial'' we will mean ``non-trivial in each degree where it can possibly be''.

Let $ E \coloneq E_2/E'_2 $ be a vector bundle over $ \Phi_{\alpha(8)} \setminus \Phi_{\alpha(7)} $.
Our goal is to show that the following map is non-trivial:
\begin{align}\label{eq:cap_9_to_8}
	- \frown e(E) \colon \overline H_*(\Phi_{\alpha(8)} \setminus \Phi_{\alpha(7)}) \to \overline H_{*-2(r_2-r_1)}(\Phi_{\alpha(8)} \setminus \Phi_{\alpha(7)}) .
\end{align}

\begin{lemma}
	There is a vector bundle $ \overline{E} $ over $ X_{8} $ such that \autoref{eq:cap_9_to_8} is non-trivial if and only if the following map is non-trivial:
	\begin{equation*}
		- \frown e(\overline{E}) \colon \overline H_*(X_{8}; \pm \Q) \to \overline H_{*-2(r_2-r_1)}(X_{8}; \pm \Q) .
	\end{equation*}
\end{lemma}

\begin{proof}
	Consider the projection map $ \pi \colon \Phi_{\alpha(8)} \setminus \Phi_{\alpha(7)} \to X_{8} $.
	Recall that its fibres are contractible, being open cones: by the construction of \cite[Section 1.2.2]{BerGorKon23}, $ \widetilde \Lambda(K) $ is the cone over $ \partial \widetilde \Lambda(K) $, which is compact by \cite[Lemma 8]{BerGorKon23}.
	Therefore, $\pi$ is a homotopy equivalence.
	Hence the pullback map is a bijection \mbox{$ \pi^* \colon \Vect^n(X_{8}) \to \Vect^n(\Phi_{\alpha(8)} \setminus \Phi_{\alpha(7)}) $} (as both base spaces are paracompact and Hausdorff).
	Moreover, taking quotients of vector bundles is compatible with pullback (by exactness).

	From this, we can define $ \overline E'_2 \subset \overline E_2 $ and $ \overline{E} = \overline E_2 / \overline E'_2 $ to be the vector bundles over $ X_{8} $ such that $ E'_2 = \pi^*(\overline E'_2), E_2 = \pi^*(\overline E_2) $ and $ E = \pi^*(\overline{E}) $ are their pullbacks via the projection map \mbox{$ \pi \colon \Phi_{\alpha(8)} \setminus \Phi_{\alpha(7)} \to X_{8} $}.

	From the fibre bundle $ \R \to \Phi_{\alpha(8)} \setminus \Phi_{\alpha(7)} \xrightarrow{\pi} X_{8} $ we obtain isomorphisms
	\[
	H^*(X_{8}) \to H^*(\Phi_{\alpha(8)} \setminus \Phi_{\alpha(7)}) \text{ and } \overline H_*(\Phi_{\alpha(8)} \setminus \Phi_{\alpha(7)}) \to \overline H_{*-1}(X_{8}; \pm \Q) .
	\]
	Moreover, they respect cap product in the following sense.

	Let  $ \iota \colon X_8 \to \Phi_{\alpha(8)} \setminus \Phi_{\alpha(7)} $ be the zero section embedding, so that $ \pi \circ \iota = Id $.
	Consider the diagram
	\begin{equation*}
		\begin{tikzcd}
			\overline H_p(\Phi_{\alpha(8)} \setminus \Phi_{\alpha(7)}) \times H^q(\Phi_{\alpha(8)} \setminus \Phi_{\alpha(7)}) \arrow[r] \arrow[d, "Th", shift right=11] & \overline H_{p-q}(\Phi_{\alpha(8)} \setminus \Phi_{\alpha(7)}) \arrow[d, "Th"] \\
			\overline H_{p-1}(X_{8}; \pm \Q) \times  H^q(X_{8}) \arrow[u, "\pi^*", shift right=15] \arrow[r]                                    & \overline H_{p-q-1}(X_{8}; \pm \Q)
		\end{tikzcd}
	\end{equation*}
	Here we used the fact that for a topological space $X$ the cap product maps as follows:
	\begin{equation*}
		- \frown - \colon \overline H_p(X, \mathscr{L}_1) \times H^q(X, \mathscr{L}_2) \to \overline H_{p-q}(X, \mathscr{L}_1 \otimes \mathscr{L}_2) .
	\end{equation*}
	We have $ \pm \Q \otimes \underline \Q = \pm \Q $, giving us the bottom morphism in the diagram above.

	Going clockwise from the top left corner we get
	\begin{equation*}
		(\iota_*(a), \pi^*(x)) \mapsto \iota_*(a) \frown \pi^*(x) \mapsto Th(\iota_*(a) \frown \pi^*(x)) .
	\end{equation*}
	Using the definition of the Thom isomorphism $ Th $ from \mbox{\cite[Proposition 10]{BerGorKon23}} and the projection formula, we find
	\begin{equation*}
		Th(\iota_*(a) \frown \pi^*(x)) = Th(\iota_*(a \frown \iota^* \pi^*(x))) = Th(\iota_*(a \frown x)) = (a \frown x) \frown u .
	\end{equation*}
	Going counter-clockwise, we get
	\begin{equation*}
		(\iota_*(a), \pi^*(x)) \mapsto (Th(\iota_*(a)), x) \mapsto Th(\iota_*(a)) \frown x = (a \frown u ) \frown x .
	\end{equation*}
	By the relation between cap and cup products,
	\begin{equation*}
		(a \frown x) \frown u = a \frown (x \smile u) = (-1)^{q} a \frown (u \smile x) = (-1)^{q} (a \frown u ) \frown x .
	\end{equation*}
	This shows that the cap product in the top row (compare to \autoref{eq:cap_9_to_8}) is non-trivial if and only if the cap product in the bottom row is non-trivial.
\end{proof}

\begin{lemma}
	The map in the lemma above is non-trivial if and only if the following map is non-trivial:
	\begin{equation*}
		- \smile e(\overline{E}) \colon H^*(X_{8}; \pm \Q) \to H^{*+2(r_2-r_1)}(X_{8}; \pm \Q) .
	\end{equation*}
\end{lemma}

\begin{proof}
	This follows from Poincar\'e duality (see \autoref{sec:map_alpha}).

	Here, Poincar\'e duality maps $ H^*(X; \mathscr{L}) \to \overline H_{\dim X - *}(X; \mathscr{L} \otimes \mathbb{O}) $, where $ \mathbb{O} $ is the orientation sheaf.
	In our case $ \mathscr{L} = \pm \Q $, and therefore $ \pm \Q \otimes \mathbb{O} = \pm \Q $, since the orientation sheaf is trivial for a complex manifold $ X $.
\end{proof}

Let $ \widetilde{X}_{8} $ be the double cover of $ X_{8} $ -- the ordered version of $ X_{8} $.
There is an action of $ S_2 $ on $ \widetilde{X}_8 $ such that \mbox{$ X_{8} = \widetilde{X}_{8} / S_2 $}.
Moreover, the action on $ \widetilde{X}_8 $ defines an action on $ H^*(\widetilde{X}_8) $.
Hence we can consider the decomposition \mbox{$ H^*(\widetilde{X}_8) = H^*(\widetilde{X}_8)^- \oplus H^*(\widetilde{X}_8)^+ $} into the invariant $ H^*(\widetilde{X}_8)^+ = H^*(\widetilde{X}_8)^{S_2} $ and anti-invariant $ H^*(\widetilde{X}_8)^- $.

Let $ \widetilde{E} $ be the pullback of $ \overline{E} $ via the covering map $ h \colon \widetilde{X}_8 \to \widetilde{X}_8 / S_2 = X_{8} $.

\begin{lemma}
	The map in the lemma above is non-trivial if and only if the following map is non-trivial:
	\begin{equation*}
		- \smile e(\widetilde{E}) \colon H^*(\widetilde{X}_8; \Q)^- \to H^{*+2(r_2-r_1)}(\widetilde{X}_8; \Q)^- .
	\end{equation*}
\end{lemma}

Since $ S_2 $ is finite and its order $ | S_2 | $ is invertible in $ \Q $, we have $ H^*(\widetilde{X}_8/S_2) = H^*(\widetilde{X}_{8})^{S_2} $.

\begin{proof}

	We now show that $ H^*(\widetilde{X}_8) = H^*(X_{8}; \pm \Q) \oplus H^*(X_{8}) $.

	Recall that $ S_2 $ acts freely on $ \widetilde{X}_8 $, and that $ h \colon \widetilde{X}_8 \to \widetilde{X}_8 / S_2 = X_{8} $ is a double cover.
	Now, there is a (piece of) long exact sequence of homotopy groups
	\begin{equation*}
		\pi_1(S_2) = 1 \to \pi_1(\widetilde{X}_8) \to \pi_1(X_{8}) \to S_2 \to 0 = \pi_0(X_{8}) .
	\end{equation*}
	Consider the regular representation $ \Q[S_2] $ of $ S_2 $.
	The composition
	\begin{equation*}
		\pi_1(X_{8}) \to S_2 \to \Aut (\Q[S_2])
	\end{equation*}
	defines a local system on $ X_{8} $.
	Moreover, it is related to the constant sheaf $ \underline{\Q}_{\widetilde{X}_8} $ on $ \widetilde{X}_8 $ as follows.

	Let $ g \in S_2 $ be non-trivial.
	We can take basis $ 1, g $ over $ \Q $ in $ \Q[ S_2 ] = \langle 1, g \rangle _{\Q} $.
	Then \mbox{$ \Q[ S_2] = R_- \oplus R_+ $}, where $ R_- = \langle 1-g \rangle _{\Q} $ and $ R_+ = \langle 1+g \rangle _{\Q} $ are irreducible components of the representation.
	The group $ S_2 $ acts as $ \pm 1 $ on $ R_- $, and on $ R_+ $ it acts trivially.
	Let $ L_-, L_+ $ be the corresponding local systems on $ X_{8} $.

	As $ h $ is a double covering map, $ R^q h_* $ vanishes for $ q > 0 $ (each fibre is just two points).
	Moreover, there is a decomposition \mbox{$ h_*\underline{\Q}_{\widetilde{X}_8} = L_- \oplus L_+ = \pm \Q \oplus \Q $}.
	Consequently, from Leray spectral sequence for $ h \colon \widetilde{X}_8 \to X_{8} $ we obtain
	\begin{equation*}
		H^*(\widetilde{X}_8; \Q) = H^*(X_{8}; R^0{h_*}\Q) = H^*(X_{8}; h_*\Q) = H^*(X_{8}; \pm \Q\oplus \Q) = H^*(X_{8}; \pm \Q) \oplus H^*(X_{8}; \Q) .
	\end{equation*}
	The isomorphism above is $ S_2 $-equivariant.
	It follows that the invariant part of $ H^*\left( \widetilde{X}_8; \Q \right) $ with respect to the action of $ S_2 $ is $ H^*(\widetilde{X}_8; \Q)^+ = H^*(\widetilde{X}_8; \Q)^{S_2} = H^*(X_{8}; \Q) $ and the anti-invariant part is $ H^*(\widetilde{X}_8; \Q)^- = H^*(X_{8}; \pm \Q) $.

	The image of $ h^* \colon H^*(X_{8}; \Q) \to H^*(\widetilde{X}_8; \Q) $ is precisely the invariant part.
	Thus $ h^* $ sends the Euler class $ e(\overline{E}) \in H^*(X_{8}; \Q) $ to $ e(\widetilde{E}) = h^*(e(\overline{E})) \in H^*(\widetilde{X}_8; \Q)^+ $.
	Moreover, for $ a \in H^*(\widetilde{X}_8; \Q)^- $, the product $ a \smile e(\widetilde{E}) $ is again in $ H^*(\widetilde{X}_8; \Q)^- $.
	Indeed, the action of $ g \in S_2 $ is an induced morphism $ g^* \colon H^*(\widetilde{X}_8; \Q) \to H^*(\widetilde{X}_8; \Q) $, so
	\[
	g^*(a \smile e(\widetilde{E})) = g^*(a) \smile g^*(e(\widetilde{E})) = (-a) \smile e(\widetilde{E}) = -(a \smile e(\widetilde{E})) .
	\]
	Hence, to prove non-triviality of the original map, it suffices to show that the map
	\begin{equation*}
		- \smile e(\widetilde{E}) \colon H^*(\widetilde{X}_8; \Q)^- \to H^{*+2(r_2-r_1)}(\widetilde{X}_8; \Q)^-
	\end{equation*}
	is non-trivial.
\end{proof}

We now compute the last map explicitly.
We show that $ \widetilde{X}_8 $ is homotopy equivalent to $ P = \PP T\C P^2 $.
As a result, $ H^*(\widetilde{X}_8; \Q) $ and $ e(\widetilde{E}) $ are represented in terms of the standard generators $ \varepsilon_0, \varepsilon_1, \varepsilon_2 $, and the multiplication is calculated explicitly.

\begin{proposition}\label{thm:wave_x8_homot_ptcp2}
	The space $ \widetilde{X}_{8} $ is homotopy equivalent to $ P = \PP T\C P^2 $.
\end{proposition}

Let us fix some notation for the rest of this paper.
Let ${ O \coloneq  (0:1:0) }$, ${ p \coloneq  (1:0:0) }$, ${ q \coloneq  (0:0:1) } $ be points of $ \C P^2 $, and let ${ l_p \coloneq \overline{Op} }$, ${ l_q \coloneq \overline{Oq} }$ be lines in $ \C P^2 $.
In homogeneous coordinates \mbox{$ (x : y : z) $}, we have ${ l_p = \{ z = 0 \} = \{ (x:y:0) \} }$ and ${ l_q = \{ x = 0 \} = \{ (0 : y : z) \} }$.
Let us also define line \mbox{$ l \coloneq \overline{pq} = \{ y = 0 \} = \{ (x : 0 : z )\} $}.

\begin{proof}
	The group $ G \coloneq \SL_3(\C) $ has a natural action on $ \widetilde{X}_8 $ induced by the natural action on $ P $ (recall that the elements of $ \widetilde X_8 $ are ordered pairs $ ((x,l_x), (y,l_y)) $, where $ (x,l_x) $ and $ (y,l_y) $ are in $ P $).

	The standard Torus subgroup $ T = G \cap \{\text{diagonal matrices}\} $ is the stabiliser of the ordered pair $ ((p,l_p), (q,l_q)) \in \widetilde{X}_{8} $.
	The orbit of $ ((p,l_p), (q,l_q)) $ is $ \widetilde{X}_{8} $.
	Indeed, an element of $ \widetilde{X}_{8} $ is completely determined by three ordered linearly independent vectors $ v_1, v_2, v_3 $ in $ \C^3 $.
	Hence the argument of \autoref{thm:transitivity} applies verbatim to show transitivity: the element $ ((p,l_p), (q,l_q)) $ maps to any other element of $ \widetilde{X}_{8} $.

	By the Orbit--Stabiliser Theorem, $ \widetilde{X}_{8} \cong \SL_3(\C) / T $.
	It follows from \autoref{rmk:homot-TB} that there is a homotopy equivalence \mbox{$ \widetilde{X}_{8} \cong \SL_3(\C) / T \sim \SL_3(\C) / B \cong P $}.
\end{proof}

\begin{proposition}\label{thm:chern_wave_e}
	The first Chern class of $ \widetilde{E} $ is
	\begin{equation*}
		c_1(\widetilde{E}) = 2\varepsilon_1 = -2(\varepsilon_0 + \varepsilon_2) .
	\end{equation*}
\end{proposition}

We want to use HVB machinery developed in \autoref{sec:hvb}.

Let $ K \coloneq \{ (p,l_p), (q,l_q) \} = \configparallelpair \in X_{8} $.
Let $ K' \coloneq \configlinepair $ be the element of $ X_{9} $ corresponding to $ K \subset K' $.

\begin{lemma}
	The fibre of the vector bundle $ \widetilde{E} \to \widetilde{X}_8 $ over a point covering $ K \in X_{8} $ as above is $ L(K) / L(K') = \langle  x^2z^2  \rangle  $.
\end{lemma}

\begin{proof}
	Since pullbacks preserve fibres,
	\begin{equation*}
		\widetilde{E}_{((p,l_p), (q,l_q))} = \overline{E}_{K} = E_{\left[K\right]} ,
	\end{equation*}
	where $ \left[K\right] $ is the principal vertex of a coherent simplex in $ \Phi_{\alpha(8)} \setminus \Phi_{\alpha(7)} $.
	But the last fibre is the quotient $ L(K) / L(K') $ of the fibres of $ E_2 $ and $ E'_2 \subset E_2 $.

	We begin by computing the generators of the quotient fibre $ L(K) / L(K') $.

	We already know
	\begin{equation*}
		L( \{ (p,l_p) \} ) = \langle  z^4 , yz^3 , y^2z^2 , y^3z , y^4 , xz^3 , xyz^2 , xy^2z , x^2z^2  \rangle  .
	\end{equation*}
	By interchanging $x$ and $z$, we obtain
	\begin{equation*}
		L( \{ (q,l_q) \} ) = \langle  x^4 , x^3y , x^2y^2 , xy^3 , y^4 , x^3z , x^2yz , xy^2z , x^2z^2  \rangle  .
	\end{equation*}
	Taking the intersection, we obtain
	\begin{equation*}
		L( \{ (p,l_p), (q,l_q) \} ) = L( \{ (p,l_p) \} ) \cap L( \{ (q,l_q) \} ) = \langle  y^4 , xy^2z , x^2z^2  \rangle  .
	\end{equation*}
	For the line $ l= \configline $, as in proof of \autoref{thm:dimensions}, we obtain
	\begin{equation*}
		L( \configline ) = \langle x^2y^2 , xy^3 , y^4 , xy^2z , y^3z , y^2z^2  \rangle .
	\end{equation*}
	To find the vector space corresponding to the union of $ l $ and the two one-point configurations $ \{ (p,l_p) \} $ and $ \{ (q,l_q) \} $, we determine the intersection
	\begin{equation*}
		L( \configlinepair ) = L( \{ (p,l_p), (q,l_q) \} ) \cap L( \configline ) = \langle  y^4 , xy^2z \rangle .
	\end{equation*}
	Taking the quotient, we obtain the quotient fibre
	\begin{equation*}
		L( \{ (p,l_p), (q,l_q) \} ) / L( \configlinepair ) = \langle  x^2z^2  \rangle .
	\end{equation*}
	Thus the fibre is generated by a single monomial $ x^2z^2 $.
\end{proof}

\begin{lemma}\label{lem:hvb_over_x8}
	The vector bundle $ \widetilde{E} \to \widetilde{X}_8 $ is a homogeneous vector bundle.
\end{lemma}

\begin{proof}
	As the pullback of a HVB under an equivariant map is again a HVB, it suffices to show that $ \overline E $ is a HVB.
	As the quotient of two HVBs is again a HVB, it suffices to show that $ \overline E_2 $ and $ \overline E'_2 $ are HVBs.

	Recall that $ E_2 = \left\{ (f,b) \in \Sigma \times \Phi_{\alpha(8)} \setminus \Phi_{\alpha(7)} \mid f \in L(\pi(b)) \right\} $.
	From the definitions of $ \overline E_2 $ and $ \pi \colon \Phi_{\alpha(8)} \setminus \Phi_{\alpha(7)} \to X_{8} $, we can deduce that $ \overline E_2 = \left\{ (f,K) \in \Sigma \times X_{8} \mid f \in L(K) \right\} $.

	From here, using \autoref{thm:wave_x8_homot_ptcp2} and the previous lemma, one checks that $ \overline E_2 $ is a HVB similarly to how it was done in the proof of \autoref{thm:euler_class_via_hvb}.
	Same for $ \overline E'_2 $.
\end{proof}

\begin{proof}{\it (\autoref{thm:chern_wave_e})}

	We now compute the corresponding character.
	Recall that the group $ G = \SL_3(\C) $ acts on the space of polynomials via $ (g . f)(v) = f(g^{-1}v) $, which induces the corresponding representation.
	Now,
	\begin{equation*}
		T(s,t) \cdot (x^2z^2) = (stx)^2(t^{-1}z)^2 = s^2x^2z^2 = \chi(s,t) \cdot x^2z^2 \implies \chi(s,t) = s^2 .
	\end{equation*}
	Next, we find the differential:
	\begin{equation*}
		\xi \coloneq \dif \chi = 2\varepsilon_1 = -2(\varepsilon_0 + \varepsilon_2) = c_1(\widetilde{E}) .
	\end{equation*}
\end{proof}

\begin{theorem}\label{thm:9to8_nontrivial}
	The map
	\begin{equation*}
		- \smile e(\widetilde{E}) \colon H^*(\widetilde{X}_8; \Q)^- \to H^{*+2(r_2-r_1)}(\widetilde{X}_8; \Q)^-
	\end{equation*}
	is non-trivial in each degree (and therefore the original mapping in \autoref{eq:cap_9_to_8} is non-trivial in each degree).
\end{theorem}

\begin{proof}
	The element $ \varepsilon_0 - \varepsilon_2 \in H^*(\widetilde{X}_8; \Q) $ is in $ H^*(\widetilde{X}_8; \Q)^- $.
	Indeed, the non-trivial element of $ S_2 $ switches $ \varepsilon_0 $ and $ \varepsilon_2 $, so the anti-symmetric part contains $ \varepsilon_0 - \varepsilon_2 $.
	See the proof of \mbox{\cite[Lemma 26]{BerGorKon23}}.

	We check that $ \varepsilon_0 - \varepsilon_2 \neq 0 $, $ (\varepsilon_0 - \varepsilon_2) \xi \neq 0 $, $ (\varepsilon_0 - \varepsilon_2)\xi^2 \neq 0 $.
	Since $ \varepsilon_0 - \varepsilon_2 \in H^*(\widetilde{X}_8; \Q)^- $ and $ (\varepsilon_0 - \varepsilon_2) \xi \in H^*(\widetilde{X}_8; \Q)^- $ generate all degrees in which the map could be non-trivial, this suffices.

	A direct computation (for example, via \cite[ideal\_membership.sage]{github_calc}) shows that $ \varepsilon_0 - \varepsilon_2 \notin I $, \mbox{$ (\varepsilon_0 - \varepsilon_2)\xi = -2 (\varepsilon_0^2 - \varepsilon_2^2) \notin I $}, and $ (\varepsilon_0 - \varepsilon_2)\xi^2 = 4 ( \varepsilon_0^3 + \varepsilon_0^2 \varepsilon_2 - \varepsilon_0 \varepsilon_2^2 - \varepsilon_2^3) \notin I $.

	It follows that the map $ - \smile e(\widetilde{E}) $ is non-trivial, and hence, by the diagram in \autoref{sec:diff_vs_leq}, the differentials from the ninth to the eighth column are non-trivial where possible.
\end{proof}

\section{Differentials from the thirteenth to the tenth column}\label{sec:diff_13_to_10}

The differentials to be computed are the two $ d^3 $ between the thirteenth and the tenth column of \autoref{table:main_sp_seq_p1},
\begin{equation*}
	d^3 \colon E^3_{13,0} \longrightarrow E^3_{10,2} \qquad \text{and} \qquad d^3 \colon E^3_{13,-2} \longrightarrow E^3_{10,0} ,
\end{equation*}
whose sources are the entries $ \Q(5) $ and $ \Q(4) $ of the thirteenth column and whose targets are the entries $ \Q(5) $ and $ \Q(4) $ of the tenth, these four entries being untouched by $ d^1 $ and $ d^2 $ (see the proof of \autoref{prop:13to10_iso}).
We show that both are isomorphisms.

Let
\begin{equation*}
	E \coloneq F_{\alpha(10)} \setminus F_{\alpha(9)} \longrightarrow \Phi_{\alpha(10)} \setminus \Phi_{\alpha(9)}
\end{equation*}
be the vector bundle of \autoref{lem:tool_vb}(1), of complex rank $ d_{10} = 1 $, and recall from the proof of \autoref{thm:spectral_groups} that $ \overline H_*(\Phi_{\alpha(10)} \setminus \Phi_{\alpha(9)}) $ is one-dimensional in the degrees $ 8 $, $ 10 $ and $ 12 $ and zero otherwise.
The argument has two parts.
In \autoref{sec:13to10_cap} we show that the cap product with the Euler class of $ E $ on $ \overline H_*(\Phi_{\alpha(10)} \setminus \Phi_{\alpha(9)}) $, the map \autoref{eq:cap_13_to_10} below, is non-trivial in each degree where it can be, namely in the degrees $ 10 $ and $ 12 $.
This is the same computation as in \autoref{sec:diff_9_to_8}, and it proceeds by the same chain of reductions.
In \autoref{sec:13to10_zero_section} we pass from that to the two differentials.
In \autoref{sec:diff_9_to_8} the corresponding passage was made by the diagram of \autoref{sec:diff_vs_leq}.
Here we instead compute the Borel--Moore homology of the open subset $ F_{\alpha(13)} \setminus F_{\alpha(9)} $ of $ \sigma $ directly, by a filtration through the zero section as in the proof of \cite[Lemma 23]{BerGorKon23}, and read the differentials off the main spectral sequence restricted to that open subset, which converges to the same groups.

\subsection{Cap product with the Euler class \texorpdfstring{$ e(E) $}{}}\label{sec:13to10_cap}

Our goal is to show that the following map is non-trivial, that is, non-trivial in the degrees $ 10 $ and $ 12 $, the only ones in which its source and its target are both non-zero:
\begin{align}\label{eq:cap_13_to_10}
	- \frown e(E) \colon \overline H_*(\Phi_{\alpha(10)} \setminus \Phi_{\alpha(9)}) \longrightarrow \overline H_{*-2}(\Phi_{\alpha(10)} \setminus \Phi_{\alpha(9)}) .
\end{align}

\begin{lemma}
	There is a vector bundle $ \overline{E} $ over $ X_{8} $ such that \autoref{eq:cap_13_to_10} is non-trivial if and only if the following map is non-trivial:
	\begin{equation*}
		- \frown e(\overline{E}) \colon \overline H_*(X_{8}; \pm \Q) \longrightarrow \overline H_{*-2}(X_{8}; \pm \Q) .
	\end{equation*}
\end{lemma}

\begin{proof}
	Consider the fibre bundle
	\begin{equation*}
		\widetilde{\Lambda}(K) \setminus  \partial \widetilde \Lambda(K) \longrightarrow \Phi_{\alpha(10)} \setminus \Phi_{\alpha(9)} \longrightarrow^{\pi} B(\C P^{2\vee}, 2) = X_{10} \ni K = \configcross{} .
	\end{equation*}
	Here, the fibre bundle projection map $ \pi $ sends a point $ x \in \Phi_{\alpha(10)} \setminus \Phi_{\alpha(9)} $ to the unique pair of lines such that either $ x $ is equal to the configuration corresponding to the two lines, or lies on an interval leading to the configuration.
	Let
	\begin{equation*}
		\overline E_{10} \coloneq \left\{ (f, K) \in \Pi \times X_{10} \ \middle | \ f \in L(K) \right\} \longrightarrow X_{10} .
	\end{equation*}
	By part $ (3) $ of \cite[Condition list 2]{BerGorKon23}, the map $ X_{10} \to \Gr(1, \Pi) $, $ K \mapsto L(K) $, is continuous.
	It can be verified that the condition holds for our family.
	So $ \overline E_{10} $ is the pullback of the tautological line bundle of $ \Gr(1, \Pi) $ along this map, hence a line bundle over $ X_{10} $ by \cite[Lemma 10]{BerGorKon23}.
	By \cite[Theorem 4]{BerGorKon23}, the bundle $ E $ is the pullback of the same tautological bundle along $ x \mapsto L(\pi(x)) $, the fibre of $ E $ over $ x $ being $ L(K) $ for $ K = \pi(x) $ the principal vertex of the simplex containing $ x $.
	Hence $ E = \pi^*(\overline E_{10}) $.

	There is also the forgetful map
	\[
	\psi \colon X_{8} \longrightarrow X_{10} = B(\C P^{2\vee}, 2) \colon \left\{ (p, l_p), (q, l_q) \right\} \longmapsto \left\{ l_p, l_q \right\} .
	\]
	It is a fibre bundle: over an unordered pair $ \{ l_1, l_2 \} $ its fibre is $ (l_1 \setminus \{ l_1 \cap l_2 \}) \times (l_2 \setminus \{ l_1 \cap l_2 \}) \simeq \C \times \C $, since by the definition of $ X_{8} $ the two points may be chosen independently on the two lines and both differ from the intersection point.
	Note that $ X_{8} $ is not the product $ B(\C P^{2\vee}, 2) \times \C^2 $: the monodromy of the unordered pair interchanges the two factors.
	The fibre is contractible, however, so $ \psi $ is a homotopy equivalence and induces the bijection \mbox{$ \psi^* \colon \Vect^n(X_{10}) \to \Vect^n(X_{8}) $}.

	Put $ \overline E \coloneq \psi^*(\overline E_{10}) $, a vector bundle over $ X_{8} $.
	The local system $ \pm \Q $ on $ X_{8} $ is the pullback along $ \psi $ of the one on $ X_{10} $, both being induced by the interchange of the two lines, so $ \psi $ identifies the two cap products.
\end{proof}

From this point, the argument proceeds identically to that in the previous section.

\begin{lemma}
	The map in the lemma above is non-trivial if and only if the following map is non-trivial:
	\begin{equation*}
		- \smile e(\overline{E}) \colon H^*(X_{8}; \pm \Q) \longrightarrow H^{*+2}(X_{8}; \pm \Q) .
	\end{equation*}
\end{lemma}

\begin{lemma}
	The map in the lemma above is non-trivial if and only if the following map is non-trivial:
	\begin{equation*}
		- \smile e(\widetilde{E}) \colon H^*(\widetilde{X}_8; \Q)^- \longrightarrow H^{*+2}(\widetilde{X}_8; \Q)^- .
	\end{equation*}
\end{lemma}

\begin{proposition}
	The space $ \widetilde{X}_{8} $ is (still) homotopy equivalent to $ P = \PP T\C P^2 $.
\end{proposition}

Let $ O, p, q, l_p, l_q $ all be as in the previous section, and let $ K \coloneq \{ (p,l_p), (q,l_q) \} = \configparallelpair \in X_{8} $ be as there.
Let $ K'' \coloneq \configcross{O} \in X_{10} $ be the configuration whose two lines are $ l_p $ and $ l_q $, so that $ \psi(K) = \{ l_p , l_q \} $ is the pair of lines underlying $ K'' $.

\begin{lemma}
	The fibre of the vector bundle $ \widetilde{E} \to \widetilde{X}_8 $ over a point covering $ K \in X_{8} $ as above is $ L(K'') = \langle  x^2z^2  \rangle  $.
\end{lemma}

\begin{proof}
	Recall that pullbacks preserve fibres.
	Now, the fibre of $ \overline E = \psi^*(\overline E_{10}) $ over $ K $ is the fibre of $ \overline E_{10} $ over $ \psi(K) $, which in turn is exactly the same as the fibre of $ E $ over the principal vertex $ [K''] $, that is, $ L(K'') $.
	As in the proof of \autoref{thm:dimensions}, $ L(K'') $ is spanned by the square of the product of the two linear forms cutting out $ l_p $ and $ l_q $, so that ${ L(K'') = \langle (zx)^2 \rangle = \langle x^2z^2 \rangle }$, in agreement with the entry $ ax^2y^2 $ in the row $ i = 10 $ of \autoref{table:configs_dim}.
	Consequently, the fibre of $ \widetilde{E} $ is also $ L(K'') = \langle  x^2z^2  \rangle  $.
\end{proof}

\begin{lemma}
	The vector bundle $ \widetilde{E} \to \widetilde{X}_8 $ is a homogeneous vector bundle.
\end{lemma}

\begin{proof}
	Using the previous lemma and the proposition before it, this follows by the same argument as in the proof of \autoref{lem:hvb_over_x8}.
\end{proof}

\begin{proposition}
	The first Chern class of $ \widetilde{E}$ is
	\begin{equation*}
		c_1(\widetilde{E}) = 2\varepsilon_1 = -2(\varepsilon_0 + \varepsilon_2) .
	\end{equation*}
\end{proposition}

\begin{theorem}\label{thm:13to10_nontrivial}
	The map
	\begin{equation*}
		- \smile e(\widetilde{E}) \colon H^*(\widetilde{X}_8; \Q)^- \longrightarrow H^{*+2}(\widetilde{X}_8; \Q)^-
	\end{equation*}
	is non-trivial (and therefore the original mapping in \autoref{eq:cap_13_to_10} is non-trivial).
\end{theorem}

\begin{proof}
	The proof of \autoref{thm:9to8_nontrivial} applies verbatim.
	Indeed, $ \widetilde X_8 $ is the same space as there, carrying the same action of $ S_2 $ and hence the same anti-invariant part $ H^*(\widetilde X_8; \Q)^- $.
	By the lemmas above $ \widetilde E $ is again a line bundle, and by the proposition above
	\begin{equation*}
		e(\widetilde E) = c_1(\widetilde E) = \xi = 2\varepsilon_1 = -2(\varepsilon_0 + \varepsilon_2) ,
	\end{equation*}
	the same class as there.
	And $ E $ is again a line bundle, so that the degree shift is the one occurring there.
	The map in question is therefore literally the map of \autoref{thm:9to8_nontrivial}, and the three non-vanishing statements verified there (i.e. $ \varepsilon_0 - \varepsilon_2 \notin I $, $ (\varepsilon_0 - \varepsilon_2)\xi \notin I $ and $ (\varepsilon_0 - \varepsilon_2)\xi^2 \notin I $) give the claim here as well.
\end{proof}

\subsection{From the cap product to the differentials}\label{sec:13to10_zero_section}

The idea is to compute the Borel--Moore homology of the open subset $ F_{\alpha(13)} \setminus F_{\alpha(9)} $ of $ \sigma $ in the degrees $ 10 $ to $ 13 $ in two ways.
Restricted to this open subset, the main spectral sequence keeps the two $ d^3 $ and has no other entries in these degrees, so both $ d^3 $ are isomorphisms once this homology vanishes.
To show that it does, we filter the open subset through its zero section, which contributes only in lower degrees.
Off the zero section only $ E $ contributes, and there the inclusion of the zero section is, under the Thom isomorphism, the cap product of \autoref{sec:13to10_cap}, whose non-triviality leaves homology only in the degrees $ 9 $ and $ 14 $.

Let $ W \coloneq F_{\alpha(13)} \setminus F_{\alpha(9)} = \sigma \setminus F_{\alpha(9)} $.
It is an open subset of $ \sigma $, since $ F_{\alpha(9)} = p^{-1}(\Phi_{\alpha(9)}) $ and $ \Phi_{\alpha(9)} $ is closed in $ \widetilde \Lambda $ (\cite[Lemma 8]{BerGorKon23}).
Let $ \sigma_0 \coloneq \{ 0 \} \times \widetilde \Lambda \subset \sigma $ be the zero section of $ p \colon \sigma \to \widetilde \Lambda $.
It lies in $ \sigma $ because $ 0 \in L(K) $ for every configuration $ K $, it is closed in $ \sigma $, being the intersection of $ \sigma $ with the closed subset $ \{ 0 \} \times \widetilde \Lambda $ of $ \Sigma \times \widetilde \Lambda $, and over each stratum $ \Phi_{\alpha(i)} \setminus \Phi_{\alpha(i-1)} $ it is the zero section of the vector bundle of \autoref{lem:tool_vb}(1), whose fibre over $ x $ is the subspace $ L(K_x) \subset \Pi $, $ K_x $ being the principal vertex of the simplex containing $ x $ (\cite[Theorem 4]{BerGorKon23}).
We write
\begin{equation*}
	W^{0} \coloneq W \cap \sigma_0 = \{ 0 \} \times ( \widetilde \Lambda \setminus \Phi_{\alpha(9)} ) \qquad \text{and} \qquad E^{\neq 0} \coloneq E \setminus \sigma_0
\end{equation*}
for the zero section of $ W $ and for the complement of the zero section in $ E = F_{\alpha(10)} \setminus F_{\alpha(9)} $.

The following two facts are instances of the naturality of the long exact sequence of a closed pair with respect to restriction to an open subset.
\begin{remark}\label{rem:filtered_ss}
	Let $ T $ be a space filtered by closed subsets $ \emptyset = T_0 \subset T_1 \subset \ldots \subset T_m = T $, and consider the Borel--Moore homology spectral sequence of this filtration, with first page \mbox{$ E^1_{p,q} = \overline H_{p+q}(T_p \setminus T_{p-1}) $}.
	\begin{enumerate}[leftmargin=*, label=(\arabic*)]
		\item The differential $ d^1 \colon \overline H_*(T_p \setminus T_{p-1}) \to \overline H_{*-1}(T_{p-1} \setminus T_{p-2}) $ is the connecting morphism of the closed pair $ ( T_p \setminus T_{p-2} , T_{p-1} \setminus T_{p-2} ) $.
		\item If $ U \subset T $ is open, then the restriction to $ U $ induces a morphism from the spectral sequence of $ T $ to that of $ U $ with the filtration $ T_p \cap U $, which on the first page is given by the restriction maps $ \overline H_*(T_p \setminus T_{p-1}) \to \overline H_*( (T_p \setminus T_{p-1}) \cap U ) $.
	\end{enumerate}
\end{remark}

Part (1) holds since, by construction, $ d^1 $ is the connecting morphism of $ (T_p, T_{p-1}) $ followed by the restriction from $ T_{p-1} $ to its open subset $ T_{p-1} \setminus T_{p-2} $, and restricting the pair $ (T_p, T_{p-1}) $ to the open subset $ T_p \setminus T_{p-2} $ of $ T_p $ turns the one into the other.

Part (2) holds since the families of restriction maps \mbox{$ \overline H_*(T_p \setminus T_{p-1}) \to \overline H_*( (T_p \setminus T_{p-1}) \cap U ) $} and \mbox{$ \overline H_*(T_p) \to \overline H_*(T_p \cap U) $} form a morphism of exact couples.
The induced morphism of spectral sequences is also compatible with the restriction $ \overline H_*(T) \to \overline H_*(U) $ on the abutments.

\begin{lemma}\label{lem:phi9_vanishing}
	$ H_k(\Phi_{\alpha(9)}) = 0 $ for $ k \geq 8 $.
\end{lemma}

\begin{proof}
	The space $ \Phi_{\alpha(9)} $ is compact and $ \Phi_{\alpha(1)} \subset \Phi_{\alpha(2)} \subset \ldots \subset \Phi_{\alpha(9)} $ is a filtration by closed subsets (\cite[Lemma 8]{BerGorKon23}), so its spectral sequence converges to \mbox{$ \overline H_*(\Phi_{\alpha(9)}) = H_*(\Phi_{\alpha(9)}) $}, and its first page consists of the groups $ \overline H_{*}(\Phi_{\alpha(p)} \setminus \Phi_{\alpha(p-1)}) $, $ 1 \leq p \leq 9 $, in the $ p $-th column.
	These groups were computed in the proof of \autoref{thm:spectral_groups} on the way to \autoref{table:main_sp_seq_p1}, of which they are the columns with the Thom shift of \autoref{lem:tool_vb_homology} removed.
	The non-zero ones among them are $ \overline H_*(\Phi_{\alpha(1)}) = H_*(\PP T\C P^2) $, concentrated in the degrees $ 0 $, $ 2 $, $ 4 $, $ 6 $, then $ \overline H_*(\Phi_{\alpha(2)} \setminus \Phi_{\alpha(1)}) $ and $ \overline H_*(\Phi_{\alpha(4)} \setminus \Phi_{\alpha(3)}) $, both concentrated in the degrees $ 3 $, $ 5 $, $ 7 $, then $ \overline H_*(\Phi_{\alpha(8)} \setminus \Phi_{\alpha(7)}) $, concentrated in the degrees $ 7 $, $ 9 $, $ 11 $, and $ \overline H_*(\Phi_{\alpha(9)} \setminus \Phi_{\alpha(8)}) $, concentrated in the degrees $ 8 $, $ 10 $, $ 12 $.
	So the only entries of the first page in total degree at least $ 8 $ lie in the columns $ 8 $ and $ 9 $.
	The differential $ d^1 $ from the column $ 9 $ to the column $ 8 $ is, by part (1) of \autoref{rem:filtered_ss}, the connecting morphism of the pair $ ( \Phi_{\alpha(9)} \setminus \Phi_{\alpha(7)} , \Phi_{\alpha(8)} \setminus \Phi_{\alpha(7)} ) $, which is the map $ \partial_{(X,Y)} $ of \autoref{sec:diff_9_to_8} and was shown there to be an isomorphism in every degree, because $ \overline H_*(\Phi_{\alpha(9)} \setminus \Phi_{\alpha(7)}) = 0 $.
	Hence the columns $ 8 $ and $ 9 $ vanish on the second page, every entry of which therefore lies in total degree at most $ 7 $, and the claim follows.
\end{proof}

\begin{lemma}\label{lem:zero_section_complement}
	$ \overline H_n(E^{\neq 0}) $ is one-dimensional for $ n = 9 $ and for $ n = 14 $, and zero for every other $ n $.
\end{lemma}

\begin{proof}
	In this proof write $ Y \coloneq \Phi_{\alpha(10)} \setminus \Phi_{\alpha(9)} $, and let $ s \colon Y \to E $ be the zero section, a closed embedding with image $ E \cap \sigma_0 $.
	The long exact sequence of the closed pair $ (E, s(Y)) $ reads
	\begin{equation*}
		\ldots \to \overline H_n(Y) \xrightarrow{s_*} \overline H_n(E) \to \overline H_n(E^{\neq 0}) \to \overline H_{n-1}(Y) \xrightarrow{s_*} \overline H_{n-1}(E) \to \ldots
	\end{equation*}
	By \cite[Proposition 10]{BerGorKon23}, applied to $ E $ as in \autoref{sec:diff_vs_leq}, the Thom isomorphism $ Th $ of $ E $ satisfies
	\begin{equation*}
		Th(s_*(a)) = a \frown e(E) \in \overline H_{n-2}(Y) \qquad \text{for all } a \in \overline H_n(Y) ,
	\end{equation*}
	so under $ Th \colon \overline H_n(E) \to \overline H_{n-2}(Y) $ the map $ s_* $ becomes the map \autoref{eq:cap_13_to_10}.
	Now $ \overline H_*(Y) $ is one-dimensional in the degrees $ 8 $, $ 10 $, $ 12 $ and zero otherwise, so $ \overline H_*(E) \cong \overline H_{*-2}(Y) $ is one-dimensional in the degrees $ 10 $, $ 12 $, $ 14 $ and zero otherwise, and by \autoref{thm:13to10_nontrivial} the map \autoref{eq:cap_13_to_10} is an isomorphism in the degrees $ 10 $ and $ 12 $, while in every other degree its source or its target is zero.
	The exact sequence gives
	\begin{equation*}
		0 \to \operatorname{coker} \left( s_* \colon \overline H_n(Y) \to \overline H_n(E) \right) \to \overline H_n(E^{\neq 0}) \to \ker \left( s_* \colon \overline H_{n-1}(Y) \to \overline H_{n-1}(E) \right) \to 0 .
	\end{equation*}
	The cokernel is non-zero only for $ n = 14 $, where $ \overline H_{14}(Y) = 0 $ and $ \overline H_{14}(E) \cong \overline H_{12}(Y) $ is one-dimensional, and the kernel is non-zero only for $ n - 1 = 8 $, where $ \overline H_8(Y) $ is one-dimensional and $ \overline H_8(E) \cong \overline H_6(Y) = 0 $.
\end{proof}

\begin{proposition}\label{prop:13to10_iso}
	The differentials $ d^3 \colon E^3_{13,0} \to E^3_{10,2} $ and $ d^3 \colon E^3_{13,-2} \to E^3_{10,0} $ of the main spectral sequence are isomorphisms.
\end{proposition}

\begin{proof}
	The first part of the proof, \emph{the zero-section filtration}, shows that $ \overline H_n(W) = 0 $ for $ 10 \leq n \leq 13 $.
	The second, \emph{the restricted main spectral sequence}, passes to the main spectral sequence restricted to $ W $, which has the same two $ d^3 $.
	In the total degrees $ 10 $ to $ 13 $ the restricted spectral sequence has only the four entries of the columns $ 10 $ and $ 13 $.
	By the first part they all vanish on the limit page, and only the two $ d^3 $ can kill them, so both are isomorphisms.

	\emph{The zero-section filtration.}
	Consider
	\begin{equation}\label{eq:zero_section_filtration}
		W^{0} \subset W^{0} \cup ( F_{\alpha(10)} \setminus F_{\alpha(9)} ) \subset W^{0} \cup ( F_{\alpha(11)} \setminus F_{\alpha(9)} ) \subset W^{0} \cup ( F_{\alpha(12)} \setminus F_{\alpha(9)} ) = W .
	\end{equation}
	Each $ F_{\alpha(i)} \setminus F_{\alpha(9)} = F_{\alpha(i)} \cap W $ is closed in $ W $, and so is $ W^{0} = W \cap \sigma_0 $, so \autoref{eq:zero_section_filtration} is a filtration of $ W $ by closed subsets.
	The last inclusion is an equality because $ d_{13} = 0 $: the vector bundle $ F_{\alpha(13)} \setminus F_{\alpha(12)} \to \Phi_{\alpha(13)} \setminus \Phi_{\alpha(12)} $ has rank zero, so $ F_{\alpha(13)} \setminus F_{\alpha(12)} $ is its own zero section and is contained in $ W^{0} $.
	The successive differences of \autoref{eq:zero_section_filtration} are $ W^{0} $, then $ E^{\neq 0} $, then the complements of the zero sections in the line bundles $ F_{\alpha(11)} \setminus F_{\alpha(10)} \to \Phi_{\alpha(11)} \setminus \Phi_{\alpha(10)} $ and $ F_{\alpha(12)} \setminus F_{\alpha(11)} \to \Phi_{\alpha(12)} \setminus \Phi_{\alpha(11)} $, of ranks $ d_{11} = d_{12} = 1 $.
	The Borel--Moore homology of the last two vanishes.
	Indeed, the proof of \autoref{thm:spectral_groups} shows that $ \overline H_*(\Phi_{\alpha(11)} \setminus \Phi_{\alpha(10)}) = 0 $ and $ \overline H_*(\Phi_{\alpha(12)} \setminus \Phi_{\alpha(11)}) = 0 $, so the total spaces of the two line bundles have vanishing Borel--Moore homology by \autoref{lem:tool_vb_homology}, and then so have the complements of their zero sections, by the long exact sequences of the closed pairs formed by a total space and its zero section.
	As for $ W^{0} \cong \widetilde \Lambda \setminus \Phi_{\alpha(9)} $, the space $ \widetilde \Lambda = \Phi_{\mathscr{N}} $ is compact, and it is contractible: every coherent simplex is a face of one with principal vertex the configuration $ \PP T\C P^2 $, the unique element of $ X_{13} $, so $ \widetilde \Lambda = \widetilde \Lambda ( \PP T\C P^2 ) $ is the cone over $ \partial \widetilde \Lambda ( \PP T\C P^2 ) $ with that vertex (\cite[Section 1.2.2]{BerGorKon23}).
	The long exact sequence of the closed pair $ ( \widetilde \Lambda , \Phi_{\alpha(9)} ) $ therefore gives $ \overline H_k(W^{0}) \cong H_{k-1}(\Phi_{\alpha(9)}) $ for $ k \geq 2 $, and this vanishes for $ k \geq 9 $ by \autoref{lem:phi9_vanishing}.
	So on the first page of the spectral sequence of \autoref{eq:zero_section_filtration} the first column is concentrated in the total degrees at most $ 8 $, the second column, equal to $\overline H_{2+q}(E^{\neq 0})$ in row $q$, is one-dimensional in the total degrees $ 9 $ and $ 14 $ and zero otherwise (\autoref{lem:zero_section_complement}), and the third and fourth columns are zero.
	Consequently
	\begin{equation*}
		\overline H_n(W) = 0 \qquad \text{for } 10 \leq n \leq 13 .
	\end{equation*}
	\emph{The restricted main spectral sequence.}
	The filtration of $ \sigma $ by the subspaces $ F_{\alpha(i)} $ restricts on the open subset $ W $ to the filtration
	\begin{equation*}
		\emptyset \subset F_{\alpha(10)} \setminus F_{\alpha(9)} \subset F_{\alpha(11)} \setminus F_{\alpha(9)} \subset F_{\alpha(12)} \setminus F_{\alpha(9)} \subset F_{\alpha(13)} \setminus F_{\alpha(9)} = W
	\end{equation*}
	by closed subsets of $ W $, the terms $ F_{\alpha(i)} \cap W $ with $ i \leq 9 $ being empty.
	Let $ E^*(W) $ denote its spectral sequence.
	Its first page consists of the columns $ 10 $ to $ 13 $ of \autoref{table:main_sp_seq_p1}, since $ F_{\alpha(p)} \setminus F_{\alpha(p-1)} \subset W $ for $ p \geq 10 $, and it converges to $ \overline H_*(W) $.
	By part (2) of \autoref{rem:filtered_ss}, the restriction to $ W $ induces a morphism of spectral sequences $ E^* \to E^*(W) $ which on the first page is the identity on the columns $ 10 $ to $ 13 $ and zero on the columns up to $ 9 $.

	In both spectral sequences the four entries at $ (13,0) $, $ (13,-2) $, $ (10,2) $ and $ (10,0) $ survive unchanged to the third page.
	Indeed, the columns $ 11 $ and $ 12 $ are zero and there is no column $ 14 $, so the entries of the column $ 13 $ are untouched by $ d^1 $ and $ d^2 $.
	The $ d^1 $ out of $ E^1_{10,q} $ lands in $ E^1_{9,q} $, which is zero for $ q = 0 $ and $ q = 2 $, the $ d^2 $ out of $ E^2_{10,q} $ lands in a subquotient of $ E^1_{8,q+1} $, which is zero for $ q + 1 = 1 $ and $ q + 1 = 3 $, and the $ d^1 $ and $ d^2 $ into the column $ 10 $ come from the zero columns $ 11 $ and $ 12 $.
	The morphism $ E^* \to E^*(W) $ commutes with the differentials and is the identity on these four entries on the third page, so the two $ d^3 $ of the proposition coincide with the corresponding differentials of $ E^*(W) $.

	In $ E^*(W) $ the entries of total degree $ 10 $, $ 11 $, $ 12 $ and $ 13 $ on the first page are exactly $ E^1_{10,0} = \Q(4) $, $ E^1_{13,-2} = \Q(4) $, $ E^1_{10,2} = \Q(5) $ and $ E^1_{13,0} = \Q(5) $, and since the columns $ 10 $ and $ 13 $ are the only non-zero ones, the only differentials of $ E^*(W) $ with a non-zero source and a non-zero target are the two $ d^3 $ between these four entries.
	As $ E^*(W) $ converges to $ \overline H_*(W) $ and $ \overline H_n(W) = 0 $ for $ 10 \leq n \leq 13 $, all four entries vanish on $ E^4(W) = E^\infty(W) $.
	Hence $ d^3 \colon E^3_{13,-2}(W) \to E^3_{10,0}(W) $ and $ d^3 \colon E^3_{13,0}(W) \to E^3_{10,2}(W) $ are surjective, and being maps between one-dimensional spaces they are isomorphisms.
\end{proof}

As a check, $ E^\infty(W) $ retains the single entry $ E^1_{10,4} = \Q(6) $, so $ \overline H_{14}(W) $ is one-dimensional, in agreement with the spectral sequence of \autoref{eq:zero_section_filtration}, whose only entry of total degree $ 14 $ is the one of \autoref{lem:zero_section_complement}, and $ \overline H_9(W) = 0 $, so the entry of total degree $ 9 $ of \autoref{lem:zero_section_complement} is killed by the $ d^1 $ into $ \overline H_8(W^{0}) \cong H_7(\Phi_{\alpha(9)}) $.

\section{Assembly of the differentials}\label{sec:assembly}

We can now complete the proof of \autoref{thm:ch2_einfty}.

\begin{proof}{\it (\autoref{thm:ch2_einfty})}
	We start from $ E^1 $, \autoref{table:main_sp_seq_p1}, and recall that $ d^r $ maps $ E^r_{p,q} $ to $ E^r_{p-r,q+r-1} $.

	On $ E^1 $ the only bidegrees in which both source and target are non-zero are \mbox{$ (2,19) \to (1,19) $} and $ (2,17) \to (1,17) $, and $ (9,7) \to (8,7) $ and $ (9,5) \to (8,5) $.
	The first two are non-trivial by \autoref{thm:2to1_nontrivial}: the first has rank $ 1 $ and the second is an isomorphism.
	The last two are non-trivial by \autoref{thm:9to8_nontrivial}, and are isomorphisms since source and target are one-dimensional.
	Passing to $ E^2 $ therefore leaves the column $ 1 $ with $ \Q(12), 2\Q(11), \Q(10) $ in the rows $ 23, 21, 19 $, the column $ 2 $ with $ \Q(8) $ in the row $ 15 $, the column $ 8 $ with $ \Q(8) $ in the row $ 9 $, and the column $ 9 $ with $ \Q(5) $ in the row $ 3 $.
	The columns $ 4 $, $ 10 $ and $ 13 $ are untouched.

	On $ E^2 $ every differential has a zero source or target: the only candidates are $ (4,q) \to (2,q+1) $ and $ (10,q) \to (8,q+1) $, and the surviving entries of the columns $ 2 $ and $ 8 $ sit in the rows $ 15 $ and $ 9 $, which are not of the form $ q + 1 $ for a non-zero $ E^2_{4,q} $, respectively $ E^2_{10,q} $.
	So $ E^3 = E^2 $.

	On $ E^3 $ the only bidegrees with a non-zero source and a non-zero target are $ (13,0) \to (10,2) $ and $ (13,-2) \to (10,0) $, and both are isomorphisms by \autoref{prop:13to10_iso}.
	They remove the rows $ 2 $ and $ 0 $ from the column $ 10 $ and the whole of the column $ 13 $.

	After that no differential has both a non-zero source and a non-zero target.
	The surviving entries are the single entries of the columns $ 2 $, $ 8 $, $ 9 $, $ 10 $ and the three entries of the column $ 4 $, together with the column $ 1 $, and no two of them are joined by any $ d^r $.
	So $ E^4 = E^{\infty} $, and the surviving entries are exactly those of \autoref{table:main_sp_seq_pinf}.
\end{proof}

\section*{Acknowledgements}

MD was partially funded by the Deutsche Forschungsgemeinschaft (DFG, German Research Foundation) under Germany's Excellence Strategy -- The Berlin Mathematics Research Center MATH+ (EXC-2046/2, project ID: 390685689).

\section*{Tool use disclosure}

AI tools (specifically, Claude models Opus and Fable) were used in the preparation of this paper for proofreading and drafting some parts of the document.
We subsequently re-checked and edited these parts, and accept full responsibility for the correctness of all results presented in this paper.

Where indicated, the computer algebra computations use SageMath \cite{sage}, with the code in \cite{github_calc}.

\printbibliography

\end{document}